\documentclass[leqno,a4paper,11pt]{report}
\usepackage{amssymb}
\usepackage{latexsym}
\usepackage{fullpage} 
\usepackage{amscd}
\usepackage{amsmath}
\usepackage{youngtab}
\usepackage{graphicx}
\usepackage{geometry}
\usepackage[latin1]{inputenc}
\usepackage{mathdots}
\usepackage{color}
\usepackage[color,matrix,arrow, curve]{xy}
\usepackage{enumerate}
\usepackage{lscape}
\usepackage{setspace}
\def\Z{\mathbb{Z}}

\def\Q{\mathbb{Q}}
\def\N{\mathbb{N}}

\def\t{\tau_\chi}

\def\End{\mathrm{End}\,}
\def\Hom{\mathrm{Hom_A}\,}

\def\Im{\mathrm{Im}}
\def\rad{\mathrm{rad}}
\def\Irr{\mathrm{Irr}}
\def\Coker{\mathrm{Coker}}
\def\Ker{\mathrm{Ker}}
\def\Ext{\mathrm{Ext}}
\def\soc{\mathrm{soc}}

\def\add{\mathrm{add}}

\def\#{\sharp}

\def\F{\mathfrak{F}}

\def\Aut{\mathrm{Aut}}

\def\bea{\begin{eqnarray*}}
\def\eea{\end{eqnarray*}}

\def\E{\mathcal{E}}

\def\Id{\mathrm{Id}}

\def\<#1,#2>{\langle\,#1,\,#2\,\rangle}
\newtheorem{Th}{Theorem} [section]
\newtheorem{Ex}[Th]{Example}
\newtheorem{Def}[Th]{Definition}
\newtheorem{Lemma}[Th]{Lemma}
\newtheorem{Cor}[Th]{Corollary}

\def\proof{\noindent\textbf{Proof:}}

\def\qed{\quad\hfill\ensuremath{\Box}}

\begin{document}


%
%

\onehalfspacing
\begin{titlepage}
\quad\vspace{0.5cm}\\
\begin{center}
{\bfseries\Huge Auslander-Reiten theory}\vspace{0.3cm}\\
{\bfseries\Huge in functorially finite}\vspace{0.3cm}\\
{\bfseries\Huge resolving subcategories}
\end{center}
\vspace{0.4in}
\begin{figure}[h!]
\centering
\includegraphics[width=4cm]{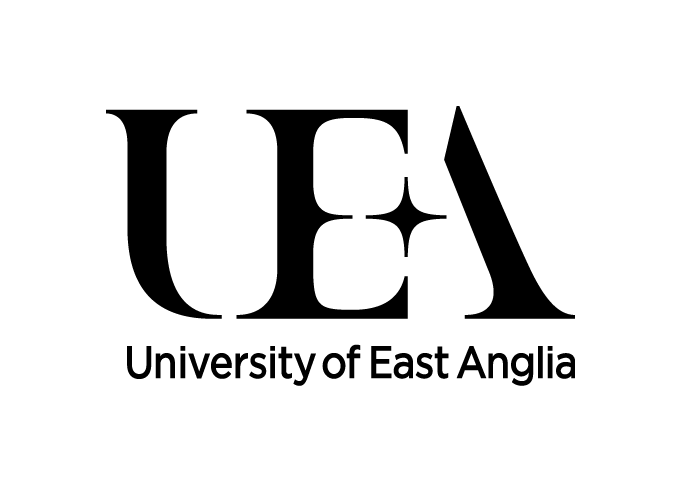}
\end{figure}
\Large
\vspace{0.6in}
\begin{center}
 {{A thesis submitted to the University of East Anglia in partial fulfilment of the requirements for the degree of Doctor of Philosophy}}
\end{center}
\vspace{0.8in}
\begin{center}
{\LARGE Matthias Krebs}\\
{\large School of Mathematics, UEA, Norwich, NR4 7TJ England}
\end{center}
\vspace{0.2in}
\begin{center}
{March 29, 2013}
\end{center}
\vfill
{\normalsize  \copyright \, This copy of the thesis has been supplied on condition that anyone
who consults it is understood to recognise that its copyright rests with
the author and that use of any information derived therefrom must be
in accordance with current UK Copyright Law. In addition, any
quotation or extract must include full attribution.}
\end{titlepage}

\chapter*{Abstract}
\pagenumbering{arabic}
\setcounter{page}{2}

We analyze the Auslander-Reiten quiver $\Gamma_\Omega$ of a functorially finite resolving subcategory $\Omega$. Chapter 1 gives a short introduction into the basic definitions and theorems of Auslander-Reiten theory in $A$-mod. We generalize these definitions and theorems in Chapter 2 and find a constant $p$ such that $l(X) \leq pl(Y)$ if there is an $\Omega$-irreducible morphism from $X$ to $Y$. This constant enables us to prove the Brauer-Thrall 1.5 conjecture for $\Omega$. Moreover, we find a connection between sectional paths in $A$-mod and irreducible morphisms in $\Omega$.
\bigskip
\\In Chapter 3 we introduce degrees of irreducible morphisms and use this notion to prove the generalization of the Happel-Preiser-Ringel theorem for $\Omega$. Finally, in Chapter 4, we analyze left stable components of $\Gamma_\Omega$ and find out that their left subgraph types are given by Dynkin diagrams if and only if $\Omega$ is finite. In the preparation of the proof we discover connected components with certain properties and name them helical components due to their shape. It turns out later that these components are the same as coray tubes. In the final section we discuss under which conditions the length of modules tends to infinity if we knit to the left in a component and give a complete description of all connected components in which this is not the case.

\chapter*{Acknowledgements}
I am greatly indebted to my supervisors Dr Vanessa Miemietz and Prof Shaun Stevens for all the assistance and many helpful discussions about my project. I would also like to thank Prof Volodymyr Mazorchuk for making it possible to come to Uppsala University alongside Dr Vanessa Miemietz for research visits. Finally, this thesis would not have been possible without my family, whose support and encouragement have been invaluable during the last couple of years.\vspace{0.2cm}\par
This work was supported by a UEA studentship and I am very grateful for the opportunity provided.

\tableofcontents

 \chapter*{Introduction}
 \label{ch:Introduction}
 \addcontentsline{toc}{chapter}{Introduction}
When the notion of almost split sequences was first introduced in the early 1970s its significance and impact on representation theory of Artin algebras was not immediately recognized. These short exact sequences contain information about how certain morphisms split and their existence was first proved in \cite{AR74}. More intensive studies of almost split sequences began after Ringel embedded them into Auslander-Reiten quivers \cite{R78}. These quivers give a complete description how any morphism between finitely generated modules splits into irreducible morphisms.
\bigskip
\\It has been proved that almost split sequences do not only exist in $A$-mod, the category of all finitely generated left A-modules over an Artin algebra $A$, but also in certain subcategories \cite{AS81}, \cite{R91}, \cite{K97}. This implies that there are also Auslander-Reiten quivers in these categories. One of the most important example of these categories is $\F(\Delta)$, the category of standard-filtered modules of a quasi-hereditary algebra. 
\bigskip
\\Since the vertices of Auslander-Reiten quivers are given by the isomorphism classes of indecomposable modules in $A$-mod or a subcategory $\chi$, these quivers are finite if and only if $A$-mod or $\chi$ is finite respectively, i.e. there are up to isomorphism only finitely many indecomposable modules in $A$-mod or $\chi$ respectively. Naturally the question arose if there are criteria on the Auslander-Reiten quiver of $A$ for the algebra to be representation finite. A partial answer to this question was given by Riedtmann by introducing the notion of tree type of stable translation quivers \cite{R80}, which allowed her to classify all selfinjective algebras of finite representation type by their Auslander-Reiten quivers \cite{R80'}, \cite{R83}, \cite{BLR81}. She proved that if an algebra is representation finite, then the tree type of each stable component of its Auslander-Reiten quiver is given by a Dynkin diagram.%
\bigskip
\\We generalize this statement and prove this generalization for Auslander-Reiten quivers of any algebra or certain subcategories that have Auslander-Reiten quivers. In order to establish this, we introduce a concept similar to the tree type which does not only work for stable components of Auslander-Reiten quivers.
\newpage
In the whole dissertation $A$ denotes an associative Artin algebra with a multiplicative identity over an algebraically closed field $K$. For convenience, we suppose that $A$ is indecomposable as this simplifies the notation significantly. Nevertheless, all results can be applied to decomposable algebras as all indecomposable direct summands of a decomposable algebra have disjoint Auslander-Reiten quivers. $X, Y$ and $Z$ usually denote modules in $A$-mod, which is the category of all finitely generated left A-modules. Simple, projective and injective modules are denoted by $S$, $P$ and $I$ respectively. If for a morphism $f : Z \to Y$ there is a module $X$ and morphisms $g: Z \to X$, $h: X \to Y$, we say $f$ factors through $g$ and $h$ and $f$ factors over $X$. The Jordan-H\"{o}lder length of a module $X$ is denoted by $l(X)$. Most definitions are taken from \cite{ARS95}. Note that for some theorems we also provide the dual statements, but only give one proof and leave the dual to the reader.
\bigskip
\\When we speak about Dynkin diagrams we always mean undirected Dynkin diagrams, which are classified below. In the infinite series $A_n$ and $D_n$ the index $n$ equals the number of vertices in the diagram.
 $$\begin{xy}
\xymatrixrowsep{0.3in}
\xymatrixcolsep{0.5in}
  \xymatrix@!0{ 
     A_n: & \bullet \ar@{-}[r]  	& 	\bullet \ar@{-}[r] 	& \cdots \ar@{-}[r]     &	\bullet \ar@{-}[r]     &   \bullet & n\geq 1 }
\end{xy}$$ 
$$\begin{xy}
\xymatrixrowsep{0.3in}
\xymatrixcolsep{0.5in}
  \xymatrix@!0{ 
     D_n: & \bullet \ar@{-}[r]  	& 	\bullet \ar@{-}[r] 	& \cdots \ar@{-}[r]     &	\bullet \ar@{-}[r]     &   \bullet & n\geq 4 \\
     	  &				&				&			&	\bullet \ar@{-}[u]	& }
\end{xy}$$ 
$$\begin{xy}
\xymatrixrowsep{0.3in}
\xymatrixcolsep{0.5in}
  \xymatrix@!0{ 
     E_6: & \bullet \ar@{-}[r]  	& 	\bullet \ar@{-}[r] 	& \bullet \ar@{-}[r]     &	\bullet \ar@{-}[r]     &   \bullet \\
     	  &				&				&	\bullet \ar@{-}[u]	& }
\end{xy}$$ 
$$\begin{xy}
\xymatrixrowsep{0.3in}
\xymatrixcolsep{0.5in}
  \xymatrix@!0{ 
E_7: & \bullet \ar@{-}[r]  	& \bullet \ar@{-}[r]  	&	\bullet \ar@{-}[r] 	& \bullet \ar@{-}[r]     &	\bullet \ar@{-}[r]     &  \bullet \\
      &				&			&				&\bullet \ar@{-}[u]	& }
\end{xy}$$ 
$$\begin{xy}
\xymatrixrowsep{0.3in}
\xymatrixcolsep{0.5in}
  \xymatrix@!0{ 
E_8: & \bullet \ar@{-}[r] & \bullet \ar@{-}[r]  & \bullet \ar@{-}[r] 	&\bullet \ar@{-}[r] & \bullet \ar@{-}[r]     &	\bullet \ar@{-}[r]     &  \bullet \\
      &			&			&			&			&\bullet \ar@{-}[u]	& }
\end{xy}$$ 
Moreover, the following are the Euclidean diagrams, which are obtained by adding one vertex to a Dynkin diagram such that the new diagram is not Dynkin anymore and every finite diagram without loops that is not a Dynkin diagram contains a Euclidean diagram. In the infinite series $A_n$ and $D_n$ the number of vertices in the diagram is $n+1$.
 $$\begin{xy}
\xymatrixrowsep{0.3in}
\xymatrixcolsep{0.5in}
  \xymatrix@!0{ 
     \widetilde{A}_n: & \bullet \ar@{-}[r]  	& 	\bullet \ar@{-}[r] 	& \cdots \ar@{-}[r]     &	\bullet \ar@{-}[r]     &   \bullet & n\geq 1 \\
      &			&		&		\bullet \ar@{-}[ull] \ar@{-}[urr] }
\end{xy}$$ 
$$\begin{xy}
\xymatrixrowsep{0.3in}
\xymatrixcolsep{0.5in}
  \xymatrix@!0{ 
     \widetilde{D}_n: & \bullet \ar@{-}[r]  	& 	\bullet \ar@{-}[r] 	& \cdots \ar@{-}[r]     &	\bullet \ar@{-}[r]     &   \bullet & n\geq 4 \\
     	  &				&	\bullet \ar@{-}[u]	&			&	\bullet \ar@{-}[u]	& }
\end{xy}$$ 
$$\begin{xy}
\xymatrixrowsep{0.3in}
\xymatrixcolsep{0.5in}
  \xymatrix@!0{ 
     \widetilde{E}_6: & \bullet \ar@{-}[r]  	& 	\bullet \ar@{-}[r] 	& \bullet \ar@{-}[r]     &	\bullet \ar@{-}[r]     &   \bullet \\
     	  &				&				&	\bullet \ar@{-}[u]	& \\
     	  &				&				&	\bullet \ar@{-}[u]	&}
\end{xy}$$ 
$$\begin{xy}
\xymatrixrowsep{0.3in}
\xymatrixcolsep{0.5in}
  \xymatrix@!0{ 
\widetilde{E}_7: & \bullet \ar@{-}[r]  	& \bullet \ar@{-}[r]  	&	\bullet \ar@{-}[r] 	& \bullet \ar@{-}[r]     &	\bullet \ar@{-}[r] &	\bullet \ar@{-}[r]     &  \bullet \\
      &				&			&				&\bullet \ar@{-}[u]	& }
\end{xy}$$ 
$$\begin{xy}
\xymatrixrowsep{0.3in}
\xymatrixcolsep{0.5in}
  \xymatrix@!0{ 
\widetilde{E}_8: & \bullet \ar@{-}[r] & \bullet \ar@{-}[r]  & \bullet \ar@{-}[r] & \bullet \ar@{-}[r]	&\bullet \ar@{-}[r] & \bullet \ar@{-}[r]     &	\bullet \ar@{-}[r]     &  \bullet \\
      &			&			&	&		&			&\bullet \ar@{-}[u]	& }
\end{xy}$$ 
Consequently, the only diagrams without loops which are not Dynkin but do not contain a Euclidean diagram are the following infinite diagrams.
 $$\begin{xy}
\xymatrixrowsep{0.3in}
\xymatrixcolsep{0.5in}
  \xymatrix@!0{ 
     A_\infty: & \bullet \ar@{-}[r]  	& 	\bullet \ar@{-}[r] 	& \bullet \ar@{-}[r]     &	\bullet \ar@{-}[r]     &   \bullet \ar@{-}[r] & \cdots }
\end{xy}$$ 
$$\begin{xy}
\xymatrixrowsep{0.3in}
\xymatrixcolsep{0.5in}
  \xymatrix@!0{ 
   D_\infty: & \bullet \ar@{-}[r]  	& 	\bullet \ar@{-}[r] 	& \bullet \ar@{-}[r]     &	\bullet \ar@{-}[r]     &   \bullet \ar@{-}[r] & \cdots \\
     	  &				&	\bullet \ar@{-}[u]	&			&	 }
\end{xy}$$ 
 $$\begin{xy}
\xymatrixrowsep{0.3in}
\xymatrixcolsep{0.5in}
  \xymatrix@!0{ 
     A_\infty^\infty: & \cdots \ar@{-}[r] &	\bullet \ar@{-}[r]  	& 	\bullet \ar@{-}[r] 	& \bullet \ar@{-}[r]     &	\bullet \ar@{-}[r]     &   \bullet \ar@{-}[r] & \cdots }
\end{xy}$$

\chapter{Introduction to Auslander-Reiten theory}

This chapter provides a short introduction to the basic definitions and most important theorems of Auslander-Reiten theory alongside with a few examples of Auslander-Reiten quivers. For further reading we recommend \cite{ARS95}.
\section{Almost split sequences}

\begin{Def}
\mbox{}
\begin{enumerate}[(a)]
\item We say a morphism $f : X \to Y$ is a \textbf{split epimorphism} if the identity on $Y$ factors through $f$, that is there is a morphism $g : Y \to X$ such that $fg = \Id_Y$.
$$\begin{xy}
  \xymatrix{ 
    &X \ar[dr]^f  	&	 \\
    Y  \ar[rr]^{\Id} \ar@{.>}[ur]^g &  	&	Y	    }
\end{xy}$$
\item Dually, a morphism $g : Y \to X$ is called a \textbf{split monomorphism} if there is a morphism $f : X \to Y$ such that $fg = Id_Y$.
$$\begin{xy}
  \xymatrix{ 
    &X \ar@{.>}[dr]^f  	&	 \\
    Y  \ar[rr]^{\Id} \ar[ur]^g &  	&	Y	    }
\end{xy}$$
\end{enumerate}
\end{Def}

\begin{Def}
\mbox{}
\begin{enumerate}[(a)]
\item A morphism $f : X \to Y$ is called \textbf{right almost split} if it is not a split epimorphism and any morphism $Z \to Y$ that is not a split epimorphism factors through $f$.
$$\begin{xy}
  \xymatrix{ 
    Z \ar@{.>}[d] \ar[dr]  	&	 \\
    X  \ar[r]^{f}    	&	Y	    }
\end{xy}$$ 
\item We dually define a morphism $g : Y \to X$ to be \textbf{left almost split} if it is not a split monomorphism and any morphism $Y \to Z$ that is not a split monomorphism factors through $g$.
$$\begin{xy}
  \xymatrix{ 
    Y  \ar[r]^g \ar[dr]  	& X \ar@{.>}[d]	 \\
    &	Z	    }
\end{xy}$$ 
\end{enumerate}
\end{Def}

For convenience, we say that a right almost split morphism that maps to a module $Y$ is a right almost split morphism for $Y$. Dually, we call a left almost split morphism $g : Y \to X$ a left almost split morphism for $Y$. An example for a right almost split morphism is the embedding $i: \rad P \to P$ for an indecomposable projective module $P$. Dually, the morphism $\Pi : I \to I/\soc I$ is left almost split for an indecomposable injective module $I$.

\begin{Lemma} \cite[V. Lemma 1.7]{ARS95}
Let $f: X \to Y$ be a morphism.
\begin{enumerate}[(a)]
 \item If $f$ is right almost split, then $Y$ is an indecomposable module. 
 \item If $f$ is left almost split, then $X$ is an indecomposable module. 
\end{enumerate}
\end{Lemma}

\begin{Def}
A short exact sequence
$$\begin{xy}
  \xymatrix{ 
      0 \ar[r]  	& 	X \ar[r]^g 	&	Y \ar[r]^f     &   Z \ar[r]	&	0	    }
\end{xy}$$
is called an \textbf{almost split sequence} if $g$ is left almost split and $f$ is right almost split.
\end{Def}

Moreover, in an almost split sequence $g$ and $f$ are \textbf{left} and \textbf{right minimal morphisms} respectively, i.e. for all $h: Y \to Y$ such that $hg = g$ or $fh = f$, $h$ is an isomorphism. Morphisms that are both right minimal and right almost split are called \textbf{minimal right almost split morphisms}, dually, a \textbf{minimal left almost split morphism} is both left minimal and left almost split. Minimal right almost split and minimal left almost split morphisms exist for every indecomposable module and their domains and codomains are unique up to isomorphism \cite[V.1]{ARS95}. This is important for the definition of Auslander-Reiten quivers in the next section.
\bigskip
\\It has been shown in \cite[V.1]{ARS95} that in an almost split sequence the modules $X$ and $Z$ are determined by each other and can be computed using certain functors. We consider a minimal projective presentation 
$$\begin{xy}
  \xymatrix{ 
      P_1 \ar[r]^f 	&	P_0 \ar[r]     &   Z \ar[r]	&	0	    }
\end{xy}$$
of Z. Let $T$ denote the contravariant functor $\Hom(\ \ ,A):$ $A$-mod $\to$ $A^{op}$-mod. Then $\Coker(T(f))$ is an $A^{op}$-module, which we call the transpose of $Z$ or $Tr(Z)$. Note that $Tr(Z)$ is zero if and only if $Z$ is projective. We then call $\tau = DTr$ the \textbf{Auslander-Reiten translation}, where $D$ denotes the usual duality $\Hom(\ \ ,K)$. 

\begin{Lemma} \cite[IV. Proposition 1.10]{ARS95}
$DTr: A$-mod $\to$ $A$-mod induces a bijection between the set of isomorphism classes of non-projective indecomposable modules and the set of isomorphism classes of non-injective indecomposable modules with $TrD$ as inverse.
\end{Lemma}

Consequently, we write $\tau^{-1}$ for $TrD$. By definition we have $\tau(P) = \tau^{-1}(I) = 0$ for a projective module $P$ and an injective module $I$ respectively. Nevertheless, for an indecomposable module $X$ and an integer $n$ we say $\tau^n(X)$ exists only if $\tau^n(X) \neq 0$.

\begin{Th} \cite[V. Theorems 1.15, 1.16]{ARS95}
\begin{enumerate}[(a)]
\item Let $Z$ be an indecomposable non-projective $A$-module. Then there is an up to isomorphism unique almost split sequence 
$$\begin{xy}
  \xymatrix{ 
      0 \ar[r]  	& 	\tau(Z) \ar[r]^g 	&	Y \ar[r]^f     &   Z \ar[r]	&	0	    }
\end{xy}$$
ending in $Z$.
\item Let $X$ be an indecomposable non-injective $A$-module. Then there is an up to isomorphism unique almost split sequence 
$$\begin{xy}
  \xymatrix{ 
      0 \ar[r]  	& 	X \ar[r]^g 	&	Y \ar[r]^f     &   \tau^{-1}(X) \ar[r]	&	0	    }
\end{xy}$$
starting in $X$.
\end{enumerate}
\end{Th}

\section{Auslander-Reiten quivers}

The Auslander-Reiten quiver of an algebra is a \textbf{biquiver}, i.e. a quiver with two disjoint sets of arrows, containing information on all almost split sequences. For its definition we need to consider a different point of view on the morphisms in an almost split sequence.

\begin{Def}
A morphism $f: X \to Y$ is called \textbf{irreducible} if $f$ is neither a split monomorphism nor a split epimorphism and, if $f = gh$ for some morphisms $h : X \to M$ and $g : M \to Y$, then either $h$ is a split monomorphism or $g$ is a split epimorphism.
\end{Def}

For example every minimal left or right almost split morphism is also an irreducible morphism. The following theorem describes the close connection between the two concepts.

\begin{Th} \cite[V. Theorem 5.3]{ARS95} \label{connectionirralmostsplit}
\begin{enumerate}[(a)]
\item Let $Z$ be an indecomposable module and $Y$ a non-zero module in $A$-mod. A morphism $g: Y \to Z$ is irreducible if and only if there exists a morphism $g' : Y' \to Z$ such that the induced morphism $(g,g') : Y \oplus Y' \to Z$ is a minimal right almost split morphism. 
\item Dually, if $X$ is an indecomposable module, then a morphism $f: X \to Y$ is irreducible if and only if there is a morphism $f':X \to Y'$ such that the induced morphism $(f,f')^T:X \to Y \oplus Y'$ is a minimal left almost split morphism.
\end{enumerate}
\end{Th}

For an almost split sequence 
$$\begin{xy}
  \xymatrix{ 
      0 \ar[r]  	& 	X \ar[r]^g 	&	Y \ar[r]^f     &   Z \ar[r]	&	0	    }
\end{xy}$$
let $Y = \oplus Y_i$ be a decomposition of $Y$ into indecomposable modules and let $\Pi_i$ denote the projection from $Y$ to $Y_i$. Furthermore, we set $f_i = f \Pi_i$ and $g_i =  \Pi_i g $ to be the induced morphisms from $X$ to $Y_i$ and from $Y_i$ to $Z$ respectively. Note that not only $f$ and $g$ are irreducible morphisms but also all $f_i$ and $g_i$. Moreover, for indecomposable modules $M$ and $N$ any irreducible morphism $h : M \to N$ can be extended to a minimal left almost split morphism and to a minimal right almost split morphism by \ref{connectionirralmostsplit}.

\begin{Th} \cite[V. Proposition 5.9]{ARS95}
Let
$$\begin{xy}
  \xymatrix{ 
      0 \ar[r]  	& 	X \ar[r]^g 	&	Y \ar[r]^f     &   Z \ar[r]	&	0	    }
\end{xy}$$
be an exact sequence. It is an almost split sequence if and only if $f$ and $g$ are both irreducible.
\end{Th}

If only either $f$ or $g$ is irreducible, then 
$$\begin{xy}
  \xymatrix{ 
      0 \ar[r]  	& 	X \ar[r]^g 	&	Y \ar[r]^f     &   Z \ar[r]	&	0	    }
\end{xy}$$
is not necessarily an almost split sequence. For example, consider a minimal right almost split morphism $f : Y \to Z$ such that $Y$ admits a non-trivial decomposition $Y_1 \oplus Y_2$ of $Y$ such that $f_1 : Y_1 \to Z$ is surjective. Then there is an exact sequence
$$\begin{xy}
  \xymatrix{ 
      0 \ar[r]  	& 	\Ker(f_1) \ar[r] 	&	Y_1 \ar[r]^{f_1}     &   Z \ar[r]	&	0	    }
\end{xy}$$
which is not an almost split sequence, but $f_1$ is clearly an irreducible morphism by \ref{connectionirralmostsplit}.
Alongside with irreducible morphisms we introduce the \textbf{radical} $\rad(X,Y)$, which is the set of all $f \in \Hom(X,Y)$ such that $gfh$ is not an isomorphism for any $h: Z \to X$ and $g: Y \to Z$ with Z indecomposable. \textbf{Powers of the radical} are defined inductively, that is an $f \in \Hom(X,Y)$ is in $\rad^n(X,Y)$ if there is a module $Z$ such that $f=gh$ for some $h \in \rad(X,Z)$ and $g \in \rad^{n-1}(Z,Y)$.
This gives rise to a sequence of submodules
$$\cdots \subset \rad^n(X,Y) \subset \rad^{n-1}(X,Y) \subset \cdots \subset \rad(X,Y) \subset \Hom(X,Y)$$
which motivates to define 
$$\rad^\infty(X,Y) = \bigcap_{n \in \N}\ \rad^n(X,Y).$$
For convenience, a morphism $f : X \to Y$ is called a \textbf{radical morphism} if $f \in \rad(X,Y)$. There is a crucial connection between the radical and irreducible morphisms.

\begin{Lemma} \cite[V. Proposition 7.3]{ARS95}
Let $f: X \to Y$ be a morphism between indecomposable modules $X, Y$. Then $f$ is irreducible if and only if $f \in \rad(X,Y) \backslash \rad^2(X,Y)$.
\end{Lemma}

The proof is also given in a more general version in Lemma \ref{irrradical}. So for $X$ and $Y$ indecomposable the elements of $\rad^n(X,Y) \backslash \rad^{n+1}(X,Y)$ are the non-zero morphisms that can be written as a sum of compositions of irreducible morphisms such that the shortest of these compositions has length $n$ and cannot be written as a composition of more than $n$ radical morphisms.
\bigskip

\begin{Def}
The \textbf{Auslander-Reiten quiver} $\Gamma_A$ of an algebra $A$ and its module category $A$-mod is a biquiver, i.e. a quiver consisting of vertices $\Gamma_0$ and two disjoint sets of arrows between them, 1-arrows $\Gamma_1$ and 2-arrows $\Gamma_2$. There is a vertex for every isomorphism class of indecomposable modules of the algebra. Given two indecomposable modules $X$ and $Y$ there is a 1-arrow from the vertex corresponding to the isomorphism class of $X$ to the vertex corresponding to the isomorphism class of $Y$ if there is an irreducible morphism from $X$ to $Y$. For each almost split sequence
$$\begin{xy}
  \xymatrix{ 
      0 \ar[r]  	& 	X \ar[r] 	&	Y \ar[r]     &   Z \ar[r]	&	0	    }
\end{xy}$$
there is a 2-arrow from the vertex corresponding to the isomorphism class of $Z$ to the vertex corresponding to the isomorphism class of $X$. In order to distinguish between different types of arrows, 2-arrows are drawn as dotted arrows.
\end{Def}

\begin{Def}\label{labeling}
Let $h : M \to N$ be an irreducible morphism between indecomposable modules and $f : X \to N$ and $g : M \to Y$ extensions to minimal right almost split and minimal left almost split morphisms respectively. We denote the number of $M$ in a sum decomposition of $X$ and the number of $N$ in a sum decomposition of $Y$ by $m$ and $n$ respectively. We then say the 1-arrow from $M$ to $N$ has \textbf{valuation} $(m,n)$. 
\end{Def}

Throughout the dissertation we usually do not distinguish between an indecomposable module, its isomorphism class or the corresponding vertices in $\Gamma_A$. There are squared brackets to mark injective and projective modules in $\Gamma_A$. For convenience, we simply say arrows for 1-arrows in $\Gamma_A$ and, if the valuation of an arrow is $(1,1)$, we say it is \textbf{trivially valuated}. Furthermore, for an indecomposable module $Y$ we call each module $X$ such that there is an arrow from $X$ to $Y$ in $\Gamma_A$ an \textbf{immediate predecessor of $\pmb{Y}$}. Dually, each module $Z$ such that there is an arrow from $Y$ to $Z$ in $\Gamma_A$ is called an \textbf{immediate successor of $\pmb{Y}$}.
\bigskip
\\If a subquiver $\Gamma$ of $\Gamma_A$ contains two modules $X$ and $Y$ such that there is a $1$-arrow or a $2$-arrow between them in $\Gamma_A$, then we define that $\Gamma$ automatically also contains this $1$-arrow or $2$-arrow respectively and the $1$-arrow in $\Gamma$ has the same valuation as in the whole Auslander-Reiten quiver. We call a non-empty subquiver $\Gamma$ of $\Gamma_A$ a \textbf{connected component} if for each $X$ in $\Gamma$ the modules in $\Gamma$ are precisely all modules $Y$ such that there is a walk between $X$ and $Y$ in $\Gamma_A$, that is $X$ and $Y$ are connected by arrows irrespective of their direction. Clearly, each indecomposable module $X$ is contained in the unique connected component of $\Gamma_A$ consisting of all modules that are connected to $X$ by a walk in the Auslander-Reiten quiver. We call this component the connected component of $X$. 
\bigskip
\\For a different point of view on arrow valuation we introduce some more notation. Recall that $X$ and $Y$ are indecomposable since they are modules in $\Gamma_A$. We denote the factor module $\rad(X,Y)/\rad^2(X,Y)$ by $\Irr(X,Y)$ and the division algebra $\End_A(X)/\rad (\End_A(X))$ by $T_X$, where in this case $\rad (\End_A(X))$ is the Jacobson radical of the algebra $\End_A(X)$. $\Irr(X,Y)$ then becomes a $T_Y$-$T_X^{op}$-bimodule and, if $X$ and $Y$ are indecomposable, a $T_Y$-$T_X^{op}$-vector space \cite[VII.1]{ARS95}.

\newpage

\begin{Th}\cite[VII. Proposition 1.3]{ARS95}
Let $X$ and $Y$ denote indecomposable $A$-modules and assume there is an irreducible morphism $f : X \to Y$. Let the valuation of the corresponding arrow be $(m,n)$. Then $m$ equals the dimension of $\Irr(X,Y)$ as a $T_X^{op}$-vector space while $n$ equals the dimension of $\Irr(X,Y)$ as a $T_Y$-vector space.
\end{Th}

\begin{Lemma}\cite[VII. Proposition 1.5]{ARS95}
Suppose there is an arrow in $\Gamma_A$ from $X$ to $Y$ with valuation $(m,n)$. If $Y$ is non-projective, then the valuation of the arrow from $\tau(Y)$ to $X$ is $(n,m)$.
\end{Lemma}

\begin{Cor}\cite[VII. Corollary 2.3]{ARS95}
Since $K$ is assumed to be algebraically closed, $m$ equals $n$ in all valuations of arrows in $\Gamma_A$.
\end{Cor}

Let $X$ and $Y$ be indecomposable and assume there is an irreducible morphism from $X$ to $Y$. If the corresponding arrow is valued $(n,n)$, we say that there are $\pmb{n}$ \textbf{arrows} from $X$ to $Y$ and we draw $n$ arrows from $X$ to $Y$ in $\Gamma_A$. In particular, if there is only one arrow from $X$ to $Y$, the arrow has trivial valuation. Consequently, if there are $n$ arrows from $X$ to $Y$ in $\Gamma_A$, then there are $n$ arrows from $\tau^k(X)$ to $\tau^k(Y)$ and from $\tau^k(Y)$ to $\tau^{k-1}(X)$ for all $k \in \Z$ such that these modules exists. Moreover, we say there are multiple arrows from $X$ to $Y$ if $n \geq 2$ and we say a connected component contains multiple arrows if and only if there are multiple arrows from $X$ to $Y$ for any modules $X$ and $Y$ in that connected component. In order to show that the module categories in the following examples are finite, we state another well-known theorem.

\begin{Th} \cite[VI. Theorem 1.4]{ARS95} \label{finitecomponentAmod}
If there is a connected component $\Gamma$ in $\Gamma_A$ such that the length of all modules in $\Gamma$ is bounded, then it is the only connected component and $A$ is representation finite.
\end{Th}

\begin{Ex}
Let $A$ be the path algebra given by the quiver
$$\begin{xy}
  \xymatrix{ 
      e_1 \ar[r]  	& 	 e_2 \ar[r] 	&	 e_3	    }
\end{xy}$$
\end{Ex}

We compute the Auslander-Reiten quiver of $A$. The indecomposable modules are $S_3 = P_3, S_2, S_1 = I_1, P_1 = I_3, P_2$ and $I_2$. Since the canonical morphisms $\rad(P) \to P$ for projective modules and $I \to I/\soc(I)$ for injective modules are irreducible, we already have given irreducible morphisms $P_3 \to P_2$, $P_2 \to P_1$, $I_3 \to I_2$ and $I_2 \to I_1$. Considering $P_2$ we notice that it is the only indecomposable module that is domain of a non-zero non-isomorphism to $S_2$. Consequently, we have the following Auslander-Reiten quiver.

$$\begin{xy}
\xymatrixrowsep{0.5in}
\xymatrixcolsep{0.5in}

  \xymatrix@!0{ 
&	& [P_1] \ar[dr] 	&	& 	\\
&[P_2 \ar[dr] \ar[ur]	&	&I_2] \ar[dr] \ar@{.>}[ll] 	&	\\
[P_3 \ar[ur]	& 	& S_2 \ar[ur] \ar@{.>}[ll]&	& I_1] \ar@{.>}[ll]	
 }
\end{xy}$$
\newpage
A slightly more complicated example shows that in some cases the Auslander-Reiten quiver contains an isomorphism class of an indecomposable module more than once. In these cases we consider the corresponding vertices to be different. This differentiation is important in Chapter 4. Nevertheless, we say a component $\Gamma$ of the Auslander-Reiten quiver is finite if its vertices correspond to only finitely many isomorphism classes of indecomposable modules. So when we draw the Auslander-Reiten quiver of an algebra of finite representation type, we make sure that every different $1$-arrow and $2$-arrow occurs at least once, i.e. there is an illustration for every almost split sequence and every irreducible morphism.

\begin{Ex}\label{exampletwistedquiver}
Let A be the path algebra
$$\begin{xy}
\xymatrix{
  e_1 \  \ar[r]^\alpha  & \ e_2 \ \ar@(ur,dr)^\beta  \\
}
\end{xy}$$
with the relation $\beta^2 = 0$. Then some isomorphism classes of indecomposable modules occur more than once in the Auslander-Reiten quiver of $A$-mod.
\end{Ex}

Since the shape of the Auslander-Reiten quiver is crucial for this example, we calculate it in more detail then usual. First of all we name all occurring modules and write down their Jordan-H\"{o}lder composition series. We have simple modules $S_1$ and $S_2$, where $S_1$ is also the injective module $I_1$. Besides that there are projective and injective modules
 $$\begin{xy}
\xymatrixrowsep{0in}
\xymatrixcolsep{-0.1in}
  \xymatrix{ 
	    &S_1    	&		&		&				&	&	&S_1					\\	
P_1=\ \ &S_2   	&\ \ \ \ \ \ \ \ \ \ P_2= \ \	&S_2	&\ \ \ \ \ \ \ \ \ \ I_2= \ \	&S_1	&	&S_2		\\	
	    &S_2    	&		&	S_2	&				&	&S_2	&
   }
\end{xy}$$
and two other modules
 $$\begin{xy}
\xymatrixrowsep{0in}
\xymatrixcolsep{-0.1in}
  \xymatrix{ 
	    &    	&					&	&	&					\\	
X=\ \ &S_1   	&\ \ \ \ \ \ \ \ \ \ Y= \ \	&S_1	&	&S_2		\\	
	    &S_2    	&		&	&S_2.	&
   }
\end{xy}$$
We start by computing $\tau(I_2)$, which is $P_2$. Clearly, $P_1$ is a direct summand of the middle term of the related almost split sequence, because the embedding $\rad(P_1) \to P_1$ is always irreducible. Consequently, $Y$ must be the other summand in order to have a short exact sequence. The projection $I_2 \to I_2/\soc(I_2) = S_1 \oplus X$ is also irreducible, so the two direct summands must be $\tau^{-1}(P_1)$ and $\tau^{-1}(Y)$ respectively. But there is no morphism from $P_1$ to $X$ that factors over a simple module, so we have $\tau^{-1}(P_1) = S_1$ by the Jordan-H\"{o}lder multiplicities of $P_1, I_2$ and $X$. It is easy to verify that $\tau(Y) = S_2$. On the other hand, by Jordan-H\"{o}lder multiplicities $S_2$ is also a direct summand of the almost split sequence from $Y$ to $X$. These calculations give rise to the Auslander-Reiten quiver, which is shown on the next page.

$$\begin{xy}
\xymatrixrowsep{0.5in}
\xymatrixcolsep{0.5in}

  \xymatrix@!0{  
&	& [P_1 \ar[dr] 	&	& I_1] \ar@{.>}[ll]	&	&	&	\\
&[P_2 \ar[dr] \ar[ur]	&	&I_2] \ar[dr] \ar[ur] \ar@{.>}[ll] 	&	&	&	&	\\
S_2\ar[dr] \ar[ur]	& 	&Y \ar[dr] \ar[ur] \ar@{.>}[ll]&	&X \ar[dr] \ar@{.>}[ll]&	&S_2 \ar[dr] \ar@{.>}[ll] &	\\ 
&X \ar[ur] &	&S_2 \ar[dr] \ar[ur] \ar@{.>}[ll] & 	&Y \ar[dr] \ar[ur] \ar@{.>}[ll]&	&X \ar@{.>}[ll]	\\ 
&	&	&	&[P_2 \ar[dr] \ar[ur]	&	&I_2] \ar[dr] \ar[ur] \ar@{.>}[ll] 	&	\\
&	&	&	&	& [P_1 \ar[ur] 	&	& I_1] \ar@{.>}[ll]	
 }
\end{xy}$$
In order to show the necessity of certain assumptions in the following chapters, we compute the Auslander-Reiten quiver of another algebra, which we refer to as the \textbf{standard example}. Its significance for this dissertation is due to the following.
\bigskip
\\If $A$ is a quasi-hereditary algebra, then there is a natural ordering of the indecomposable projective modules $P_1, \ldots, P_n$. We then define $\Delta_i$ to be the largest factor module of $P_i$ such that for all simple blocks $S_j = P_j/(\rad(P_j)$ in the Jordan-H\"{o}lder series of $\Delta_i$ we have $j \leq i$. This gives a collection of modules $\Delta_1, \ldots, \Delta_n$, which are called the \textbf{standard modules} of $A$. Then a module $M$ in $A$-mod such that there is a filtration $0 = M_0 \subset M_1  \subset \cdots \subset M_k = M$ with $M_j/M_{j-1} \cong \Delta_{i_j}$ for all $j = 1, \ldots, k$ and some $1 \leq i_j \leq k$ is called \textbf{standard-filtered} and $\F(\Delta)$ consists of all modules with this property. It has been shown that $\F(\Delta)$ is a functorially finite resolving subcategory \cite{R91}. Dually, we define the \textbf{costandard module} $\nabla_i$ as the largest submodule of $I_i$ such that for all simple blocks $S_j$ in $\nabla_i$ we have $j \leq i$.

\begin{Ex}\label{standardexample}
Let $A$ be the quasi-hereditary path algebra given by the quiver
$$\begin{xy}
  \xymatrix{ 
e_1 \ar@<1ex>[r]^{f_1}  &e_2	\ar@<1ex>[l]^{g_1} \ar@<1ex>[r]^{f_2}  &e_3	\ar@<1ex>[l]^{g_2} \ar@<1ex>[r]^{f_3}  &	e_4	\ar@<1ex>[l]^{g_3}  }
\end{xy}
$$ 
with relations $f_1g_1 = g_2f_2$, $f_2g_2 = g_3f_3$ and $f_3g_3 = f_2f_1 = f_3f_2 = g_1g_2 = g_2g_3 = 0$.
\end{Ex}

The structure of $_AA$, which is $A$ considered as a left $A$-module, is given by
 $$\begin{xy}
\xymatrixrowsep{0in}
\xymatrixcolsep{-0.1in}
  \xymatrix{ 
    	S_1    	&		&	&S_2	&	&		&	&S_3	&	&		&S_4 	\\	
	S_2    	&\ \oplus \	&S_1	&	&S_3	&\ \oplus \	&S_2	&	&S_4	&\ \oplus \	&\ S_3.	\\	
	S_1    	&		&	&S_2	&	&		&	&S_3	&	&		&		
    }
\end{xy}$$
For convenience, we give all modules proper names. Later Theorem \ref{finitecomponentAmod} is used to see that the modules mentioned are in fact all indecomposable modules up to isomorphism. $S_1, S_2, S_3$ and $S_4$ naturally denote the simple modules. The projective and injective modules are given by
 
 $$\begin{xy}
\xymatrixrowsep{0in}
\xymatrixcolsep{-0.1in}
  \xymatrix{ 
	    &S_1    	&					&	&S_2	&	\\	
P_1=I_1=\ \ &S_2   	&\ \ \ \ \ \ \ \ \ \ P_2=I_2= \ \	&S_1	&	&S_3	\\	
	    &S_1    	&					&	&S_2	&	
   }
\end{xy}$$
 $$\begin{xy}
\xymatrixrowsep{0in}
\xymatrixcolsep{-0.1in}
  \xymatrix{ 
		&	&S_3	&	&				&S_4	&				&S_3	\\	
P_3=I_3= \ \	&S_2	&	&S_4	& \ \ \ \ \ \ \ \ \ \ P_4= \ \	&S_3	&\ \ \ \ \ \ \ \ \ \ I_4= \ \	&S_4.	\\	
	       	&	&S_3	&	&				&	&				&
    }
\end{xy}$$
The standard modules and costandard modules are
 $$\begin{xy}
\xymatrixrowsep{0in}
\xymatrixcolsep{-0.1in}
  \xymatrix{ 
			&	&					&S_2	&					&S_3		\\	
\Delta_1= \nabla_1= \ \	&S_1	&\ \ \ \ \ \ \ \ \ \ \Delta_2 = \ \	&S_1	& \ \ \ \ \ \ \ \ \ \ \Delta_3= \ \	&S_2		
    }
\end{xy}$$
 $$\begin{xy}
\xymatrixrowsep{0in}
\xymatrixcolsep{-0.1in}
  \xymatrix{ 
		&S_1	&					&S_2	\\	
\nabla_2= \ \	&S_2	&\ \ \ \ \ \ \ \ \ \ \nabla_3 = \ \	&S_3,		
    }
\end{xy}$$
where $\Delta_4 = P_4$ and $\nabla_4 = I_4$. Finally we name all other modules.
 $$\begin{xy}
\xymatrixrowsep{0in}
\xymatrixcolsep{-0.1in}
  \xymatrix{ 
	&	&S_2	&	&S_4	&				&S_1	&	&S_3	&	\\	
X = \ \	&S_1	&	&S_3	&	&\ \ \ \ \ \ \ \ \ \ \ Y = \ \	&	&S_2	&	&S_4		
    }
\end{xy}$$
 $$\begin{xy}
\xymatrixrowsep{0in}
\xymatrixcolsep{-0.1in}
  \xymatrix{ 
		&	&S_2	&	&					&	&S_3	&	\\	
M_2 = \ \	&S_1	&	&S_3	&\ \ \ \ \ \ \ \ \ \ \ M_3 = \ \	&S_2	&	&S_4		
    }
\end{xy}$$
 $$\begin{xy}
\xymatrixrowsep{0in}
\xymatrixcolsep{-0.1in}
  \xymatrix{ 
	  &S_1	&	&S_3	&					&S_2	&	&S_4		\\	
N_2 = \ \ &	&S_2	&	&\ \ \ \ \ \ \ \ \ \ \ N_3 = \ \	&	&S_3	&		
    }
\end{xy}$$
We start the computation of the Auslander-Reiten quiver of $A$-mod by observing that $P_1$ and $P_3$ are also injective. It follows that there are almost split sequences
$$\begin{xy}
  \xymatrix{ 
      0 \ar[r]  	& 	\Delta_2 \ar[r] 	&	P_1 \oplus S_2 \ar[r]     &   \nabla_2 \ar[r]	&	0	    }
\end{xy}$$
$$\begin{xy}
  \xymatrix{ 
      0 \ar[r]  	& 	N_3 \ar[r] 	&	S_2 \oplus S_4 \oplus P_3 \ar[r]     &   M_3 \ar[r]	&	0	    }
\end{xy}$$
by \cite[V. Proposition 5.5]{ARS95}. We connect these almost split sequences at the module $S_2$ and obtain the whole Auslander-Reiten quiver by knitting to either left or right.

$$\begin{xy}
\xymatrixrowsep{0.5in}
\xymatrixcolsep{0.5in}

  \xymatrix@!0{  
&	&		&	&					& [P_1] \ar[dr] 	&	&	&	&	\\
&	&[P_4 \ar[dr]	&	&\Delta_2 \ar[dr] \ar[ur] \ar@{.>}[ll]	&		&\nabla_2 \ar[dr] \ar@{.>}[ll]	&	&I_4] \ar[dr] \ar@{.>}[ll] & \\
&S_3\ar[dr] \ar[ur]	& &X \ar[dr] \ar[ur] \ar@{.>}[ll]&	&S_2 \ar[dr] \ar[ur] \ar@{.>}[ll]&	&Y \ar[dr] \ar[ur] \ar@{.>}[ll] &	&S_3  \ar@{.>}[ll]\\ 
N_2 \ar[ddr] \ar[dr] \ar[ur] &	&M_2 \ar[dr] \ar[ur] \ar@{.>}[ll]&	&N_3 \ar[ddr] \ar[dr] \ar[ur] \ar@{.>}[ll]& &M_3 \ar[dr] \ar[ur] \ar@{.>}[ll]&	&N_2 \ar[dr] \ar[ur] \ar@{.>}[ll] & \\
&S_1 \ar[ur]&  &\nabla_3 \ar[ur] \ar@{.>}[ll]&	&S_4 \ar[ur] \ar@{.>}[ll]&	&\Delta_3 \ar[ur] \ar@{.>}[ll] &	&S_1  \ar@{.>}[ll]\\ 
&[P_2] \ar[uur]	&	&	&	& [P_3] \ar[uur] 	&	&	&	&	
 }
\end{xy}$$

\chapter{Irreducible morphisms in functorially finite resolving subcategories}

\section{Functorially finite resolving subcategories}

As mentioned in the introduction, almost split sequences do not only exist in $A$-mod but also occur in subcategories of $A$-mod. Let $\chi$ be a subcategory of $A$-mod. By \textbf{subcategory} we always mean a full subcategory closed under isomorphisms, direct sums and summands. 

\begin{Def}
\mbox{}
\begin{enumerate}[(a)]
 \item For modules $X$ and $Y$ in $\chi$ a morphism $f : X \to Y$ is called \textbf{right almost split in} $\pmb{\chi}$ if it is not a split epimorphism and any morphism $Z \to Y$ with $Z$ in $\chi$ that is not a split epimorphism factors through $f$.
$$\begin{xy}
  \xymatrix{ 
    Z \ar@{.>}[d] \ar[dr]  	&	 \\
    X  \ar[r]^{f}    	&	Y	    }
\end{xy}$$ 
\item We dually define a morphism $g : Y \to X$ to be \textbf{left almost split in} $\pmb{\chi}$ if it is not a split monomorphism and any morphism $Y \to Z$ with $Z$ in $\chi$ that is not a split monomorphism factors through $g$.
$$\begin{xy}
  \xymatrix{ 
    Y  \ar[r]^g \ar[dr]  	& X \ar@{.>}[d]	 \\
    &	Z	    }
\end{xy}$$ 
\end{enumerate}

\end{Def}

Similarly to $A$-mod, we say that a right almost split morphism $f: X \to Y$ in $\chi$ is a right almost split morphism for $Y$ in $\chi$. Dually, we call a left almost split morphism $g : Y \to X$ in $\chi$ a left almost split morphism for $Y$ in $\chi$ .

\begin{Def}
For modules $X, Y, Z$ in $\chi$ a short exact sequence
$$\begin{xy}
  \xymatrix{ 
      0 \ar[r]  	& 	X \ar[r]^g 	&	Y \ar[r]^f     &   Z \ar[r]	&	0	    }
\end{xy}$$
is called an \textbf{almost split sequence in }$\pmb{\chi}$ if $g$ is minimal left almost split in $\chi$ and $f$ is minimal right almost split in $\chi$.
\end{Def}

If there are up to isomorphism unique almost split sequences in $\chi$ ending in every indecomposable non-projective module in $\chi$, then we clearly also have an Auslander-Reiten quiver for $\chi$. One type of subcategory that we are particularly interested in with this property are functorially finite resolving subcategories. The motivating example for these categories is the category of all modules with a standard filtration of a quasi-hereditary algebra or in general categories obtained by a generalized cotilting module. For an arbitrary indecomposable non-projective module $X$ in $\chi$, the module $\tau(X)$ is not necessarily in $\chi$. Therefore, we introduce approximations, which give rise to an operation on the indecomposable modules in $\chi$ whose properties in $\chi$ are similar to those of the Auslander-Reiten translation in $A$-mod. We start by recalling the definitions.

\begin{Def}\mbox{}
\begin{enumerate}[(a)]
 \item
A \textbf{right $\pmb{\chi}$-approximation} of a module $Y$ is a morphism $f_Y: X_Y \to Y$, where $X_Y$ is in $\chi$, such that for all $Z$ in $\chi$ every morphism $g: Z \to Y$ factors through $f_Y$.
$$\begin{xy}
  \xymatrix{ 
    Z \ar[dr]^g \ar@{.>}[d]  	&	 \\
    X_Y  \ar[r]^{f_Y}  	&	Y	    }
\end{xy}$$
 \item Dually, a \textbf{left $\pmb{\chi}$-approximation} of a module $Y$ is a morphism $f^Y: Y \to X^Y$, where $X^Y$ is in $\chi$, such that for all $Z$ in $\chi$ every morphism $g : Y \to Z$ factors through $f^Y$.
$$\begin{xy}
  \xymatrix{ 
        	& Z	 \\
    Y  \ar[r]^{f^Y}\ar[ur]^g  	&	X^Y \ar@{.>}[u].	    }
\end{xy}$$
\end{enumerate}

\end{Def}

\begin{Def}
A subcategory $\chi$ of $A$-mod is called 
\begin{enumerate}[(a)]
 \item \textbf{contravariantly finite} if right $\chi$-approximations exist for all modules in $A$-mod,
 \item \textbf{covariantly finite} if left $\chi$-approximations exist for all modules in $A$-mod,
 \item \textbf{functorially finite} if it is both contravariantly and covariantly finite.
\end{enumerate}
\end{Def}

If a right $\chi$-approximation is also a right minimal morphism, we call it a \textbf{minimal right $\pmb{\chi}$-approximation}. Dually, if a left $\chi$-approximation is also a left minimal morphism, we call it a \textbf{minimal left $\pmb{\chi}$-approximation}. Unless otherwise specified we assume all right and left $\chi$-approximations to be right minimal and left minimal respectively. By definition the domain $X_Y$ of a minimal right $\chi$-approximation and the codomain of a minimal left $\chi$-approximation are unique up to isomorphism.
\begin{Def}
Let $Y$ be an arbitrary module and let $X$ be a module in $\chi$. 
\begin{enumerate}[(a)]
 \item A $\pmb{\chi}$\textbf{-section} is a morphism $f: X \to Y$ such that for all modules $Z$ in $\chi$ and morphisms $h : X \to Z$ and $g:Z \to Y$ with $f = gh$, $h$ must be a split monomorphism. 
\item A $\pmb{\chi}$\textbf{-contraction} is a morphism $f: Y \to X$ such that for all modules $Z$ in $\chi$ and morphisms $h : Z \to X$ and $g:Y \to Z$ with $f = hg$, $h$ must be a split epimorphism. 
\end{enumerate}
\end{Def}

The relation between $\chi$-sections and minimal right $\chi$-approximations corresponds to the relation between irreducible and minimal right almost split morphisms. This means a morphism $f: X \to Y$ is a $\chi$-section if and only if there are a module $X'$ in $\chi$ and a morphism $f': X' \to Y$ such that $(f,f') : X \oplus X' \to Y$ is a minimal right $\chi$-approximation \cite[Proposition 2.2]{KP03}. We call the morphism $(f,f')$ an \textbf{extension of $\pmb{f}$} to a minimal right $\chi$-approximation.

\begin{Def}
\mbox{}
\begin{enumerate}[(a)]
 \item A module $I$ in $\chi$ is called \textbf{Ext-injective in }$\pmb{\chi}$ if $\Ext_A^i(M,I) = 0$ for all $M$ in $\chi$ and all $i \geq 1$. 
 \item Dually, a module $P$ in $\chi$ is called \textbf{Ext-projective in }$\pmb{\chi}$ if $\Ext_A^i(P,M) = 0$ for all $M$ in $\chi$ and all $i \geq 1$.
\end{enumerate}
\end{Def}

If we deal with only one subcategory, we simply say $I$ is Ext-injective or $P$ is Ext-projective.


\begin{Def}
A subcategory $\chi$ of $A$-mod is called \textbf{closed under extensions} if for all $X,Z$ in $\chi$ the existence of a short exact sequence
$$\begin{xy}
  \xymatrix{ 
        0 \ar[r]   &	X \ar[r]		&	Y  \ar[r]		&	Z \ar[r]		&	0	 }
\end{xy}$$
implies that $Y$ is in $\chi$.
\end{Def}

\begin{Def}
Let $\Omega$ be a subcategory of $A$-mod. We call $\Omega$ \textbf{resolving} if 
\begin{enumerate}[(a)]
\item it is closed under extensions,
\item $_AA$ is in $\Omega$ and
\item it is closed under kernels of epimorphisms.
\end{enumerate}
\end{Def}


\begin{Def}
An $A$-module $T$ is called a \textbf{generalized cotilting module} if
\begin{enumerate}[(a)]
\item its injective dimension is finite, that is there is an exact sequence
$$\begin{xy}
  \xymatrix{ 
        0 \ar[r]   &	T \ar[r]		&	I_0 \ar[r]	&	I_1 \ar[r]	&	\cdots \ar[r]	&	I_n \ar[r] &	0,	 }
\end{xy}$$
where all $I_j$ are injective $A$-modules,
\item $\Ext_A^i(T,T) = 0$ for all $i \geq 1$ and
\item there is a finite resolution
$$\begin{xy}
  \xymatrix{ 
        0 \ar[r]   &	T_n \ar[r]		&	\cdots  \ar[r]	&	T_1 \ar[r]	&	T_0 \ar[r]&	I \ar[r] &	0	 }
\end{xy}$$
for each injective module $I$ where $T_i$ in $\add(T)$, i.e. each $T_i$ is a sum of direct summands of $T$.
\end{enumerate}

\end{Def}

It has been shown by Auslander and Reiten in \cite{AR91} that for a generalized cotilting module $T$ the category $^\bot T = \{\Ext_A^i(-,T) = 0\}$ is contravariantly finite and resolving. Moreover, if $A$ has finite global dimension, then there is a bijection between contravariantly finite resolving subcategories $^\bot T$ and isomorphism classes of basic cotilting modules $T$ \cite[Corollary 2.4]{R07}, where basic means that no indecomposable summand of $T$ occurs more than once in $T$.
\bigskip
\\Let $\chi$ denote a functorially finite subcategory of $A$-mod that is closed under extensions. Let $Y$ be an indecomposable module in $\chi$ that is not Ext-projective, we then can decompose the domain of a right $\chi$-approximation $X_{\tau(Y)} = M \oplus I$ where $M \neq 0$ is indecomposable and not Ext-injective \cite[Theorem 2.3]{K97} whereas $I$ is Ext-injective \cite[Theorem 5.3]{KP03}. The module $M$ is then denoted by $\tau_{\chi}(Y)$. For left $\chi$-approximations the dual statement holds, i.e. if $Y$ is not Ext-injective, we can decompose the codomain of a left $\chi$-approximation $X^{\tau^{-1}(Y)} = N \oplus P$ where $N \neq 0$ is indecomposable and not Ext-projective while $P$ is Ext-projective. Analogously, the module $N$ is denoted by $\tau_{\chi}^{-1}(Y)$. We call $\tau_{\chi}$ the \textbf{relative Auslander-Reiten translation of} $\pmb{\chi}$ and $\tau_{\chi}(Y)$ the relative Auslander-Reiten translate of a module $Y$ in $\chi$. It follows from these results that almost split sequences 
in $\chi$ have similar properties as they do in $A$-mod.

\begin{Lemma}\cite[Proposition 2.7]{K97}\label{almostsplitsequenceinchi}
\mbox{}
\begin{enumerate}[(a)]
 \item Let
$$\begin{xy}
  \xymatrix{ 
        0 \ar[r]   &	\Ker(f) \ar[r]^(.6)g		&	Y  \ar[r]^f		&	Z \ar[r]		&	0	 }
\end{xy}$$
be an exact sequence in $\chi$. If $f$ is a minimal right almost split morphism in $\chi$, then $\Ker(f) \cong \t(Z)$ and $g$ is a minimal left almost split morphism in $\chi$.
 \item Let
$$\begin{xy}
  \xymatrix{ 
        0 \ar[r]   &	X \ar[r]^g		&	Y  \ar[r]^(.35)f		&	\Coker(g) \ar[r]		&	0	 }
\end{xy}$$
be an exact sequence in $\chi$. If $g$ is a minimal left almost split morphism in $\chi$, then $\Coker(g) \cong \t^{-1}(X)$ and $f$ is a minimal right almost split morphism in $\chi$.
\end{enumerate}
\end{Lemma}

Throughout the dissertation $\chi$ denotes a functorially finite subcategory of $A$-mod that is closed under extensions while $\Omega =$ $^\bot T$ denotes a functorially finite resolving subcategory of $\chi$ that is generated by a generalized cotilting module $T$. In particular, $\chi$ contains all projective modules since they are contained in $\Omega$. Most of the upcoming results still hold if we just demand $\chi$ and $\Omega$ to be contravariantly finite, but all examples that are discussed are functorially finite subcategories. Moreover, we want to deal with Auslander-Reiten quivers in which we are able to knit both to the left and the right. 
\bigskip
\\By definition resolving subcategories contain all projective modules and these are precisely the Ext-projective modules. Hence if $\chi$ contains an Ext-projective module $P$ which is not projective, then $P$ is not in $\Omega$. An arbitrary resolving subcategory usually does not contain all injective modules. In the case of $^\bot T$ the Ext-injective modules are all modules in $\add(T)$ by construction. In particular, in $\Omega$ the number of indecomposable non-isomorphic projective modules coincides with the number of indecomposable non-isomorphic Ext-injective modules by \cite[Proposition 1.4]{R07}. Note that $A$-mod itself is clearly a functorially finite resolving subcategory.

\begin{Lemma}\label{approxdiagram}
Let
$$\begin{xy}
  \xymatrix{ 
        0 \ar[r]   &	N \ar[r]^f		&	M  \ar[r]^g		&	Z \ar[r]		&	0	 }
\end{xy}$$
be a short exact sequence and let $Z$ be an indecomposable module in a functorially finite resolving subcategory $\Omega$. Then there is a commutative diagram 
 $$\begin{xy}
  \xymatrix{ 
    0 \ar[r]   &   	X_N  \ar[r]^i \ar[d]_{f_N} &	X_M \ar[r]^{gf_M} \ar[d]_{f_M} 	&    	Z \ar[r]\ar@{=}[d]  	&	0	\\  
    0 \ar[r]   &	N \ar[r]^f		&	M  \ar[r]^g		&	Z \ar[r]		&	0	 }
\end{xy}$$
such that the top row splits if and only if the bottom row splits.
\end{Lemma}
\proof
\\Let $f_M$ be a right $\Omega$-approximation of $M$, then the composition $gf_M$ is surjective as a composition of surjective maps. We denote its kernel by $X_N$ and prove that the natural morphism from $X_N$ to $N$ is a right $\Omega$-approximation. As $0 = gf_Mi$, we know that $\Im(f_Mi) \subset \Ker(g) \cong N$ embeds naturally into $N$. We call this morphism $f_N$ and show that it is a right $\Omega$-approximation of $N$. Let $X$ be a module in $\Omega$ and $h : X \to N$ a morphism, then the composition $fh$ factors through $f_M$, so we have $fh = f_Mh'$ for some $h' : X \to X_M$. Moreover, $gf_Mh' = gfh = 0$, so $h'$ maps to the Kernel of $gf_M$, which equals $X_N$. In other words, there is a $h'': X \to X_N$ such that $h' = ih''$ and, in particular, $fh = f_Mih'' = ff_Nh''$. As $f$ is injective, this implies $h = f_Nh''$, which proves that $f_N$ is a right $\Omega$-approximation.
\bigskip
\\Suppose now that the top row splits, then $Z$ is a direct summand of $X_M$ and $gf_M$ is a split epimorphism. In particular, $g$ is a split epimorphism and the bottom row also splits. Lastly, we assume that the bottom row splits, then $Z$ is a direct summand of $M$ and a minimal right $\Omega$-approximation is given by the identity on $Z$. Hence the top row also splits, which completes the proof.
\qed
\bigskip
\\Note that if the bottom row is an almost split sequence in $A$-mod, then the diagram can be reduced to
 $$\begin{xy}
  \xymatrix{ 
    0 \ar[r]   &   	X  \ar[r]^{i|_X} \ar[d]_{f_N|_X} 	&	Y \ar[r]^{gf_M|_Y} \ar[d]_{f_M|_Y} 	&    	Z \ar[r]\ar@{=}[d]  	&	0	\\  
    0 \ar[r]   &	N \ar[r]_f				&	M  \ar[r]_g				&	Z \ar[r]		&	0	 }
\end{xy}$$
where $X = \tau_\Omega(Z)$ is the unique direct summand of $X_N = X \oplus I$ that is not Ext-injective and $Y \cong X_M/i(I)$.
\bigskip
\\We now recall a result on contravariantly finite subcategories that we apply to $\Omega$ frequently throughout the dissertation.

\begin{Lemma}\cite[Proposition 3.3]{AR91} \label{Extiso}
\begin{enumerate}[(a)]
\item Let $f_Y : X_Y \to Y$ be a minimal right $\Omega$-approximation of $Y$. Then $f_Y$ induces an isomorphism
$$\Ext_A^i(X,X_Y) \cong \Ext_A^i(X,Y)$$
for all $i \geq 1$ and $X$ in $\Omega$.
\item Dually, if $f^X : X \to X^X$ is a minimal left $\Omega$-approximation of $X$, then $f^X$ induces an isomorphism
$$\Ext_A^i(X^X,Y) \cong \Ext_A^i(X,Y)$$
for all $i \geq 1$ and $Y$ in $\Omega$.
\end{enumerate}
\end{Lemma}

Clearly, the short exact sequences corresponding to the isomorphism in $(a)$ for $i = 1$ are given by the diagram of Lemma \ref{approxdiagram}. In particular, the middle term of the sequence corresponding to $\Ext_A^1(X,X_Y)$ is the domain of a minimal right $\Omega$-approximation of the middle term of the sequence corresponding to $\Ext_A^1(X,Y)$. More generally, the following holds for a functorially finite subcategory that is closed under extensions.

\begin{Lemma}\cite[Lemma 2.1]{K97}\label{functormono}
Let $N$ be a module and let $f_N : X_N \to N$ and $f^N: N \to X^N$ denote minimal right and left $\chi$-approximations respectively. Then we have the following.
\begin{enumerate}[(a)]
\item $\Ext_A^1(\ \ ,f_N)|_\chi : \Ext_A^1(\ \ , X_N)|_\chi \to \Ext_A^1(\ \ , N)|_\chi$ is a monomorphism of functors.
\item $\Ext_A^1(f^N,\ \ )|_\chi : \Ext_A^1(X^N,\ \ )|_\chi \to \Ext_A^1(N, \ \ )|_\chi$ is a monomorphism of functors.
\end{enumerate}
\end{Lemma}

\begin{Lemma}\cite[Lemma 4.3]{KP03}\label{functormonoExtinj}
Let $X$ be a module in $\chi$, $N$ an arbitrary module and $f : X \to N$ a morphism that induces a monomorphism of contravariant functors
$$\Ext_A^1(\ \ ,f)|_\chi : \Ext_A^1(\ \ , X)|_\chi \to \Ext_A^1(\ \ , N)|_\chi.$$ 
\begin{enumerate}[(a)]

 \item If $\Ext_A^1(\tau^{-1}(N),X) = 0$, then $X$ is Ext-injective in $\chi$.
 \item If $\tau^{-1}(N)$ is in $\chi$, then $\Ext_A^1(\tau^{-1}(N),X) = 0$ if and only if $X$ is Ext-injective in $\chi$.
\end{enumerate}
Dually, let $f : N \to X$ be a morphism that induces a monomorphism of covariant functors
$$\Ext_A^1(f,\ \ )|_\chi : \Ext_A^1(X,\ \ )|_\chi \to \Ext_A^1(N, \ \ )|_\chi.$$ 
\begin{enumerate}[(a)]
\setcounter{enumi}{2}
 \item If $\Ext_A^1(X, \tau(N)) = 0$, then $X$ is Ext-projective in $\chi$.
 \item If $\tau(N)$ is in $\chi$, then $\Ext_A^1(X, \tau(N)) = 0$ if and only if $X$ is Ext-projective in $\chi$.
\end{enumerate}

\end{Lemma}

In order to show that the relative Auslander-Reiten translation in $\Omega$ can be defined analogously via the relative Auslander-Reiten translation in $\chi$, we prove the following statements.

\newpage
\begin{Th}\label{Omegatranslateinchi}

Let $Z$ be a module in $\Omega$.
\begin{enumerate}[(a)]
\item If $Z$ is not projective, then the domain of a minimal right $\Omega$-approximation $f_{\t(Z)} : X_{\t(Z)} \to \t(Z)$ can be decomposed $X_{\t(Z)} = I \oplus \tau_\Omega(Z)$ where $I$ is Ext-injective in $\Omega$.
 \item If $Z$ is not Ext-injective in $\Omega$, then the codomain of a minimal left $\Omega$-approximation $f^{\t^{-1}(Z)} : X^{\t^{-1}(Z)} \to \t^{-1}(Z)$ can be decomposed $X^{\t^{-1}(Z)} = P \oplus \tau_\Omega^{-1}(Z)$ where $P$ is projective.
\end{enumerate}

\end{Th}

\proof
\\Firstly, note that $Z$ is not projective implies it is not Ext-projective in $\Omega$ and hence cannot be Ext-projective in $\chi$, so $\t(Z)$ exists. We prove statement $(a)$ and statement $(b)$ follows by duality. 
\bigskip
\\Let $f_{\tau(Z)} : X_{\tau(Z)} \to \tau(Z)$ denote a minimal right $\chi$-approximation of $\tau(Z)$ and consider a minimal right $\Omega$-approximation $g_{\tau(Z)} : Y_{\tau(Z)} \to \tau(Z)$, which is also an $\Omega$-section. Then there is a factorization $g_{\tau(Z)} = f_{\tau(Z)} g'$ for some morphism $g' : Y_{\tau(Z)} \to X_{\tau(Z)}$ as $Y_{\tau(Z)}$ is in $\chi$. On the other hand, there is also a factorization $g' = f_{\t(Z)}h$ for some morphism $h : Y_{\tau(Z)} \to X_{X_{\tau(Z)}}$, where $f_{X_{\tau(Z)}}: X_{X_{\tau(Z)}} \to X_{\tau(Z)}$ is a minimal right $\Omega$-approximation. Consequently, as $g_{\tau(Z)} = f_{\tau(Z)} f_{\t(Z)}h$ is an $\Omega$-section, $h$ must be a split monomorphism and $Y_{\tau(Z)}$ is a direct summand of $X_{X_{\tau(Z)}}$.
$$\begin{xy}
  \xymatrixrowsep{0.6in}
\xymatrixcolsep{1.2in}

  \xymatrix@!0{ 
	&	Y_{\tau(Z)} \ar[dr]^{g_{\tau(Z)}} \ar@{.>}[d]^{g'} \ar@{.>}[ld]_h	&    		\\  
   	X_{X_{\tau(Z)}}  \ar[r]^(.6){f_{X_{\tau(Z)}}}			&	X_{\tau(Z)}  \ar[r]^{f_{\tau(Z)}}			&	\tau(Z) }
\end{xy}$$
Hence we get a chain of monomorphisms
 $$\begin{xy}
  \xymatrix{ 
\Ext_A^1(Z, Y_{\tau(Z)}) \ar@{^{(}->}[r]	&	\Ext_A^1(Z, X_{X_{\tau(Z)}})	\ar@{^{(}->}[r] &  \Ext_A^1(Z, X_{\tau(Z)}) \ar@{^{(}->}[r] &	\Ext_A^1(Z, \tau(Z))  		 }
\end{xy}$$
of $\End_A(Z)^{op}$-modules by Lemma \ref{functormono}, which is in fact a chain of isomorphisms by Lemma \ref{Extiso}. In particular, $\Ext_A^1(Z, X_{X_{\tau(Z)}}) \cong \Ext_A^1(Z, Y_{\tau(Z)}) \cong \Ext_A^1(Z, \tau_\Omega(Z))$ as $\tau_\Omega(Z)$ is the only direct summand of $Y_\tau(Z)$ that is not Ext-injective. Hence we can decompose $X_{X_{\tau(Z)}} = I \oplus \tau_\Omega(Z)$ and obtain $\Ext_A^1(Z, I) = 0$. Note that the restriction of $f = f_{X_{\tau(Z)}}|_I : I \to X_{\tau(Z)}$ still induces a monomorphism of contravariant functors $\Ext_A^1(\ \ ,f)|_\Omega : \Ext_A^1(\ \ , I)|_\Omega \to \Ext_A^1(\ \ , X_{\tau(Z)})|_\Omega$ and hence $I$ is Ext-injective in $\Omega$ by Lemma \ref{functormonoExtinj}.
\bigskip
\\It remains to show that $\tau_\Omega(Z)$ is actually a direct summand of $X_{\t(Z)}$ if $X_{\tau(Z)} = J \oplus \t(Z)$ denotes the decomposition of $X_{\tau(Z)}$ into its Ext-injective part and the indecomposable relative Auslander-Reiten translate of $Z$ in $\chi$. Suppose $\tau_\Omega(Z)$ is not a direct summand of $X_{\t(Z)}$ then we have 
$$\Ext_A^1(Z, \tau_\Omega(Z)) \cong \Ext_A^1(Z, X_J) \cong \Ext_A^1(Z, J) = 0,$$
\newpage
which cannot be true. Consequently, $\tau_\Omega(Z)$ is a direct summand of the domain of $f_{\t(Z)} : X_{\t(Z)} \to \t(Z)$ and all other direct summands of $X_{\t(Z)}$ must be Ext-injective as they are also direct summands of $I$. \qed

\begin{Lemma}\label{approxinequality}
If $M$ is module in $A$-mod and $X_M$ denotes the domain of a minimal right $\Omega$-approximation of $M$, then we have 
$$l(X_M) \leq l(M)s$$
where $s = \max\{l(X_{S})| S$ a simple module in $A$-mod$\}$
\end{Lemma}

\proof
\\We prove the statement by induction on $l(M)$, but first we recall a result from \cite{AR91}. Let
 $$\begin{xy}
  \xymatrix{ 
    0 \ar[r]   &	X \ar[r]				&	Y  \ar[r]^f				&	Z \ar[r]		&	0	 }
\end{xy}$$
be an exact sequence in $A$-mod. Moreover, let $f_X : X_X \to X$ and $f_Z : X_Z \to Z$ denote minimal right $\Omega$-approximations of $X$ and $Z$ respectively and let $X_Z \times_Z Y$ denote the pullback of $f_Z$ and $f$. We then obtain a short exact sequence
 $$\begin{xy}
  \xymatrix{ 
    0 \ar[r]   &	X \ar[r]				&	X_Z \times_Z Y \ar[r] 		&	X_Z \ar[r]		&	0	 }
\end{xy}$$
and by Lemma \ref{Extiso} we have $\Ext_A^1(X_Z, X) \cong \Ext_A^1(X_Z, X_X)$ from which we obtain a commutative diagram 
 $$\begin{xy}
  \xymatrix{ 
    0 \ar[r]   &   	X_X  \ar[r] \ar[d]_{f_X} 	&	Y_Y \ar[r] \ar[d]_h 	&    	X_Z \ar[r] \ar@{=}[d]  	&	0	\\
     0 \ar[r]   &	X \ar[r] \ar@{=}[d]		&	X_Z \times_Z Y \ar[d]_g	\ar[r] 	&	X_Z \ar[r] \ar[d]_{f_Z}		&	0	\\
    0 \ar[r]   &	X \ar[r]			&	Y  \ar[r]^f			&	Z \ar[r]			&	0.	 }
\end{xy}$$
By \cite[Proposition 3.6]{AR91} the composition $gh = g_Y: Y_Y \to Y$ is a right $\Omega$-approximation of $Y$, which is not necessarily right minimal. If $l(M) = 1$, then $M$ is a simple $A$-module and it follows that 
$$l(X_M) \leq \max\{l(X_{S})| S \text{ a simple module in } A \text{-mod}\} = l(M)s.$$ 
Suppose now the statement is true for modules of length $n-1$ and $l(M) = n$. Let $S$ be a simple module that embeds into $M$, we then obtain a commutative diagram
 $$\begin{xy}
  \xymatrix{ 
    0 \ar[r]   &   	X_S  \ar[r] \ar[d]_{f_S} 	&	Y_M \ar[r] \ar[d]_{g_M} 	&    	X_{M/S} \ar[r] \ar[d]_{f_{M/S}}  	&	0	\\  
    0 \ar[r]   &	S \ar[r]			&	M  \ar[r]^f			&	M/S \ar[r]			&	0.	 }
\end{xy}$$
By induction we get 
$$l(Y_M) = l(X_S) + l(X_{M/S})  \leq l(S)s + l(M/S)s = l(M)s.$$
In particular, if $f_M : X_Y \to Y$ denotes a minimal right $\Omega$-approximation, then
$$l(X_M) \leq l(Y_M) \leq l(M)s.$$\qed
\bigskip
\\It follows from this proof that if $f_M$ in the commutative diagram of Lemma \ref{approxdiagram} is a minimal right $\Omega$-approximation of $M$, then $f_N$ is also right minimal. Note that in general even if $f_X$ and $f_Z$ are minimal right $\Omega$-approximations, the morphism $f_Y$ we obtain by the construction in the previous proof is not necessarily a minimal right $\Omega$-approximation as the standard example shows.

\begin{Ex}
Consider the standard example \ref{standardexample} and let $\Omega = \F(\Delta)$ be the category of standard-filtered modules. Then the almost split sequence 
 $$\begin{xy}
  \xymatrix{ 
    0 \ar[r]   &	M_2 \ar[r]			&	X \oplus \nabla_3  \ar[r]			&	N_3 \ar[r]		&	0.	 }
\end{xy}$$
gives rise to a commutative diagram with an approximation of $X \oplus \nabla_3$ that is not minimal.
\end{Ex}

Let us first recall the Auslander-Reiten quiver and mark all standard-filtered modules red.

$$\begin{xy}
\xymatrixrowsep{0.5in}
\xymatrixcolsep{0.5in}

  \xymatrix@!0{ 
&	&		&	&					& {\color{red}[P_1]} \ar[dr] 	&	&	&	&	\\
&	&[{\color{red}P_4} \ar[dr]	&	&{\color{red}\Delta_2} \ar[dr] \ar[ur] \ar@{.>}[ll]	&		&\nabla_2 \ar[dr] \ar@{.>}[ll]	&	&I_4] \ar[dr] \ar@{.>}[ll] & \\
&S_3\ar[dr] \ar[ur]	& &{\color{red}X} \ar[dr] \ar[ur] \ar@{.>}[ll]&	&S_2 \ar[dr] \ar[ur] \ar@{.>}[ll]&	&Y \ar[dr] \ar[ur] \ar@{.>}[ll] &	&S_3  \ar@{.>}[ll]\\ 
{\color{red}N_2} \ar[ddr] \ar[dr] \ar[ur] &	&M_2 \ar[dr] \ar[ur] \ar@{.>}[ll]&	&N_3 \ar[ddr] \ar[dr] \ar[ur] \ar@{.>}[ll]& &M_3 \ar[dr] \ar[ur] \ar@{.>}[ll]&	&{\color{red}N_2} \ar[dr] \ar[ur] \ar@{.>}[ll] & \\
&{\color{red}S_1} \ar[ur]&  &\nabla_3 \ar[ur] \ar@{.>}[ll]&	&S_4 \ar[ur] \ar@{.>}[ll]&	&{\color{red}\Delta_3} \ar[ur] \ar@{.>}[ll] &	&{\color{red}S_1}  \ar@{.>}[ll]\\ 
&{\color{red}[P_2]} \ar[uur]	&	&	&	& {\color{red}[P_3]} \ar[uur] 	&	&	&	&	
 }
\end{xy}$$

It is now easy to see that $f_{M_2} : P_2 \oplus S_1 \to M_2$, $f_{X \oplus \nabla_3} : X \oplus P_2 \to X \oplus \nabla_3$ and $f_{N_3} : X \to N_3$
are minimal right $\Omega$-approximations of the modules in the given almost split sequence. On the other hand, the commutative diagram we obtain by the construction in the proof Lemma \ref{approxinequality} is the following.
$$\begin{xy}
  \xymatrix{ 
0 \ar[r]   &P_2 \oplus S_1  \ar[r] \ar[d]_{f_{M_2}}	&X \oplus P_2 \oplus S_1 \ar[r] \ar[d] 	&    	X \ar[r] \ar[d]_{f_{N_3}}  	&	0	\\  
    0 \ar[r]   &	M_2 \ar[r]			&	X \oplus \nabla_3  \ar[r]			&	N_3 \ar[r]		&	0.	 }
\end{xy}$$
Clearly the right $\Omega$-approximation of $X \oplus \nabla_3$ is not minimal as the direct summand $S_1$ is redundant.

\section{Irreducible morphisms in subcategories}

\begin{Def}
Let $X,Y$ be modules in $\chi$. We call a morphism $f: X \to Y$ \textbf{irreducible in }$\pmb{\chi}$ or $\pmb{\chi}$\textbf{-irreducible}, if $f$ is neither a split monomorphism nor a split epimorphism and, whenever we have a $Z$ in $\chi$ such that there are morphisms $h: X \to Z$ and $ g: Z \to Y$ that satisfy $f = gh$, then either $h$ is a split monomorphism or $g$ a split epimorphism.
\end{Def}

The connection between irreducible and minimal right almost split morphisms is the same as for $A$-mod.

\begin{Th}\label{irreducibleminimalrightalmostsplit}
\mbox{}
\begin{enumerate}[(a)]
 \item Let $Z$ be an indecomposable module in $\chi$ and $Y$ a non-zero module in $\chi$. Then a morphism $g: Y \to Z$ is irreducible in $\chi$ if and only if there exists a morphism $g' : Y' \to Z$ such that the induced morphism $(g,g') : Y \oplus Y' \to Z$ is a minimal right almost split morphism in $\chi$.
 \item  Dually, if $X$ is an indecomposable module in $\chi$, then a morphism $f: X \to Y$ is $\chi$-irreducible if and only if there is a morphism $f':X \to Y'$ such that the induced morphism $(f,f')^T:X \to Y \oplus Y'$ is a minimal left almost split morphism in $\chi$.
\end{enumerate}

\end{Th}

\proof
\\Assume first that $g: Y \to Z$ is irreducible in $\chi$ and let $h: M \to Z$ be a minimal right almost split morphism in $\chi$. Since $g$ is not a split epimorphism, it factors through $h$, i.e. $g = h \varphi$ where $\varphi$ is a split monomorphism by irreducibility of $g$. We set $Y' = \Coker(\varphi)$ and hence obtain $M \cong Y \oplus Y'$. Moreover, $g' = h|_{Y'} : Y' \to Z$ is a morphism such that $(g,g') : Y \oplus Y' \to Z$ is a minimal right almost split morphism.
\bigskip
\\Suppose now that $h : M \to Z$ is a minimal right almost split morphism in $\chi$. Furthermore, let $M = Y \oplus Y'$ with $Y$ non-zero and let $g = h|_Y: Y \to Z$. Assume that $g = st$ for some morphisms $t : Y \to N$ and $s : N \to Z$ such that $s$ is not a split epimorphism. Since $h$ is a right almost split morphism, there exist a morphism $(u,v)^T: N \to Y \oplus Y'$ such that $s = (g,g')(u,v)^T$. We obtain the following commutative diagram
 $$\begin{xy}
\xymatrixrowsep{0.6in}
\xymatrixcolsep{1.2in}

  \xymatrix@!0{ 
      Y \oplus Y'\ar[r]^{\left( \begin{smallmatrix} t&0\\ 0&\Id_{Y'} \end{smallmatrix} \right)} \ar[dr]_{(g,g')}	&	N \oplus Y' \ar[r]^{\left( \begin{smallmatrix} u&0\\ v&\Id_{Y'} \end{smallmatrix} \right)} \ar[d]^{(s,g')}  	&	Y \oplus Y' \ar[dl]^{(g,g')} \\
        	&	Z 	& }
\end{xy}$$

Hence the composition ${\left( \begin{smallmatrix} ut&0\\ vt&\Id_{Y'} \end{smallmatrix} \right)}$ is an isomorphism as $(g,g')$ is right minimal. Thus $ut : Y \to Y$ is an isomorphism, showing that $t$ is a split monomorphism. It follows that $g$ is $\chi$-irreducible. \qed

\bigskip

Analogously to $A$-mod, one can define the $\pmb{\chi}$\textbf{-radical} of $\Hom(X,Y)$ for modules $X, Y$ in $\chi$, i.e. 
$\rad_\chi(X,Y)$ is the set of all $f \in \Hom(X,Y)$ such that $gfh$ is not an isomorphism for any $h: Z \to X$ and $g: Y \to Z$ with Z indecomposable in $\chi$. Again, we define \textbf{powers of the $\pmb{\chi}$-radical} inductively, so an $f \in \Hom(X,Y)$ is in $\rad_\chi^n(X,Y)$ if there is a $Z$ in $\chi$ such that $f=gh$ for some $h \in \rad_\chi(X,Z)$, $g \in \rad_\chi^{n-1}(Z,Y)$. Moreover, we also set $\rad_\chi^\infty(X,Y) = \bigcap_{n \in \N}\ \rad_\chi^n(X,Y)$ and $\Irr_\chi(X,Y) = \rad_\chi(X,Y)/\rad_\chi^2(X,Y)$. Note that we can swap the order of $g$ and $h$ in the definition of power of the $\chi$-radical, that is $\rad_\chi^n(X,Y)$ is also the set of all $f \in \Hom(X,Y)$ such that there is a $Z$ in $\chi$ such that $f=gh$ for some $h \in \rad_\chi^{n-1}(X,Z)$, $g \in \rad_\chi(Z,Y)$.

\begin{Lemma} \label{irrradical}
Let $f: X \to Y$ be a morphism between indecomposable modules $X, Y$ in $\chi$. Then $f$ is an irreducible morphism in $\chi$ if and only if $f \in \rad_\chi(X,Y) \backslash \rad_\chi^2(X,Y)$.
\end{Lemma}

\proof
\\We first show for $X$ an indecomposable and $M$ an arbitrary module in $\chi$ that $\rad_\chi(X,M)$ is the set of morphisms from $X$ to $M$ that are not split monomorphisms whereas $\rad_\chi(M,X)$ is the set of morphisms from $M$ to $X$ that are not split epimorphisms.
\bigskip
\\Let $\phi : X \to M$ be a morphism. If $\phi$ is a split monomorphism, we have $\phi'\phi = \Id_X$ and trivially get an isomorphism $\phi'\phi \Id_X$. If there are morphisms $h: Z \to X$ and $g: M \to Z$ with $Z$ indecomposable in $\chi$ and $g\phi h$ an isomorphism, it is clear that $hg\phi$
is an isomorphism on $\Im(h)$ and hence $X$ decomposes into $\Im(h) \oplus \Ker(g\phi)$. Since $X$ is indecomposable, we then know that $h$ must be an isomorphism and $\phi$ a split monomorphism.
\bigskip
\\Now we consider $\psi : M \to X$. If $\psi$ is a split epimorphism, we obtain an isomorphism $\Id_X \psi \psi'$. If we have morphisms $h: Z \to M$ and 
$g: X \to Z$ with $Z$ indecomposable in $\chi$ and $g\psi h$ an isomorphism, we know that $\psi hg$ is an isomorphism on $\Im(\psi h)$ and $X$ decomposes into $\Im(\psi h) \oplus \Ker(g)$, i.e. $\Ker(g) = 0$, $g$ is an isomorphism and $\psi$ is a split epimorphism.
\bigskip
\\Let $f$ be irreducible, then it is neither a split monomorphism nor a split epimorphism and, therefore, in $\rad_\chi(X,Y)$. If $f \in \rad_\chi^2(X,Y)$, we have a decomposition $f=gh$ with $g \in \rad_\chi(X,Z)$ and $h \in \rad_\chi(Z,Y)$, thus neither $g$ nor $h$ splits, contradicting $f$ to be irreducible. Hence $f$ is not in $\rad_\chi^2(X,Y)$.
\bigskip
\\On the other hand, let $f \in \rad\chi(X,Y) \backslash \rad_\chi^2(X,Y)$ and let $f = gh$ be a factorization of $f$ in $\chi$. Then either $g$ or $h$ is not in some $\chi$-radical and, therefore, a split monomorphism or split epimorphism respectively as $X$ and $Y$ are indecomposable. \qed
\bigskip
\\Note that we have $\rad_\chi(X,Y) = \rad(X,Y)$ for $X, Y$ in $\chi$, but usually $\rad_\chi^n(X,Y) \subsetneq \rad^n(X,Y)$ for $n > 1$.

\begin{Def}
The \textbf{Auslander-Reiten quiver} of $\chi$ is a biquiver, i.e. a quiver consisting of vertices $\Gamma_0$ and two disjoint sets of arrows between them, 1-arrows $\Gamma_1$ and 2-arrows $\Gamma_2$. There is a vertex for every isomorphism class of indecomposable modules of $\chi$. Given two indecomposable modules $X$ and $Y$ in $\chi$ there is a 1-arrow from the vertex corresponding to the isomorphism class of $X$ to the vertex corresponding to the isomorphism class of $Y$ if there is an irreducible morphism from $X$ to $Y$. For each almost split sequence
$$\begin{xy}
  \xymatrix{ 
      0 \ar[r]  	& 	X \ar[r] 	&	Y \ar[r]     &   Z \ar[r]	&	0	    }
\end{xy}$$
there is a 2-arrow from the vertex corresponding to the isomorphism class of $Z$ to the vertex corresponding to the isomorphism class of $X$. To distinguish between different types of arrows, 2-arrows are drawn as dotted arrows.
\end{Def}

\begin{Def}\label{sublabeling}
Let $h : M \to N$ be an irreducible morphism in $\chi$ between indecomposable modules and $f : X \to N$ and $g : M \to Y$ extensions to minimal right almost split and minimal left almost split morphisms in $\chi$ respectively. We denote the number of copies of $M$ in a sum decomposition of $X$ and the number of copies of $N$ in a sum decomposition of $Y$ by $n$ and $m$ respectively. As in $A$-mod, the arrow from $X$ to $Y$ in the Auslander-Reiten quiver of $\chi$ has \textbf{valuation} $(n,m)$.
\end{Def}

Throughout the dissertation we denote the Auslander-Reiten quiver of $\chi$ by $\Gamma_\chi$. Just as for Auslander-Reiten quivers of $A$, we mostly do not distinguish between an indecomposable module, its isomorphism class or a corresponding vertex in $\Gamma_\chi$. There are squared brackets to mark Ext-injective and Ext-projective modules in $\Gamma_\chi$. For convenience, we say arrows for 1-arrows and, if the valuation of an arrow is $(1,1)$, we say it is \textbf{trivially valuated}. 
\bigskip
\\If a subquiver $\Gamma$ of $\Gamma_\chi$ contains two modules $X$ and $Y$ such that there is a $1$-arrow or a $2$-arrow between them, then $\Gamma$ also contains this $1$-arrow or $2$-arrow respectively and the $1$-arrow in $\Gamma$ has the same valuation as in the whole Auslander-Reiten quiver. We call a non-empty subquiver $\Gamma$ of $\Gamma_\chi$ a \textbf{connected component} of $\chi$ if for each $X$ in $\Gamma$ the modules in $\Gamma$ are precisely all modules $Y$ such that there is a walk between $X$ and $Y$ in $\Gamma_\chi$. Clearly, each indecomposable module $X$ in $\chi$ is contained in a unique connected component of $\Gamma_\chi$, which we call the connected component of $X$ in $\Gamma_\chi$. Analogously to $A$-mod, we say a connected component is finite if it contains only finitely many indecomposable non-isomorphic modules.
\bigskip
\\Note that also in the Auslander-Reiten quiver of $\chi$ it frequently happens that the same indecomposable module occurs more than once. In this case we again differentiate between the corresponding vertices of this module in $\Gamma_\chi$. Having defined valuation for arrows in $\chi$, we can generalize a fact about valuation of arrows in $A$-mod.

\begin{Th}\label{radbasis}
Let $X$ and $Y$ be indecomposable modules in $\chi$ such that there is an irreducible morphism from $X$ to $Y$. Suppose the valuation of the 1-arrow from $X$ to $Y$ is $(n,m)$, i.e. there is a module $M$ such that $X$ is not a direct summand of $M$ and a minimal right almost split morphism $f : nX \oplus M \to Y$  in $\chi$. Let $f_1, \ldots, f_n$ name the morphisms obtained by restricting $f$ to the different copies of $X$ in $nX \oplus M$. Then $\{\overline{f_1}, \ldots, \overline{f_n}\}$ is a basis for $\Irr_\chi(X,Y)$ as a $T_X^{op}$ and $T_Y$-vector space. In particular, $n = m$.
\end{Th}

\proof
\\The induced morphism $(f_1, \ldots, f_n)  : nX \to Y$ is irreducible as it can be extended to a minimal right almost split morphism in $\chi$. Suppose $\{\overline{f_1}, \ldots, \overline{f_n}\}$ is linearly dependent, i.e. there are endomorphisms $\alpha_i : X \to X$ such that $\overline{f_1 \alpha_1} + \overline{f_2 \alpha_2} + \cdots + \overline{f_n \alpha_n} = 0$ with at least one $\overline{\alpha_i}$ being non-zero. This relation implies that $f_1 \alpha_1 + f_2 \alpha_2 + \cdots + f_n \alpha_n$ is in $\rad_\chi^2(X,Y)$. 
\newpage
However, since at least one $\overline{\alpha_i}$ is non-zero, $\alpha_i$ must be an isomorphism and the morphism $\alpha = (\alpha_i)^T : X \to nX$ is a split monomorphism. Therefore, there is a morphism $\beta : nX \to X$ such that $\beta\alpha = \Id_X$ and $nX = (n-1)X \oplus X = \Ker(\beta) \oplus \Im(\alpha)$. It follows that we can extend $\alpha$ to an isomorphism $\varphi = (\alpha, \alpha') : X \oplus (n-1)X \to nX$. This gives rise to a new minimal right almost split morphism $f \varphi$. The diagram below shows how an arbitrary morphism in $\chi$ mapping to $Y$ factors over $f \varphi$ using that $f$ is a minimal right almost split morphism.
$$\begin{xy}
  \xymatrix{ 
    & Z \ar@{.>}[d]^g \ar@{.>}[dl]_{\varphi^{-1} g} \ar[dr]  	&	 \\
  nX  \ar[r]^\varphi & nX  \ar[r]^{f}    	&	Y	    }
\end{xy}$$
Since we have extended the composition $f \alpha = f_1 \alpha_1 + f_2 \alpha_2 + \cdots + f_n \alpha_n$ to a minimal right almost split morphism, it must be irreducible in $\chi$ by \ref{irreducibleminimalrightalmostsplit}, contradicting the fact that $f_1 \alpha_1 + f_2 \alpha_2 + \cdots + f_n \alpha_n$ is in $\rad_\chi^2(X,Y)$. Consequently, $\{\overline{f_1}, \ldots, \overline{f_n}\}$ is linearly independent.
\bigskip
\\Now we show that $\{\overline{f_1}, \ldots, \overline{f_n}\}$ span $\Irr_\chi(X,Y)$. Let $g : X \to Y$ be irreducible in $\chi$ and $\overline{g}$ its coset in $\Irr_\chi(X,Y)$. Since $g$ is not a split epimorphism, there is a morphism $h: X \to nX \oplus M$ such that $g = fh$. We can decompose this morphism into $g = f|_{nX}ph + f|_{M}qh$ where $p$ and $q$ name the projections from $nX \oplus M$ to $nX$ and $M$ respectively. Clearly, $qh$ is in $\rad_\chi(X,M)$ as $M$ does not contain $X$ as a direct summand. Thus we get $f|_{M}qh \in \rad_\chi^2(X,Y)$ and $\overline{g} = \overline{f|_{nX}ph}$. Using that $ph$ maps from $X$ to $nX$ we decompose it into a set of endomorphisms $\{\alpha_i\} \subset \End_A(X)$ such that $\overline{g} = \sum \overline{f_i\alpha_i}$.
\bigskip
\\Since $A$ is a $K$-algebra over an algebraically closed field, we have $T_X^{op} \cong T_Y \cong K$, hence the statement follows immediately for $T_Y$. As a trivial consequence of that, we obtain $n = m$.
\qed
\bigskip
\\This theorem enables us to generalize another well-established result.
\begin{Cor}\label{valuationisinvariantoftranslation}
Let $X$ and $Y$ be indecomposable modules in $\chi$ such that there is an arrow from $X$ to $Y$ in the Auslander-Reiten quiver of $\chi$ with valuation $(n,n)$. If $Y$ is not Ext-projective, then the arrow from $\t(Y)$ to $X$ is also valued $(n,n)$.
\end{Cor}

\proof
\\By our assumptions $n$ is the number of direct summand isomorphic to $X$ in the domain of a minimal right almost split morphism $f: M \to Y$ in $\chi$. Consequently, by Lemma \ref{almostsplitsequenceinchi}, $n$ is also the number of direct summands isomorphic to $X$ in the codomain of a minimal left almost split morphism $g: \t(Y) \to M$. It follows from Theorem \ref{radbasis} that the valuation of the corresponding arrow is indeed $(n,n)$. \qed
\bigskip
\\Since any arrow from $X$ to $Y$ in the Auslander-Reiten quiver of $\chi$ is valued $(n,n)$ for some $n \in \N$, it again makes sense to say there are $\pmb{n}$ \textbf{arrows} from $X$ to $Y$. Moreover, we draw $n$ arrows from $X$ to $Y$ in the Auslander-Reiten quiver of $\chi$. We say there are multiple arrows from $X$ to $Y$ if $n \geq 2$ and a connected component contains multiple arrows if and only if there are multiple arrows from $X$ to $Y$ for some modules $X$ and $Y$ in that connected component.
\bigskip
\\Let us recall the well-known lemma by Harada and Sai, which clearly can be applied to irreducible morphisms in $\chi$. Therefore, we can generalize a criterion for the Auslander-Reiten quiver of an Artin algebra to consist of only one component to a criterion for the Auslander-Reiten quiver of $\chi$ to have this property. 

\begin{Lemma}[Harada-Sai Lemma]\cite[VI. Lemma 1.2]{ARS95} \label{HaradaSai}
Let $X_0, \ldots,$ $X_{2^n-1}$ be indecomposable modules such that $l(X_i) \leq n$ for $i = 0, \ldots, 2^n-1$. Suppose $f_i :X_i \to X_{i+1}$ are non-isomorphisms, then $f_{2^n-1} \cdots f_1 = 0$.
\end{Lemma}

\begin{Th}\label{finitecomponent}
Let $A$ be an indecomposable algebra. If there is a component $\Gamma$ in the Auslander-Reiten quiver of $\chi$ such that the length of all modules in $\Gamma$ is bounded, then it is the only component. 
\end{Th}

\proof
\\Let $n$ be the positive integer such that $l(X) \leq n$ for all $X$ in $\Gamma$. Let $X$ be in $\Gamma$ and $Y$ be an arbitrary indecomposable module in $\chi$ such that $\rad_\chi^{2^n}(X,Y) \neq 0$. We then know there must be a chain of at least $2^n$ non-isomorphisms from $X$ to $Y$. This contradicts Lemma \ref{HaradaSai}, so we have $\rad_\chi^{2^n}(X,-) = \rad_\chi^{2^n}(-,X) = 0$. In particular, we have $\Hom(X,Y) = \Hom(Y,X) = 0$ for $Y$ not in $\Gamma$. 
\bigskip
\\Since $\chi$ contains a functorially finite resolving subcategory $\Omega$, it follows that it, in particular, contains all projective modules. But there is a projective module $P$ such that $\Hom(P,X) \neq 0$ and hence $P$ is in $\Gamma$. As projective modules of an indecomposable algebra are always connected, we conclude that all projective modules are in $\Gamma$. Moreover, for every indecomposable module $Y$ there is a non-zero homomorphism from a projective module to $Y$, so every indecomposable module must be contained in $\Gamma$. 
\qed
\bigskip
\\Note that this theorem does not hold for arbitrary functorially finite subcategories that are closed under extensions, as the standard example shows. 
\begin{Ex}
Let $A$ be the algebra of the form
 $$\begin{xy}
\xymatrixrowsep{0in}
\xymatrixcolsep{-0.1in}
  \xymatrix{ 
    	S_1    	&		&	&S_2	&	&		&	&S_3	&	&		&S_4	\\	
	S_2    	&\ \oplus \	&S_1	&	&S_3	&\ \oplus \	&S_2	&	&S_4	&\ \oplus \	&S_3	\\	
	S_1    	&		&	&S_2	&	&		&	&S_3	&	&		&		
    }
\end{xy}$$
which we have analyzed in Example \ref{standardexample}. Then there is a functorially finite subcategory $\chi$ of $A$-mod that is finite and closed under extensions, but its Auslander-Reiten quiver has 2 connected components.
\end{Ex}

We have already seen that $A$ is representation finite in Example \ref{standardexample}. Hence every subcategory of $A$-mod is functorially finite. We consider the subcategory $\chi = \add(P_1 \oplus I_4)$ that consists of all modules whose indecomposable direct summands are isomorphic to $P_1$ or $I_4$. Their Jordan-H\"{o}lder composition series 
 $$\begin{xy}
\xymatrixrowsep{0in}
\xymatrixcolsep{-0.1in}
  \xymatrix{ 
	&S_1&	&S_3	\\	
P_1 = I_1 = \ \ & S_2 &\ \ \ \ \ \ \ \ \ \ I_4= \ \	&S_4.\\
	  &S_1&	&
    }
\end{xy}$$
show that $\Hom(P_1,I_4) = \Hom(I_4,P_1) = 0$. Moreover, since $P_1$ is also injective, $\chi$ is closed under extensions, but its Auslander-Reiten quiver consists of two components that contain one indecomposable module each.
\bigskip
\\Moreover, with the help of the lemma of Harada and Sai we can also prove a generalization of the first Brauer-Thrall conjecture.

\begin{Cor}
The subcategory $\chi$ is finite if and only if the Jordan-H\"{o}lder length of its indecomposable modules is bounded.
\end{Cor}

\proof
\\Let $n$ be a positive integer such $l(X) \leq n$ for all indecomposable modules $X$ in $\chi$. It follows that $\rad_\chi^{2^n}(X,Y) = \rad_\chi^\infty(X,Y) = 0$ for all indecomposable modules $X,Y$ in $\chi$ by Lemma \ref{HaradaSai}. Moreover, for every module $X$ there is a projective module $P$ and a path in $\Gamma_\chi$ from $P$ to $X$ of length shorter than $2^n$. Since there are only finitely many projective modules and finally many paths of length shorter than $2^n$, there can only be finitely many indecomposable modules in $\chi$ up to isomorphism. The converse is trivial. \qed
\bigskip
\\The second Brauer-Thrall conjecture, that is an Artin algebra $A$ is representation infinite if and only if there are infinitely many positive integers $n_1, n_2, \ldots$ such that for each $i \in \N$ there are infinitely many non-isomorphic indecomposable modules of Jordan-H\"{o}lder length $n_i$, has been proved for Artin algebras over an algebraically closed field. The original proof can be found in \cite{NR75}, a sketch of the proof in English is given in \cite{R80*}. The problem is still open for algebras over arbitrary fields or rings, but it has been established in \cite{S80} that if there is one $n \in \N$ such that there are infinitely many non-isomorphic indecomposable modules of length $n$, then there are infinitely many $n_i \in \N$ with that property. This is nothing but the induction step for a proof of the second Brauer-Thrall conjecture and is often referred to as the one and a half Brauer-Thrall conjecture. We generalize the aforementioned result using the same techniques, but we need to find a boundary 
for the length of an indecomposable module $X$ in terms of the length of an indecomposable module $Y$ if there is an irreducible morphism from $X$ to $Y$ in $\Omega$. The first lemma is the dual of \cite[Lemma 3.1]{S80}, which fits better into our context.

\begin{Lemma}\label{irreducibleinequality}
Let $X$ and $Y$ be indecomposable modules and assume $f: X \to Y$ is an irreducible morphism. Then
$$|l(X) - l(Y)| \leq l(Y) \cdot m^2$$
where $m = \max\{l(_AA),l(A_A)\}$.
\end{Lemma}

\proof
\\If $Y$ is projective, then $f$ is a monomorphism and $l(X) < l(Y)$, so obviously the stated inequation holds in this case. Suppose now that $Y$ is not projective, then there is an almost split sequence
 $$\begin{xy}
  \xymatrix{ 
    0 \ar[r]   &	DTr(Y) \ar[r]	&	X \oplus X'  \ar[r]		&	Y \ar[r]		&	0.	 }
\end{xy}$$
This gives us $l(DTr(Y)) + l(Y) = l(X) + l(X')$ and, therefore, we obtain $|l(X) - l(Y)| \leq \max\{l(Y),l(DTr(Y))\}$. Moreover, we have $l(DTr(Y)) = l(Tr(Y))$ since $D$ is a duality. So let $P_1 \to P_0 \to Y$ be a minimal projective resolution, then we obtain a morphism $\Hom(f,A) : \Hom(P_0,A) \to \Hom(P_1,A)$ such that $Tr(Y) = \Coker(\Hom(f,A))$. Therefore, we have
$$l(DTr(Y)) = l(Tr(Y)) \leq l(\Hom(P_1,A)).$$
Since $\Hom(P_1,A)$ is a projective $A^{op}$-module, we know that $l(\Hom(P_1,A)) \leq l(A_A)n$ where $n$ denotes the number of indecomposable direct summands of $P_1$. Clearly, $n \leq l(P_0)$ and the number of direct summands of $P_0$ is at most $l(Y)$ and hence $n \leq l(P_0) \leq l(Y)l(_AA)$ as the Jordan-H\"{o}lder length of each direct summand of $P_0$ is less than or equal to $l(_AA)$. Overall we obtain $l(DTr(Y)) \leq l(Y)m^2$ and
$$|l(X) - l(Y)| \leq l(Y)m^2$$
\qed

\begin{Lemma}\label{irreducibleboundary}
Let $X$ and $Y$ be indecomposable modules in $\Omega$ and suppose there is an irreducible morphism $f : X \to Y$ in $\Omega$. Then 
$$|l(X) - l(Y)| \leq l(Y) p $$ where $p = (s(1 + m^2) - 1)$, $m = \max\{l(_AA),l(A_A)\}$ and $s = \max\{l(X_{S_i})\}$ for all simple modules $S_i$. 
\end{Lemma}
\proof
\\Since a minimal right almost split morphism in $\Omega$ is given by a minimal right $\Omega$-approxima-tion of the domain of a minimal right almost split morphism in $A$-mod, there is a module $M$ in $A$-mod such that there is an irreducible morphism $g: M \to Y$ and $X$ is a direct summand of the domain $X_M$ of a right $\Omega$-approximation of $M$. We then know that 
$$l(X) \leq l(M)s$$
by Lemma \ref{approxinequality}, where $s$ is the maximal Jordan-H\"{o}lder length of a minimal right $\Omega$-approximation of a simple $A$-module. Note that $s$ is always finite as $A$ is an Artin algebra.
$$|l(X) - l(Y)| \leq l(Y)(s(1 + m^2) - 1)$$
clearly holds if $l(X) \leq l(Y)$. So assume now $l(X) > l(Y)$, then
$$|l(X) - l(Y)| \leq |l(M)s - l(Y)| = |l(M)s -l(Y)s + l(Y)s - l(Y)|$$
$$\leq |l(M) - l(Y)|s + |l(Y)s - l(Y)| \leq l(Y)s m^2 + l(Y)(s - 1)  = l(Y)p,$$
where we use the inequation obtained in Lemma \ref{irreducibleinequality}. \qed
\newpage
This inequation allows us to prove the one and a half Brauer-Thrall conjecture for a functorially finite resolving subcategory $\Omega$. The proof closely follows \cite[Theorem 3.2]{S80}.

\begin{Th} \label{BrauerThrall1.5}
Let $\Omega$ be a functorially finite resolving subcategory of $A$-mod such that there exist $\mathfrak{N} \geq \mathfrak{N}_0$ non-isomorphic indecomposable modules of length $n$, where $\mathfrak{N}_0$ is the cardinality of a countable set. Then there are infinitely many positive integers $n_i$ with $\mathfrak{N}$ non-isomorphic indecomposable modules of length $n_i$.
\end{Th}

\proof
\\Let $\{M_i\}_{i \in \mathfrak{N}}$ be the set of non-isomorphic, indecomposable modules of length $n$. We prove the statement by showing that for each integer $n' > n$ there is an integer $k \geq n'$ such that there are $\mathfrak{N}$ indecomposable, non-isomorphic modules of length $k$.
\bigskip
\\Suppose there is an integer $n' > n$ such that for all $m \geq n'$ the number of modules of length $m$ is strictly less than $\mathfrak{N}$ and fix such an $m$. We know that there is a simple module $S$ such that $S$ is a direct summand of $\soc(M_i)$ for at least $\mathfrak{N}$ modules $M_i$ and, therefore, there are non-zero morphisms $f_i : S \to M_i$ for these modules. It follows that there exist $\mathfrak{N}$ non-zero morphisms $g_i : X_S \to M_i$.
\bigskip
\\Let us first assume that $\mathfrak{N}$ of the morphisms $g_i$ are not sums of compositions of irreducible morphisms in $\Omega$, that is there are $\mathfrak{N}$ morphisms $g_i$ in $\rad_\Omega^\infty(X_S, M_i)$. We then get $\mathfrak{N}$ chains 
$$\begin{xy}
  \xymatrix{ 
  \ldots \ar[r] & M_{i,2^m+1} \ar[r] & M_{i,2^m} \ar[r] & \ldots \ar[r] & M_{i,1} \ar[r] & M_{i,0} \ar[r] & M_i	 }
\end{xy}$$
of arbitrary length of indecomposable modules $M_{i,j}$ and $\Omega$-irreducible morphisms. Then by Lemma \ref{HaradaSai} there is an $m' \in \N$ such that in each of these chain in the last $m'$ steps there is a module of length greater or equal to $m$. Let $M$ be such a module, then there is only a finite number of modules $M_i$ such that $M$ occurs in the last $m'$ steps of the chains associated to the modules $M_i$ as $\Gamma_\Omega$ is a locally finite Auslander-Reiten quiver. Furthermore, the length of modules in the last $m'$ steps of each chain associated to an $M_i$ is bounded by $(np)^{m'}$ where $p$ is the constant defined in Lemma \ref{irreducibleboundary}. It follows that there are at least $\mathfrak{N}$ non-isomorphic indecomposable modules with length in an interval $[m, (np)^{m'}]$. So there exists a $k \in \N$, $m \leq k \leq (np)^{m'}$ such that there are $\mathfrak{N}$ non-isomorphic indecomposable modules of length $k$, which is a contradiction.
\bigskip
\\Suppose now that for $\mathfrak{N}$ of the $M_i$ there is a chain of $\Omega$-irreducible morphisms with non-zero composition starting in $X_S$ and ending in $M_i$. Since there are only finitely many arrows ending in any module in an Auslander-Reiten quiver, we know that only finitely many of the chains can be shorter than $m'$ and the same argument as in the first case can be applied to the other chains.
\bigskip
\\So in both cases we get that there is no bound on the integers $k$ such that there exist $\mathfrak{N}$ non-isomorphic indecomposable modules of length $k$, which completes the proof of the theorem. \qed
\newpage
\section{Decompositions of morphisms}\label{decomposition}
The purpose of this section is to create the general setup for the subsequent section. When an arbitrary morphism is decomposed into a sum of compositions of irreducible morphisms between indecomposable modules, this decomposition is not unique in general. In particular, not even the aforementioned indecomposable modules are determined by the morphism. In order to give a description when at least these modules are unique, we need to establish a setting in which it makes sense to consider decompositions. First of all we generalize a result from \cite[V 7.4]{ARS95}.

\begin{Lemma}\label{decomp}
Let $X$ and $Y$ be indecomposable modules in $\chi$ and let $f \in \rad_\chi^n(X,Y)$ with $n \geq 2$. Then we have the following.
\begin{enumerate}[(a)]
 \item There exist an integer $s \geq 1$, indecomposable modules $Z_1, \ldots , Z_s$, morphisms $f_i \in \rad_\chi(X, Z_i)$ and morphisms $g_i :Z_i \to Y$ with each $g_i$ a sum of compositions of $n-1$ irreducible morphisms between indecomposable modules such that $f = \sum_{i=1}^s g_if_i$.
\item If $f \in \rad_\chi^n(X,Y) \backslash \rad_\chi^{n+1}(X,Y)$ then at least one of the $f_i$ in $(a)$ is irreducible and $f = u + v$, where $u$ is a non-zero sum of compositions of $n$ irreducible morphisms between indecomposable modules and $v \in \rad_\chi^{n+1}(X,Y)$.
\item Dually, there exist an integer $t \geq 1$, indecomposable modules $W_1, \ldots , W_t$, morphisms $f_i : X \to W_i$ and morphisms $g_i \in \rad_\chi(W_i, Y)$ with each $f_i$ a sum of compositions of $n-1$ irreducible morphisms between indecomposable modules such that $f = \sum_{i=1}^t g_if_i$.
\item If $f \in \rad_\chi^n(X,Y) \backslash \rad_\chi^{n+1}(X,Y)$ then at least one of the $g_i$ in $(c)$ is irreducible and $f = u + v$, where $u$ is a non-zero sum of compositions of $n$ irreducible morphisms between indecomposable modules and $v \in \rad_\chi^{n+1}(X,Y)$.
\end{enumerate}
\end{Lemma}

\proof
\\We proof all statements by induction on $n$. For $n = 2$ let $g : Z \to Y$ be a minimal right almost split morphism in $\chi$. The decomposition of $Z$ into indecomposable modules $Z \cong \bigoplus_{i=1}^s Z_i$ induces morphisms $g_i : Z_i \to Y$. Since $f \in \rad_\chi^2(X,Y)$, it is not a split epimorphism and there is a morphism $f' : X \to Z$ such that $f = gf'$. Let $f_i: X \to Z_i$ be the morphisms induced by $f'$ and the decomposition $Z \cong \bigoplus_{i=1}^s Z_i$. Then we get that $f_i \in \rad_\chi(X,Z_i)$, the morphisms $g_i$ are irreducible and $f = \sum_{i=1}^s g_if_i$. Moreover, if $f \notin \rad_\chi^3(X,Y)$, then not all $f_i$ are in $\rad_\chi^2(X, Z_i)$. Hence for at least one $i \in \{1, \ldots, s\}$ we have that $f_i$ is irreducible. We obtain the decomposition $f = u + v$ by setting $u = \sum g_{i}f_{i}$, where we sum over all $i$ such that $f_{i}$ is irreducible. This establishes the claim for $n=2$.
\bigskip
\\Assume now that $f \in \rad_\chi^n(X,Y)$ for $n \geq 3$. Let $g' : Z \to Y$ be a minimal right almost split morphism in $\chi$ with $Z \cong \bigoplus_{i=1}^t Z_i$ and induced morphisms $g_i' : Z_i \to Y$. As $f$ is in $\rad_\chi^n(X,Y)$, there is a morphism $f' \in \rad_\chi^{n-1}(X,Z)$ such that $f = g' f'$. Furthermore, we consider the induced morphisms $f_i' : X \to Z_i$ and rewrite $f = \sum_{i=1}^t g_i' f_i'$. Each $f_i'$ is in $\rad_\chi^{n-1}(X,Z_i)$ as it maps to a distinct direct summand of $Z$. By induction there are indecomposable modules $Z_{ij}$ for $j = 1, \ldots, s_i$ and morphisms $f_{ij} \in \rad_\chi(X, Z_{ij})$ and $h_{ij} :Z_{ij} \to Z_i$ which are sums of compositions of $n-2$ irreducible morphisms between indecomposable modules such that $f_i' = \sum_{j=1}^{s_i} h_{ij}f_{ij}$. Then the composed morphism $g_{ij} = g_i' h_{ij}$ is clearly a sum of compositions of $n-1$ irreducible morphisms between indecomposable modules and we have that 
$$f = \sum_{i=1}^t \sum_{j=1}^{s_i} g_{ij}f_{ij}.$$ 
In order to prove $(b)$, observe that when $f \notin \rad_\chi^{n+1}$ then not all $f_{ij}$ can be in $\rad_\chi^2(X,Z_{ij})$. This shows that at least one $f_{ij}$ is irreducible. We obtain the decomposition $f = u + v$ by setting $u = \sum g_{ij}f_{ij}$, where we sum over all $ij$ such that $f_{ij}$ is irreducible.
\qed

\bigskip

Let $X$ and $Y$ be indecomposable modules in $\chi$ such that there is a $\chi$-irreducible morphism from $X$ to $Y$. Suppose the corresponding arrow in $\Gamma_\chi$ is valued $(n,n)$, i.e. there are $n$ arrows from $X$ to $Y$. By Theorem \ref{radbasis} there are irreducible morphisms $f_1, \ldots , f_n \in \rad_\chi(X,Y) \backslash \rad_\chi^2(X,Y)$ such that $\{ \overline{f_1}, \ldots, \overline{f_n} \}$ is a basis of $\Irr_\chi(X,Y)$. We label each of the $n$ arrows with one morphism $f_i$ and fix this labeling for the chapter.
\bigskip
\\By doing this we have not only chosen a morphism for each arrow in $\Gamma_\chi$, but also a minimal right almost split morphism and a minimal left almost split morphism for each module. After making this choice it is possible to consider a path in $\Gamma_\chi$ as a morphism that is a composition of irreducible morphisms between indecomposable modules; therefore, we do not distinguish between a path $\gamma$ and a composition of irreducible morphisms between indecomposable modules $g$ if $\gamma = g \varphi$ for an isomorphism $\varphi$. In particular, we can consider the morphisms $g_i$ appearing in the first part of Lemma \ref{decomp} as a sum of paths, because they are all obtained from minimal right almost split morphisms by construction. Applying this labeling and the choice of minimal right almost split morphisms to \ref{decomp} we can prove the existence of decompositions given by our fixed labeling.

\begin{Lemma}\label{pathdecomp}
Let $X$ and $Y$ be indecomposable modules in $\chi$ and let $f$ be a morphism in $\rad_\chi^n(X,Y) \backslash \rad_\chi^{n+1}(X,Y)$. Then for each $m \geq n$ $f$ can be written as $f = v + \sum_{i=n}^{m-1}u_i$, where $u_i$ is a sum of paths of length $i$ and $v \in \rad_\chi^m(X,Y)$.
\end{Lemma}
\proof
\\We have already seen that $f = \sum_{i=1}^s g_if_i$, where each $g_i$ is a sum of paths of length $n-1$ and $f_i \in \rad_\chi(X, Z_i)$. Without loss of generality, let $f_1, \ldots, f_k$ be the only irreducible morphisms of the $f_i$. So each $f_i$ for $i = 1, \ldots, k$ can be written as $f_i = h_i \varphi_i$, where $h_i $ is the minimal right almost split morphism mapping to $Z_i$ given by the labeling of the Auslander-Reiten quiver and $\varphi_i$ is a split monomorphism. We go into more detail and write
$$f_i = \sum_{j=1}^{k_i} h_{ij} \varphi_{ij} + \sum_{j=k_i+1}^{n_i} h_{ij} \varphi_{ij}$$ 
where $\varphi_{ij}: X \to X$ are isomorphisms for $j = 1, \ldots, k_i$ and $\varphi_{ij}$ radical morphisms for $j= k_i+1, \ldots, n_i$ which might be zero. Note that $h_{ij}$ are the irreducible morphisms between indecomposable modules forming $h_i$, i.e. $h_i = (h_{i1}, \ldots, h_{in_i})$. By construction $g_i h_{ij} \varphi_{ij}$ is a sum of paths from $X$ to $Y$ for all $j = 1, \ldots, k_i$, where $g_i h_{ij} \varphi_{ij} \in \rad_\chi^{n+1}(X,Y)$ if $j \geq k_i+1$. Consequently, we set $u_n = \sum_{i=1}^k \sum_{j=1}^{k_i} g_i h_{ij} \varphi_{ij}$. We continue inductively with $f- u_n \in \rad_\chi^{n+1}(X,Y)$ until we have the designated decomposition. \qed 
\newpage
To make sure that for an arbitrary morphism this inductive method terminates at some point we need to prove another Lemma. The proof of the first statement closely follows \cite[V 7.2]{ARS95}.

\begin{Lemma}
Let $X$ and $Y$ be indecomposable modules in $\chi$. Then there is an $m \in \N$ such that $\rad_\chi^m(X,Y) = \rad_\chi^\infty(X,Y)$. In particular, there are only finitely many non-zero paths from $X$ to $Y$ in the Auslander-Reiten quiver of $\chi$.
\end{Lemma}

\proof
\\Both $X$ and $Y$ are finitely generated modules, so $\Hom(X,Y)$ is a finitely generated $K$-module and hence of finite length. Therefore, the descending chain
$$\Hom(X,Y) \supset \rad_\chi(X,Y) \supset \cdots \supset \rad_\chi^m(X,Y) \supset \cdots$$ 
becomes stable and there is an $m \in \N$ such that 
$$\rad_\chi^\infty(X,Y) = \bigcap_{n \in \N}\ \rad_\chi^n(X,Y) = \rad_\chi^m(X,Y).$$
This implies that there cannot be non-zero paths of length greater or equal to $m$ from $X$ to $Y$ in the Auslander-Reiten quiver of $\chi$. Since for each module there is only a finite number of immediate predecessors in an Auslander-Reiten quiver, there are clearly at most finitely many paths of each length from $X$ to $Y$ and overall there are only finitely many paths from $X$ to $Y$. \qed
\bigskip
\\The entire setup suggests the following definition.

\begin{Def}
Let $f : X \to Y$ be a morphism between indecomposable modules in $\chi$ that is in $\rad_\chi^n(X,Y) \backslash \rad_\chi^{n+1}(X,Y)$. A \textbf{decomposition} of $f$ is a sum of morphisms $$f = v + \sum_{i=n}^m u_i$$ such that $v \in \rad_\chi^\infty(X,Y)$ and $u_i$ is a sum of paths of length $i$, where the single paths are given by the construction in Lemma \ref{decomp} and, therefore, equal the labeling of the Auslander-Reiten quiver up to an isomorphism of $X$. 
\end{Def}

Without loss of generality, we assume that any sum of paths occurring in a decomposition is non-zero.

\begin{Def}
Let $\gamma$ be a path of length $k$ from $X$ to $Y$ in the Auslander-Reiten quiver of $\chi$. We then say that a decomposition $f = v + \sum_{i=n}^m u_i$ contains $\gamma$ if $\gamma$ occurs as a summand of $u_k$, i.e. $\gamma$ equals a path in $u_k$ up to an isomorphism of $X$.
\end{Def}

\newpage

\section{Uniqueness of sectional paths}\label{secunique}

Recall that a \textbf{sectional path}
 $$\begin{xy}
  \xymatrix{ 
      X_0 \ar[r]^{f_1} 	&	X_1 \ar[r]^{f_2}     &   \cdots	 \ar[r]	&	X_{n-1} \ar[r]^{f_n}	&	X_n	    }
\end{xy}$$
in the Auslander-Reiten quiver of $A$-mod contains no module $X_i$ such that $\tau(X_i) = X_{i-2}$.

\begin{Th}\label{secTh}
Let
 $$\begin{xy}
  \xymatrix{ 
      X_0 \ar[r]^{f_1} 	&	X_1 \ar[r]^{f_2}     &   \cdots	 \ar[r]	&	X_{n-1} \ar[r]^{f_n}	&	X_n	    }
\end{xy}$$
be a sequence of irreducible morphisms between indecomposable A-modules, i.e. a path in the Auslander-Reiten quiver of $A$-mod. Let $f = f_n \cdots f_1 + f_n\phi$, where $\phi$ is a morphism from $X_0$ to $X_{n-1}$ such that there exists a decomposition of $\phi$ that does not contain the path $\gamma = f_{n-1} \cdots f_1$. Suppose that $f$ factors through a morphism $g: Y \to X_n$ for some module $Y$ such that $(f_n,g): X_{n-1} \oplus Y \to X_n$ is irreducible. Then the sequence is not sectional and $f$ factors over $\tau(X_n)$. 
\end{Th}

\proof
\\Note that $f=0$ becomes a special case of this theorem by setting $Y=0$ and $g=0$, because $(f_n,0): X_{n-1} \oplus 0 \to X_n$ is clearly irreducible. The proof is roughly based on \cite[VII. Lemma 2.5]{ARS95} and is done by induction on $n$. Let $n=2$ and $f = f_2(f_1 + \phi) = gh$. We consider
$$\begin{xy}
\xymatrixcolsep{0.6in}
  \xymatrix{ 
      X_0 \ar[r]^(.45){(f_1 + \phi, -h)^T} 	&	X_1 \oplus Y \ar[r]^{(f_2, g)}     &   	X_2	    }
\end{xy}$$
and conclude $(f_2, g)(f_1 + \phi,-h)^T = 0$. If $(f_2, g)$ is a monomorphism, we have $f_1 + \phi = 0$ and, in particular, $\overline{f_1} = \overline{\phi}$ in $\rad(X_0,X_1)/\rad^2(X_0,X_1)$, which is a contradiction to our assumption that there is a decomposition of $\phi$ not containing the arrow labeled $f_1$. Therefore, $(f_2, g)$ is a non-injective irreducible morphism and hence an epimorphism. Consequently, $X_2$ is not projective and there is an almost split sequence
$$\begin{xy}
\xymatrixcolsep{0.6in}
  \xymatrix{ 
      0 \ar[r]  	& 	\tau(X_2) \ar[r]^(.4){(f_2', g', t')^T} 	&	X_1 \oplus Y \oplus Z \ar[r]^(.6){(f_2, g ,t)}     &   X_2 \ar[r]	&	0.	    }
\end{xy}$$
Due to $\Im(f_1 + \phi, -h) \subset \Ker(f_2, g) \subset \Ker(f_2, g, t) = \Im(f_2', g', t') \cong \tau(X_2)$ there is a natural morphism $h' : X_0 \to \tau(X_2)$ such that the diagram
$$\begin{xy}
\xymatrixcolsep{0.6in}
  \xymatrix{ 
        	& 	 	&	X_0 \ar[d]^{(f_1 + \phi, -h,0)^T} \ar[dl]_{h'}    &   	&		  \\
      0 \ar[r]  	& 	\tau(X_2) \ar[r]^{(f_2', g', t')^T} 	&	X_1 \oplus Y \oplus Z \ar[r]^(.6){(f_2, g ,t)}     &   X_2 \ar[r]	&	0	    }
\end{xy}$$
commutes. This shows the second statement; moreover, we have $f_1 + \phi = f_2'h'$. Suppose now that $f_1 + \phi$ is not irreducible, i.e. $f_1 + \phi \in \rad^2(X_0,X_1)$ and $\overline{f_1} = \overline{\phi}$ in $\rad(X_0,X_1)/\rad^2(X_0,X_1)$. It follows that every decomposition of $\phi$ contains the arrow labeled $f_1$, which contradicts our assumptions. Thus both $f_1 + \phi$ and $f_2'$ are irreducible, so $h'$ is a split monomorphism and hence an isomorphism, because $\tau(X_2)$ is indecomposable.
\bigskip
\\We now assume the claim holds for $n - 1$. Let $\gamma = f_{n-1}\cdots f_1$ and $f_n(\gamma + \phi) = gh$ for some $h : X_0 \to Y$ and $g : Y \to X_n$ such that $(f_n,g): X_{n-1} \oplus Y \to X_n$ is irreducible. Then the composition $(f_n, g)(\gamma + \phi, -h)^T$ is zero. If $(f_n,g)$ is an epimorphism, then $X_n$ is not projective. There is an irreducible morphism $f_n': \tau(X_n) \to X_{n-1}$ and we obtain a morphism $h' : X_0 \to \tau(X_n)$, in the same way as for $n = 2$, such that $\gamma + \phi = f_n'h'$. Hence we have shown that $f$ factors over $\tau(X_n)$.
$$\begin{xy}
\xymatrixcolsep{0.6in}
  \xymatrix{ 
        	& 	 	&	X_0 \ar[d]^{(\gamma + \phi, -h,0)^T} \ar[dl]_{h'}    &   	&		  \\
      0 \ar[r]  	& 	\tau(X_n) \ar[r]^{(f_n', g', t')^T} 	&	X_{n-1} \oplus Y \oplus Z \ar[r]^(.6){(f_n, g ,t)}     &   X_n \ar[r]	&	0	    }
\end{xy}$$
If $X_{n-2} \cong \tau(X_n)$, we are done, if $X_{n-2} \ncong \tau(X_n)$, we take a closer look at $\phi$. By our assumption there is a decomposition of $\phi = v + \sum_{i=k}^m u_i$ that does not contain $\gamma$, where $k$ and $m$ denote the minimal and maximal length of paths in this decomposition respectively. We arbitrarily choose $v_1, v_2, v_3$ such that $v = (f_{n-1}, \phi', f_n')(v_1, v_2, v_3)^T$, where $\phi': M \to X_{n-1}$ is irreducible such that $(f_{n-1}, \phi', f_n')$ is the minimal right almost split morphism of $X_{n-1}$ given by the labeling of the Auslander-Reiten quiver of $A$-mod. Furthermore, we divide the $u_i$ in the same way, namely $\sum_{i=k}^m u_i = (f_{n-1}, \phi', f_n')(\phi_1, \phi_2, \phi_3)^T$. Note that $\phi_1$ clearly has a decomposition not containing $f_{n-2} \cdots f_1$ by construction. This implies that $\phi_1 + v_1$ has a decomposition not containing $f_{n-2} \cdots f_1$ since $v_1 \in \rad_\chi^\infty(X_0,X_{n-2})$. Moreover, we have $\phi = f_{n-1}(\phi_1 + v_1) + \phi'(\phi_2 +v_2) + f_n'(\phi_3 +v_3)$.
\bigskip
\\Inserting this equation into $\gamma + \phi = f_n'h'$ we obtain $\gamma + f_{n-1}(\phi_1 + v_1) = f_n'(h' - \phi_3 - v_3) - \phi'(\phi_2 + v_2)$. Since $(f_{n-1}, \phi', f_n'): X_{n-2} \oplus M \oplus \tau X_n \to X_{n-1}$ is irreducible as a minimal right almost split morphism, we can apply the induction hypothesis on $\gamma + f_{n-1}(\phi_1 + v_1)$. 
$$\begin{xy}
\xymatrixrowsep{0.5in}
\xymatrixcolsep{0.7in}

  \xymatrix@!0{
        X_0 \ar[dr]^{f_1} &				& 	 	&	  &   	&		  \\
		&\ar@{.}[dr]	& &	&   	&		  \\
		&	& \ar[dr]^{f_{n-2}}&	&   	&		  \\
							& 	& 	& X_{n-2} \ar[dr]^{f_{n-1}}	   &   	&		  \\
X_0 \ar[urrr]^{\phi_1 + v_1} \ar[rrr]^{\phi_2 + v_2} \ar[drrr]^{\phi_3 + v_3}	& &	 	& M \ar[r]^{\phi' }	 &X_{n-1} \ar[dr]^{f_n}	&	  \\
	&		& 			&\tau(X_n) \ar[r]^(.6){t'} \ar[ur]^{f_n'} \ar[dr]^{g'} 	&	Z \ar[r]^t     &   X_n 	 \\
X_0 \ar[urrr]^{h'}	& 	&	&						 		&	Y \ar[ur]^g     &    	 }
\end{xy}$$

On the other hand, if $(f_n,g)$ is injective, then $\gamma + \phi$ is zero and we consider a decomposition of $\phi = v + \sum_{i=k}^m u_i$ that does not contain $\gamma$. We can write $v = (f_{n-1}, \phi')(v_1, v_2)^T$, where $(f_{n-1}, \phi') : X_{n-2} \oplus M \to X_{n-1}$ is the minimal right almost split morphism induced by the labeling of the Auslander-Reiten quiver. Furthermore, we split up the the given sum of paths in the same way as before, $\sum_{i=k}^m u_i = (f_{n-1}, \phi')(\phi_1, \phi_2)^T$. By construction there is a decomposition of $\phi_1 + v_1$ not containing $f_{n-2} \cdots f_1$ and we can apply the induction hypothesis on $\gamma + f_{n-1}(\phi_1 + v_1)$ as it factors through $\phi'$. This completes the proof.
\qed
\bigskip
\\The same proof holds, with a minor restriction, for functorially finite subcategories. Therefore, we first need to generalize the definition of sectional paths in the natural way.

\begin{Def}
Let 
$$\begin{xy}
  \xymatrix{ 
      X_0 \ar[r]^{f_1} 	&	X_1 \ar[r]^{f_2}     &   \cdots	 \ar[r]	&	X_{n-1} \ar[r]^{f_n}	&	X_n	    }
\end{xy}$$
be a path in the Auslander-Reiten quiver of $\chi$. The path is called \textbf{sectional} if there is no $i$ such that $X_{i+2} = \tau_{\chi}(X_i)$.
\end{Def}

\begin{Th}\label{secThsub}
Let
 $$\begin{xy}
  \xymatrix{ 
      X_0 \ar[r]^{f_1} 	&	X_1 \ar[r]^{f_2}     &   \cdots	 \ar[r]	&	X_{n-1} \ar[r]^{f_n}	&	X_n	    }
\end{xy}$$
be a sequence of irreducible morphisms in $\chi$ between indecomposable modules such that $X_i$ is not Ext-projective for $i\geq 1$. Let $f = f_n \cdots f_1 + f_n\phi$, where $\phi$ is a morphism from $X_0$ to $X_{n-1}$ such that there exist a decomposition of $\phi$ that does not contain the path $\gamma = f_{n-1} \cdots f_1$. Suppose $f$ factors through a morphism $g: Y \to X_n$  for some module $Y$ in $\chi$ such that $(f_n,g): X_{n-1} \oplus Y \to X_n$ is irreducible $\chi$. Then the sequence is not sectional in $\chi$ and $f$ factors over $\tau_\chi(X_n)$.
\end{Th}

\proof
\\The proof works just as for $A$-mod except that we do not have to show that the $X_i$ are not Ext-projective.
\qed
\bigskip
\\In the following, there are several statements for a functorially finite subcategory $\chi$ and a sectional path 
 $$\begin{xy}
  \xymatrix{ 
      X_0 \ar[r]^{f_1} 	&	X_1 \ar[r]^{f_2}     &   \cdots	 \ar[r]	&	X_{n-1} \ar[r]^{f_n}	&	X_n	    }
\end{xy}$$
in its Auslander-Reiten quiver such that $X_i$ is not Ext-projective for $i \geq 1$. Note that if $\chi$ is $A$-mod, the restriction to non-projective modules is not necessary, because the proofs of these statements are always based on the previous theorems. 

\begin{Th}\label{unique}
 Let
 $$\begin{xy}
  \xymatrix{ 
      X_0 \ar[r]^{f_1} 	&	X_1 \ar[r]^{f_2}     &   \cdots	 \ar[r]	&	X_{n-1} \ar[r]^{f_n}	&	X_n	    }
\end{xy}$$
be a sectional path in $\chi$ such that $X_i$ is not Ext-projective for $i \geq 1$. Let $f$ be a morphism that has a decomposition containing $\gamma = f_n \cdots f_1$. Then every decomposition of $f$ contains $\gamma$.
\end{Th}

\proof
\\Suppose that there is a decomposition of $f$ not containing $\gamma$. We rewrite this decomposition as $f = g_1 + g_2$, where $g_1$ factors through $f_n$ and $g_2$ factors through a morphism $g : Y \to X_n$ such that $(f_n, g) : X_{n-1} \oplus Y \to X_n$ is the minimal right almost split morphism in $\chi$ given by the labeling of the Auslander-Reiten quiver. Let $f' = f - \gamma\varphi$, where $\varphi$ is the isomorphism of the path $\gamma$ in the decomposition of $f$ containing $\gamma$. Consequently, there is a decomposition of $f'$ not containing $\gamma$. Moreover, we decompose this decomposition of $f'$ in the same way as before, i.e. $f' = h_1 + h_2$, where $h_1$ factors through $f_n$ and $h_2$ factors through $g$. So we have $\gamma\varphi + h_1 - g_1 = g_2 - h_2$ and, consequently, $\gamma + (h_1 - g_1)\varphi^{-1} = (g_2 - h_2)\varphi^{-1}$. By construction $\gamma$ is not one of the paths contained in $(h_1 - g_1)\varphi^{-1}$, so $\gamma$ is not sectional by Theorem \ref{secThsub}, 
which contradicts the assumptions. Hence every decomposition of $f$ must contain $\gamma$. \qed
\bigskip
\\Note that different non-sectional paths can occur in distinct decompositions, even the number of paths in a decomposition may depend on the choice of an isomorphism in another path as the following example shows.

\begin{Ex}
Let $A$ be the path algebra of the quiver 
$$\begin{xy}
\xymatrix{
  e_1 \ \  \ar@<2pt>[r]^f   & \ \ e_2 \ \  \ar@<2pt>[l]^g \\
}
\end{xy}$$
with the relation $gfg = 0$. We set $P_1 = Ae_1$ and $P_2 = Ae_2$ and note that $P_2 = \rad(P_1)$. Then the decomposition of the embedding $i : P_2 \to P_1$ depends on the choice of isomorphisms in $P_2$.
\end{Ex}

The Jordan-H\"{o}lder composition series of the non-simple indecomposable modules are
 $$\begin{xy}
\xymatrixrowsep{0in}
\xymatrixcolsep{-0.1in}
  \xymatrix{ 
	    &S_1    	&				&S_2	&				&S_1		\\	
P_1=I_2=\ \ &S_2   	&\ \ \ \ \ \ \ \ \ \ P_2 = \ \	&S_1	&\ \ \ \ \ \ \ \ \ \ I_1 = \ \	&S_2		\\	
	    &S_1    	&				&S_2	&				&S_1		\\
	    &S_2    	&				&	&	
   }
\end{xy}$$
 $$\begin{xy}
\xymatrixrowsep{0in}
\xymatrixcolsep{-0.1in}
  \xymatrix{ 
	    &S_1    	&					&S_2		\\	
I_1/S_1=\ \ &S_2   	&\ \ \ \ \ \ \ \ \ \ P_2/S_2 = \ \	&S_1.	
   }
\end{xy}$$

We label the Auslander-Reiten quiver of $A$ with the natural embedding and projection morphisms. Then there is a decomposition of $i$ just containing the arrow from $P_2$ to $P_1$. Since this arrow is also a sectional path, it is contained in any decomposition. If we choose the isomorphism corresponding to this path to be $e_2 \mapsto e_2 + fg$ instead of the identity, then the decomposition of $i$ contains the arrow from $P_2$ to $P_1$ and the path from $P_2$ to $P_1$ that factors over $S_2$.

$$\begin{xy}
\xymatrixrowsep{0.5in}
\xymatrixcolsep{0.5in}

  \xymatrix@!0{  
& 	&	& [P_1] \ar@*{[red]}[dr]	&	& 	&	&[P_1] \ar[dr]\\
&	&[P_2 \ar[dr]\color{red} \ar[ur]	&	&I_1] \ar[dr] \ar@{.>}[ll] &	&[P_2 \ar[ur] \ar[dr]&	&I_1] \ar@{.>}[ll]	\\
& I_1/S_1 \ar[ur] \ar[dr] &	&P_2/S_2 \ar[ur] \ar[dr] \ar@{.>}[ll]	&	&I_1/S_1 \ar[dr] \ar[ur] \ar@{.>}[ll]&	&P_2/S_2 \ar[ur] \ar@{.>}[ll]	\\
S_2 \ar[ur]& 	&S_1 \ar[ur] \ar@{.>}[ll]&	&S_2 \ar[ur] \ar@{.>}[ll]&	&S_1 \ar@{.>}[ll] \ar[ur]	
 }
\end{xy}$$

\begin{Cor}\label{submodule}
Let  
$$\begin{xy}
   \xymatrix{ 
       X_0 \ar[r]^{f_1} 	&	X_1 \ar[r]^{f_2}     &   \cdots	 \ar[r]	&	X_{n-1} \ar[r]^{f_n}	&	X_n	    }
\end{xy}$$
be a sectional path in the Auslander-Reiten quiver of $\chi$ such that $X_i$ is not Ext-projective for $i \geq 1$. If a morphism $f : X_0 \to X_n$ has a decomposition containing $f_n \cdots f_1$ and $f$ factors over some module $Y$ in $\chi$, then $X_i$ is a direct summand of $Y$ for some $i \in \{0, \ldots, n\}$.
\end{Cor}

\proof
\\Assume $X_i$ is not a direct summand of $Y$, then there is a decomposition of $f$ not containing $f_n \cdots f_1$, which clearly contradicts Theorem \ref{unique}. \qed

\begin{Cor}\label{non-zero}
Let $f: X_0 \to X_n$ be a morphism between indecomposable modules in $\chi$ that has a decomposition containing a sectional path 
 $$\begin{xy}
  \xymatrix{ 
      X_0 \ar[r]^{f_1} 	&	X_1 \ar[r]^{f_2}     &   \cdots	 \ar[r]	&	X_{n-1} \ar[r]^{f_n}	&	X_n	    }
\end{xy}$$
such that $X_i$ is not Ext-projective for $i \geq 1$. Then $f$ is not in $\rad_\chi^{n+1}(X_0,X_n)$; in particular, $f$ is non-zero. 
\end{Cor}

\proof
\\By Theorem \ref{unique} every decomposition of $f$ contains $f_n \cdots f_1$. This path can clearly only be decomposed into at most $n$ radical morphisms and, therefore, $f \notin \rad_\chi^{n+1}(X_0,X_n)$. \qed
\bigskip
\\Note that $f$ can be any sum of distinct paths such that at least one path is sectional.

\begin{Cor}\label{cycle}
If there is a sectional cycle in the Auslander-Reiten quiver of $\chi$, then there is at least one Ext-projective and one Ext-injective module on it.
\end{Cor}

\proof
\\As a reminder, a sectional cycle in an Auslander-Reiten quiver is a sectional path 
$$\begin{xy}
  \xymatrix{ 
      X_0 \ar[r]^{f_1} 	&	X_1 \ar[r]^{f_2}     &   \cdots	 \ar[r]	&	X_{n-1} \ar[r]^(.4){f_n}	&	X_n = X_0     }
\end{xy}$$
from some indecomposable $X_0$ in $\chi$ to itself such that the composition of the path with itself is again sectional. Let $f = f_n \cdots f_1$, then by the lemma of Harada and Sai there is an integer $k$ such that $f^k = 0$. Suppose none of the $X_i$ is Ext-projective, then $f^k \neq 0$ by Corollary \ref{non-zero}, which is a contradiction. If we assume that none of the $X_i$ is Ext-injective, we can apply the same arguments to 
$$\begin{xy}
  \xymatrix{ 
      \tau_\chi^{-1}(X_0) \ar[r]^{f_1} 	&	\tau_\chi^{-1}(X_1) \ar[r]^(.6){f_2}     &   \cdots	 \ar[r]	&	\tau_\chi^{-1}(X_{n-1}) \ar[r]^(.4){f_n}	&	\tau_\chi^{-1}(X_n) = \tau_\chi^{-1}(X_0),     }
\end{xy}$$
which completes the proof.
\qed

\newpage
\section{Sectional paths and subcategories}\label{secsub}

Starting in this section we consider how compositions of irreducible morphisms in $\chi$ behave in a functorially finite resolving subcategory of $\chi$. Recall that $\Omega$ denotes such a subcategory of $\chi$ in the whole dissertation. In order to be able to find relations between $\chi$ and $\Omega$, we have to require that they have the same Ext-projective modules. But clearly in $\Omega$ every Ext-projective module is projective as $\Omega$ contains all projective $A$-modules and is closed under kernels of epimorphisms. From now on let $\chi$ be a functorially finite subcategory such that all Ext-projective modules in $\chi$ are projective. In particular, for all results in this chapter $\chi$ can be specialized to equal $A$-mod. Firstly, we observe another easy consequence from the previous section.

\begin{Cor}\label{extension}
Let $M$ be an indecomposable module in $\chi$. If there is an indecomposable module $X$ in $\Omega$ and there is a sectional path $\gamma$ from $X$ to $M$ in $\Gamma_\chi$ such that no other module besides $X$ on this path is in $\Omega$, then $\gamma$ is an $\Omega$-section and can be extended to a minimal right $\Omega$-approximation of $M$.
\end{Cor}

\proof
\\Let $\gamma$ be given by
$$\begin{xy}
   \xymatrix{ 
       X = X_0 \ar[r]^{f_1} 	&	X_1 \ar[r]^{f_2}     &   \cdots	 \ar[r]	&	X_{n-1} \ar[r]^{f_n}	&	X_n = M.	    }
\end{xy}$$
Moreover, let $Y$ be a module in $\Omega$ and $g: X \to Y$ and $f: Y \to M$ morphisms such that $\gamma = fg$. Since no module on the path is in $\Omega$, the modules $X_i$ are, in particular, not projective for $i \geq 1$. Therefore, every module that occurs in a factorization of $\gamma$ contains a direct summand $X_i$ by Corollary \ref{submodule}. Consequently, $X_0 = X$ is a direct summand of $Y$ and, as every decomposition of $fg$ contains $\gamma$, $g$ is a split monomorphism. \qed

\begin{Th}\label{sectionalirreducible}
Let $X, Y$ be modules in $\Omega$ such that $Y$ is not projective and let
 $$\begin{xy}
  \xymatrix{ 
    X =  X_0 \ar[r]^{f_0} 	&	X_1 \ar[r]^{f_1}     &   \cdots	 \ar[r]^{f_{n-1}}	&	X_{n-1} \ar[r]^{f_n}	&	X_n = Y	    }
\end{xy}$$
be a sectional path in $\chi$ such that $X_i$ is not in $\Omega$ for $i = 1, \ldots ,n-1$ . Then every morphism $f$ whose decompositions contain $\gamma = f_n \cdots f_0$ is irreducible in $\Omega$. Moreover, if we denote all sectional paths from $X$ to $Y$ in $\chi$ such that all modules along these paths are not in $\Omega$ by $\gamma_1, \ldots, \gamma_n$, then their cosets $\{\overline{\gamma_1}, \ldots, \overline{\gamma_n}\}$ in $\Irr_\Omega(X,Y)$ are linearly independent.
\end{Th}

\proof
\\Let $Z$ be a module in $\Omega$ and $g: X \to Z$ and $h: Z \to Y$ morphisms such that $f = hg$. Since no module on the path is in $\Omega$, the modules $X_i$ are, in particular, not projective for $i \geq 1$. Therefore, every module that occurs in a factorization of $f$ contains a direct summand $X_i$ by Corollary \ref{submodule}. Hence $X$ or $Y$ is a direct summand of $Z$ and either $g$ is a split monomorphism or $h$ is a split epimorphism.
\bigskip
\\Suppose $\{\overline{\gamma_1}, \ldots, \overline{\gamma_n}\}$ is a linearly dependent set in $\Irr_\Omega(X,Y)$. This means we can rewrite $\overline{\gamma_1} = \sum_{i=2}^n \overline{\gamma_ih_i}$ for some $h_i \in \End_A(X)$, i.e. $\gamma_1 - \sum_{i=2}^n \gamma_ih_i \in \rad_\Omega^2(X,Y)$. But $\gamma_1 - \sum_{i=2}^n \gamma_ih_i$ is a morphism clearly admitting a decomposition containing $\gamma_1$, hence by Theorem \ref{unique} every decomposition of $\gamma_1 - \sum_{i=2}^n \gamma_ih_i$ contains $\gamma_1$. By the first part of the theorem we conclude that $\gamma_1 - \sum_{i=2}^n \gamma_ih_i$ must be irreducible in $\Omega$, which is a contradiction. Thus the cosets in $\Irr_\Omega(X,Y)$ of sectional paths in $\chi$ such that no module along the path is in $\Omega$ are linearly independent.
\qed
\bigskip
\\Later we see that the converse statement, i.e. every non-sectional path between $X$ and $Y$ gives rise to a reducible morphism in $\Omega$, is not always true, but we provide conditions under which it holds. Note that if $\chi$ equals $A$-mod in the last and most the following theorems, the statements are still true even if $Y$ is projective. This follows from the fact that these theorems are based on Theorem \ref{secTh} and Theorem \ref{secThsub}.

\begin{Lemma}\label{exact}
Given two short exact sequences
$$\begin{xy}
\xymatrixcolsep{0.5in}
  \xymatrix{ 
    0 \ar[r]   &   X \ar[r]^(.4){(f_1,f_2)} 	&	U \oplus Y \ar[r]^(.55){(g_1,g_2)}    &    Z \ar[r]  	&	0	    }
\end{xy}
$$
and
$$
\begin{xy}
\xymatrixcolsep{0.5in}
  \xymatrix{ 
    0 \ar[r]   &   U \ar[r]^(.4){(g_1,g_3)} 	&	Z \oplus V \ar[r]^(.55){(h_1,h_2)}    &    W \ar[r]  	&	0,	    }
\end{xy}$$
then there is a short exact sequence
$$\begin{xy}
\xymatrixcolsep{0.5in}
  \xymatrix{ 
    0 \ar[r]   &   X \ar[r]^(.4){(g_3f_1,f_2)} 	&	V \oplus Y \ar[r]^(.55){(-h_2,h_1g_2)}    &    W \ar[r]  	&	0.	    }
\end{xy}$$
\end{Lemma}

\proof
\\Let $(g_3f_1,f_2)(x) = 0$ for some $x \in X$, then, in particular, $f_2(x) = 0$. We assume that $f_1(x) \neq 0$ and deduce $g_1f_1(x) \neq 0$ from $g_3f_1(x) = 0$ using that $(g_1,g_3)$ is injective.
On the other hand, we know that $g_1f_1(x) = - g_2f_2(x)$, contradicting $f_2(x) = 0$. Thus $f_1(x) = 0$ and, consequently, $x = 0$ by injectivity of $(f_1, f_2)$. Hence $(g_3f_1,f_2)$ is a monomorphism.
\bigskip
\\For some $w \in W$ there is an element $z + v \in Z \oplus V$ such that $h_1(z) + h_2(v) = w$. Repeating the procedure for $z$ we get $g_2(y) + g_1(u) = z$. Using $h_1g_1(u) = - h_2g_3(u)$ we know that $h_1g_2(y) -  h_2(g_3(u) - v) = c$. Therefore, $(-h_2,h_1g_2)$ is surjective.
\bigskip
\\We now prove that $\Im(g_3f_1,f_2) \subset \Ker(-h_2,h_1g_2)$. For an $x \in X$ we get
$h_1g_2f_2(x) - h_2g_3f_1(x) = h_1g_2f_2(x) + h_1g_1f_1(x) = h_1g_2f_2(x) - h_1g_2f_2(x) = 0$
using the exactness of the given sequences once each.
\bigskip
\\In the last step of the proof we verify that in fact $\Im(g_3f_1,f_2) = \Ker(-h_2,h_1g_2)$ using the Jordan-H\"{o}lder length of its modules. We have $l(V) = l(U) + l(W) - l(Z)$ and $l(Y) = l(X) + l(Z) - l(U)$. Thus 
$$ l(V \oplus Y) = l(X) + l(W)$$
and the sequence
$$\begin{xy}
\xymatrixcolsep{0.5in}
  \xymatrix{ 
    0 \ar[r]   &   X \ar[r]^(.4){(g_3f_1,f_2)} 	&	V \oplus Y \ar[r]^(.55){(-h_2,h_1g_2)}    &    W \ar[r]  	&	0	    }
\end{xy}$$
is exact.
\qed
\newpage
The new sequence splits if and only if both of the given sequences split. This is easily seen because $(-h_2,h_1g_2)$ and $(g_3f_1,f_2)$ are radical morphisms if $(h_1,h_2)$ and $(f_1,f_2)$ are radical morphisms respectively.

\begin{Lemma}\label{sequence}
Let 
$$\begin{xy}
   \xymatrix{ 
       X_0 \ar[r]^{f_1} 	&	X_1 \ar[r]^{f_2}     &   \cdots	 \ar[r]	&	X_{n-1} \ar[r]^{f_n}	&	X_n	    }
\end{xy}$$
be a sectional path in $\chi$ such that $X_i$ is not projective for $i = 0, \ldots, n$ and let $\alpha_i$ name the almost split sequences 
$$\begin{xy}
  \xymatrix{ 
    0 \ar[r]   &   \tau_\chi(X_i) \ar[r] 	&	M_i \oplus \tau_\chi(X_{i+1}) \oplus X_{i-1} \ar[r]   &    X_i  \ar[r]	&	0,	    }
\end{xy}$$
where we require that $M = \bigoplus_{i=0}^n M_i \neq 0$ and define $\tau_\chi(X_{n+1})$ and $X_{-1}$ to be the zero module. Then there is a non-split exact sequence
$$\begin{xy}
  \xymatrix{ 
    0 \ar[r]   &   \tau_\chi(X_0) \ar[r] 	&	M \ar[r]^f   &    X_n  \ar[r]	&	0	    }
\end{xy}$$
such that for all modules $X$ in $\chi$ that do not contain a direct summand isomorphic to $X_i$ for $i = 0, \ldots, n$, all morphisms $g : X \to X_n$ factor through $f$.
\end{Lemma}

\proof
\\We prove the lemma by induction on $n$. For $n=0$ both existence and the factorization property follow immediately from the fact that $\alpha_0$ is an almost split sequence in $\chi$. Suppose the statement has been proved for $n-1$, then by induction there is a non-split short exact sequence 
$$\begin{xy}
  \xymatrix{ 
    0 \ar[r]   &   \tau_\chi(X_0) \ar[r] 	&	\bigoplus\limits_{i=0}^{n-1} M_i \oplus \tau_\chi(X_n) \ar[r]^(.65){(h,h')}   &    X_{n-1}  \ar[r]	&	0.	    }
\end{xy}$$
Applying Lemma \ref{exact} to this sequence and $\alpha_n$ we obtain another non-split short exact sequence
$$\begin{xy}
  \xymatrix{ 
    0 \ar[r]   &   \tau_\chi(X_0) \ar[r] 	&	\bigoplus\limits_{i=0}^{n} M_i  \ar[r]^(.6)f   &    X_n  \ar[r]	&	0.	    }
\end{xy}$$
Let $X$ be an indecomposable module in $\chi$ not isomorphic to $X_i$ for $i = 0, \ldots, n$ and let $g : X \to X_n$ be a morphism. Then there is a morphism $(g_n, g_n')^T: X \to X_{n-1} \oplus M_n$ such that $g =  (f_n,f_n')(g_n, g_n')^T$, where $(f_n,f_n'): X_{n-1} \oplus M_n \to X_n$ is a minimal right almost split morphism in $\chi$. By induction there is a factorization $g_n = (h,h')(s,s')^T$ with $(s,s') : X \to \bigoplus_{i=0}^{n-1} M_i \oplus \tau_\chi(X_n)$ and $(h,h') : \bigoplus_{i=0}^{n-1} M_i \oplus \tau_\chi(X_n) \to X_{n-1}$. 
\bigskip
\\If $M_n = 0$, then the composition $f_nh's'$ is zero and by construction $f$ equals the composition $f_n(h,h')$, which gives us $g = f_n(h,h')(s,s')^T =  f(s,s')^T$ as stated. On the other hand, if $M_n \neq 0$, then there is a morphism $f': \tau_\chi(X_n) \to M_n$ such that $f_nh' = f_n'f'$ and we have $g = (f_n,f_n')(g_n, g_n')^T = f_ng_n + f_n'g_n' = f_nhs + f_nh's' + f_n'g_n' = f_nhs + f_n'f's' + f_n'g_n' = f_nhs + f_n'(f's' + g_n') = f(s, f's' + g_n')^T$, where $f = (f_nh, f_n') : M \to X_n$.

$$\begin{xy}
\xymatrixrowsep{0.75in}
\xymatrixcolsep{0.75in}
  \xymatrix@!0{ 
\tau_\chi(X_0) \ar[r] \ar [dr]	&M_0 \ar[r]	&        X_0 \ar[dr]^{f_1}				& 	 	&	  &   	&		  \\
	&\tau_\chi(X_1)	\ar[dr] \ar[ur]&	&\ar@{.}[dr]  &	&   	&		  \\
&	&	\ar@{.}[dr]&	&\ar[dr]^{f_{n-2}}	&	&	&	\\
&	&	& \ar[dr]& 	 	& X_{n-2} \ar[dr]^{f_{n-1}}	   &   	&		  \\
&	&	&	&  \tau_\chi(X_{n-1}) \ar[dr] \ar[r] \ar[ur]	 	& M_{n-1} \ar[r]	   	&X_{n-1} \ar[dr]^{f_n}   	&	 \\
&	&	&	& 	& \tau_\chi(X_n)  \ar[r]^{f'} \ar[ur]^{h'}	&	M_n \ar[r]^{f_n'}     &   X_n 	}
\end{xy}$$
This proves the statement for an indecomposable module $X$ not isomorphic to $X_i$ for $i = 0, \ldots, n$. Let $X$ now be any module in $\chi$ such that none of its direct summands is isomorphic to $X_i$ for $i = 0, \ldots, n$ and let $g : Y \to X_n$ be a morphism. Clearly, all induced morphisms from a direct summand of $Y$ to $X_n$ factor through $f$, which gives rise to a factorization of $g$ through $f$, which completes the proof.
\qed

\begin{Th}\label{irrfactor}
Let $X, Y$ be indecomposable modules in $\Omega$ such that $X$ is not Ext-injective in $\Omega$, $Y$ is not projective and let $f : X \to Y$ be an $\Omega$-irreducible morphism. Then $f$ does not factor over $\t(Y)$ non-trivially, i.e. if there are morphisms $h: X \to \t(Y)$ and $g: \t(Y) \to Y$ such that $f = gh$, then $g : \t(Y) \to Y$ is an isomorphism.
\end{Th}

\proof
\\Suppose there is a factorization $f = gh$ where $h: X \to \t(Y)$ and $g : \t(Y) \to Y$. By definition $h$ factors through the minimal right $\Omega$-approximation $h_{\t(Y)} : X_{\t(Y)} \to \t(Y)$, so we have $f = gh_{\t(Y)}h'$ for some $h': X \to X_{\t(Y)}$. Since $f$ is irreducible, either $gh_{\t(Y)}$ is a split epimorphism or $h'$ is a split monomorphism. But if $gh_{\t(Y)}$ is a split epimorphism, then, in particular, $g$ is an isomorphism and the factorization $f = gh$ is trivial. Suppose now $h'$ is a split monomorphism. It follows that $X \cong \tau_\Omega(Y)$ as it is the only direct summand of $X_{\t(Y)}$ that is not Ext-injective by Theorem \ref{Omegatranslateinchi}. This gives rise to a commutative diagram
$$\begin{xy}
  \xymatrix{ 
    0 \ar[r]	&	X \ar[r]^\varphi \ar[d]^h	&	M' \ar[r]\ar[d]    	&    	Y \ar[r]\ar@{=}[d]  	&	0	\\  
    0 \ar[r]    &	\tau_\chi(Y) \ar[r]^{g'}		&	M  \ar[r]	&	Y \ar[r]	&	0,	 }  
\end{xy}$$
by Lemma \ref{approxdiagram}. Clearly, $g'h$ is an $\Omega$-section by Corollary \ref{extension}, so $\varphi$ must be a split monomorphism, which is a contradiction as by construction the upper row is the almost split sequence in $\Omega$ ending in $Y$.
\qed
\bigskip
\\Now we show that if a connected component $\Gamma$ of $\Gamma_\chi$ is in some sense large, then it does not contain $\Omega$-irreducible morphisms given by non-sectional paths in $\Gamma$.

\begin{Def}
Let $X, Y$ be indecomposable modules in $\Omega$ in a connected component $\Gamma$ of $\Gamma_\chi$. Suppose there is an immediate successor $Z$ of $X$ such that either there is a sectional path
$$\begin{xy}
  \xymatrix{ 
    Z =  Z_0 \ar[r] 	&	Z_1 = \t^{-1}(X) \ar[r]     &   \cdots	 \ar[r]	&	Z_{k-1} \ar[r]	&	Z_{k} = Y	    }
\end{xy}$$ 
from $Z$ to $Y$ such that each $Z_i$ is not Ext-projective for $i = 1, \ldots, k$ and where we define $l = 0$ or there exist an $l > 0$ such that there are two paths 
 $$\begin{xy}
  \xymatrix{ 
    Z =  Z_0 \ar[r] 	&	Z_1 \ar[r]     &   \cdots	 \ar[r]	&	Z_{k+l-1} \ar[r]	&	Z_{k+l} = Y	    }
\end{xy}$$
and
 $$\begin{xy}
  \xymatrix{ 
    Z =  Y_0 \ar[r] 	&	Y_1 \ar[r]     &   \cdots	 \ar[r]	&	Y_{k+l-1} \ar[r]	&	Y_{k+l} = Y	    }
\end{xy}.$$
from $Z$ to $Y$ with the following properties.
\begin{enumerate}
 \item Both paths are not sectional in precisely one module, i.e. we have $\t(Z_{l+1})$ $= Z_{l-1}$ and $\t(Y_{k+1}) = Y_{k-1}$, while the induced paths from $Z$ to $Z_l$, from $Z_l$ to $Y$, from $Z$ to $Y_k$ and from $Y_k$ to $Y$ are all sectional. 
\item $\t(Y_1) = X$, $Z_1 \neq Y_1$ and $Z_{k+l-1} \neq Y_{k+l-1}$.
\item $Y_1, \ldots, Y_k$ are not projective. 
\item If $l \leq k$ we have $Y_i = \t^l(Z_{i+2l})$ for $i = 0, \ldots, k - l$ and $Y_{k-i} = \t^i(Y_{k+i})$, $Z_{l-i} = \t^i(Z_{l+i})$ for $i = 0, \ldots, l$. 
\item If $k \leq l$ we have $Z_i = \t^k(Y_{i+2k})$ for $i = 0, \ldots, l - k$ and $Y_{k-i} = \t^i(Y_{k+i})$, $Z_{l-i} = \t^i(Z_{l+i})$ for $i = 0, \ldots, k$. 
\end{enumerate}
We then say that $\Gamma$ is \textbf{large between $\pmb{X}$ and $\pmb{Y}$}. 
\\Moreover, we call the following the \textbf{inner modules of $\pmb{X}$ and $\pmb{Y}$}:
\begin{enumerate}[(a)]
\item If $l \leq k$, the inner modules are $\t^j(Z_{i+2l})$ for $i = 0, \ldots, k - l$, $j= 0, \ldots, l$, 
\\$\t^j(Y_{k+i})$ for $i = 0, \ldots, l$, $j=0, \ldots, i$ and 
\\$\t^j(Z_{l+i})$ for $i = 1, \ldots, l$, $j=0, \ldots, i-1$. 
\item If $k \leq l$, the inner modules are $\t^j(Y_{i+2k})$ for $i = 0, \ldots, l - k$, $j= 0, \ldots, k$, 
\\$\t^j(Y_{k+i})$ for $i = 0, \ldots, k$, $j=0, \ldots, i$ and 
\\$\t^j(Z_{l+i})$ for $i = 1, \ldots, k$, $j=0, \ldots, i-1$. 
\end{enumerate}
\end{Def}

The following is an example for a large component between $X$ and $Y$ for $k = 5$ and $l = 2$. Note that in general in a large component the number of middle terms for an almost split sequence does not have to be $2$. Inner modules of the given paths are marked in red.

$$\begin{xy}
\xymatrixrowsep{0.5in}
\xymatrixcolsep{0.5in}

  \xymatrix@!0{ 
&	&	&		&	&					& {\color{red}Y_5} \ar[dr] 	&	&	\\
&	&	&	&	&{\color{red}Y_4} \ar[dr] \ar[ur] 	&		&{\color{red}Y_6} \ar[dr] \ar@{.>}[ll]	&	\\
&	&	& 	&{\color{red}Y_3} \ar[dr] \ar[ur] &	&{\color{red}\t(Y)} \ar[dr] \ar[ur] \ar@{.>}[ll]&	&{\color{red}Y}  \ar@{.>}[ll] 	\\ 
&	&	&{\color{red}Y_2} \ar[dr] \ar[ur]&	&{\color{red}\t(Z_6)} \ar[dr] \ar[ur] \ar@{.>}[ll]& &{\color{red}Z_6} \ar[ur] \ar@{.>}[ll]&	\\
X \ar[dr]&	&{\color{red}Y_1} \ar[ur] \ar[dr] \ar@{.>}[ll]&  &{\color{red}\t(Z_5)} \ar[dr] \ar[ur] \ar@{.>}[ll]&	&{\color{red}Z_5} \ar[ur] \ar@{.>}[ll]&	&	\\ 
&Z \ar[dr] \ar[ur]	&	&{\color{red}\t(Z_4)} \ar[ur] \ar[dr] \ar@{.>}[ll]	&	&{\color{red}Z_4} \ar[ur] \ar@{.>}[ll]& 	&	&	\\
& 	&Z_1 \ar[dr] \ar[ur]	&	&{\color{red}Z_3} \ar[ur]\ar@{.>}[ll]	&	& 	&	&	\\
&	&	&Z_2 \ar[ur]	&	&	& 	&	&	
 }
\end{xy}$$

Note that a component $\Gamma$ may be large between $X$ and $Y$ using different paths in the Auslander-Reiten quiver. 

\begin{Ex}
Let $A$ be the path algebra of 
$$\begin{xy}
  \xymatrix{ 
      \stackrel{1}{\cdot} \ar[r] \ar@/_1pc/[rr]  	& 	 \stackrel{2}{\cdot} \ar[r] 	&	 \stackrel{3}{\cdot}	    }
\end{xy}$$
\end{Ex}
The preprojective component of its Auslander-Reiten quiver is illustrated on the next page. It is easy to see that there are $\frac{n}{2}$ or $\frac{n+1}{2}$ different combinations of paths that satisfy the definition of largeness between $P_3$ and $Y_n$ for $n$ even and odd respectively. Furthermore, the inner modules of a large component are not unique either, but they are uniquely determined by the paths from $Z_l$ and $Y_k$ to $Y$. Let us denote these paths by $\gamma_l$ and $\gamma_k$ respectively. We then say that $\Gamma$ is large between $X$ and $Y$ \textbf{with respect to $\pmb{\gamma_l}$ and $\pmb{\gamma_k}$} and refer to inner modules of $\gamma_l$ and $\gamma_k$ if necessary.

$$\begin{xy}
\xymatrixrowsep{0.5in}
\xymatrixcolsep{0.5in}

  \xymatrix@!0{ 
&	&	&		&	&					&  	&	&	\\
&	&	&[P_2 \ar[dr] \ar[ur]	&		&Y_2 \ar[ur] \ar[dr]\ar@{.>}[ll]	&					&  	&	&	\\
&	&[P_3 \ar[dr] \ar[ur]	&	&Y_1 \ar[ur] \ar[dr]\ar@{.>}[ll]	&	&Y_4 \ar[dr] \ar[ur]\ar@{.>}[ll] 	&	&	&	\\
&	&	&[P_1 \ar[ur]\ar[dr]	& 	&Y_3 \ar[dr] \ar[ur]\ar@{.>}[ll] &	&Y_6 \ar[dr] \ar[ur] \ar@{.>}[ll]&	& 	\\ 
&	&[P_2 \ar[ur]\ar[dr]	&	&Y_2 \ar[dr] \ar[ur]\ar@{.>}[ll]&	&Y_5 \ar[dr] \ar[ur] \ar@{.>}[ll]& &Y_8 \ar[ur] \ar[dr] \ar@{.>}[ll]&	\\
&	[P_3 \ar[dr] \ar[ur]&	&Y_1 \ar[ur] \ar[dr] \ar@{.>}[ll]&  &Y_4 \ar[dr] \ar[ur] \ar@{.>}[ll]&	&Y_7 \ar[ur] \ar[dr] \ar@{.>}[ll]&	&	\\ 
&	&[P_1 \ar[dr] \ar[ur]	&	&Y_3 \ar[ur] \ar[dr] \ar@{.>}[ll]	&	&Y_6 \ar[ur] \ar[dr] \ar@{.>}[ll]& 	&	&	\\
&[P_1 \ar[dr] \ar[ur]	& 	&Y_2 \ar[dr] \ar[ur]\ar@{.>}[ll]	&	&Y_5 \ar[ur] \ar[dr] \ar@{.>}[ll]	&	& 	&	&	\\
[P_2 \ar[dr] \ar[ur]&	&Y_1 \ar[ur] \ar[dr] \ar@{.>}[ll]	&	&Y_4 \ar[ur]\ar[dr]\ar@{.>}[ll]	&	&	& 	&	& \\
&	&	&		&	&					&  	&	&	
 }
\end{xy}$$

\begin{Lemma}\label{large}
Let $X$ and $Y$ be indecomposable modules in $\chi$ such that their connected component $\Gamma$ is large between them. Then there is a non-split short exact sequence 
$$\begin{xy}
  \xymatrix{ 
    0 \ar[r]   &   X \ar[r] 	&	M \ar[r]    &    Y \ar[r]  	&	0	    }
\end{xy}$$
and every morphism $f : N \to Y$ such that $N$ does not contain a direct summand isomorphic to an inner module of $X$ and $Y$ factors over $M$.
\end{Lemma}
\proof
\\We prove the lemma by induction on $l$ as defined for large components. For $l = 0$ we are precisely in the situation of Lemma \ref{sequence} and there is nothing left to prove. Suppose now that $l > 0$ and both statements are true for $l-1$. Since $l > 0$, there are paths 
 $$\begin{xy}
  \xymatrix{ 
  \gamma_l: &  Z =  Z_0 \ar[r] 	&	Z_1 \ar[r]     &   \cdots	 \ar[r]	&	Z_{k+l-1} \ar[r]	&	Z_{k+l} = Y,	    }
\end{xy}$$
 $$\begin{xy}
  \xymatrix{ 
   \gamma_k: & Z =  Y_0 \ar[r] 	&	Y_1 \ar[r]     &   \cdots	 \ar[r]	&	Y_{k+l-1} \ar[r]	&	Y_{k+l} = Y	    }
\end{xy}$$
such that $\Gamma$ is large between $X$ and $Y$ with respect to $\gamma_l$ and $\gamma_k$. By Theorem \ref{sequence} there are $l$ short exact sequences of the form
$$\begin{xy}
  \xymatrix{ 
    0 \ar[r]   &   Z_i \ar[r] 	&	Z_{i+1} \oplus Y_{k+i} \oplus M_i \ar[r]   &    Y_{k+1 + i} \ar[r]  	&	0	    }
\end{xy}$$
for $i = 0, \ldots, l-1$ 
and one short exact sequence 
$$\begin{xy}
  \xymatrix{ 
    0 \ar[r]   &   X \ar[r] 	&	Z_0 \oplus M_{-1} \ar[r]^(.6){(f_0,f_{-1})}    &    Y_k \ar[r]  	&	0	    }
\end{xy}.$$
The connected component of $X$ and $Y$ is also large between $Z$ and $Y$ either since
$$\begin{xy}
  \xymatrix{ 
  Z_1 \ar[r]     &   \cdots	 \ar[r]	&	Z_{k+l-1} \ar[r]	&	Z_{k+l} = Y	    }
\end{xy}$$
is a sectional path for $l = 1$ or with respect to the paths 
$$\begin{xy}
  \xymatrix{ 
  Z_l \ar[r]     &   \cdots	 \ar[r]	&	Z_{k+l-1} \ar[r]	&	Z_{k+l} = Y	    }
\end{xy}$$
and
$$\begin{xy}
  \xymatrix{ 
Y_{k+1} \ar[r] &  Y_{k+2} \ar[r]	& \cdots	 \ar[r]	&	Y_{k+l-1} \ar[r]	&	Y_{k+l} = Y	    }
\end{xy}$$ 
for $l \geq 2$. So by inductive construction there is a short exact sequence
$$\begin{xy}
  \xymatrix{ 
    0 \ar[r]   &  Z \ar[r] 	&	Y_k \oplus Z_l \oplus \bigoplus\limits_{i=0}^{l-1} M_i \ar[r]^(.72){(g_k,g)}    &    Y \ar[r]  	&	0.	    }
\end{xy}$$
with $g_k : Y_k \to Y$ and $g : Z_l \oplus \bigoplus_{i=0}^{l-1} M_i \to   Y$. We apply Lemma \ref{exact} and obtain a short exact sequence 
$$\begin{xy}
  \xymatrix{ 
    0 \ar[r]   &  X \ar[r] 	&	Z_l \oplus \bigoplus\limits_{i=-1}^{l-1} M_i \ar[r]   &    Y \ar[r]  	&	0,	    }
\end{xy}$$
where we define $M = Z_l \oplus \bigoplus\limits_{i=-1}^{l-1} M_i$. Let $N$ be a module such that none of its direct summands is isomorphic to an inner module of $X$ and $Y$, in particular, none of its direct summands is isomorphic to an inner module of $Z$ and $Y$. Then, by induction, every morphism $f : N \to Y$ has a factorization $f = (g_k,g)(g'_k,g')^T$ for some $g'_k : N \to Y_k$ and $g' : N \to Z_l \oplus \bigoplus_{i=0}^{l-1} M_i$. By Lemma \ref{sequence} $g'_k$ factors through $(f_0,f_{-1}) : Z_0 \oplus M_{-1} \to Y_k$ since $Y_1, \ldots, Y_k$ are inner modules of $X$ and $Y$, i.e. there is an $(f_0',f_{-1}')^T: N \to Z_0 \oplus M_{-1}$ such that $g'_k = (f_0,f_{-1})(f_0',f_{-1}')^T$. Moreover, by construction $g_kf_0 : Z_0 \to Y$ factors through $g$, i.e. $g_kf_0 = gh$ for some $h : Z_0 \to Z_l \oplus \bigoplus_{i=0}^{l-1} M_i$. Composing these morphisms we have a factorization $f = g_kf_{-1}f_{-1}' + g(g' + hf_0')$, where $(g_kf_{-1},g)$ is a morphism mapping from $M$ to $Y$ which completes the proof. \qed
\bigskip
\\For the following theorem we define the distance in $\Gamma\chi$ from an indecomposable module $X$ in $\chi$ to another indecomposable module $Y$ in $\chi$. Suppose there are non-negative integers $l$ and $n$ with $n + l$ minimal such that there is a sectional path 
 $$\begin{xy}
  \xymatrix{ 
\t^{-l}(X) = Y_0 \ar[r] 	&	Y_1 \ar[r]     &   \cdots \ar[r] 	&	Y_{n-1} \ar[r]	&	Y_n = Y	    }
\end{xy}$$
Then $n+l$ is the \textbf{distance from $\pmb{X}$ to $\pmb{Y}$} in the Auslander-Reiten quiver of $\chi$.

\begin{Th}\label{largetheorem}
Let $X$ and $Y$ be indecomposable modules in $\Omega$ that are in the same connected component of $\Gamma_\chi$ and let $l$ be a positive and $n$ a non-negative integer with $l + n$ minimal such that there is a sectional path
 $$\begin{xy}
  \xymatrix{ 
  \gamma:	& \t^{-l}(X) = Y_0 \ar[r] 	&	Y_1 \ar[r]     &   \cdots \ar[r] 	&	Y_{n-1} \ar[r]	&	Y_n = Y.	    }
\end{xy}$$
Moreover, we assume that $l$ is also minimal, i.e. for all non-negative integers $l'$ and $n'$ with $l + n = l' + n'$ such that there is a sectional path
 $$\begin{xy}
  \xymatrix{ 
 \t^{-l'}(X) = Y_0' \ar[r] 	&	Y_1' \ar[r]     &   \cdots \ar[r] 	&	Y_{n'-1}' \ar[r]	&	Y_{n'}' = Y	    }
\end{xy}$$
we have $l' \geq l$. Suppose there are sectional paths
 $$\begin{xy}
  \xymatrix{ 
    \gamma_1:	&   Y_{-t} \ar[r] & Y_{-t+1} \ar[r] &\cdots \ar[r] 	&	 Y_{-1} \ar[r]	&	Y_0 = \t^{-l}(X),   }
\end{xy}$$
$$\begin{xy}
  \xymatrix{ 
   \gamma_2:	&   X_{-s} \ar[r] & X_{-s+1} \ar[r] &\cdots \ar[r] 	&	 X_{-1} \ar[r]	&	X_0 = Y    }
\end{xy}$$
 such that
\begin{enumerate}
 \item $\gamma \gamma_1$ is sectional,
 \item $X_{-1}$ is not isomorphic to $Y_{n-1}$,
 \item Either $t = l$ and $s = l - 1$ or $s = l$ and $t = l - 1$ and
 \item $X_{j}$ is not projective for $j = 1-l, \ldots, 0$ and $Y_{j}$ is not projective for $j = 1-l, \ldots, n$.
\end{enumerate}
Then there is no $\Omega$-irreducible morphism from $X$ to $Y$.
\end{Th}

\proof
\\We consider the following diagram for the case that $s = l$ to visualize the assumptions.
$$\begin{xy}
\xymatrixrowsep{0.38in}
\xymatrixcolsep{0.4in}

  \xymatrix@!0{ 
&	&	&	&	&	&X_{-l} \ar[dr]	&	&					&  	&	&	\\
&	&	&	&	&	&	&X_{-l+1} \ar[dr]\ar@{.>}[ll]	&	&					&  	&	&	\\
&	&	&	&	&	&	&	& \ar@{.}[dr]	&					&  	&	&	\\
&	&	&	&	&	&	&	&	& \ar[dr] 	&	&	&	\\
&	&	&	&	&	&	& 	&	&	&X_{-1} \ar[dr]  \ar@{.>}[ll]&	& 	\\ 
&	&	&	&	&	&	&	&	& 	&	&Y \ar@{.>}[ll]&	\\
&	&	&	&	&	&	&	&	&	&Y_{n-1} \ar[ur]  \ar@{.>}[ll]&	&	\\ 
&	&	&	&	&	&	&	&	& \ar[ur] & 	&	&	\\
&	&	&	&	& 	&	&	&\ar@{.}[ur] &	& 	&	&	\\
&	&	&	&	&	&	&Y_1 \ar[ur]\ar@{.>}[ll]	&	&	& 	&	& \\
&	&	&	&	&	&Y_0 = \t^{-l}(X) \ar[ur]\ar@{.>}[ll]	&	&	\\	
&	&	&	&	&Y_{-1} \ar[ur]\ar@{.>}[ll]	&	&	& 	&	& \\
& 	&	&	&\ar[ur]&	&	& 	&	& \\
&	&	&\ar@{.}[ur]	&	&	& 	&	& \\
&	&Y_{-l+1} \ar[ur]\ar@{.>}[ll]	&	&	&	& 	&	& 
}
\end{xy}$$
Clearly, $X$ and $Y$ are in the same connected component $\Gamma$ of $\Gamma_\chi$ as there is a sectional path from $\t^{-l}(X)$ to $Y$. If $\Gamma$ is large between $X$ and $Y$ with respect to $\gamma \gamma_1$ and $\gamma_2$, then by Lemma \ref{large} there is a module $M$ in $\Omega$ such that all morphisms $Z \to Y$ factor over $M$ if no inner module of $X$ and $Y$ is isomorphic to a direct summand of $Z$. By minimality of $n + l$ we know that $X$ cannot be an inner module as the distance from an inner module of $\gamma \gamma_1$ and $\gamma_2$ to $Y$ is always strictly smaller than $n + l$ by definition.
\bigskip
\\Suppose $N$ is an indecomposable direct summand of $M$ such that there is a sectional path whose distance to $Y$ is given by $l' + n'$. By construction in Lemma \ref{large} we have $n' + l' \leq n + l$ and $l' < l$. Since we have chosen $n + l$ and $l$ to be minimal, $X$ is not a direct summand of $M$.
\bigskip
\\Suppose now $\Gamma$ is not large between $X$ and $Y$ with respect to $\gamma \gamma_1$ and $\gamma_2$, i.e. for some $l' \geq 1$ and some module $Y' = Y_i$ or $Y' = X_j$ we have $\t^{l'}(Y') = P$ is projective. Let $n'$ denote the length of the implicit path from $Y'$ to $Y$ given by either $\gamma \gamma_1$ or $\gamma_2$. Without loss of generality, we can assume that $P$ is the projective module obtained in this way such that $n' + l'$ is minimal. Clearly, $n' + l' < l+n$, so by Lemma \ref{large} there is a module $M$ in $\Omega$ such that all morphisms $Z \to Y$ factor over $M$ if no inner module of $P$ and $Y$ is isomorphic to a direct summand of $M$. Since we have chosen $n + l$ to be minimal, $X$ can neither be an inner module of $P$ and $Y$ nor a direct summand of $M$.
\bigskip
\\In both cases $M$ is in $\Omega$ as it is closed under extensions and all morphisms $f: X \to Y$ have a factorization $f = gh$ where $h : X \to M$ and $g: M \to Y$ are radical morphisms. Hence the existence of an irreducible morphism $f : X \to Y$ is impossible. \qed

\begin{Th}\label{sec}
Let $X, Y$ be indecomposable modules in $\Omega$ such that $X$ is not Ext-injective in $\Omega$, $Y$ is not projective and $f : X \to Y$ an $\Omega$-irreducible morphism given by a path
 $$\begin{xy}
  \xymatrix{ 
 \gamma: &   X =  X_0 \ar[r]^{f_1} 	&	X_1 \ar[r]^{f_2}     &   \cdots	 \ar[r]^{f_{n-1}}	&	X_{n-1} \ar[r]^{f_n}	&	X_n = Y	    }
\end{xy}.$$
Then $\gamma$ is sectional if one of the following conditions hold.
\begin{enumerate}[(a)]
 \item Every $X_i$ has at most $2$ predecessors.
 \item Every domain of each right $\Omega$-approximation of $\t(X_i)$ has at most one direct summand that is not Ext-injective in $\Omega$.
\end{enumerate}
\end{Th}

\proof
\\In the whole proof we assume that $\gamma$ is not sectional and find various contradictions. Therefore, let $j \leq n - 2$ be the largest integer such that $X_j = \t(X_{j+2})$. Firstly, suppose that every $X_i$ has at most $2$ predecessors. Clearly, we have a sectional path from $X_{j+2}$ to $X_{n-1}$ such that all modules on this path have exactly $2$ predecessors. By Theorem \ref{sequence} there is a short exact sequence 
$$\begin{xy}
  \xymatrix{ 
    0 \ar[r]   &   X_j \ar[r] 	&	X_{j+1} \oplus \t(Y) \ar[r]    &    X_{n-1} \ar[r]  	&	0	    }
\end{xy}$$
where the induced morphisms are given by $f_{j+1} : X_j \to X_{j+1}$ and $f_{n-1} \cdots f_{j+2}: X_{j+1} \to X_{n-1}$. Hence $f$ factors over $\t(Y)$ and, consequently, cannot be irreducible by Theorem \ref{irrfactor}.
\bigskip
\\Let us now suppose every domain of each right $\Omega$-approximation of $\t(X_i)$ has at most one direct summand that is not Ext-injective in $\Omega$, then, in particular, the domain of the approximation $f_{X_j} : X_{X_j} \to X_j$ has only one direct summand that is not Ext-injective. Since $f_{j} \cdots f_1 : X \to X_j$ is an $\Omega$-section, we know that $X$ must be the aforementioned only direct summand of $X_{X_j}$ that is not Ext-injective in $\Omega$ and, without loss of generality, we can assume that $f_{X_j}$ is an extension of $f_{j} \cdots f_1$ to a minimal right $\Omega$-approximation. By Theorem \ref{sequence} there is a non-split short exact sequence 
 $$\begin{xy}
  \xymatrix{ 
    0 \ar[r]   &	X_j \ar[r]		&	X_{j+1} \oplus M  \ar[r]		&	Y \ar[r]	&	0	 }
\end{xy}$$
such that the induced morphism from $X_j$ to $X_{j+1}$ is given by $f_{j+1}$. By Lemma \ref{approxdiagram} there is a commutative diagram
 $$\begin{xy}
  \xymatrix{ 
    0 \ar[r]   &   	X_{X_j}  \ar[r]^(.35)h\ar[d]_{f_{X_j}} 	&X_{X_{j+1}} \oplus X_M \ar[r]\ar[d]^g 	&Y \ar[r]\ar@{=}[d]  	&	0	\\  
    0 \ar[r]   &	X_j \ar[r]		&	X_{j+1} \oplus M  \ar[r]		&	Y \ar[r]	&	0	 }
\end{xy}$$
such that that the upper row does not split. Since $X$ is the only direct summand of $X_{X_j}$ that is not Ext-injective, it follows that $h|_X$ cannot be a split monomorphism. On the other hand, by construction we have that $gh|_X$ equals $f_{j+1} \cdots f_1$, which is an $\Omega$-section. Since $\Omega$ is closed under extensions, $M$ is in $\Omega$. Therefore, $h|_X$ must be a split monomorphism, which is a contradiction.
\qed


\begin{Cor}\label{seccontra}
Let $X, Y$ be indecomposable modules in $\Omega$ such that $X$ is not Ext-injective in $\Omega$ and $Y$ is not projective. Moreover, we suppose that $X$ and $Y$ are in the same connected component $\Gamma$ of $\Gamma_\Omega$ and $f : X \to Y$ is an $\Omega$-irreducible morphism given by a path
 $$\begin{xy}
  \xymatrix{ 
  \gamma:&  X =  X_0 \ar[r]^{f_1} 	&	X_1 \ar[r]^{f_2}     &   \cdots	 \ar[r]^{f_{n-1}}	&	X_{n-1} \ar[r]^{f_n}	&	X_n = Y	    }
\end{xy}.$$
If $\gamma$ is not sectional, then the following statements hold.
\begin{enumerate}[(a)]
 \item There is an $X_i$ with at least $3$ predecessors.
 \item The domain of the $\Omega$-approximation of the greatest $j$ such that $\t(X_{j+2}) = X_j$ has at least two direct summands that are not Ext-injective in $\Omega$.
 \item $\Gamma$ is not large between $X$ and $Y$.
\end{enumerate}

\end{Cor}
\proof
\\The statements are just the contrapositives of Theorem \ref{sec} and Theorem \ref{largetheorem}.
\qed

\begin{Cor}
Let $f : X \to Y$ be an $\Omega$-irreducible morphism between indecomposable modules, where $X$ is not Ext-injective whereas $Y$ is not projective. If $f = v + \sum u_i$ denotes a decomposition such that $v \in \rad_\chi^\infty(X,Y)$ and $v \in \rad_\Omega^2(X,Y)$, then $X$ and $Y$ are in the same connected component $\Gamma$ of $\Gamma_\chi$ and the decomposition contains a path
 $$\begin{xy}
  \xymatrix{ 
  \gamma:&  X =  X_0 \ar[r]^{f_1} 	&	X_1 \ar[r]^{f_2}     &   \cdots	 \ar[r]^{f_{n-1}}	&	X_{n-1} \ar[r]^{f_n}	&	X_n = Y	    }
\end{xy}$$
such that $f_n \cdots f_1$ is irreducible. Moreover, $\gamma$ is either sectional or it is non-sectional and the following holds
\begin{enumerate}[(a)]
 \item There is an $X_i$ with at least $3$ predecessors.
 \item The domain of the $\Omega$-approximation of the greatest $j$ such that $\t(X_{j+2}) = X_j$ has at least two direct summands that are not Ext-injective in $\Omega$.
 \item $\Gamma$ is not large between $X$ and $Y$.
\end{enumerate}

\end{Cor}

\proof
\\Since $f$ is irreducible by assumption and $v \in \rad_\Omega^2(X,Y)$, we conclude that $\sum u_i \notin \rad_\Omega^2(X,Y)$. In particular, it must contain an individual path 
 $$\begin{xy}
  \xymatrix{ 
  \gamma:&  X =  X_0 \ar[r]^{f_1} 	&	X_1 \ar[r]^{f_2}     &   \cdots	 \ar[r]^{f_{n-1}}	&	X_{n-1} \ar[r]^{f_n}	&	X_n = Y	    }
\end{xy}$$
such that $f_n \cdots f_1 \notin \rad_\Omega^2(X,Y)$, hence $f_n \cdots f_1$ is $\Omega$-irreducible. All other statements follow immediately from Corollary \ref{seccontra}.
\qed
\bigskip
\\The previous corollary suggests a way of labeling arrows in a functorially finite resolving subcategory $\Omega$. We first label the Auslander-Reiten quiver of $\chi$ as described in Section \ref{decomposition}. Let $X$ and $Y$ be modules in $\Omega$ such that $Y$ is not projective and there is a sectional path
 $$\begin{xy}
  \xymatrix{ 
      X = X_0 \ar[r]^{f_1} 	&	X_1 \ar[r]^{f_2}     &   \cdots	 \ar[r]	&	X_{n-1} \ar[r]^{f_n}	&	X_n = Y	    }
\end{xy}$$
in $\chi$ such that no module along this path is in $\Omega$. Then we label an arrow from $X$ to $Y$ with $f_n \cdots f_1$. Doing this for all modules and sectional paths we have labeled all arrows in $\Omega$ but those who either end in projective modules and those that have to be labeled either with a morphism in $\rad_\chi^\infty(X,Y)$ or a morphism whose decomposition does not contain a sectional path. We fix this induced labeling for the rest of the chapter. 

\begin{Th}
Let $f : X \to Y$ be a morphism between indecomposable modules in $\Omega$, that is given by a path
$$\begin{xy}
  \xymatrix{ 
     \gamma_\chi: & X = Y_0 \ar[r]^{g_1} 	&	Y_1 \ar[r]^{g_2}     &   \cdots	 \ar[r]	&	Y_{m-1} \ar[r]^{g_m}	&	Y_m = Y	    }
\end{xy}$$ 
in $\chi$ and a path
 $$\begin{xy}
  \xymatrix{ 
       \gamma_\Omega: &     X = X_0 \ar[r]^{f_1} 	&	X_1 \ar[r]^{f_2}     &   \cdots	 \ar[r]	&	X_{n-1} \ar[r]^{f_n}	&	X_n = Y	    }
\end{xy}$$
in $\Omega$.
\begin{enumerate}[(a)]
 \item If $\gamma_\chi$ is sectional and all $Y_i$ are non-projective for $i \geq 1$, then $\gamma_\Omega$ is also sectional. 
 \item Suppose $\gamma_\Omega$ is sectional and all $X_i$ are neither projective for $i \geq 1$ nor Ext-injective in $\Omega$ for $i < n$. Furthermore, we assume that every module in the connected component of $\Gamma_\chi$ containing $X$ and $Y$ has at most 2 immediate predecessors or the domain of every right $\Omega$-approximation of $Y_i$ has at most one direct summand that is not Ext-injective in $\Omega$. Then $\gamma_\chi$ is sectional. 
\end{enumerate}
\end{Th}

\proof
\\Let $\{i_0, i_1, \ldots, i_n\} \subset \{0, \ldots, m\}$ be the set indexing all modules on $\gamma_\chi$ in $\Omega$. Note that in both in $(a)$ and $(b)$ these modules are exactly all modules on $\gamma_\Omega$ by Theorem \ref{unique}.
\bigskip
\\In order to prove the first statement, let $\gamma_\chi$ be sectional and suppose  
 $$\begin{xy}
  \xymatrix{ 
      X = X_0 \ar[r]^{f_1} 	&	X_1 \ar[r]^{f_2}     &   \cdots	 \ar[r]	&	X_{n-1} \ar[r]^{f_n}	&	X_n = Y	    }
\end{xy}$$
is not sectional, then there is a $j$ such that $X_j = \tau_\Omega (X_{j+2})$. Let $M \oplus X_{j+1}$ denote the whole middle term of the almost split sequence
$$\begin{xy}
   \xymatrix{ 
     0 \ar[r]   &   X_j \ar[r] 	&	M \oplus X_{j+1} \ar[r]    &    X_{j+2} \ar[r]  	&	0.	    }
 \end{xy}
 $$ 
The module $M$ is in $\Omega$ as it is closed under extensions and, consequently, $M$ is also a module in $\chi$. Note that $M \neq 0$ since $f$ is non-zero by Corollary \ref{non-zero} and, in particular, the morphism from $X_j$ to $X_{j+2}$ induced by $f$ is non-zero. Therefore, $f$ factors over $M$, so there is a decomposition of $f$ not containing $\gamma_\chi$, which contradicts Theorem \ref{unique}. Hence
$$\begin{xy}
  \xymatrix{ 
      X = X_0 \ar[r]^{f_1} 	&	X_1 \ar[r]^{f_2}     &   \cdots	 \ar[r]	&	X_{n-1} \ar[r]^{f_n}	&	X_n = Y	    }
\end{xy}$$
is sectional, which proves $(a)$.
\bigskip
\\By the assumptions that we have made in $(b)$ it is clear that we can apply Theorem \ref{sec} to each path from $Y_{i_j} = X_j$ to $Y_{i_{j+1}} = X_{j+1}$. Consequently, these paths must be sectional and, if we suppose $\gamma_\chi$ to be non-sectional, there is a $Y_{i_j}$ in $\Omega$ such that $\t(Y_{i_j+1}) = Y_{i_j-1}$. Let us first suppose that the domain of the right $\Omega$-approximation of each $Y_i$ has only one direct summand that is not Ext-injective. Then $X_{j-1}$ is the only direct summand of $X_{Y_{i_j-1}}$ that is not Ext-injective and we can rewrite $X_{Y_{i_j-1}} = X_{j-1} \oplus I$. By Theorem \ref{sequence} there is a non-split short exact sequence 
 $$\begin{xy}
  \xymatrix{ 
    0 \ar[r]   &	Y_{i_j-1} \ar[r]		&	X_j \oplus M  \ar[r]		&	X_{j+1} \ar[r]	&	0	 }
\end{xy}$$
such that the induced morphism from $X_j$ to $X_{j+1}$ is given by $f_{j+1}$. By Lemma \ref{approxdiagram} there is a commutative diagram
 $$\begin{xy}
 \xymatrixcolsep{0.5in}
  \xymatrix{ 
    0 \ar[r]   &X_{j-1} \oplus I  \ar[r]^{\left( \begin{smallmatrix} f_j&\ \\ \  & \ \end{smallmatrix} \right)}\ar[d]_{f_{Y_{i_j-1}}} 	&X_j \oplus X_M \ar[r] \ar[d]    	&X_{j+1} \ar[r]\ar@{=}[d]  	&	0	\\  
    0 \ar[r]   &	Y_{i_j-1} \ar[r]		&	X_j \oplus M  \ar[r]^{(f_{j+1}, \ )}		&	X_{j+1} \ar[r]	&	0.	 }
\end{xy}$$
Since $X_j$ is in $\Omega$, its minimal right $\Omega$-approximation must be an isomorphism. Hence, without loss of generality, we can assume that the same irreducible morphism $f_{j+1} : X_j \to X_{j+1}$ also occurs in the top row after reducing it to 
 $$\begin{xy}
 \xymatrixcolsep{0.5in}
  \xymatrix{ 
    0 \ar[r]   &X_{j-1}  \ar[r]^{(f_j, \ )^T} 	&X_j \oplus N \ar[r]^{(f_{j+1}, \ )}   	&X_{j+1} \ar[r]  	&	0.	 }
\end{xy}$$
The morphism $(f_{j+1}, \ ): X_j \oplus N \to X_{j+1}$ factors through a minimal right almost split morphism which we obtain by an extension of $f_{j+1}$, so we can extend the diagram by the almost split sequence ending in $X_{j+1}$.
 $$\begin{xy}
\xymatrixcolsep{0.5in}
  \xymatrix{ 
0 \ar[r]   &   	\tau_\Omega(X_{j+1})  \ar[r]^{(g, \ )^T} 	&	X_j \oplus M' \ar[r]^{(f_{j+1}, \ )}   &    	X_{j+1} \ar[r]\ar@{=}[d]  	& 0\\	
0 \ar[r]   &   	X_{j-1} \ar[u]^f  \ar[r]^(.4){(f_j, \ )^T} &  X_j \oplus N \ar[u]^{\left( \begin{smallmatrix} \Id & 0 \\ \  & \ \end{smallmatrix} \right)} \ar[r]^{(f_{j+1}, \ )}      	&    	X_{j+1} \ar[r] 	&	0	 }
\end{xy}$$
The irreducible morphism $f_j$ factors over the indecomposable module $\tau_\Omega(X_{j+1})$, it follows that either $f$ or $g$ is an isomorphism. In the latter case the top row would split, which is a contradiction as it is an almost split sequence. In the former we have $X_{j-1} \cong \tau_\Omega(X_{j+1})$, which is a contradiction to the assumptions made in $(b)$. Thus $\gamma_\chi$ must be a sectional path.
\bigskip
\\Let us now suppose that every module in the connected component of $\Gamma_\chi$ containing $X$ and $Y$ has at most 2 immediate predecessors. Then, as the composition $f_{j+1}f_j = g_{i_{j+1}} \cdots g_{i_{j-1}+1}$ is non-zero, the component in $\chi$ between $Y_{i_j-1} = X_{j-1}$ and $Y_{i_j+1} = X_{j+1}$ is large and there is a short exact sequence 
 $$\begin{xy}
  \xymatrix{ 
    0 \ar[r]   &	X_{j-1} \ar[r]		&	X_j \oplus M  \ar[r]		&	X_{j+1} \ar[r]	&	0	 }
\end{xy}$$
by Lemma \ref{large} where $M$ is an indecomposable module. By our assumptions this is not an almost split sequence because $\gamma_\Omega$ is sectional. By construction the morphism $f: M \to X_{j+1}$ is given by a sectional path which factors through an irreducible morphism $f': M' \to X_{j+1}$, where $M'$ must be a module on the sectional path from $M$ to $X_{j+1}$ by Theorem \ref{unique}. But then $(f',f_{j+1}) : M' \oplus X_j \to X_{j+1}$ is a minimal right almost split morphism in $\Omega$ whose kernel is an indecomposable module isomorphic to some $Y_k$, $i_{j-1} < k < i_j$. This contradicts $f_j$ to be irreducible since $\Omega$ is closed under kernels of surjective morphisms and $f_j$ factors over $Y_k$ as the following example diagram shows.

$$\begin{xy}
\xymatrixrowsep{0.5in}
\xymatrixcolsep{0.5in}

  \xymatrix@!0{  
&	&	&		&	&M \ar[dr]			&  	&	&	\\
&	&	&		&\ar[ur] \ar[dr]	&			& M'\ar[dr]\ar@{.>}[ll] 	&	&	\\
&	&	& \ar[ur] \ar[dr]	&	& \ar[dr] \ar[ur]\ar@{.>}[ll] 	&	& \ar[dr]\ar@{.>}[ll]	&	\\
&	& \ar[ur]\ar[dr]	& 	& \ar[dr] \ar[ur]\ar@{.>}[ll] &	& \ar[dr] \ar[ur] \ar@{.>}[ll]&	& X_{j+1}\ar@{.>}[ll]	\\ 
& \ar[ur]\ar[dr]	&	& \ar[dr] \ar[ur]\ar@{.>}[ll]&	& \ar[dr] \ar[ur] \ar@{.>}[ll]& & \ar[ur] \ar@{.>}[ll]&	\\
X_{j-1} \ar[dr] \ar[ur]&	& \ar[ur] \ar[dr] \ar@{.>}[ll]&  & \ar[dr] \ar[ur] \ar@{.>}[ll]&	& \ar[ur]  \ar@{.>}[ll]&	&	\\ 
& Y_k\ar[dr] \ar[ur]	&	& \ar[ur] \ar[dr] \ar@{.>}[ll]	&	& \ar[ur]  \ar@{.>}[ll]& 	&	&	\\
& 	& \ar[dr] \ar[ur]	&	& \ar[ur]  \ar@{.>}[ll]	&	& 	&	&	\\
&	&	&X_j \ar[ur]	&	&	& 	&	&	
 }
\end{xy}$$

\section{Applications and Examples}\label{examples}

The main purpose of the theorems in Section \ref{secsub} is to derive $\Gamma_\Omega$ from the Auslander-Reiten quiver of a larger category such as $A$-mod. If the larger quiver is finite, the procedure always works and one can get subquivers without computation of approximations. However, in infinite cases the theorems are not always applicable. Certain subcategories of the standard example show that the restriction to not Ext-injective and non-projective modules is necessary.
\bigskip
\begin{Ex}\label{standardexamplesubcategories}
We derive the Auslander-Reiten quiver of $\F(\Delta)$ from the Auslander-Reiten quiver of $A$-mod of the quasi-hereditary algebra $A$ whose Jordan-H\"{o}lder composition series as a left $A$-module is given by
 $$\begin{xy}
\xymatrixrowsep{0in}
\xymatrixcolsep{-0.1in}
  \xymatrix{ 
	  &	S_1    	&		&	&S_2	&	&		&	&S_3	&	&		&S_4 \	\\	
_AA = \ \ &	S_2    	&\ \oplus \	&S_1	&	&S_3	&\ \oplus \	&S_2	&	&S_4	&\ \oplus \	&S_3.	\\	
	  &	S_1    	&		&	&S_2	&	&		&	&S_3	&	&		& }	
 \end{xy}$$
\end{Ex}

Recall the Auslander-Reiten quiver of $A$-mod that has already been calculated in Example \ref{standardexample} and modules in $\F(\Delta)$ are again marked red.

$$\begin{xy}
\xymatrixrowsep{0.5in}
\xymatrixcolsep{0.5in}

  \xymatrix@!0{ 
&	&		&	&					& {\color{red}[P_1]} \ar[dr] 	&	&	&	&	\\
&	&[{\color{red}P_4} \ar[dr]	&	&{\color{red}\Delta_2} \ar[dr] \ar[ur] \ar@{.>}[ll]	&		&\nabla_2 \ar[dr] \ar@{.>}[ll]	&	&I_4] \ar[dr] \ar@{.>}[ll] & \\
&S_3\ar[dr] \ar[ur]	& &{\color{red}X} \ar[dr] \ar[ur] \ar@{.>}[ll]&	&S_2 \ar[dr] \ar[ur] \ar@{.>}[ll]&	&Y \ar[dr] \ar[ur] \ar@{.>}[ll] &	&S_3  \ar@{.>}[ll]\\ 
{\color{red}N_2} \ar[ddr] \ar[dr] \ar[ur] &	&M_2 \ar[dr] \ar[ur] \ar@{.>}[ll]&	&N_3 \ar[ddr] \ar[dr] \ar[ur] \ar@{.>}[ll]& &M_3 \ar[dr] \ar[ur] \ar@{.>}[ll]&	&{\color{red}N_2} \ar[dr] \ar[ur] \ar@{.>}[ll] & \\
&{\color{red}S_1} \ar[ur]&  &\nabla_3 \ar[ur] \ar@{.>}[ll]&	&S_4 \ar[ur] \ar@{.>}[ll]&	&{\color{red}\Delta_3} \ar[ur] \ar@{.>}[ll] &	&{\color{red}S_1}  \ar@{.>}[ll]\\ 
&{\color{red}[P_2]} \ar[uur]	&	&	&	& {\color{red}[P_3]} \ar[uur] 	&	&	&	&	
 }
\end{xy}$$
There are sectional paths from $X$ to $P_3$, $\Delta_2$ to $\Delta_3$, $P_1$ to $N_2$, $P_3$ to $\Delta_3$, $N_2$ to $P_4$, $S_1$ to $X$ and $P_2$ to $X$. Since these are sectional paths in $A$-mod, we can apply Theorem \ref{sectionalirreducible} to all these paths including those that end in a projective module and obtain that they all give rise to irreducible morphisms in $\F(\Delta)$. We can now easily derive the Auslander-Reiten quiver of $\F(\Delta)$.
$$\begin{xy}
\xymatrixrowsep{0.5in}
\xymatrixcolsep{0.5in}

  \xymatrix@!0{ 
		&	&					& [P_1] \ar[dr] &				&[P_4 \ar[dr]			&		\\
[P_4 \ar[dr]	&	&\Delta_2 \ar[dr] \ar[ur] \ar@{.>}[ll]	&	&N_2 \ar[ddr] \ar[dr] \ar[ur] \ar@{.>}[ll]&				&X \ar@{.>}[ll]	\\
		&X \ar[dr] \ar[ur] &				&\Delta_3 \ar[ur] \ar@{.>}[ll]		&	& S_1] \ar[ur] \ar@{.>}[ll]	&	\\ 
&	&[P_3] \ar[ur]	&	&	& [P_2] \ar[uur] 	&		
 }
\end{xy}$$
In order to prove that certain restrictions to non-projective modules in Section \ref{secunique} are necessary, we continue by computing the Auslander-Reiten quiver of the functorially finite resolving subcategory $^\bot T = \{\Ext_A^i(-,T) = 0\}$ for $T = P_1 \oplus P_2 \oplus P_3 \oplus \Delta_2$, which is clearly a generalized cotilting module. The indecomposable modules in this category are $\Delta_2, \Delta_3, P_1, P_2, P_3, P_4$ and $X$. The Auslander-Reiten quiver is easily obtained in the same manner as before.
$$\begin{xy}
\xymatrixrowsep{0.5in}
\xymatrixcolsep{0.5in}

  \xymatrix@!0{ 
	&		&					& 	&[P_1]	\\
	&[P_4 \ar[dr]	&	&\Delta_2] \ar[dr] \ar[ur] \ar@{.>}[ll]	&	\\
\Delta_3 \ar[ur] \ar[dr]&		&X \ar[dr] \ar[ur] \ar@{.>}[ll] &	&\Delta_3 \ar@{.>}[ll]	\\ 
	&[P_2] \ar[ur]	&	& [P_3] \ar[ur] 	&		\\
[P_1] \ar[ur] &		&	&	&
 }
\end{xy}$$

Considering the Jordan-H\"{o}lder decompositions of $\Delta_2$ and $P_2$ it is obvious that there is only one morphism from $\Delta_2$ to $P_2$ up to scalar multiplication. This morphism factors over both $P_1$ and $\Delta_3$ and can clearly be given by two different sectional paths. This contradicts neither Theorem \ref{unique} nor Corollary \ref{submodule} as $P_2$ is projective.
\bigskip
\\In the Auslander-Reiten quiver of $^\bot T$ there is a sectional path $\gamma$ from $P_1$ to itself. This shows that sectional cycles exist if they contain both an Ext-injective and Ext-projective module. As $\gamma^2$ is still a sectional path but $\gamma^2 = 0$, we can see that Corollary \ref{non-zero} does not hold for Ext-projective modules.
\bigskip
\\If we consider another generalized cotilting module $T' = P_1 \oplus P_2 \oplus P_3 \oplus \Delta_3$, we obtain another subcategory $^\bot T'$ containing the indecomposable modules $P_1, P_2, P_3, P_4$ and $\Delta_3$. Hence $^\bot T'$ is a subcategory of $^\bot T$. Although it does not follow from previous corollaries, it is easy to see that there is an irreducible morphism $f : P_2 \to P_3$ in $^\bot T'$ since it only factors over $X$ in $^\bot T$.
\bigskip
\\The morphism $f$ considered in $A$-mod clearly corresponds to a non-sectional path factoring over $N_3$, $X$ and $M_2=\tau(N_3)$. All domains of the right $^\bot T'$-approximations of these modules have at most one direct summand which is not Ext-injective in $^\bot T'$. This shows that Theorem \ref{sec} is not true when not restricted to modules that are not projective and Ext-injective respectively.
\bigskip
\\One of the main questions arising is whether there are irreducible morphisms in $\Omega$ given by non-sectional paths. We can answer this question considering the following example.

\begin{Ex}
Let $A$ be the path algebra of the quiver of Dynkin type $D_5$ with the following orientation.
$$\begin{xy}
\xymatrixrowsep{0.1in}
\xymatrix{
      &		&	&e_4	\\
e_1 \ar[r]& e_2 \ar[r] & e_3 \ar[ur] \ar[dr] \\
      &		&	&e_5}
   
\end{xy}$$

Then there is a functorially finite resolving subcategory $\Omega$ such that there is an $\Omega$-irreducible morphism from $P_2$ to $S_2$.
\end{Ex}

As usual, the projective, injective and simple modules are named $P_i, I_i$ and $S_i$ respectively. Note that $S_4 = P_4$, $S_5 = P_5$ and $I_1 = S_1$. The Auslander-Reiten quiver of $A$ is given by
$$\begin{xy}
\xymatrixrowsep{0.5in}
\xymatrixcolsep{0.5in}

  \xymatrix@!0{ 
& & & [P_1 \ar[dr] & & S_3 \ar[dr] \ar@{.>}[ll]	& &S_2 \ar[dr] \ar@{.>}[ll]& &I_1] \ar@{.>}[ll] &	&	& \\
& & [P_2 \ar[dr] \ar[ur] & & \tau^2(I_2) \ar[dr] \ar[ur] \ar@{.>}[ll]	& &\tau(I_2) \ar[dr]\ar[ur] \ar@{.>}[ll]& &I_2] \ar[ur] \ar@{.>}[ll]&  \\
& [P_3 \ar[ddr]\ar[dr] \ar[ur] & & \tau^2(I_3) \ar[ddr]\ar[dr] \ar[ur] \ar@{.>}[ll]	& &\tau(I_3) \ar[ddr]\ar[dr] \ar[ur] \ar@{.>}[ll]& &I_3] \ar[ur] \ar@{.>}[ll]	&	& \\
[P_4 \ar[ur] &	& \tau^2(I_4) \ar[ur] \ar@{.>}[ll]	&	&\tau(I_4) \ar[ur] \ar@{.>}[ll]&  	&I_4] \ar[ur] \ar@{.>}[ll]	&	&	& \\ 
[P_5 \ar[uur] &	& \tau^2(I_5) \ar[uur] \ar@{.>}[ll]	&	&\tau(I_5) \ar[uur] \ar@{.>}[ll]&  	&I_5]. \ar[uur] \ar@{.>}[ll]	&	&	& }
\end{xy}$$

For a better understanding, we provide the Jordan-H\"{o}lder composition series of the non-simple, indecomposable modules.

 $$\begin{xy}
\xymatrixrowsep{0in}
\xymatrixcolsep{-0.1in}
  \xymatrix{ 
	 &   	&S_1   	&	&				&	&	&	&				&	&	&	\\	
	 &	&S_2   	&	&				&	&S_2	&	&				&	&	&	\\
P_1=\ \	 &	&S_3    &	&\ \ \ \ \ \ \ \ \ \ P_2= \ \	&	&S_3	&	&\ \ \ \ \ \ \ \ \ \ P_3= \ \	&	&S_3	&	\\		
	 &S_4   &	&S_5	&				&S_4	&	&S_5	&				&S_4	&	&S_5	   }
\end{xy}$$

 $$\begin{xy}
\xymatrixrowsep{0in}
\xymatrixcolsep{-0.1in}
  \xymatrix{ 
	 &S_1	&				&S_1	&				&S_1	&				&S_1		\\	
I_5=\ \	 &S_2  	&\ \ \ \ \ \ \ \ \ \ I_4= \ \	&S_2	&\ \ \ \ \ \ \ \ \ \ I_3= \ \	&S_2	& \ \ \ \ \ \ \ \ \ \ I_2= \ \ 	&S_2		\\
	 &S_3   &				&S_3	&				&S_3	&				&		\\
	 &S_5   &				&S_4	&				&	&	  			&		}
\end{xy}$$

 $$\begin{xy}
\xymatrixrowsep{0in}
\xymatrixcolsep{-0.1in}
  \xymatrix{ 
		 &   	&S_1   	&	&					&	&S_1	&	&		\\	
		 &	&S_2   	&	&					&	&S_2^2	&	&		\\
\tau^2(I_2)=\ \	 &	&S_3^2  &	&\ \ \ \ \ \ \ \ \ \ \tau(I_3)= \ \	&	&S_3^2	&	&		\\	
		 &S_4   &	&S_5	&					&S_4	&	&S_5	&		   }
\end{xy}$$

 $$\begin{xy}
\xymatrixrowsep{0in}
\xymatrixcolsep{-0.1in}
  \xymatrix{ 
		 &	&S_2   	&	&					&S_2	&					&S_2	\\
\tau^2(I_3)=\ \	 &	&S_3^2  &	&\ \ \ \ \ \ \ \ \ \ \tau(I_4)= \ \	&S_3	&\ \ \ \ \ \ \ \ \ \ \tau(I_5)= \ \	&S_3		\\			 	 &S_4   &	&S_5	&					&S_5	&					&S_4	   }
\end{xy}$$

$$\begin{xy}
\xymatrixrowsep{0in}
\xymatrixcolsep{-0.1in}
  \xymatrix{ 
\tau(I_2)=\ \	&S_2  &\ \ \ \ \ \ \ \ \ \ \tau^2(I_4)= \ \ &S_3	&\ \ \ \ \ \ \ \ \ \ \tau^2(I_5)= \ \	&S_3		\\			 	 		&S_3  &					    &S_4	&					&S_5	   }
\end{xy}$$

It can easily be checked that $T = P_1 \oplus P_4 \oplus P_5 \oplus S_2 \oplus I_2$ is a generalized cotilting module and the functorially finite resolving subcategory $\Omega = $ $^\bot T$ contains $P_2, P_3$ and $I_1$ in addition to the direct summands of $T$. The morphism $f$ given by the path
$$\begin{xy}

  \xymatrix@!{ 
P_2 \ar[r]& \tau^2(I_3) \ar[r]	 & \tau(I_5) \ar[r]	& \tau(I_3) \ar[r]	& \tau(I_2) \ar[r]	& S_2		}
\end{xy}$$
does not factor over $\tau^2(I_2)$. Therefore, the only other decomposition of $f$ is given by the path 
$$\begin{xy}

  \xymatrix@!{ 
P_2 \ar[r]& \tau^2(I_3) \ar[r]	 & \tau(I_4) \ar[r]	& \tau(I_3) \ar[r]	& \tau(I_2) \ar[r]	& S_2		}
\end{xy}$$
and hence $f$ is $\Omega$-irreducible. The Auslander-Reiten quiver of $\Omega$ also contains an irreducible morphism from $P_1$ to $I_2$ and looks as follows.

$$\begin{xy}
\xymatrixrowsep{0.5in}
\xymatrixcolsep{0.5in}

  \xymatrix@!0{ 
& 	& 			& [P_1] \ar[dr]  & 				&	 		\\
& 	& [P_2 \ar[dr] \ar[ur]  & 		& I_2] \ar[dr] \ar@{.>}[ll]	&	 		\\
& [P_3  \ar[ur]  & 		& S_2 \ar[ur] \ar@{.>}[ll]& 			&S_1] \ar@{.>}[ll]	\\
[P_4] \ar[ur] &	 &		&			&			&	 		\\ 
[P_5] \ar[uur]&	 & 		&			&			&	 		}
\end{xy}$$

Moreover, let us observe what Corollary \ref{seccontra} means for this example. All decompositions of $f$ factor over $\tau(I_3)$, which is a module with $3$ predecessors and must occur by the aforementioned corollary. Furthermore, the approximation of $\tau^2(I_3)$ is a morphism $f_{\tau^2(I_3)}: P_3 \oplus P_2 \to \tau^2(I_3)$, so its domain contains $2$ non-isomorphic indecomposable direct summands that are not Ext-injective.

\chapter{Degrees of irreducible morphisms in functorially finite subcategories}

Degrees of irreducible morphisms have been established for irreducible morphism in $A$-mod by Liu in \cite{L92}. In this chapter we generalize this notion for functorially finite subcategories and see that irreducible morphisms and paths behave similarly to $A$-mod when Ext-projective modules are not involved. We give a proof for a generalized version of the Happel-Preiser-Ringel theorem, which says that infinite connected components of $\Gamma_A$ that contain a $\tau$-periodic module are stable tubes \cite{HPR80}. Moreover, we analyze the shape of left stable connected components that contain an oriented cycle but no periodic modules. The whole chapter closely follows \cite{L92} and \cite{L93}. In order to do this, we introduce some additional notation. 
\bigskip
\\Let $Z$ be an indecomposable module in $\chi$ that is not Ext-projective, we then denote by $\E_\chi(Z)$ the almost split sequence in $\chi$ ending in $Z$,
$$\begin{xy}
  \xymatrix{ 
    0 \ar[r]   &   \t(Z) \ar[r] 	&	Y \ar[r]    &    Z \ar[r]  	&	0.	    }
\end{xy}$$
Moreover, we denote the almost split sequence 
$$\begin{xy}
  \xymatrix{ 
    0 \ar[r]   &   X \ar[r] 	&	Y \ar[r]    &    \t^{-1}(X) \ar[r]  	&	0	    }
\end{xy}$$
by $\E'_\chi(X)$ if $\t^{-1}(X)$ exists for a given indecomposable module $X$ in $\chi$. The module $Y$ is called the \textbf{middle term} of $\E_\chi(Z)$ and $\E'_\chi(X)$.
\bigskip
\\Let $f: X \to Y$ and $g: Y \to Z$ be irreducible morphisms in $\chi$. We say that a pair $\{f, g\}$ is a \textbf{component} of an almost split sequence if either there is an almost split sequence
$$\begin{xy}
  \xymatrix{ 
    0 \ar[r]   &   X \ar[r]^g 	&	Y \ar[r]^f    &    Z \ar[r]  	&	0	    }
\end{xy}$$
in $\chi$ or if there are $\chi$-irreducible morphisms $f' : X \to Y'$ and $g' : Y' \to Z$ such that 
$$\begin{xy}
  \xymatrix{ 
    0 \ar[r]   &   X \ar[r]^{(g,g')^T} 	&	Y \oplus Y' \ar[r]^{(f,f')}    &    Z \ar[r]  	&	0	    }
\end{xy}$$
is an almost split sequence in $\chi$.

\section{Definitions and basic properties}

\begin{Def}
Let $f: X \to Y$ be an irreducible morphism in $\chi$. 
\begin{enumerate}[(a)]
 \item We define the \textbf{left degree} $d_\chi^l(f)$ of $f$ in $\chi$ to be the least integer $n$ such that there is a module $Z$ and a morphism $g \in \rad_\chi^n (Z, X) \backslash \rad_\chi^{n+1}(Z,X)$ such that $fg \in \rad_\chi^{n+2}(Z,Y)$. If none such integer exists, we define $d_\chi^l(f)$ to be $\infty$. 
 \item Dually, we define the \textbf{right degree} $d_\chi^r(f)$ of $f$ in $\chi$ to be the least integer $n$ such that there is a module $Z$ and a morphism $g \in \rad_\chi^n (Y, Z)\backslash \rad_\chi^{n+1}(Y,Z)$ such that $gf \in \rad_\chi^{n+2}(X,Z)$. If none such integer exists, we define $d_\chi^r(f)$ to be $\infty$. 
\end{enumerate}
\end{Def}

\begin{Lemma}\label{degree1}
\mbox{}
\begin{enumerate}[(a)]
 \item Let $n \geq 1$ be an integer and let $p : X \to Y$ and $f : Y \to Z$ be morphisms in $\chi$, where $f$ is $\chi$-irreducible and $Z$ is indecomposable and not Ext-projective. If $p \notin \rad_\chi^{n+1}(X,Y)$, $fp \in \rad_\chi^{n+2}(X,Z)$ and 
$$\begin{xy}
  \xymatrix{ 
    0 \ar[r]   &   \tau_\chi(Z) \ar[r]^{(g,g')^T} 	&	Y \oplus Y' \ar[r]^(.6){(f,f')}    &    Z \ar[r]  	&	0	    }
\end{xy}$$
is an almost split sequence in $\chi$, then there exists a morphism $q : X \to \tau_\chi(Z)$ such that $q \notin \rad_\chi^{n}(X,\tau_\chi(Z))$, $p+gq \in \rad_\chi^{n+1}(X,Y)$ and $g'q \in \rad_\chi^{n+1}(X,Y')$. 
\item Dually, let $n \geq 1$ be an integer and let $g : X \to Y$ and $p : Y \to Z$ morphisms in $\chi$, where $g$ is $\chi$-irreducible and $X$ is indecomposable and not Ext-injective. If $p \notin \rad_\chi^{n+1}(Y,Z)$, $pg \in \rad_\chi^{n+2}(X,Z)$ and 
$$\begin{xy}
  \xymatrix{ 
    0 \ar[r]   &   X \ar[r]^(.45){(g,g')^T} 	&	Y \oplus Y' \ar[r]^{(f,f')}    &    \tau_\chi^{-1}(X) \ar[r]  	&	0	    }
\end{xy}$$
is an almost split sequence in $\chi$, then there exists a morphism $q : \tau_\chi^{-1}(X) \to Z$ such that $q \notin \rad_\chi^{n}(\tau_\chi^{-1}(X),Z)$, $p+qf \in \rad_\chi^{n+1}(Y,Z)$ and $qf' \in \rad_\chi^{n+1}(Y',Z)$. 
\end{enumerate}

\end{Lemma}

\proof
\\There is a factorization $fp = ts$ with $s \in \rad_\chi^{n+1}(X,W)$ and $t \in \rad_\chi(W,Z)$. The morphism $t$ factors through $(f,f')$ as the latter is a minimal right almost split morphism. If $t = (f,f')(u,u')^T$, then $(f,f')(us-p,u's)^T = 0$. Since $\Im(us-p,u's)^T \subset \Im(g,g')^T \cong \tau_\chi(Z)$, there exists a morphism $q: X \to \tau_\chi(Z)$ such that $(us-p, u's)^T = (g,g')^Tq$, which is equivalent to $(p + gq, g'q)^T = (us, u's)^T \in \rad_\chi^{n+1}(X,Y \oplus Y')$ as $s \in \rad_\chi^{n+1}(X,W)$. Thus $p + gq \in \rad_\chi^{n+1}(X,Y)$ and $g'q \in \rad_\chi^{n+1}(X,Y')$. Moreover, as $p \notin \rad_\chi^{n+1}(X,Y)$ and $p + gq \in \rad_\chi^{n+1}(X,Y)$, we know that $gq \notin \rad_\chi^{n+1}(X,Y)$ and hence $q \notin \rad_\chi^{n}(X,\tau_\chi(Z))$. \qed

\begin{Cor}\label{degree1Cor}
\mbox{}
\begin{enumerate}[(a)]
 \item Let $f: Y \to Z$ be an irreducible morphism in $\chi$ with finite left degree in $\chi$, where $Z$ is an indecomposable module that is not Ext-projective. If $ Y \oplus Y'$ is a direct summand of the whole middle term of $\E_\chi(Z)$ with $Y' \neq 0$, then there is an irreducible morphism $g' : \tau_\chi(Z) \to Y'$ with $d_\chi^l(g') < d_\chi^l(f)$. Consequently, if $d_\chi^l(f) = 1$, then $f$ is a surjective minimal right almost split morphism.
 \item Dually, let $g: X \to Y$ be an irreducible morphism in $\chi$ with finite right degree in $\chi$, where $X$ is an indecomposable module that is not Ext-injective. If $ Y \oplus Y'$ is a direct summand of the whole middle term of $\E'(X)$ with $Y' \neq 0$, then there is an irreducible morphism $f' : Y' \to \tau_\chi^{-1}(X)$ with $d_\chi^r(f') < d_\chi^r(g)$. Consequently, if $d_\chi^r(g) = 1$, then $g$ is an injective minimal left almost split morphism.
\end{enumerate}
\end{Cor}
\proof
\\Suppose that $d_\chi^l(f) = n$, i.e. there is a morphism $p \in \rad_\chi^{n}(X,Y) \backslash \rad_\chi^{n+1}(X,Y)$ such that $fp \in \rad_\chi^{n+2}(X,Z)$. Since $Y \oplus Y'$ is a direct summand of the middle term of $\E_\chi(Z)$, we know there are irreducible morphisms $(g,g')^T : \t(Z) \to Y \oplus Y'$ and $(f, f'): Y \oplus Y' \to Z$ such that $\{(g,g')^T,(f,f')\}$ is a component of $\E_\chi(Z)$. If $(g,g',g'')^T : \t(Z) \to Y \oplus Y' \oplus Y''$ denotes an extension of $(g,g')^T$ to a minimal left almost split morphism, then there is a morphism $q \notin \rad_\chi^{n}(X,\tau_\chi(Z))$ such that $(g',g'')^Tq \in \rad_\chi^{n+1}(X,Y' \oplus Y'')$ by Lemma \ref{degree1}. Consequently, we also have $g'q \in \rad_\chi^{n+1}(X,Y')$, which implies that $d_\chi^l(g') \leq n-1 < d_\chi^l(f)$. \qed

\begin{Def}
A path 
$$\begin{xy}
  \xymatrix{ 
     X_0 \ar[r] 	&	X_1 \ar[r]     &   \cdots	 \ar[r]	&	X_{n-1} \ar[r]	&	X_n 	    }
\end{xy}$$
in $\Gamma_\chi$ is said to be \textbf{pre-sectional in $\pmb{\chi}$} if for all $i = 2, \ldots, n$ such that $\tau_\chi(X_i)$ exists, $X_{i-2} \oplus \tau_\chi(X_i)$ is a direct summand of the domain of a minimal right almost split morphism mapping to $X_{i-1}$.
\end{Def}

Equivalently, one can define that the path is pre-sectional in $\chi$ if for all $i = 2, \ldots, n$ such that $\tau_\chi^{-1}(X_{i-2})$ exists, $\tau_\chi^{-1}(X_{i-2}) \oplus X_i$ is a direct summand of the codomain of a minimal left almost split morphism mapping from $X_{i-1}$. Clearly, every sectional path is pre-sectional.

\begin{Lemma}
Let
$$\begin{xy}
  \xymatrix{ 
     X_0 \ar[r] 	&	X_1 \ar[r]     &   \cdots	 \ar[r]	&	X_{n-1} \ar[r]	&	X_n 	    }
\end{xy}$$
be a pre-sectional path in $\Gamma_\chi$ and let $m$ be an integer. Then the following paths are also pre-sectional if they are defined.
\begin{enumerate}[(a)]
 \item $\begin{xy}
  \xymatrix{ 
     \tau_\chi^n(X_n) \ar[r] 	&	\tau_\chi^{n-1}(X_{n-1}) \ar[r]     &   \cdots	 \ar[r]	&	\tau_\chi(X_1) \ar[r]	&	X_0 	    }
\end{xy}$
 \item $\begin{xy}
  \xymatrix{ 
     X_n \ar[r] 	&	\tau_\chi^{-1}(X_{n-1}) \ar[r]     &   \cdots	 \ar[r]	&	\tau_\chi^{1-n}(X_1) \ar[r]	&	\tau_\chi^{-n}(X_0)    }
\end{xy}$
\item $\begin{xy}
  \xymatrix{ 
     \tau_\chi^m(X_0) \ar[r] 	&	\tau_\chi^m(X_1) \ar[r]     &   \cdots	 \ar[r]	&	\tau_\chi^m(X_{n-1}) \ar[r]	&	\tau_\chi^m(X_n)&   }
\end{xy}$
\end{enumerate}

\end{Lemma}

\proof
\\The lemma follows immediately from Corollary \ref{valuationisinvariantoftranslation}. \qed

\newpage

\begin{Lemma}\label{lowerdegree}
\mbox{}
\begin{enumerate}[(a)]
 \item Let $f : X \to Y$ be an irreducible morphism in $\chi$ such that $Y$ is indecomposable and $d_\chi^l(f) < \infty$. Let 
$$\begin{xy}
  \xymatrix{ 
     X_n \ar[r] 	&	X_{n-1} \ar[r]     &   \cdots	 \ar[r]	&	X_1 \ar[r]	&	X_0 = Y 	    }
\end{xy}$$
be a pre-sectional path in $\Gamma_\chi$ such that no $X_i$ is Ext-projective and $n \geq 1$. If $X \oplus X_1$ is a direct summand of the whole middle term of $\E_\chi(Y)$, then for each $1 \leq i \leq n$ there is an irreducible morphism $f_i: \tau_\chi(X_{i-1}) \to X_i$ such that
$$d_\chi^l(f_n) < d_\chi^l(f_{n-1}) < \cdots < d_\chi^l(f_1) < d_\chi^l(f).$$
Consequently, $n < d_\chi^l(f)$.
 \item Dually, suppose $g : X \to Y$ be an irreducible morphism in $\chi$ such that $X$ is indecomposable and $d_\chi^r(g) < \infty$. Let 
$$\begin{xy}
  \xymatrix{ 
    X = Y_0 \ar[r] 	&	Y_1 \ar[r]     &   \cdots	 \ar[r]	&	Y_{n-1} \ar[r]	&	Y_n 	    }
\end{xy}$$
be a pre-sectional path in $\Gamma_\chi$ such that no $Y_i$ is Ext-injective and $n \geq 1$. If $Y \oplus Y_1$ is a direct summand of the whole middle term of $\E'_\chi(X)$, then for each $1 \leq i \leq n$ there is an irreducible morphism $g_i: Y_i \to \tau_\chi^{-1}(Y_{i-1})$ such that
$$d_\chi^r(g_n) < d_\chi^r(g_{n-1}) < \cdots < d_\chi^r(g_1) < d_\chi^r(g).$$
Consequently, $n < d_\chi^r(g)$.
\end{enumerate}
\end{Lemma}

\proof
\\Let $X \oplus X_1$ be a direct summand of the whole middle term of $\E_\chi(Y)$. Then by Corollary \ref{degree1Cor} there is an irreducible morphism $f_1: \tau_\chi(Y) \to X_1$ with $d_\chi^l(f_1) < d_\chi^l(f).$ Suppose that $1 \leq m < n$ and for each $i$ such that $1 \leq i \leq m$ there is an irreducible morphism $f_i: \tau_\chi(X_{i-1}) \to X_i$ such that $d_\chi^l(f_m) < \cdots < d_\chi^l(f_1) < d_\chi^l(f).$ By definition of pre-sectional paths $\tau_\chi(X_{m-1}) \oplus X_{m+1}$ is a direct summand of the whole middle term of $\E_\chi(X_m)$. We can apply Corollary \ref{degree1Cor} again to find an irreducible morphism $f_{m+1}: \tau_\chi(X_{m}) \to X_{m+1}$ such that $d_\chi^l(f_{m+1}) < d_\chi^l(f_m)$. Hence the result follows inductively. \qed

\begin{Cor}\label{infdegree}
\mbox{}
\begin{enumerate}[(a)]
 \item Let $Y$ be indecomposable and let $f : X \to Y$ be an irreducible morphism in $\chi$. Let 
$$\begin{xy}
  \xymatrix{ 
     \cdots \ar[r] &	X_n \ar[r] 	&	X_{n-1} \ar[r]     &   \cdots	 \ar[r]	&	X_1 \ar[r]	&	X_0 = Y 	    }
\end{xy}$$
be an infinite pre-sectional path in $\Gamma_\chi$ such that no $X_i$ is Ext-projective and $X \oplus X_1$ is a direct summand of $\E_\chi(Y)$. Then $f$ has infinite left degree in $\chi$.

\newpage
 \item Dually, let $g : X \to Y$ be an irreducible morphism in $\chi$ for some indecomposable module $X$. Let 
$$\begin{xy}
  \xymatrix{ 
     X = Y_0 \ar[r] 	&	Y_1 \ar[r]     &   \cdots	 \ar[r]	&	Y_{n-1} \ar[r]	&	Y_n \ar[r]	&	\cdots 	    }
\end{xy}$$
be an infinite pre-sectional path in $\Gamma_\chi$ such that no $Y_i$ is Ext-injective and $Y \oplus Y_1$ is a direct summand of $\E'_\chi(X)$. Then $g$ has infinite right degree in $\chi$.
\end{enumerate}

\end{Cor}

\proof
\\As the path is pre-sectional, we obtain an infinite collection of irreducible morphisms $f_i: \tau_\chi(X_{i-1}) \to X_i$ such that 
$$\cdots < d_\chi^l(f_n) < d_\chi^l(f_{n-1}) < \cdots < d_\chi^l(f_1) < d_\chi^l(f).$$
by Lemma \ref{lowerdegree}. \qed

\begin{Def}
Let $Y$ and $Y'$ be indecomposable modules in $\chi$ that are not Ext-projective and let $(f,f') : \t(Y) \oplus \t(Y') \to X$ and $(g,g')^T : X \to Y \oplus Y'$ be irreducible morphisms. If $\{f,g\}$ and $\{f',g'\}$ are components of $\E_\chi(Y)$ and $\E_\chi(Y')$ respectively, we say that $(f,f')$ is a \textbf{left neighbour} of $(g,g')^T$ and $(g,g')^T$ is a \textbf{right neighbour} of $(f,f')$.
\end{Def}

\begin{Lemma}\label{neighbour}
Let $f : X \to Y$ be an irreducible morphism in $\chi$.
\begin{enumerate}[(a)]
 \item If $f$ has finite left degree and $Y = Y_1 \oplus Y_2$ such that both $Y_i$ are indecomposable and not Ext-projective, then $f$ has a left neighbour 
 $$g = (g_1, g_2): \tau_\chi(Y_1) \oplus \tau_\chi(Y_2) \to X$$
with $d_\chi^l(g) < d_\chi^l(f)$
 \item If $f$ has finite right degree and $X = X_1 \oplus X_2$ such that both $X_i$ are indecomposable and not Ext-injective, then $f$ has a right neighbour 
 $$g = (g_1, g_2)^T: Y \to \tau_\chi^{-1}(X_1) \oplus \tau_\chi^{-1}(X_2)$$
with $d_\chi^r(g) < d_\chi^r(f)$.

\end{enumerate}

\end{Lemma}

\proof
\\Assume that $d_\chi^l(f) = n$. Then there is a $p \in \rad_\chi^n (Z,X) \backslash \rad_\chi^{n+1} (Z,X)$ such that $fp \in \rad_\chi^{n+2} (Z,Y)$. Let 
 $$g = (g_1, g_2) : \tau_\chi(Y_1) \oplus \tau_\chi(Y_2) \to X$$
be a left neighbour of $f$. By definition $\{g_i,f_i\}$ is a component of $\E_\chi(Y_i)$. So by Lemma \ref{degree1} there are $q_i \notin \rad_\chi^n (Z, \tau_\chi(Y_i))$ such that $p + g_iq_i \in \rad_\chi^{n+1} (Z, X)$. Consequently, $g(q_1, -q_2)^T = g_1q_1 - g_2q_2 = p + g_1q_1 - (p + g_2q_2) \in \rad_\chi^{n+1} (Z, X)$. Since $(q_1, q_2)^T \notin \rad_\chi^n (Z, \tau_\chi(Y_1) \oplus \tau_\chi(Y_2))$, we obtain $d_\chi^l(g) \leq n -1 < d_\chi^l(f)$, which completes the proof. \qed

\newpage

\begin{Def}\label{defstable}
Let $X$ be an indecomposable module in $\chi$. 
\begin{enumerate}[(a)]
\item We call $X$ \textbf{left stable} in $\chi$ if $\t^n(X)$ exists for all $n \geq 0$.
\item Dually, $X$ is called \textbf{right stable} in $\chi$ if $\t^n(X)$ exists for all $n \leq 0$.
\item $X$ is called \textbf{stable} in $\chi$ if it is both left and right stable.
\item If there is a positive integer $n \geq 1$ such that $\t^n(X) = X$ we call $X$ $\pmb{\t}$\textbf{-periodic} or, for convenience, periodic.
\end{enumerate}
Moreover, we call a subquiver $\Gamma$ of the Auslander-Reiten quiver of $\chi$ left stable, right stable, stable or $\t$-periodic if all their modules are left stable, right stable, stable or $\t$-periodic respectively and, if $X$ is in $\Gamma$, then the module $\t^n(X)$ is also $\Gamma$ for all $n \in \Z$ such that $\t^n(X)$ exists.
\end{Def}

Note that a module that is stable in a subcategory does not have to be stable in neither a larger nor a smaller subcategory. Consider the standard example \ref{standardexample} and let $\chi = \F(\Delta)$ be the category of standard-filtered modules and $\Omega$ = $^\bot T$ the Ext-orthogonal of the generalized cotilting module $T = P_1 \oplus P_2 \oplus P_3 \oplus \Delta_2$. In \ref{standardexamplesubcategories} we have shown that $\Omega$ is a subcategory of $\chi$ just as in our general setup. The indecomposable module $X$ is stable and periodic in $A$-mod as $\tau^4(X) = X$, but $X$ is neither left or right stable in $\chi$ as $\t^3(X)$ is projective and $\t^{-3}(X)$ is Ext-injective in $\chi$. On the other hand, in $\Omega$ we have $\tau_\Omega^2(X) = X$, hence $X$ is $\tau_\Omega$-periodic and, consequently, stable in $\Omega$.

\begin{Lemma}\label{finitedegreevaluation}
Let $X \to Y$ be an arrow in $\Gamma_\chi$ with valuation $(a,a)$. 
If $X$ and $Y$ are left stable in $\chi$ and there is an irreducible morphism $f: X \to Y$ with finite left degree or if $X$ and $Y$ are right stable in $\chi$ and there is an irreducible morphism $f: X \to Y$ with finite right degree, then $a = 1$.
\end{Lemma}

\proof
\\Suppose $X$ and $Y$ are left stable and $f$ has finite left degree, hence we prove the statement by induction on $d_\chi^l(f)$. If $d_\chi^l(f) = 1$, then $f$ must be a minimal right almost split morphism by Corollary \ref{degree1Cor}, so $a = 1$. Assume the statement is true for $d_\chi^l(f) < m$. If $d_\chi^l(f) = m$ and $a > 1$, then $X \oplus X$ is a direct summand of the whole middle term of $\E_\chi(Y)$. We use Corollary \ref{degree1Cor} to obtain a morphism $g : \tau_\chi(Y) \to X$ with $d_\chi^l(g) < m$. By induction the arrow from $\tau_\chi(Y)$ to $X$ must be valued $(1,1)$. On the other hand, the same arrow must be valued $(a,a)$ by Corollary \ref{valuationisinvariantoftranslation}, which is a contradiction. The proof for the other condition works dually. \qed

\begin{Cor}
Let $X \to Y$ be an arrow in $\Gamma_\chi$ and let $f : X \to Y$ and $g : X \to Y$ be irreducible morphisms in $\chi$. 
\begin{enumerate}[(a)]
 \item If the valuation of the arrow is $(1,1)$ or $X$ and $Y$ are left stable in $\chi$, then $d_\chi^l(f) = d_\chi^l(g)$.
 \item Dually, if the valuation of the arrow is $(1,1)$ or $X$ and $Y$ are right stable in $\chi$, then $d_\chi^r(f) = d_\chi^r(g)$.
\end{enumerate}

\end{Cor}

\proof
\\First we assume that $X$ and $Y$ are left stable in $\chi$. Let $(a,a)$ be the valuation of the arrow from $X$ to $Y$, if $a > 1$, then $d_\chi^l(f) = d_\chi^l(g) = \infty$ by Lemma \ref{finitedegreevaluation}. So let $a = 1$ and, without loss of generality, suppose $d_\chi^l(f) < \infty$, i.e. there is a $p \in \rad_\chi^n (Z,X) \backslash \rad_\chi^{n+1} (Z,X)$ such that $fp \in \rad_\chi^{n+2} (Z,Y)$. As $ a = 1$, there is an isomorphism $\varphi : X \to X$ such that $f - g \varphi \in \rad_\chi^2 (X,Y)$. Therefore, $g \varphi p = fp - (f - g \varphi)p \in \rad_\chi^{n+2} (Z,Y)$ and hence $d_\chi^l(g) \leq d_\chi^l(f) < \infty$ as $\varphi p \notin \rad_\chi^{n+1}(Z,X)$ . We similarly obtain $d_\chi^l(f) \leq d_\chi^l(g)$, which completes the proof in this case. \qed
\bigskip
\\Now the following definition makes sense.

\begin{Def}
Let $X \to Y$ be an arrow in $\Gamma_\chi$ and $f: X \to Y$ an irreducible morphism. 
\begin{enumerate}[(a)]
 \item If the valuation of the arrow is $(1,1)$ or $X$ and $Y$ are left stable in $\chi$, we then define the \textbf{left degree of the arrow }$\pmb{X \to Y}$ to be $d_\chi^l(f)$.
 \item Dually, if the valuation of the arrow is $(1,1)$ or $X$ and $Y$ are right stable in $\chi$, we then define the \textbf{right degree of the arrow }$\pmb{X \to Y}$ to be $d_\chi^r(f)$.
\end{enumerate}
\end{Def}

\begin{Ex}
We compute the left degrees of all arrows of the Auslander-Reiten quiver of $A$-mod and $\F(\Delta)$ of the algebra discussed in the standard example.
\end{Ex}

Note that for this algebra all arrows that are given by an irreducible monomorphisms have left degree $\infty$. The almost split sequences $\E(\Delta_2)$, $\E(I_4)$, $\E(\Delta_3)$, $\E(S_4)$, $\E(\nabla_3)$ and $\E(S_1)$ have an indecomposable middle term, hence the arrows from these middle terms to $\Delta_2, I_4, \Delta_3, S_4, \nabla_3$ and $S_1$ respectively must have left degree $1$. Let $ f : X \to Y$ be an irreducible morphism between indecomposable $A$-modules such that we have not obtained the left degree of the corresponding arrow yet. It is not hard to see that in this case if
$$\begin{xy}
  \xymatrix{ 
    X_n \ar[r]^{f_n} 	&	X_{n-1} \ar[r]^{f_{n-1}}     &   \cdots	 \ar[r]^{f_2}	&	X_1 \ar[r]^{f_1}	&	X_0 = X 	    }
\end{xy}$$
is the shortest sectional path such that $X_1 = \tau(Y)$ and $ff_1 \cdots f_n = 0$, then $d_\chi^l(f) = n$ and corresponding arrow also has left degree $n$. This completes the calculation and gives us the following quiver with arrows labeled with their left degrees. 
$$\begin{xy}
\xymatrixrowsep{0.5in}
\xymatrixcolsep{0.5in}

  \xymatrix@!0{  
&	&		&	&					& [P_1] \ar[dr]_4 	&	&	&	&	\\
&	&[P_4 \ar[dr]_\infty&	&\Delta_2 \ar[dr]_3 \ar[ur]^\infty \ar@{.>}[ll]&	&\nabla_2 \ar[dr]_\infty \ar@{.>}[ll]&	&I_4] \ar[dr]_3 \ar@{.>}[ll] & \\
&S_3\ar[dr]_\infty \ar[ur]^\infty& &X \ar[dr]_2 \ar[ur]^1 \ar@{.>}[ll]&	&S_2 \ar[dr]_\infty \ar[ur]^\infty \ar@{.>}[ll]&	&Y \ar[dr]_2 \ar[ur]^1 \ar@{.>}[ll] &	&S_3  \ar@{.>}[ll]\\ 
N_2 \ar[ddr]_(.6)\infty \ar[dr]^1 \ar[ur]^2 &	&M_2 \ar[dr]_1 \ar[ur]^\infty \ar@{.>}[ll]&	&N_3 \ar[ddr]_(.6)\infty \ar[dr]^1 \ar[ur]^2 \ar@{.>}[ll]& &M_3 \ar[dr]_1 \ar[ur]^\infty \ar@{.>}[ll]&	&N_2 \ar[dr]_1 \ar[ur]^2 \ar@{.>}[ll] & \\
&S_1 \ar[ur]^\infty&  &\nabla_3 \ar[ur]^\infty \ar@{.>}[ll]&	&S_4 \ar[ur]^\infty \ar@{.>}[ll]&	&\Delta_3 \ar[ur]^\infty \ar@{.>}[ll] &	&S_1  \ar@{.>}[ll]\\ 
&[P_2] \ar[uur]_(.4)4	&	&	&	& [P_3] \ar[uur]_(.4)4 	&	&	&	&	
 }
\end{xy}$$     

In $\F(\Delta)$, just as in $A$-mod, all irreducible monomorphisms have infinite left degree. Moreover, all $\F(\Delta)$-irreducible morphisms between irreducible modules in $\F(\Delta)$ that are minimal right almost split morphisms have left degree $1$. We obtain the left degrees of the remaining arrows similarly to $A$-mod, except for the arrow given by the irreducible morphism $f : P_2 \to X$. Composing this morphism with a sectional path in $\chi$ of length $n$ starting in some module $M$ always gives a morphism in $\rad_{\F(\Delta)}^{n+1} (M, X) \backslash \rad_{\F(\Delta)}^{n+2} (M, X)$. On the other hand, if we compose $f$ with the up to a scalar unique non-zero morphism $g \in \rad_{\F(\Delta)}^3 (\Delta_2, P_2) \backslash \rad_{\F(\Delta)}^4 (\Delta_2, P_2)$, which can be considered as a non-sectional path, we obtain $fg = 0$. Hence the left degree of the arrow from $P_2$ to $X$ is $3$ and we have finished the computations.

$$\begin{xy}
\xymatrixrowsep{0.5in}
\xymatrixcolsep{0.5in}

  \xymatrix@!0{ 
		&	&					& [P_1] \ar[dr]_3 &				&[P_4 \ar[dr]_\infty			&		\\
[P_4 \ar[dr]_\infty	&	&\Delta_2 \ar[dr]_2 \ar[ur]^\infty \ar@{.>}[ll]	&	&N_2 \ar[ddr]_(.6)\infty \ar[dr]^1 \ar[ur]^1 \ar@{.>}[ll]&				&X \ar@{.>}[ll]	\\
		&X \ar[dr]_1 \ar[ur]^1 &				&\Delta_3 \ar[ur]^\infty \ar@{.>}[ll]		&	& S_1] \ar[ur]^\infty \ar@{.>}[ll]	&	\\ 
&	&[P_3] \ar[ur]^2	&	&	& [P_2] \ar[uur]_(.4)3 	&		
 }
\end{xy}$$

Contrary to $A$-mod, where arrows pointing at a projective module always have infinite left degree, the arrow in the Auslander-Reiten quiver of $\F(\Delta)$ from $N_2$ to $P_4$ has left degree $1$.

\begin{Lemma}\label{vmtinfinity}
\mbox{}
\begin{enumerate}[(a)]
 \item Let  
$$\begin{xy}
  \xymatrix{ 
     \cdots \ar[r] &	X_{i+1} \ar[r] 	&	X_i \ar[r]     &   \cdots	 \ar[r]	&	X_1 \ar[r]	&	X_0	    }
\end{xy}$$
be an infinite pre-sectional path in $\Gamma_\chi$ with all $X_i$ left stable. If there is some integer $n \geq 0$ such that the almost split sequence $\E_\chi(X_n)$ has three left stable middle terms or the arrow $X_{n+1} \to X_n$ has non-trivial valuation, then all arrows $\tau_\chi^j(X_{i+1}) \to \tau_\chi^j(X_{i})$ and $\tau_\chi^{j+1}(X_{i}) \to \tau_\chi^j(X_{i+1})$ with $j\geq 0$ and $i > n+1$ have infinite left degree.
\item Let  
$$\begin{xy}
  \xymatrix{ 
     Y_0 \ar[r]  &	Y_1 \ar[r] 	&	  \cdots	 \ar[r]	&	Y_i \ar[r]	&	\cdots	    }
\end{xy}$$
be an infinite pre-sectional path in $\Gamma_\chi$ with all $Y_i$ right stable. If there is some integer $n \geq 0$ such that the almost split sequence $\E'_\chi(Y_n)$ has three right stable middle terms or the arrow $Y_n \to Y_{n+1}$ has non-trivial valuation, then all arrows $\tau_\chi^j(Y_i) \to \tau_\chi^j(Y_{i+1})$ and $\tau_\chi^j(Y_{i+1}) \to \tau_\chi^{j-1}(Y_i)$ with $j\leq 0$ and $i > n+1$ have infinite right degree.
\end{enumerate}

\end{Lemma}

\proof
\\As all $X_i$ are left stable, it is sufficient to prove that all arrows $X_{i+1} \to X_i$ and $\tau(X_i) \to X_{i+1}$ with $i > n + 1$ have infinite left degree. Since $\tau_\chi(X_i) \oplus X_{i+2}$ is a direct summand of the whole middle term of $\E_\chi(X_{i+1})$ by definition, all arrows $\tau_\chi(X_i) \to X_{i+1}$ with $i\geq 0$ have infinite left degree by Corollary \ref{infdegree}.
\newpage
Suppose that $\E_\chi(X_n)$ has three left stable middle terms $X_{n+1}, Y_n$ and $Z_n$ and there is an irreducible morphism $f : X_{i+1} \to X_i$ with finite left degree in $\chi$. We can find a pre-sectional path 
$$\begin{xy}
  \xymatrix{ 
     \cdots \ar[r] &	\tau_\chi^{i-n}(X_n) \ar[r] 	&	\tau_\chi^{i-n-1}(X_{n+1}) \ar[r]     &   \cdots	 \ar[r]	&	\tau_\chi(X_{i-1}) \ar[r]	&	X_i 	    }
\end{xy}$$
in $\Gamma_\chi$ such that $X_{i+1} \oplus \tau_\chi(X_{i-1})$ is a direct summand of the whole middle term of $\E_\chi(X_i)$. By Lemma \ref{lowerdegree} there exists an irreducible morphism $f': \tau_\chi^{i-n}(X_{n+1}) \to \tau_\chi^{i-n}(X_n)$ such that $d_\chi^l(f') < d_\chi^l(f)$. Since $X_{n+1}, Y_n$ and $Z_n$ are left stable, the whole middle term of $\E_\chi(\tau_\chi^{i-n}(X_n))$ has a direct summand $\tau_\chi^{i-n}(X_{n+1}) \oplus \tau_\chi^{i-n}(Y_n) \oplus \tau_\chi^{i-n}(Z_n)$. There is an irreducible morphism $h : \tau_\chi^{i-n+1}(X_n) \to \tau_\chi^{i-n}(Y_n) \oplus \tau_\chi^{i-n}(Z_n)$ with $d_\chi^l(h) < d_\chi^l(f')$ by Corollary \ref{degree1Cor} and hence another irreducible morphism $h' : \tau_\chi^{i-n+1}(Y_n) \oplus \tau_\chi^{i-n+1}(Z_n) \to \tau_\chi^{i-n+1}(X_n)$ such that $d_\chi^l(h') < d_\chi^l(h)$ by Lemma \ref{neighbour}. On the other hand, there is a pre-sectional path 
$$\begin{xy}
  \xymatrix{ 
     \cdots \ar[r] &	\tau_\chi^{i-n+1}(X_{n+k}) \ar[r] 	&	\tau_\chi^{i-n+1}(X_{n+k-1}) \ar[r]  & \cdots    }
\end{xy}$$
$$\begin{xy}
  \xymatrix{ 
 \cdots	 \ar[r]	&	\tau_\chi^{i-n+1}(X_{n+1}) \ar[r]	&	\tau_\chi^{i-n+1}(X_n) 	  &  }
\end{xy}$$
such that $\tau_\chi^{i-n+1}(X_{n+1}) \oplus \tau_\chi^{i-n+1}(Y_n) \oplus \tau_\chi^{i-n+1}(Z_n)$ is a direct summand of the whole middle term of $\E_\chi(\tau_\chi^{i-n+1}(X_n))$. Therefore, by Corollary \ref{infdegree}, $h'$ has infinite left degree, which is a contradiction.
\bigskip
\\Suppose now there is an arrow $X_{n+1} \to X_n$  valued $(a,a)$ with $a > 1$. Then for each integer $i > n$ there is a pre-sectional path
$$\begin{xy}
  \xymatrix{ 
     \cdots \ar[r] &	\tau_\chi^{i-n}(X_{n+k}) \ar[r] 	& \cdots \ar[r] &	\tau_\chi^{i-n}(X_{n+1}) \ar[r]     &  \tau_\chi^{i-n}(X_n) \ar[r] &   }
\end{xy}$$
$$\begin{xy}
  \xymatrix{ 
\ar[r] & \tau_\chi^{i-n-1}(X_{n+1}) \ar[r] & \cdots	 \ar[r]	&	\tau_\chi(X_{i-1}) \ar[r]	&	X_i &	    }
\end{xy}$$
in $\Gamma_\chi$ such that $\t(X_{i-1}) \oplus X_{i+1}$ is a summand of the whole middle term of $\E_\chi(X_i)$. Thus by Corollary \ref{infdegree} the arrow $X_{i+1} \to X_i$ has infinite left degree. \qed

\begin{Lemma}\label{A_inftydegree}
Let
$$\begin{xy}
  \xymatrix{ 
     \cdots \ar[r] &	X_{-n} \ar[r] 	& \cdots \ar[r] &	X_{-1} \ar[r]     &  X_0 \ar[r] &  X_1 \ar[r] & \cdots \ar[r] & X_n \ar[r] & \cdots	    }
\end{xy}$$ 
be a be-infinite pre-sectional path in $\Gamma_\chi$.
\begin{enumerate}[(a)]
 \item If all $X_i$ are left stable, then all arrows $\tau_\chi^j(X_i) \to \tau_\chi^j(X_{i+1})$ and $\tau_\chi^{j+1}(X_{i+1})$ $\to \tau_\chi^j(X_i)$ with $j \geq 0$ and $i \in \Z$ have infinite left degree.
 \item Dually, if all $X_i$ are right stable, then all arrows $\tau_\chi^j(X_i) \to \tau_\chi^j(X_{i+1})$ and $\tau_\chi^{j+1}(X_{i+1}) \to \tau_\chi^j(X_i)$ with $j \leq 0$ and $i \in \Z$ have infinite right degree.
\end{enumerate}

\end{Lemma}

\proof
\\Assume that all $X_i$ are left stable. There are infinite pre-sectional paths
$$\begin{xy}
  \xymatrix{ 
     \cdots \ar[r] &	\tau_\chi^j(X_{i-2}) \ar[r] 	& \tau_\chi^j(X_{i-1}) \ar[r]     &\tau_\chi^j(X_i)	 \ar[r] & \tau_\chi^j(X_{i+1})   }
\end{xy}$$ 
and
$$\begin{xy}
  \xymatrix{ 
     \cdots \ar[r] &	\tau_\chi^{j+2}(X_{i+3}) \ar[r] 	& \tau_\chi^{j+1}(X_{i+2}) \ar[r]     &\tau_\chi^j(X_{i+1})	    }
\end{xy}$$ 
in $\Gamma_\chi$ for each $j \geq 0, i \in \Z$ such that $\tau_\chi^{j+1}(X_{i+1}) \oplus \tau_\chi^j(X_{i-1})$ and $\tau_\chi^{j+1}(X_{i+2}) \oplus \tau_\chi^j(X_i)$ are direct summands of the whole middle terms of $\E_\chi(\tau_\chi^j(X_i))$ and $\E_\chi(\tau_\chi^j(X_{i+1}))$, respectively. Then the statement follows immediately from Corollary \ref{infdegree}. \qed

\begin{Lemma}\label{pre-sectionalnon-zero}
Let 
$$\begin{xy}
  \xymatrix{ 
   X_0 \ar[r] &	X_1 \ar[r] 	& \cdots \ar[r] &	X_{n-1} \ar[r]     &  X_n	    }
\end{xy}$$ 
be a pre-sectional path in $\Gamma_\chi$ such that $X_i$ is not Ext-projective for $i = 1, \ldots, n$ or not Ext-injective for $i = 0, \ldots, n-1$. Then there are $\chi$-irreducible morphisms $f_i: X_{i-1} \to X_i$ such that $f_n \cdots f_1 \notin \rad_\chi^{n+1}(X_0,X_n)$. In particular, $f_n \cdots f_1$ is non-zero.
\end{Lemma}

\proof
\\Suppose none of the $X_i$ is Ext-projective for $i = 1, \ldots, n$. For convenience, we define $\t(X_{n+1}) = 0$, so that $X_{i-1} \oplus \t(X_{i+1})$ is always a direct summand of the domain of a minimal right almost split morphism mapping to $X_i$. We use induction to prove the following statement from which we conclude the lemma.
\bigskip
\\For each $i \in \{1, \ldots, n\}$ there is an irreducible morphism $(f_i,g_i) : X_{i-1} \oplus \t(X_{i+1}) \to X_i$ such that $f_i \cdots f_1 + g_ip_{i-1} \notin \rad_\chi^{i+1}(X_0,X_i)$ for every morphism $p_{i-1} : X_0 \to \t(X_{i+1})$ in $\chi$.
\bigskip
\\For $i=1$ we choose $(f_1,g_1) : X_0 \oplus \t(X_2) \to X_1$ to be the restriction of a minimal right almost split morphism mapping to $X_1$. Let $p_0: X_0 \to \t(X_2)$ be an arbitrary morphism. If $X_0 \ncong \t(X_2)$, it is clear that $f_1 + g_1p_0 \notin \rad_\chi^2(X_0,X_1)$. If $X_0 \cong \t(X_2)$, then $f_1$ and $g_1$ are linearly independent in $\Irr_\chi(X_0,X_1)$ as a $T_{X_0}^{op}$-module by Theorem \ref{radbasis}, hence $f_1 + g_1p_0 \notin \rad_\chi^2(X_0,X_1)$.
\bigskip
\\Suppose now $1 < i < n$ and we have found a $\chi$-irreducible morphism $(f_i,g_i)$ as required. We may choose $p_i = 0$ and obtain $f_i \cdots f_1 \notin \rad_\chi^{i+1}(X_0, X_i)$. Since $X_i \oplus \t(X_{i+2})$ is a direct summand of the domain of a minimal right almost split morphism mapping to $X_{i+1}$, we know there are $\chi$-irreducible morphisms $(g_i,h_i)^T : \t(X_{i+1}) \to X_i \oplus \t(X_{i+2})$ and $(f_{i+1},g_{i+1}) : X_i \oplus \t(X_{i+2}) \to X_{i+1}$ such that $\{(g_i,h_i)^T,(f_{i+1},g_{i+1})\}$ is a component of the almost split sequence $\E_\chi(X_{i+1})$. Suppose there is a $p_i : X_0 \to \t(X_{i+2})$ such that $f_{i+1} \cdots f_1 + g_{i+1}p_i \in \rad_\chi^{i+2}(X_0, X_{i+1})$. Then by Lemma \ref{degree1} there is a $p_{i-1} : X_0 \to \t(X_{i+1})$ such that $(f_i \ldots f_1,p_i) + (g_i,h_i)p_{i-1} \in \rad_\chi^{i+1}(X_0, X_i \oplus \t(X_{i+2}))$. In particular, $f_i \ldots f_1 + g_ip_{i-1} \in \rad_\chi^{i+1}(X_0, X_i)$, which contradicts our inductive assumption. The case where 
none of 
the $X_i$ is Ext-injective for $i = 0, \ldots, n-1$ can be proved dually. \qed

\begin{Cor}
Let
$$\begin{xy}
  \xymatrix{ 
   X_0 \ar[r]^{f_1} &	X_1 \ar[r]^{f_2} 	& \cdots \ar[r]^{f_{n-1}} &	X_{n-1} \ar[r]^{f_n}     &  X_n	    }
\end{xy}$$ 
be a pre-sectional path in $\Gamma_\chi$ such that $X_i$ is not Ext-projective for $i = 1, \ldots, n$ or not Ext-injective for $i = 0, \ldots, n-1$. Then for each given constant $b$ there are less than $2^b$ integers $i \in \{0, \ldots, n\}$ such that $l(X_i) \leq b$.
\end{Cor}
\proof
\\Suppose that each $X_i$ is not Ext-projective for $i = 1, \ldots, n$. Without loss of generality, we assume that $f_n \cdots f_1 \neq 0$ by Lemma \ref{pre-sectionalnon-zero}. Suppose the statement is not true, i.e. there is a positive integer $b$ such that there are at least $2^b$ modules $X_{i_1}, X_{i_2}, \ldots, X_{i_{2^b}}$ with $i_j < i_{j+1}$ for $j = 1, \ldots, 2^b-1 $ and $l(X_{i_j}) \leq b$ for $j = 1, \ldots, 2^b$. We set $g_j = f_{i_{j+1}} \cdots f_{i_j+1} : X_{i_j} \to X_{i_{j+1}}$ for $j = 1, \ldots, 2^b-1$, so the composition $g_{2^b-1} \cdots g_1 : X_{i_1} \to X_{i_{2^b}}$ is zero by Lemma \ref{HaradaSai}. By construction this is a contradiction to $f_n \cdots f_1 \neq 0$. Again, the other case is proved dually. \qed

\begin{Cor}
Every pre-sectional cycle in $\Gamma_\chi$ contains at least one Ext-projective and one Ext-injective module.
\end{Cor}
\proof
\\A pre-sectional cycle in an Auslander-Reiten quiver is a pre-sectional path 
$$\begin{xy}
  \xymatrix{ 
      X_0 \ar[r]^{f_1} 	&	X_1 \ar[r]^{f_2}     &   \cdots	 \ar[r]	&	X_{n-1} \ar[r]^{f_n}	&	X_n = X_0     }
\end{xy}$$
from some indecomposable $X_0$ in $\chi$ to itself such that the composition of the path with itself is again pre-sectional. Let $f = f_n \cdots f_1$, then by the lemma of Harada and Sai there is an integer $k$ such that $f^k = 0$. Suppose none of the $X_i$ is Ext-projective, then $f^k \neq 0$ by Lemma \ref{pre-sectionalnon-zero}, which is a contradiction. If we assume that none of the $X_i$ is Ext-injective, we can apply the same arguments to 
$$\begin{xy}
  \xymatrix{ 
      \tau_\chi^{-1}(X_0) \ar[r]^{f_1} 	&	\tau_\chi^{-1}(X_1) \ar[r]^(.55){f_2}     &   \cdots	 \ar[r]	&	\tau_\chi^{-1}(X_{n-1}) \ar[r]^(.4){f_n}	&	\tau_\chi^{-1}(X_n) = \tau_\chi^{-1}(X_0).     }
\end{xy}$$
\qed

\section{The Happel-Preiser-Ringel theorem for subcategories}
Recall that an \textbf{oriented cycle} is a path in an Auslander-Reiten quiver that starts and ends in the same module.
\begin{Lemma}\label{cycledegree}
Every oriented cycle in $\Gamma_\chi$ contains both an arrow of finite left degree and an arrow of finite right degree.
\end{Lemma}

\proof
\\Let 
$$\begin{xy}
  \xymatrix{ 
     X_0 \ar[r]^{f_1} &	X_1 \ar[r]^{f_2} 	& \cdots \ar[r] & X_{n-1} \ar[r]^{f_n}    & X_n = X_0	    }
\end{xy}$$ 
be an oriented cycle in $\Gamma_\chi$. By Lemma \ref{HaradaSai} there is a $k \geq 1$ such that $(f_n \cdots f_1)^k=0$, which completes the proof. \qed

\begin{Th}
 
Let $\Gamma$ be a connected component of $\Gamma_\chi$ that is left stable or right stable and assume every almost split sequence in $\Gamma$ has at least 2 middle terms. Then for every path 
$$\begin{xy}
  \xymatrix{ 
     X_0 \ar[r]^{f_1} &	X_1 \ar[r]^{f_2} 	& \cdots \ar[r] & X_{n-1} \ar[r]^{f_n}     & X_n	    }
\end{xy}$$ 
in $\Gamma$ we have $f_n \cdots f_1 \notin \rad_\chi^{n+1}(X_0,X_n)$. In particular, $f_n \cdots f_1$ is non-zero.
\end{Th}

\proof
\\Suppose there are no Ext-projective modules in $\Gamma$. We prove by induction that for each $n \in \N$ the left degree of each irreducible morphism $f: X \to Y$ between indecomposable modules $X$ and $Y$ in $\Gamma$ is greater then $n$. Clearly, $f$ cannot be a surjective right almost split morphism because $\E_\chi(Y)$ has at least two middle terms. Hence $d_\chi^l(f) > 1$ by \ref{degree1Cor}.
\bigskip
\\We assume that the statement is true for $n - 1$, but false for $n$, i.e. there exists an irreducible morphism $f: X \to Y$ such that $d_\chi^l(f) = n$. Since every almost split sequence has at least two middle terms, there is a module $X'$ such that $X \oplus X'$ is a direct summand of the whole middle term of $\E_\chi(Y)$. So by Corollary \ref{degree1Cor} there is an irreducible morphism $f' : \tau_\chi(Y) \to X'$ with $d_\chi^l(f') < n$, which is a contradiction.
\bigskip
\\We have proved that every arrow in $\Gamma$ has infinite left degree and hence every path 
$$\begin{xy}
  \xymatrix{ 
     X_0 \ar[r]^{f_1} &	X_1 \ar[r]^{f_2} 	& \cdots \ar[r] & X_{n-1} \ar[r]^{f_n}     & X_n	    }
\end{xy}$$ 
in $\Gamma$ of length $n$ is not in $\rad_\chi^{n+1}(X_0, X_n)$. If $\Gamma$ does not contain any Ext-injective modules, the theorem can be proved dually. \qed

\begin{Lemma}\label{pathexistence}
Let $\Gamma$ be a connected stable subquiver of $\Gamma_\chi$ and let $X$ and $Y$ be modules in $\Gamma$. Then the following holds.
\begin{enumerate}[(a)]
\item There is a path in $\Gamma$ from $X$ to $\tau_\chi^n(Y)$ for some $n \in \Z$.
\item If there is a path from $X$ to $Y$ in $\Gamma$, then either $X = \tau_\chi^n(Y)$ for some $n \geq 1$ or there is a sectional path in $\Gamma$ from $X$ to $\tau_\chi^n(Y)$ for some $n \geq 0$.
\end{enumerate}
\end{Lemma}


\proof
\\Since $\Gamma$ is connected, there is a walk 
$$\begin{xy}
  \xymatrix{ 
    X = X_0 \ar@{-}[r] &	X_1 \ar@{-}[r] 	& \cdots \ar@{-}[r] & X_{s-1} \ar@{-}[r]     & X_s	= Y    }
\end{xy}$$ 
in $\Gamma$. We prove the first statement by induction on $s$. For $s = 1$ the statement is trivial. Suppose now $s > 0$ and there exists a path in $\Gamma$ from $X$ to $\t^n(X_{s-1})$ for some integer $n$. Since there is an edge between $X_{s-1}$ and $Y$, there is either an arrow $X_{s-1} \to Y$ or $Y \to X_{s-1}$. In the first case there is also an arrow from $\t^n(X_{s-1}) \to \t^n(Y)$ as $\Gamma$ is stable and we obtain a path from $X$ to $\t^n(Y)$ in $\Gamma$. Otherwise there is an arrow $\t^n(X_{s-1}) \to \t^{n-1}(Y)$ and a path from $X$ to $\t^{n-1}(Y)$ again using that $\Gamma$ is stable.
\bigskip 
\\In order to prove the second statement, let
$$\begin{xy}
  \xymatrix{ 
   X =  X_0 \ar[r] &	X_1 \ar[r] 	& \cdots \ar[r] & X_{s-1} \ar[r]     & X_s = Y	    }
\end{xy}$$ 
be a path in $\Gamma$. Again the statement is trivial for $s = 1$. Assume $s > 0$, then by the inductive hypothesis we either have $X = \t^n(X_{s-1})$ for some $n \geq 1$ or there is a sectional path 
$$\begin{xy}
  \xymatrix{ 
   X =  Y_0 \ar[r] &	Y_1 \ar[r] 	& \cdots \ar[r] & Y_{k-1} \ar[r]     & Y_k = \t^m(X_{s-1})  	    }
\end{xy}$$ 
for some $m \geq 0$. In the first case $X \to \t^n(Y)$ is a sectional path. In the second case if $Y_{k-1} = \t^{m+1}(Y)$, then either $X = \t^{m+1}(Y)$ with $m + 1 > 0$ or 
$$\begin{xy}
  \xymatrix{ 
   X =  Y_0 \ar[r] &	Y_1 \ar[r] 	& \cdots \ar[r] & Y_{k-1} = \t^{m+1}(Y)  	    }
\end{xy}$$ 
is a sectional path. If $Y_{k-1} \neq \t^{m+1}(Y)$, then 
$$\begin{xy}
  \xymatrix{ 
   X =  Y_0 \ar[r] &	Y_1 \ar[r] 	& \cdots \ar[r] & Y_{k-1} \ar[r]     & Y_k \ar[r] & \t^m(Y)  	    }
\end{xy}$$ 
is a sectional path. \qed
\bigskip
\\In order to proof the main theorem of this section, we need another definition and lemma.

\begin{Def}
Let $X$ and $Y$ be indecomposable modules in $\chi$. Then the set $\{\t^n(X) | n \in \Z\}$ is called the \textbf{$\pmb{\t}$-orbit} of $X$. 
We call a $\t$-orbit left stable, right stable, stable or periodic if all the modules it contains are left stable, right stable, stable or periodic respectively. A $\t$-orbit that is neither left stable nor right stable is a finite set and hence called a finite $\t$-orbit. Moreover, we say two orbits $\{\t^n(X) | n \in \Z\}$ and $\{\t^n(Y) | n \in \Z\}$ are \textbf{adjacent} if there are integers $i,j \in \Z$ such that there is an irreducible morphism $f: \t^i(X) \to \t^j(Y)$ or $g: \t^j(Y) \to \t^i(X)$.
\end{Def}

It follows directly from the definitions that if one module $X$ is left stable, right stable, stable or periodic, then all modules in the $\t$-orbit of $X$ are left stable, right stable, stable or periodic respectively. Note that although periodic $\t$-orbits only contain a finite number of non-isomorphic indecomposable modules, we do not call them finite $\t$-orbits.

\begin{Lemma}\label{orbit}
Let $\Gamma_\chi$ be the Auslander-Reiten quiver of $\chi$ with a module $X$ in a periodic $\t$-orbit, i.e. there is an $n \in \N$ such that $\t^n(X) = X$. Then any $\t$-orbit adjacent to the $\t$-orbit of $X$ is either periodic or finite.
\end{Lemma}

\proof
\\Suppose there is an $\t$-orbit $\{\t^n(Y) | n \in \Z\}$ that is adjacent to the $\t$-orbit of $X$ which is neither periodic nor finite. Without loss of generality, we assume that this $\t$-orbit contains a module $Y$ such that there is an arrow from $X$ to $Y$ in the Auslander-Reiten quiver of $\chi$ and that $\t^m(Y)$ exists for all $m \geq 0$. By our assumption we know that $\tau_\chi^m(Y) \neq Y$ for all $m \geq 1$. Since both orbits are left stable, we conclude that there is an arrow from $\t^n(X) = X$ to $\t^n(Y)$. Proceeding like this it is easy to see that there are arrows from $X$ to $Y, \t^n(Y), \t^{2n}(Y), \ldots$, i.e. the almost split sequence $\E_\chi'(X)$ has infinitely many middle terms, which is impossible.
\qed

\begin{Def}
A stable subquiver $\Gamma$ of $\Gamma_\chi$ is called a \textbf{stable tube} if every module in $\Gamma$ is $\t$-periodic and it contains an infinite sectional path 
$$\begin{xy}
  \xymatrix{ 
     \cdots \ar[r] &	X_n \ar[r] 	&	X_{n-1} \ar[r]     &   \cdots	 \ar[r]	&	X_1 \ar[r]	&	X_0 = Y 	    }
\end{xy}$$
such that each orbit in $\Gamma$ is generated by a unique module $X_i$.
\end{Def}

\newpage

\begin{Th}[Happel-Preiser-Ringel Theorem]\label{HPRT}
Let $\Gamma$ be an infinite connected stable subquiver of $\Gamma_\chi$. If there is a $\tau_\chi$-periodic module in $\Gamma$, then $\Gamma$ is a stable tube.
\end{Th}

\proof
\\Since $\Gamma$ is stable and connected, if one module in $\Gamma$ is $\tau_\chi$-periodic, then every module in $\Gamma$ is $\tau_\chi$-periodic by Lemma \ref{orbit}. Thus there must be infinitely many $\t$-orbits in $\chi$. Fix a module $X_0$ in $\Gamma$, then for every module $X$ that is not in the same $\t$-orbit as $X_0$ there is a sectional path in $\Gamma$ from $\tau_\chi^n(X)$ to $X_0$ for some $n \in \Z$ by Lemma \ref{pathexistence}. Consequently, there are arbitrarily long sectional paths in $\Gamma$ ending in $X_0$ and hence there exists an infinite sectional path
$$\begin{xy}
  \xymatrix{ 
 \gamma:	&   \cdots \ar[r] &	X_{i+1} \ar[r] 	&	X_i \ar[r]     &   \cdots	 \ar[r]	&	X_1 \ar[r]	&	X_0	    }
\end{xy}$$
in $\Gamma$. For each integer $i \geq 0$, there is an oriented cycle
$$\begin{xy}
  \xymatrix{ 
     X_{i+1} = \tau_\chi^{n_i}(X_{i+1}) \ar[r] &	\tau_\chi^{n_i}(X_{i}) \ar[r] 	&	\tau_\chi^{n_i-1}(X_{i+1}) \ar[r]     &   \cdots	   }
\end{xy}$$
$$\begin{xy}
  \xymatrix{ 
        \cdots	 \ar[r]	&	\tau_\chi(X_{i}) \ar[r]	&	X_{i+1}	    }
\end{xy}$$
in $\Gamma$. By Lemma \ref{cycledegree} there is an integer $j_i \geq 0$ such that $\tau_\chi^{j_i}(X_{i+1}) \to \tau_\chi^{j_i}(X_{i})$ or $\tau_\chi^{j_i+1}(X_{i}) \to \tau_\chi^{j_i}(X_{i+1})$ has finite left degree. As a consequence, every arrow in $\gamma$ has trivial valuation and $X_i$ has two immediate predecessors $X_{i+1}$ and $\tau_\chi(X_{i-1})$ and two immediate successors $X_{i-1}$ and $\tau_\chi^{-1}(X_{i+1})$ for all integers $i \geq 1 $ by Lemma \ref{vmtinfinity}. Moreover, $\gamma$ cannot be contained in a bi-infinite sectional path
$$\begin{xy}
  \xymatrix{ 
     \cdots \ar[r] &	X_{-n} \ar[r] 	& \cdots \ar[r] &	X_{-1} \ar[r]     &  X_0 \ar[r] &  X_1 \ar[r] & \cdots \ar[r] & X_n \ar[r] & \cdots	    }
\end{xy}$$ 
by Lemma \ref{A_inftydegree}. Without loss of generality, we can assume that $\gamma$ is a maximal sectional path, i.e. $X_0$ has only one immediate predecessor $X_1$ and one immediate successor $\tau_\chi^{-1}(X_1)$. Therefore, $\gamma$ contains at least one module of every $\t$-orbit and all $X_i$ have the same $\tau_\chi$-periodicity.
\bigskip
\\It remains to prove that the $X_i$ belong to pairwise different $\tau_\chi$-orbits. So let $X_m$ and $X_j$ be in the same $\t$-orbit with $m < j$ and $j - m$ minimal. Then $X_{j-1}$ and $X_{m-1}$ also belong to the same $\t$-orbit as every module has at most two successors and we have chosen $j - m$ to be minimal. We continue inductively and obtain that $X_0$ and $X_{j-m}$ are contained in the same $\t$-orbit. Since these modules only have 1 predecessor and successor respectively, $X_1$ and $X_{j-m-1}$ must belong to the same $\t$-orbit as well, which is a contradiction to the minimality of $j - m$. \qed


\begin{Lemma}\label{sectionalfinite}
Let $\Gamma$ be a connected stable subquiver of $\Gamma_\chi$ and let $X$ be a module in $\Gamma$. If there is a sectional path in $\Gamma$ from $X \to \tau_\chi^n(X)$ for some $n \geq 1$, then $\Gamma$ is a finite subquiver of $\tau_\chi$-periodic modules.
\end{Lemma}

\proof
\\Let
$$\begin{xy}
  \xymatrix{ 
    X = X_0 \ar[r] &	X_1 \ar[r] 	&	\cdots  \ar[r]     &   X_{k-1} \ar[r]	&	X_k = \tau_\chi^n(X) 	    }
\end{xy}$$
be a sectional path in $\Gamma$. Without loss of generality, we assume the path to be minimal, i.e. if $k \geq 3$ then $X_{k-1} \neq \tau_\chi^m (X_1)$ for all $m \geq 0$. We construct a bi-infinite path
$$\begin{xy}
  \xymatrix{ 
    \gamma: & \cdots \ar[r] &  \tau_\chi^{-2n}(X_{k-1}) \ar[r] &	\tau_\chi^{-n}(X_0) \ar[r] 	& \tau_\chi^{-n}(X_1) \ar[r] &	\cdots      }
\end{xy}$$
$$\begin{xy}
  \xymatrix{ 
     \cdots \ar[r] &  \tau_\chi^{-n}(X_{k-1}) \ar[r] & X_0 \ar[r] & X_1 \ar[r] & \cdots       }
\end{xy}$$
$$\begin{xy}
  \xymatrix{ 
 \cdots \ar[r] &  X_{k-1} \ar[r] & 	\tau_\chi^n(X_0) \ar[r]	&	\tau_\chi^n(X_1) \ar[r]	&	\cdots     }
\end{xy}$$

in $\Gamma$. Note that there is also an oriented cycle in $\Gamma$ as follows
$$\begin{xy}
  \xymatrix{ 
      X_0 \ar[r]& X_1 \ar[r] & \cdots \ar[r] & X_{k-1} \ar[r] &	X_k \ar[r] 	& \tau_\chi^{-1}(X_{k-1}) \ar[r] & \tau_\chi^{-1}(X_k) \ar[r] &	\cdots   }
\end{xy}$$
$$\begin{xy}
  \xymatrix{ 
    	\cdots \ar[r]	&	\tau_\chi^{-n}(X_{k-1}) \ar[r]	&	\tau_\chi^{-n}(X_k)  = X_0  }
\end{xy}$$
By Lemma \ref{cycledegree} there are some $0 \leq i \leq k - 1$ and $0 \leq j \leq n $ such that $\tau_\chi^{-j}(X_i) \to \tau_\chi^{-j}(X_{i+1})$ or $\tau_\chi^{-j}(X_{i+1}) \to \tau_\chi^{-j-1}(X_i)$ has finite left degree. Then by Lemma \ref{infdegree} $\gamma$ cannot be sectional, which gives us $1 \leq k \leq 2$. If $k = 1$, then there is a non-sectional path
$$\begin{xy}
\xymatrixrowsep{0in}
  \xymatrix{ 
X_0 \ar[r] & \t^n(X_0) \ar[r] & \t^{2n}(X_0)   }
\end{xy}$$
and, consequently, $X_0 = \t^{2n+1}(X_0)$. Similarly, if $k = 2$, then the path
$$\begin{xy}
\xymatrixrowsep{0in}
  \xymatrix{ 
X_1 \ar[r] & \t^n(X_0) \ar[r] & \t^{n}(X_1)   }
\end{xy}$$
is not sectional and we have $X_1 = \t^{n+1}(X_1)$. In both cases there are periodic modules in $\Gamma$. If $\Gamma$ is infinite, then it is a stable tube by Theorem \ref{HPRT} and there is no sectional path $$\begin{xy}
  \xymatrix{ 
    X = X_0 \ar[r] &	X_1 \ar[r] 	&	\cdots  \ar[r]     &   X_{k-1} \ar[r]	&	X_k = \tau_\chi^n(X). 	    }
\end{xy}$$
Thus $\Gamma$ must be a finite subquiver consisting of $\tau_\chi$-periodic modules.\qed

\begin{Th}\label{orientedperiodic}
Let $\Gamma$ be a connected stable subquiver of $\Gamma_\chi$. Then $\Gamma$ contains an oriented cycle if and only if it consists of $\t$-periodic modules.
\end{Th}

\proof
\\Clearly, $\Gamma$ contains an oriented cycle if it consists of $\tau_\chi$-periodic modules. On the other hand, suppose $\Gamma$ contains a non-trivial oriented cycle from $X$ to itself. Then, by Lemma \ref{pathexistence}, there is either an $n \geq 1$ such that $X = \tau_\chi^n X$ or there is a sectional path from $X$ to $\t^n(X)$ where $n \geq 0$. Hence $\Gamma$ contains $\t$-periodic modules by \ref{sectionalfinite}. \qed

\begin{Cor}
Let $X$ be a module in $\Gamma_\chi$ such that there is an oriented cycle 
$$\begin{xy}
  \xymatrix{ 
    X = X_0 \ar[r] &	X_1 \ar[r] 	&	\cdots  \ar[r]     &   X_{n-1} \ar[r]	&	X_n = X 	    }
\end{xy}$$
such that all $X_i$ are stable modules. Then all $X_i$ are $\t$-periodic modules.
\end{Cor}

\proof
\\We obtain the result by applying Theorem \ref{orientedperiodic} to the connected stable subquiver generated by the $\t$-orbits of the $X_i$. \qed

\section{Global degrees of irreducible morphisms}

\begin{Def}
Let $X \to Y$ be an arrow in $\Gamma_\chi$.
\begin{enumerate}[(a)]
\item Let $X$ and $Y$ be left stable modules. We then define the \textbf{global left degree} in $\chi$ of the arrow $X \to Y$ to be the minimum of left degrees in $\chi$ of all arrows $\t^{n+1}(Y) \to \t^n(X)$ and $\t^n(X) \to \t^n(Y)$ for all $n \geq 0$.
\item Dually, if $X$ and $Y$ are right stable modules, we define the \textbf{global right degree} in $\chi$ of the arrow $X \to Y$ to be the minimum of right degrees in $\chi$ of all arrows $\t^{n}(Y) \to \t^{n-1}(X)$ and $\t^n(X) \to \t^n(Y)$ for all $n \leq 0$.
\end{enumerate}
\end{Def}

These definitions allow us to rephrase Lemmas \ref{vmtinfinity} and \ref{A_inftydegree} as follows.

\begin{Lemma}\label{vmtinfinityglobal}
\mbox{}
\begin{enumerate}[(a)]
 \item Let  
$$\begin{xy}
  \xymatrix{ 
     \cdots \ar[r] &	X_{i+1} \ar[r] 	&	X_i \ar[r]     &   \cdots	 \ar[r]	&	X_1 \ar[r]	&	X_0	    }
\end{xy}$$
be a sectional path in $\Gamma_\chi$ with all $X_i$ left stable. If the path contains infinitely many arrows with finite global left degree in $\chi$, then for each integer $i \geq 0$ the arrow $X_{i+1} \to X_i$ has trivial valuation and the module $X_i$ has at most two left stable immediate predecessors in $\Gamma_\chi$.
\item Let  
$$\begin{xy}
  \xymatrix{ 
     Y_0 \ar[r] &	Y_1 \ar[r] 	&	\cdots	 \ar[r]	&	Y_i \ar[r]	&	Y_{i+1} \ar[r] &\cdots 	    }
\end{xy}$$
be a sectional path in $\Gamma_\chi$ with all $Y_i$ right stable. If the path contains infinitely many arrows with finite global right degree in $\chi$, then for each integer $i \geq 0$ the arrow $Y_i \to Y_{i+1}$ has trivial valuation and the module $Y_i$ has at most two right stable immediate successors in $\Gamma_\chi$.
\end{enumerate}

\end{Lemma}

\begin{Lemma}\label{A_inftydegreeglobal}
Let
$$\begin{xy}
  \xymatrix{ 
     \cdots \ar[r] &	X_{-n} \ar[r] 	& \cdots \ar[r] &	X_{-1} \ar[r]     &  X_0 \ar[r] &  X_1 \ar[r] & \cdots \ar[r] & X_n \ar[r] & \cdots	    }
\end{xy}$$ 
be a bi-infinite sectional path in $\Gamma_\chi$. 
\begin{enumerate}[(a)]
\item If all modules $X_i$ are left stable in $\chi$ then all arrows $X_{i+1} \to X_i$ have infinite global left degree in $\chi$.
\item Dually, if all modules $X_i$ are right stable in $\chi$, then all arrows $X_{i+1} \to X_i$ have infinite global right degree in $\chi$.
\end{enumerate}

\end{Lemma}

The following lemma is an easy consequence of Corollary \ref{infdegree}.

\newpage

\begin{Lemma}\label{twosectionalpathsinfiniteglobaldegree}
\mbox{}
\begin{enumerate}[(a)]
 \item Let 
$$\begin{xy}
  \xymatrix{ 
     \cdots \ar[r] &	X_{i+1} \ar[r] 	&	X_i \ar[r]     &   \cdots	 \ar[r]	&	X_1 \ar[r]	&	X_0	    }
\end{xy}$$
and
$$\begin{xy}
  \xymatrix{ 
     \cdots \ar[r] &	Y_{i+1} \ar[r] 	&	Y_i \ar[r]     &   \cdots	 \ar[r]	&	Y_1 \ar[r]	&	Y_0	    }
\end{xy}$$
be infinite sectional paths in $\Gamma_\chi$ containing only left stable modules. If $X_0 = Y_0$ and $X_1 \neq Y_1$, then both arrows $X_{i+1} \to X_i$ and $Y_{i+1} \to Y_i$ have infinite global left degree in $\chi$ for each $i \geq 0$ where $Z = X_0 = Y_0$.
\item Let 
$$\begin{xy}
  \xymatrix{ 
    X_0 \ar[r] &	X_1 \ar[r] 	&	\cdots	 \ar[r]	&	X_i \ar[r]	&	X_{i+1} \ar[r] &\cdots 	    }
\end{xy}$$
and
$$\begin{xy}
  \xymatrix{ 
    Y_0 \ar[r] &	Y_1 \ar[r] 	&	\cdots	 \ar[r]	&	Y_i \ar[r]	&	Y_{i+1} \ar[r] &\cdots 	    }
\end{xy}$$
be infinite sectional paths in $\Gamma_\chi$ containing only right stable modules. If $X_0 = Y_0$ and $X_1 \neq Y_1$, then both arrows $X_i \to X_{i+1}$ and $Y_i \to Y_{i+1}$ have infinite global right degree in $\chi$ for each $i \geq 0$.
\end{enumerate}

\end{Lemma}

\proof
\\Clearly, it is sufficient to prove the statement for an arrow $X_i \to X_{i-1}$. For any arrow $\t^n(X_i) \to \t^n(X_{i-1})$ we can construct an infinite sectional path
$$\begin{xy}
  \xymatrix{ 
     \cdots \ar[r]& \t^{n+i}(Y_2) \ar[r] &	\t^{n+i}(Y_1) \ar[r] &	\t^{n+i}(Z) \ar[r] 	&	\t^{n+i-1}(X_1) \ar[r]     &   \cdots	     }
\end{xy}$$
$$\begin{xy}
  \xymatrix{ 
      \cdots \ar[r] &	 \t^{n+2}(X_{i-2})	\ar[r]	&	\t^{n+1}(X_{i-1}) \ar[r]	&	\t^n(X_{i}).	    }
\end{xy}$$
It then follows from Corollary \ref{infdegree} that the arrow $\t^n(X_i) \to \t^n(X_{i-1})$ has infinite left degree in $\chi$. Similarly, we use the sectional path 
$$\begin{xy}
  \xymatrix{ 
     \cdots \ar[r]& \t^{n}(X_{i+2}) \ar[r] &	\t^{n}(X_{i+1}) \ar[r] 	&	\t^{n}(X_i) \ar[r]     &  \t^{n}(X_{i-1}) 	     }
\end{xy}$$
to verify that $\t^n(X_{i-1}) \to \t^{n-1}(X_i)$ also has infinite left degree in $\chi$ again by Corollary \ref{infdegree}. Hence by definition the arrow $X_i \to X_{i-1}$ has infinite global left degree in $\chi$.
\qed

\begin{Lemma}\label{leftpathexistence}
\mbox{}
\begin{enumerate}[(a)]
 \item Let $\Gamma$ be a connected left stable subquiver of $\Gamma_\chi$ and let $X$ and $Y$ be modules in $\Gamma$. If there is a non-trivial path from $X$ to $Y$ in $\Gamma$, then either $X = \tau_\chi^n(Y)$ for some $n \geq 1$ or there is a sectional path in $\Gamma$ from $X$ to $\tau_\chi^n(Y)$ for some $n \geq 0$.
 \item  Dually, if $\Gamma$ is a connected right stable subquiver of $\Gamma_\chi$ and there is a non-trivial path from $X$ to $Y$ in $\Gamma$, then either $X = \tau_\chi^n(Y)$ for some $n \geq 0$ or there is a sectional path in $\Gamma$ from $\tau_\chi^n(X)$ to $Y$ for some $n \leq 0$.
\end{enumerate}

\end{Lemma}

\newpage

\proof
\\We prove the Lemma by induction on $s$. Let
$$\begin{xy}
  \xymatrix{ 
   X =  X_0 \ar[r] &	X_1 \ar[r] 	& \cdots \ar[r] & X_{s-1} \ar[r]     & X_s = Y	    }
\end{xy}$$ 
be a path in $\Gamma$. The statement is trivial for $s = 1$. Assume $s > 0$, then by the inductive hypothesis we either have $X = \tau_\chi^n(X_{s-1})$ for some $n \geq 1$ or there is a sectional path 
$$\begin{xy}
  \xymatrix{ 
   X =  Y_0 \ar[r] &	Y_1 \ar[r] 	& \cdots \ar[r] & Y_{k-1} \ar[r]     & Y_k = \tau_\chi^m(X_{s-1})  	    }
\end{xy}$$ 
for some $m \geq 0$. In the first case $X \to \tau_\chi^n(Y)$ is a sectional path. In the second case if $Y_{k-1} = \tau_\chi^{m+1}(Y)$, then either $X = \tau_\chi^{m+1}(Y)$ with $m + 1 > 0$ or 
$$\begin{xy}
  \xymatrix{ 
   X =  Y_0 \ar[r] &	Y_1 \ar[r] 	& \cdots \ar[r] & Y_{k-1} = \tau_\chi^{m+1}(Y)  	    }
\end{xy}$$ 
is a sectional path. If $Y_{k-1} \neq \tau_\chi^{m+1}(Y)$, then 
$$\begin{xy}
  \xymatrix{ 
   X =  Y_0 \ar[r] &	Y_1 \ar[r] 	& \cdots \ar[r] & Y_{k-1} \ar[r]     & Y_k \ar[r] & \tau_\chi^m(Y)  	    }
\end{xy}$$ 
is a sectional path. \qed

\begin{Lemma}\label{infinitesectionalpathinleftstablecomponent}
\mbox{}
\begin{enumerate}[(a)]
 \item Let $\Gamma$ be a left stable component of $\Gamma_\chi$ containing no $\t$-periodic module. If there is an oriented cycle in $\Gamma$, then there is an infinite sectional path
$$\begin{xy}
  \xymatrix{ 
     \cdots \ar[r]& \t^{2r}(X_2) \ar[r] &	\t^{2r}(X_1) \ar[r] &	\t^r(X_s) \ar[r] 	& \cdots	     }
\end{xy}$$
$$\begin{xy}
  \xymatrix{ 
   \cdots \ar[r] & \t^r(X_2) \ar[r] &  \t^r(X_1) \ar[r] & X_s \ar[r] &  \cdots	\ar[r]	&	X_2 \ar[r]	&	X_1	    }
\end{xy}$$
in $\Gamma$ containing infinitely many arrows of finite global left degree in $\chi$, where $r > s$, the $X_i$ belong to pairwise different $\t$-orbits and at least one of the $X_i$ is not stable.
 \item Dually, let $\Gamma$ be a right stable component of $\Gamma_\chi$ containing no $\t$-periodic module. If there is an oriented cycle in $\Gamma$, then there is an infinite sectional path
$$\begin{xy}
  \xymatrix{ 
    X_1 \ar[r]& X_2 \ar[r] & \cdots \ar[r] &	X_s \ar[r] &	\t^{-r}(X_1) \ar[r] 	& \t^{-r}(X_2) \ar[r] & \cdots	     }
\end{xy}$$
$$\begin{xy}
  \xymatrix{ 
   \cdots \ar[r] & \t^{-r}(X_s) \ar[r] &  \t^{-2r}(X_1) \ar[r] & \t^{-2r}(X_2) \ar[r] &  \cdots	    }
\end{xy}$$
in $\Gamma$ containing infinitely many arrows of finite global right degree in $\chi$, where $r > s$, the $X_i$ belong to pairwise different $\t$-orbits and at least one of the $X_i$ is not stable.
\end{enumerate}
\end{Lemma}

\proof
\\Suppose that $\Gamma$ contains an oriented cycle from $X$ to itself. Then there is a sectional path from $X$ to $\t^r(X)$ for some $r \geq 0$ by Lemma \ref{leftpathexistence}. Let $\gamma$ denote a sectional path
$$\begin{xy}
  \xymatrix{ 
    Y_1 \ar[r]& Y_2 \ar[r] & \cdots \ar[r] &	Y_{s-1} \ar[r] &	Y_s \ar[r] 	& Y_{s+1}	     }
\end{xy}$$
in $\Gamma$ of minimal positive length such that $Y_{s+1} = \t^n(Y_1)$ for some $n \geq 0$. We prove that $n > 0$. Suppose that $n = 0$, i.e. we have $Y_1 = Y_{s+1}$. By Corollary \ref{cycle} we know that if the composition 
$$\begin{xy}
  \xymatrix{ 
    Y_1 \ar[r]& Y_2 \ar[r] & \cdots \ar[r] &	Y_{s-1} \ar[r] &	Y_s \ar[r] 	& Y_1 \ar[r] & \cdots  }
\end{xy}$$
$$\begin{xy}
  \xymatrix{ 
\cdots \ar[r] 	& Y_{s-1} \ar[r] &	Y_s \ar[r] 	& Y_1 \ar[r] & \cdots  	  }
\end{xy}$$
was sectional, it would contain both an Ext-projective and an Ext-injective module. But by our assumptions $\Gamma$ is left stable and hence does not contain Ext-projective modules. It follows that the composition of $\gamma$ with itself is not sectional and, therefore, $\t(Y_2) = Y_s$. For $s = 1$ this composition is
$$\begin{xy}
  \xymatrix{ 
    Y_1 \ar[r]& Y_1 \ar[r] & Y_1 \ar[r] & \cdots     }
\end{xy}$$
while for $s = 2$ we trivially obtain $\t(Y_2) = Y_2$. So in both cases we have $\t(Y_s) = Y_s$, which contradicts our assumption that there is no $\t$-periodic module in $\Gamma$. It necessarily follows that $s \geq 3$, but then the path 
$$\begin{xy}
  \xymatrix{ 
    Y_2 \ar[r]& Y_3 \ar[r] & \cdots \ar[r] &	Y_s \ar[r] &	     }
\end{xy}$$
is of positive length and $\t(Y_2) = Y_s$, which contradicts $\gamma$ to be of minimal length with this property. We conclude that $n$ is strictly greater than $0$.
\bigskip
\\The next step is to show that all modules $Y_1, \ldots, Y_s$ belong to pairwise different $\t$-orbits. Since $\gamma$ is of minimal length, if $1 \leq i < j \leq s$, then $Y_j \neq \t^p(Y_i)$ for any $p \geq 0$. Hence if $Y_i$ an $Y_j$ belong to the same $\t$-orbit, then $Y_i = \t^q(Y_j)$ for some $q > 0$. By the minimality of $\gamma$ we know that $Y_s \neq \t^{n+1}(Y_2)$ and thus the path
$$\begin{xy}
  \xymatrix{ 
    Y_j \ar[r]& Y_{j+1} \ar[r] & \cdots \ar[r] &	Y_s \ar[r] &	\t^n(Y_1) \ar[r] 	& \t^n(Y_2) \ar[r] & \cdots  }
\end{xy}$$
$$\begin{xy}
  \xymatrix{ 
\cdots \ar[r] 	& \t^n(Y_{i-1}) \ar[r] & \t^n(Y_i) = \t^{n+q}(Y_j) 	  }
\end{xy}$$
must be sectional. However, it is of length $s - j + i < s$, which contradicts $\gamma$ to be of minimal length. Consequently, each $Y_i$ belongs to a different $\t$-orbit.
\bigskip
\\Now let $r = s + n$ and $X_i = \t^{i-1}(Y_i)$ for $i = 1, \ldots, s$. We then have $r > s$, the $X_i$ belong to pairwise different $\t$-orbits and $$\begin{xy}
  \xymatrix{ 
     \cdots \ar[r]& \t^{2r}(X_2) \ar[r] &	\t^{2r}(X_1) \ar[r] &	\t^r(X_s) \ar[r] 	& \cdots	     }
\end{xy}$$
$$\begin{xy}
  \xymatrix{ 
   \cdots \ar[r] & \t^r(X_2) \ar[r] &  \t^r(X_1) \ar[r] & X_s \ar[r] &  \cdots	\ar[r]	&	X_2 \ar[r]	&	X_1	    }
\end{xy}$$
is an infinite sectional path in $\Gamma$. For every $j \geq 0$ denote the sectional path
$$\begin{xy}
  \xymatrix{ 
   \t^{(j+1)r}(X_1) \ar[r] & \t^{jr}(X_s) \ar[r] &  \cdots	\ar[r]	&	\t^{jr}(X_2) \ar[r]	&	\t^{jr}(X_1)	    }
\end{xy}$$
by $\gamma_j$. Moreover, for each $j$ there are paths

$$\begin{xy}
  \xymatrix{ 
   \t^{jr+s}(X_1) \ar[r] & \t^{jr + s - 1}(X_2) \ar[r] &  \cdots	\ar[r]	&	\t^{jr + 1}(X_s) \ar[r]	&	\t^{(j + 1)r}(X_1)	    }
\end{xy}$$
and, since $r > s$,
$$\begin{xy}
\xymatrix{ 
\t^{(j+1)r}(X_1) \ar[r] & \t^{(j+1)r - 1}(X_2) \ar[r] & \t^{(j+1)r - 1}(X_1) \ar[r] &  \cdots	}
\end{xy}$$
$$\begin{xy}
\xymatrix{ 
\cdots 		\ar[r] &	\t^{jr + s + 1}(X_1) \ar[r]	&\t^{jr + s}(X_2) \ar[r] &\t^{jr + s}(X_1).	    }
\end{xy}$$
The last two paths form an oriented cycle through $\t^{(j+1)r}(X_1)$ that does not contain an Ext-projective module. By Lemma \ref{cycledegree} an oriented cycle in $\chi$ always contains an arrow of finite left degree in $\chi$, hence $\gamma_j$ always contains an arrow of finite global left degree in $\chi$, which can easily be seen in the picture on the next page.
\bigskip
\\Since the sectional path 
$$\begin{xy}
  \xymatrix{ 
     \cdots \ar[r]& \t^{2r}(X_2) \ar[r] &	\t^{2r}(X_1) \ar[r] &	\t^r(X_s) \ar[r] 	& \cdots	     }
\end{xy}$$
$$\begin{xy}
  \xymatrix{ 
   \cdots \ar[r] & \t^r(X_2) \ar[r] &  \t^r(X_1) \ar[r] & X_s \ar[r] &  \cdots	\ar[r]	&	X_2 \ar[r]	&	X_1	    }
\end{xy}$$
is nothing but the composition $ \gamma_0 \gamma_1 \cdots \gamma_{j-1}\gamma_j\gamma_{j+1} \cdots$, it contains infinitely many arrows of finite global left degree in $\chi$.
\bigskip
\\In the last step we show that at least one of the $X_i$ is not stable. Suppose all $X_i$ are stable, then there is an infinite sectional path
$$\begin{xy}
  \xymatrix{ 
   \cdots \ar[r] & \t^r(X_2) \ar[r] &  \t^r(X_1) \ar[r] & X_s \ar[r] &  \cdots			}
\end{xy}$$
$$\begin{xy}
  \xymatrix{
  \cdots \ar[r] &	X_2 \ar[r]	&	X_1 \ar[r] & \t^{-r}(X_s) \ar[r]	& \cdots	    }
\end{xy}$$
$$\begin{xy}
  \xymatrix{ 
     \cdots \ar[r]& \t^{-r}(X_2) \ar[r] &	\t^{-r}(X_1) \ar[r] &	\t^{-2r}(X_s) \ar[r] 	& \cdots	     }
\end{xy}$$
which only contains arrows of infinite global left degree in $\chi$ by Lemma \ref{A_inftydegreeglobal}. This clearly contradicts that 
$$\begin{xy}
  \xymatrix{ 
   \cdots \ar[r] & \t^r(X_2) \ar[r] &  \t^r(X_1) \ar[r] & X_s \ar[r] &  \cdots	\ar[r]	&	X_2 \ar[r]	&	X_1	    }
\end{xy}$$
contains arrows of finite global left degree in $\chi$, which we have proved previously in this lemma. \qed
\begin{landscape}

$$\begin{xy}
\xymatrixrowsep{0.6in}
\xymatrixcolsep{0.6in}

  \xymatrix@!0{  
&	&	&	&	&	&	&	&	&\t^{(j+1)r}(X_1) \ar[dr]	&  	&	&	\\
&	&	&	&	&	&	&	&\t^{jr + 1}(X_s)\ar[ur]	 &	& \t^{jr}(X_s)\ar[dr]\ar@{.>}[ll] 	&	&	\\
&	&	&	&	&	&	& \ar[ur]	&	&  	&	& \ar@{.}[dr]	&	\\
&	&	&	&	&	&\ar@{.}[ur]	&	&	&	&	&	& \ar[dr]	&	& \\
&\t^{(j + 1)r - 1}(X_2)	&	&\cdots \ar@{.>}[ll]	&	& \t^{jr+s-1}(X_2)\ar[ur]\ar@{.>}[ll]& 	&  	&	&	&	& 	&	& \t^{jr}(X_2)\ar[dr]	\\ 
\t^{(j + 1)r}(X_1) \ar[ur]&	&\cdots \ar@{.>}[ll]	&	&\t^{jr+s}(X_1) \ar[ur]\ar@{.>}[ll]	&	& 	&	& 	& 	&	&	& 	&	 & \t^{jr}(X_1)
 }
\end{xy}$$

\end{landscape}

\begin{Th}\label{leftstablecomponentstructure}
Let $\Gamma$ be a left or right stable component of $\Gamma_\chi$ that contains an oriented cycle but does not contain a $\t$-periodic module. Then
\begin{enumerate}[(a)]
 \item $\Gamma$ is not a stable component of $\Gamma_\chi$, 
 \item $\Gamma$ contains only finitely many $\t$-orbits,
 \item every arrow in $\Gamma$ has trivial valuation and
 \item each module in $\Gamma$ has at most two immediate predecessors and at most two immediate successors in $\Gamma$.
\end{enumerate}

\end{Th}

\proof
\\As usual, we prove the theorem only in the case where $\Gamma$ is left stable. Since there is an oriented cycle in $\Gamma$ and no $\t$-periodic module, there is an infinite sectional path 
$$\begin{xy}
  \xymatrix{ 
     \cdots \ar[r]& \t^{2r}(X_2) \ar[r] &	\t^{2r}(X_1) \ar[r] &	\t^r(X_s) \ar[r] 	& \cdots	     }
\end{xy}$$
$$\begin{xy}
  \xymatrix{ 
   \cdots \ar[r] & \t^r(X_2) \ar[r] &  \t^r(X_1) \ar[r] & X_s \ar[r] &  \cdots	\ar[r]	&	X_2 \ar[r]	&	X_1	    }
\end{xy}$$
in $\Gamma$ by Lemma \ref{infinitesectionalpathinleftstablecomponent}, which we denote by $\gamma$. Furthermore, by the same lemma, $\gamma$ contains infinitely many arrows of finite global left degree in $\chi$, all $X_i$ belong to pairwise different $\t$-orbits and at least one of the $X_i$ is not stable. Thus $\Gamma$ is not a stable component. In addition, each arrow in $\gamma$ has trivial valuation and each module in $\gamma$ has at most two immediate predecessors in $\Gamma$ by Lemma \ref{vmtinfinityglobal}.
\bigskip
\\We prove that in fact each $\t$-orbit in $\Gamma$ is generated by some $X_i$ in $\gamma$. Since $\Gamma$ is left stable, it is sufficient to show that each immediate predecessor of a module in a $\t$-orbit generated by some $X_i$ is again in a $\t$-orbit of some $X_j$. First note that each module in $\gamma$ clearly belongs to the $\t$-orbit of some module $X_i$ and has at most two immediate predecessors in $\Gamma$. Hence every module in $\gamma$ other than $X_1$ has precisely two immediate predecessors in $\Gamma$, each of which belongs to the $\t$-orbit of one of the $X_i$. Let $Y$ be an immediate predecessor of $X_1$, then either $\t^r(Y) = \t^r(X_2)$ or $\t^r(Y) = \t(X_s)$ as $\t^r(X_1)$ has two immediate predecessors. In both cases $\t^r(Y)$ and $Y$ belong to the $\t$-orbit of one of the $X_i$.
\bigskip
\\Now assume that $Z$ is an immediate predecessor of $\t^p(X_j)$ in $\Gamma$ for some $p \geq 0$ and $1 \leq j \leq s$. Then $\t^{pr-p}(Z)$ is an immediate predecessor of $\t^{pr}(X_j)$ in $\Gamma$. Since $\t^{pr}(X_j)$ is a module in $\gamma$, the modules $\t^{pr-p}(Z)$ and $Z$ belong either to the $\t$-orbit of $X_{j-1}$ or the orbit of $X_{j+1}$, where $X_0 = X_s$ and $X_1 = X_{s+1}$. Consequently, every orbit in $\Gamma$ is generated by some $X_i$ for $1 \leq i \leq s$.
\bigskip
\\So if $Z$ now denotes an arbitrary module in $\Gamma$, then we have $Z = \t^m(X_i)$ for some $m \in \Z$ and $1 \leq i \leq s$. It follows that there exists some $p > 0$ such that $\t^p(Z)$ is a module in $\gamma$. Therefore, $Z$ has at most two immediate predecessors in $\Gamma$ since $\Gamma$ is left stable. Furthermore, all arrows ending in $Z$ have trivial valuation and so have all arrows in the whole component $\Gamma$, which completes the proof. \qed

\newpage

\begin{Cor}\label{oneinfinitesectionalpathpermodule}
\mbox{}
\begin{enumerate}[(a)]
\item For each module $Z$ in a left stable component $\Gamma$ of $\Gamma_\chi$ that contains an oriented cycle but no $\t$-periodic module there is at most one infinite sectional path 
$$\begin{xy}
  \xymatrix{ 
   \cdots \ar[r] & Z_j \ar[r] &    \cdots	\ar[r]	&	Z_1 \ar[r]	&	Z_0 = Z	    }
\end{xy}$$
ending in $Z$.
\item Dually, if $\Gamma$ is a right stable component of $\Gamma_\chi$ containing an oriented cycle but no $\t$-periodic module, then for every module $Z$ in $\Gamma$ there is at most one infinite sectional path 
$$\begin{xy}
  \xymatrix{ 
   Z = Z_0 \ar[r] & Z_1 \ar[r] &    \cdots	\ar[r]	&	Z_j \ar[r]	&	\cdots	    }
\end{xy}$$
starting in $Z$.
\end{enumerate}
\end{Cor}

\proof
\\We know by Lemma \ref{infinitesectionalpathinleftstablecomponent} that $\Gamma$ contains an infinite sectional path
$$\begin{xy}
  \xymatrix{ 
     \cdots \ar[r]& \t^{2r}(X_2) \ar[r] &	\t^{2r}(X_1) \ar[r] &	\t^r(X_s) \ar[r] 	& \cdots	     }
\end{xy}$$
$$\begin{xy}
  \xymatrix{ 
   \cdots \ar[r] & \t^r(X_2) \ar[r] &  \t^r(X_1) \ar[r] & X_s \ar[r] &  \cdots	\ar[r]	&	X_2 \ar[r]	&	X_1,	    }
\end{xy}$$
and every module $Z$ in $\Gamma$ equals $\t^m(X_j)$ for some $m \in \Z$ and $1 \leq j \leq s$ by Theorem \ref{leftstablecomponentstructure}. Moreover, we can find a non-negative integers $n$ and $p$ such that $\t^n(Z) = \t^{pr}(X_j)$. So if we have an infinite sectional path 
$$\begin{xy}
  \xymatrix{ 
   \cdots \ar[r] & Z_j \ar[r] &    \cdots	\ar[r]	&	Z_1 \ar[r]	&	Z_0 = Z	    }
\end{xy}$$
ending in $Z$, we also obtain an infinite sectional path 
$$\begin{xy}
  \xymatrix{ 
   \cdots \ar[r] & \t^n(Z_j) \ar[r] &    \cdots	\ar[r]	&	\t^n(Z_1) \ar[r]	&	\t^n(Z_0)	    }
\end{xy}$$ 
ending in $\t^n(Z) = \t^{pr}(X_j)$. We show that this path equals
$$\begin{xy}
   \xymatrix{ 
      \cdots \ar[r]& \t^{p(r+1)}(X_2) \ar[r] &	\t^{p(r+1)}(X_1) \ar[r] &	\t^{pr}(X_s) \ar[r] 	& \cdots	     }
 \end{xy}$$
$$\begin{xy}
   \xymatrix{ 
    \cdots \ar[r] & \t^{pr}(X_{j+1}) \ar[r] &  \t^{pr}(X_j).  }
\end{xy}$$
Suppose they are not equal and let $k$ be the lowest positive integer such that $\t^n(Z_k) =  \t^{qr}(X_i)$ but $\t^n(Z_{k+1}) \neq  \t^{qr}(X_{i+1})$ for some $q \geq p$ and $1 \leq i \leq s$, where we define $X_{s+1} = \t^r(X_1)$. Then there are two infinite sectional paths ending in $\t^n(Z_k)$ and by Lemma \ref{twosectionalpathsinfiniteglobaldegree} all arrows in 
$$\begin{xy}
   \xymatrix{ 
      \cdots \ar[r]& \t^{q(r+1)}(X_2) \ar[r] &	\t^{q(r+1)}(X_1) \ar[r] &	\t^{qr}(X_s) \ar[r] 	& \cdots	     }
 \end{xy}$$
$$\begin{xy}
   \xymatrix{ 
    \cdots \ar[r] & \t^{qr}(X_{i+1}) \ar[r] &  \t^{qr}(X_i).  }
\end{xy}$$
have infinite global left degree in $\chi$. This contradicts the previous result from Lemma \ref{infinitesectionalpathinleftstablecomponent} that
\newpage
$$\begin{xy}
  \xymatrix{ 
     \cdots \ar[r]& \t^{2r}(X_2) \ar[r] &	\t^{2r}(X_1) \ar[r] &	\t^r(X_s) \ar[r] 	& \cdots	     }
\end{xy}$$
$$\begin{xy}
  \xymatrix{ 
   \cdots \ar[r] & \t^r(X_2) \ar[r] &  \t^r(X_1) \ar[r] & X_s \ar[r] &  \cdots	\ar[r]	&	X_2 \ar[r]	&	X_1	    }
\end{xy}$$
contains infinitely many arrows of finite global left degree in $\chi$. Hence 
$$\begin{xy}
  \xymatrix{ 
   \cdots \ar[r] & \t^n(Z_j) \ar[r] &    \cdots	\ar[r]	&	\t^n(Z_1) \ar[r]	&	\t^n(Z_0)	    }
\end{xy}$$ 
must equal
$$\begin{xy}
   \xymatrix{ 
      \cdots \ar[r]& \t^{p(r+1)}(X_2) \ar[r] &	\t^{p(r+1)}(X_1) \ar[r] &	\t^{pr}(X_s) \ar[r] 	& \cdots	     }
 \end{xy}$$
$$\begin{xy}
   \xymatrix{ 
    \cdots \ar[r] & \t^{pr}(X_{j+1}) \ar[r] &  \t^{pr}(X_j)  }
\end{xy}$$
and 
$$\begin{xy}
  \xymatrix{ 
   \cdots \ar[r] & Z_j \ar[r] &    \cdots	\ar[r]	&	Z_1 \ar[r]	&	Z_0 = Z	    }
\end{xy}$$
is the only infinite sectional path ending in $Z$. \qed

\begin{Cor}\label{doubleinfinitesectionalpathsdontexist}
Let $\Gamma$ be a left or right stable component that contains an oriented cycle but no $\t$-periodic module. Then $\Gamma$ does not contain a bi-infinite sectional path 
$$\begin{xy}
  \xymatrix{ 
   \cdots \ar[r] & Z_j \ar[r] &    \cdots	\ar[r]	&	Z_1 \ar[r]	&	Z_0 \ar[r]	&	Z_{-1} \ar[r]  & \cdots \ar[r] & Z_{-j} \ar[r] & \cdots	    }
\end{xy}$$
 \end{Cor}

\proof
\\We prove the statement in the case that $\Gamma$ is left stable, the other follows by duality. Let us assume there is an infinite sectional path
$$\begin{xy}
  \xymatrix{ 
   \cdots \ar[r] & Z_j \ar[r] &    \cdots	\ar[r]	&	Z_1 \ar[r]	&	Z_0 \ar[r]	&	Z_{-1} \ar[r]  & \cdots \ar[r] & Z_{-j} \ar[r] & \cdots	    }
\end{xy}$$
then we have two different infinite sectional paths 
$$\begin{xy}
  \xymatrix{ 
   \cdots \ar[r] & Z_j \ar[r] &    \cdots	\ar[r]	&	Z_1 \ar[r]	&	Z_0	    }
\end{xy}$$
and
$$\begin{xy}
  \xymatrix{ 
   \cdots \ar[r] & \t^j(Z_{-j}) \ar[r] &    \cdots	\ar[r]	&	\t^2(Z_{-2}) \ar[r]	&	\t(Z_{-1}) \ar[r]	&	Z_0	    }
\end{xy}$$
ending in $Z_0$, which clearly contradicts Corollary \ref{oneinfinitesectionalpathpermodule}. \qed

\chapter{Shapes of Auslander-Reiten quivers of functorially finite resolving subcategories}

In this chapter we find criteria for Auslander-Reiten quivers of resolving subcategories to be finite or infinite. Therefore, we introduce two different concepts of assigning a graph to a connected component of an Auslander-Reiten quiver, that coincide in some cases. However, some crucial examples show why both ideas are necessary.

\section{Left stable components of Auslander-Reiten quivers}

Let $\Gamma_\Omega$ be the Auslander-Reiten quiver of a functorially finite resolving subcategory $\Omega$. The subquivers $\Gamma_l$, $\Gamma_r$ and $\Gamma_s$ of $\Gamma_\Omega$ containing all left stable, right stable and stable $\tau_\Omega$-orbits are called the stable, left stable and right stable Auslander-Reiten quivers of $\Omega$ respectively. Connected components of $\Gamma_l$, $\Gamma_r$ and $\Gamma_s$ are defined analogously to connected components of the whole Auslander-Reiten quiver. Note that if a connected component $\Gamma$ of either $\Gamma_l$ or $\Gamma_r$ contains a $\tau_\Omega$-periodic module, then by Lemma \ref{orbit} every module in $\Gamma$ is $\tau_\Omega$-periodic and $\Gamma$ is a connected component of $\Gamma_s$.

\begin{Def}
Let $\Gamma$ be a subquiver of $\Gamma_\Omega$ and let $\Sigma$ be a connected subquiver of $\Gamma$. We call $\Sigma$ a \textbf{sectional subgraph} if all paths of length two contained in $\Sigma$ are sectional. $\Sigma$ is called a \textbf{full sectional subgraph} of $\Gamma$ if any connected subquiver $\Sigma'$ of $\Gamma$ such that $\Sigma \subsetneq \Sigma'$ is not a sectional subgraph. The undirected graph $\overline{\Sigma}$ associated to $\Sigma$ is called the \textbf{type of $\pmb{\Sigma}$} and for a vertex $X$ in $\Sigma$ the corresponding vertex in $\overline{\Sigma}$ is denoted by $\overline{X}$.
\end{Def}

Naturally, we call a sectional subgraph and its type finite if it consists of only finitely many vertices. Note that the same module can occur several times in an Auslander-Reiten quiver and also in a sectional subgraph. In particular, there can be two arrows in a sectional subgraph whose composition is a non-sectional path if they are not adjacent in the subgraph, as the following example shows.
\newpage
\begin{Ex}
Let A be the path algebra of
$$\begin{xy}
\xymatrix{
  e_1 \  \ar[r]^\alpha  & \ e_2 \ \ar@(ur,dr)^\beta  \\
}
\end{xy}$$
with the relation $\beta^2 = 0$. Then the Auslander-Reiten quiver of $A$-mod contains a sectional path which contains two arrows whose composition is not sectional.
\end{Ex}

Recall from Example \ref{exampletwistedquiver} that the Auslander-Reiten quiver of $A$-mod looks as follows.
$$\begin{xy}
\xymatrixrowsep{0.5in}
\xymatrixcolsep{0.5in}

  \xymatrix@!0{  
&	& [P_1 \ar[dr] 	&	& I_1] \ar@{.>}[ll]	&	&	&	\\
&[P_2 \ar[dr] \ar[ur]	&	&I_4] \ar[dr] \ar[ur] \ar@{.>}[ll] 	&	&	&	&	\\
S_2\ar[dr] \ar[ur]	& 	&Y \ar[dr] \ar[ur] \ar@{.>}[ll]&	&X \ar[dr] \ar@{.>}[ll]&	&S_2 \ar[dr] \ar@{.>}[ll] &	\\ 
&X \ar[ur] &	&S_2 \ar[dr] \ar[ur] \ar@{.>}[ll] & 	&Y \ar[dr] \ar[ur] \ar@{.>}[ll]&	&X \ar@{.>}[ll]	\\ 
&	&	&	&[P_2 \ar[dr] \ar[ur]	&	&I_4] \ar[dr] \ar[ur] \ar@{.>}[ll] 	&	\\
&	&	&	&	& [P_1 \ar[ur] 	&	& I_1] \ar@{.>}[ll]	
 }
\end{xy}$$

We denote the sectional path

$$\begin{xy}
  \xymatrix{ 
P_2 \ar[r]^f & Y \ar[r]	& S_2 \ar[r]^g	& P_2
 }
\end{xy}$$
by $\gamma$. Although the composition $fg$ is not sectional, $\gamma$ is a sectional subgraph since the arrows of $f$ and $g$ are not adjacent and do not form a subpath of $\gamma$. Note that $\gamma$ is not a full sectional subgraph as the arrow from $P_2$ to $P_1$ can be added to $\gamma$ on both sides. In the following example we see that the differentiation between vertices that are given by the same module is very important in order to assign a graph to a connected component of $\Gamma_l$. 

\begin{Ex}
Let $A$ be the hereditary path algebra given by 
$$\begin{xy}
\xymatrix{
  e_1  \ar[r] \ar@/^1pc/[rr] &  e_2 \ar[r] &  e_3  .
}
\end{xy}$$
In its preinjective component $\Gamma$ there are sectional subgraphs of type $A_\infty^\infty$ that contain non-consecutive arrows that form a non-sectional path.
\end{Ex}

The injective modules are given by the following Jordan-H\"{o}lder compositions.
 $$\begin{xy}
\xymatrixrowsep{0in}
\xymatrixcolsep{-0.1in}
  \xymatrix{  
	    &    	&		&		&				&	&	&S_1					\\	
I_1=\ \ &S_1   	&\ \ \ \ \ \ \ \ \ \ I_2= \ \	&S_1	&\ \ \ \ \ \ \ \ \ \ I_3= \ \	&S_1	&	&S_2		\\	
	    &    	&		&	S_2	&				&	&S_3	&
   }
\end{xy}$$
This information is already enough to determine the structure of the preinjective component $\Gamma$.

$$\begin{xy}
\xymatrixrowsep{0.5in}
\xymatrixcolsep{0.5in}

  \xymatrix@!0{ 
&	& \ar[dr]	&  	&\ar[dr]	& 	& \ar[dr]	&	&	\\
&	\cdots &	& \tau_\Omega^2(I_1) \ar[dr] \ar[ur] \ar@{.>}[ll] 	&	& \tau_\Omega(I_1) \ar[dr] \ar[ur] \ar@{.>}[ll]	&	&I_1] \ar@{.>}[ll]	&	\\
\cdots&		&\tau_\Omega^2(I_3) \ar[dr] \ar[ur] \ar@{.>}[ll]	&	&\tau_\Omega(I_3) \ar[dr] \ar[ur] \ar@{.>}[ll] 	&	&I_3] \ar[dr] \ar[ur] \ar@{.>}[ll]	&	&	\\
&\cdots	& 	&\tau_\Omega^2(I_2) \ar[dr] \ar[ur] \ar@{.>}[ll]&	&\tau_\Omega(I_2) \ar[dr]\ar[ur] \ar@{.>}[ll]&	&I_2] \ar[dr] \ar@{.>}[ll] &	\\ 
&	&\cdots &	&\tau_\Omega^2(I_1) \ar[dr] \ar[ur] \ar@{.>}[ll] & 	&\tau_\Omega(I_1) \ar[dr] \ar[ur] \ar@{.>}[ll]&	&I_1] \ar@{.>}[ll]	\\ 
&\cdots&	&\tau_\Omega^2(I_3) \ar[dr] \ar[ur] \ar@{.>}[ll]	&	&\tau_\Omega(I_3) \ar[dr] \ar[ur] \ar@{.>}[ll]	&	&I_3] \ar[dr] \ar[ur] \ar@{.>}[ll] 	&	\\
&	&\cdots	&	&\tau_\Omega^2(I_2) \ar[dr] \ar[ur] \ar@{.>}[ll]	&	& \tau_\Omega(I_2) \ar[dr] \ar[ur] \ar@{.>}[ll] 	&	& I_2] \ar[dr] \ar@{.>}[ll]	\\
&	&	& \ar[ur]&						&\ar[ur]	&							&\ar[ur]&	& 
 }
\end{xy}$$
We denote the subgraph 
 $$\begin{xy}
  \xymatrix{ 
\cdots \ar[r]& \tau_\Omega^n(I_3) \ar[r]	&\tau_\Omega^n(I_2) \ar[r]     &\tau_\Omega^n(I_1)\ar[r]& \tau_\Omega^{n-1}(I_3) \ar[r]	&\cdots 	    }
\end{xy}$$
 $$\begin{xy}
  \xymatrix{ 
\cdots \ar[r] &   I_3	 \ar[r]	&	I_2 \ar[r]&	I_1 	  & I_3 \ar[l]	& \tau_\Omega(I_2) \ar[l] & \tau_\Omega^2(I_1) \ar[l] &	\cdots \ar[l]  }
\end{xy}$$
 $$\begin{xy}
  \xymatrix{ 
\cdots  &   \tau_\Omega^n(I_3)	 \ar[l]	&\tau_\Omega^{n+1}(I_2) \ar[l]	&\tau_\Omega^{n+2}(I_1) \ar[l] & \tau_\Omega^{n+2}(I_3) \ar[l] & \cdots \ar[l] }
\end{xy}$$
by $\Sigma$ and note that its type is $A_\infty^\infty$. It clearly contains arrows that, if we compose them, form a non-sectional path, but since these arrows are not adjacent in $\Sigma$ the subgraph is a full sectional subgraph. However, $\Gamma$ also contains a full sectional subgraph
 $$\begin{xy}
  \xymatrix{ 
     \cdots \ar[r]&	I_1 & I_3 \ar[r] \ar[l] & I_2 \ar[r]& 	I_1     &   I_3	 \ar[r] \ar[l]	&	I_2 \ar[r]	&	\cdots 	    }
\end{xy}$$
of type $A_\infty^\infty$ which does not contain arrows whose composition is non-sectional. 
\bigskip
\\Suppose now we would identify vertices of an Auslander-Reiten quiver with each other if they are given by the same module. Then we would have full sectional subgraphs 
 $$\begin{xy}
  \xymatrix{ 
     I_1 & I_3 \ar[r] \ar[l] & I_2 \ar[r]& 	I_1  }
\end{xy}$$
and 
 $$\begin{xy}
  \xymatrix{ 
\cdots \ar[r] &	\tau_\Omega^n(I_2) \ar[r]     &\tau_\Omega^n(I_1)\ar[r]& \tau_\Omega^{n-1}(I_3) \ar[r]& \cdots \ar[r] &   I_3	 \ar[r]	&	I_2 \ar[r]&	I_1,}
\end{xy}$$
which would be of type $\widetilde{A}_2$ and $A_\infty$ respectively. In particular, $\Gamma$ would contain full sectional subgraphs of different types, which we avoid in our setup as the subsequent statements show.

\newpage

\begin{Lemma}\label{adjacentorbits}
Let $X$ be a vertex in a full sectional subgraph $\Sigma$ of a connected component $\Gamma$ in the Auslander-Reiten quiver of $\Gamma_\Omega$. Then the number of vertices adjacent to $\overline{X}$ in $\overline{\Sigma}$ equals the number of immediate predecessors of $X$ in $\Gamma$ plus the number of projective successors of $X$ in $\Gamma$. Furthermore, let $Y$ be an immediate predecessor or immediate successor of $X$ in $\Sigma$, the number of undirected edges between $\overline{X}$ and $\overline{Y}$ equals the number of arrows from $Y$ to $X$ or $X$ to $Y$ in $\Gamma$ respectively.
\end{Lemma}
\proof
\\We prove the second statement first. Let $X$ be a vertex in $\Sigma$ and, without loss of generality, let $Y$ be an immediate predecessor of $X$ in $\Sigma$. Since $\Sigma$ is full, it must necessarily contain all arrows from $Y$ to $X$, and so the number of undirected edges between $\overline{X}$ and $\overline{Y}$ coincides with the number of arrows from $Y$ to $X$.
\bigskip
\\It is easy to see that every projective successor of $X$ in $\Gamma$ must be contained in $\Sigma$ as $\Sigma$ is a full sectional subgraph. Assume there is a predecessor $Y$ of $X$ in $\Gamma$ that is not contained in $\Sigma$. Then, as $\Sigma$ is full, $\tau_\Omega^{-1}(Y)$ exists and all arrows from $X$ to $\tau_\Omega^{-1}(Y)$ are contained in $\Sigma$. On the other hand, if both $Y$ and $\tau_\Omega^{-1}(Y)$ are contained in $\Sigma$, then $\Sigma$ also contains a path of length two from $Y$ to $\tau_\Omega^{-1}(Y)$, which is a contradiction. \qed

\begin{Th}\label{subgraphtype}
\mbox{}
\begin{enumerate}[(a)]
\item Let $\Gamma$ be a connected component of the left stable Auslander-Reiten quiver $\Gamma_l$ and let $\Sigma$ and $\Sigma'$ be two full sectional subgraphs of $\Gamma$ such that for each Ext-injective module $I$ in $\Gamma$ there is no path from $I$ to any vertex in $\Sigma$ or $\Sigma'$. Then $\Sigma$ and $\Sigma'$ are of the same type. 
\item If for every full sectional subgraph $\Sigma$ of $\Gamma$ there is an Ext-injective module $I$ and a path from $I$ to some vertex of $\Sigma$, then $\Gamma$ contains an oriented cycle but no $\tau_\Omega$-periodic module and $\Sigma$ is of type $A_\infty$.
\item Dually, let $\Gamma$ denote a connected component of $\Gamma_r$, then two full sectional subgraphs $\Sigma$ and $\Sigma'$ have the same type if there is no path from $\Sigma$ or $\Sigma'$ to a projective module in $\Gamma$.
\item If for every full sectional subgraph $\Sigma$ of $\Gamma$ there is a projective module $P$ and a path from some vertex of $\Sigma$ to $P$, then $\Gamma$ contains an oriented cycle but no $\tau_\Omega$-periodic module and $\Sigma$ is of type $A_\infty$.
\end{enumerate}
\end{Th}

\proof
\\Suppose $\Sigma$ and $\Sigma'$ are full sectional subgraphs in $\Gamma$ such that there is no path from an Ext-injective module $I$ to $\Sigma$ or $\Sigma'$. Due to this property and since $\Gamma$ is a connected component of $\Gamma_l$, we can choose modules $X$ in $\Sigma$, $X'$ in $\Sigma'$ and an integer $m \in \Z$ such that there is a sectional path
 $$\begin{xy}
  \xymatrix{ 
     X = X_0 \ar[r] 	&	X_1 \ar[r]     &   \cdots	 \ar[r]	&	X_{n-1} \ar[r]	&	X_n = \tau_\Omega^m(X')	    }
\end{xy}$$
of minimal length $n$. If $n$ does not equal zero, then $X_1$ is not a vertex of $\Sigma$ by minimality of $n$. Consider the subgraph $\Sigma_1$ obtained by adding the module $X_1$ and the arrow from $X$ to $X_1$ to the full sectional subgraph $\Sigma$. By definition of full sectional subgraphs there must be a path in $\Sigma_1$ that is not sectional and ends in $X_1$. Hence the path also contains $\tau_\Omega(X_1)$, which then must be a vertex of $\Sigma$ and there is a sectional path from $\tau_\Omega(X_1)$ to $\tau_\Omega^{m+1}(X')$, which contradicts $n$ to be minimal. Hence we can find modules $X$ and $X'$ such that $X = \tau_\Omega^m(X')$.
\bigskip
\\Without loss of generality, we can assume that $m \geq 0$. Since $\Gamma$ is left stable, the numbers of vertices and undirected edges adjacent to $\overline{X'}$ and $\overline{X}$ equal the numbers of immediate predecessors of $X'$ and $X$ and arrows from immediate predecessors to $X'$ and $X$ respectively by Lemma \ref{adjacentorbits}. Therefore, as $m \geq 0$ and $\Gamma$ is left stable, the number of vertices adjacent to $\overline{X}$ is greater than or equal to the number of vertices adjacent to $\overline{X'}$. Suppose this number is strictly higher, then there is an immediate predecessor $Y$ of $X$ such that $\tau_\Omega^{-m}(Y)$ is not an immediate predecessor of $X'$. Consequently, $\tau_\Omega^{-m}(Y)$ does not exist, the orbit of $Y$ must be Ext-injective and there is a path from some Ext-injective module $I$ to $\Sigma'$, which is a contradiction. So let $Y'$ be an immediate predecessor or immediate successor of $X'$ in $\Sigma'$ and without loss of generality, let $Y = \tau_\Omega^m(Y')$ be an adjacent vertex to $X$ in $\Sigma$. Then the numbers of undirected edges between $\overline{X'}$ and $\overline{Y'}$ and between $\overline{X}$ and $\overline{Y}$ coincide by Corollary \ref{valuationisinvariantoftranslation}, which completes the proof of statement $(a)$.
\bigskip
\\Assume now that $\Gamma$ is a connected component of $\Gamma_l$ with the property that for every full sectional subgraph $\Sigma$ there is an Ext-injective module $I$ such that there is a path from $I$ to some vertex in $\Sigma$. Since every module $X$ is contained in some full sectional subgraph, it follows that for every module $X$ in $\Gamma$ there is an Ext-injective module $I$ such that there is a path from $I$ to $X$ as $\Gamma$ is left stable. In particular, $\Gamma$ contains an Ext-injective module $I_0$ and there is an Ext-injective module $I_1$ with a non-trivial path to $I_0$. Furthermore, there is also an Ext-injective module $I_2$ with a non-trivial path to $I_1$ which can be extended to a path from $I_2$ to $I_0$.  We continue inductively and obtain an infinite collection of paths $I_j \to I_i$ for all $j > i \geq 0$. Since there are up to isomorphism only finitely many indecomposable Ext-injective modules in $\Omega$, there is an oriented cycle in $\Gamma$.
\bigskip
\\Suppose $\Gamma$ contains a $\tau_\Omega$-periodic module $X$. Then all adjacent orbits are either $\tau_\Omega$-periodic or finite by Lemma \ref{orbit}. Since $\Gamma$ is left stable, it cannot contain finite orbits and hence every orbit must be $\tau_\Omega$-periodic. This is a contradiction as $\Gamma$ contains at least one Ext-injective module.
\bigskip
\\Now by Lemma \ref{infinitesectionalpathinleftstablecomponent} $\Gamma$ contains an infinite sectional path of type $A_\infty$ and by Theorem \ref{leftstablecomponentstructure} every module in $\Gamma$ has at most two immediate predecessors. Hence the type of a full sectional subgraph is either $A_\infty$ or $A_\infty^\infty$. Suppose $\Sigma$ is a full sectional subgraph of type $A_\infty^\infty$ and let $X$ be any module in $\Sigma$. Since $\Gamma$ is left stable, we can then construct two different sectional paths ending in $X$ in the way the following example suggests. 

$$\begin{xy}
\xymatrixrowsep{0.5in}
\xymatrixcolsep{0.5in}

  \xymatrix@!0{ 
&	&	&\ar@{.}[dl]	&	&	&  \ar@{.}[dr]	&	&	&	\\
&	&	\ar[dl]	&  	&	&	&	&\ar[dr]&	\\
& Y_2 \ar[dl]  	&	& 	&	&	&	&	&\tau_\Omega(Y_2) \ar[dr]	\\
Y_1 \ar[dr]&	&	& 	&	&	&	&	&	&Y_1 \ar[dr]	\\
	&X	&	&	&	&	&	&	&	&	&X\\ 
X_1\ar[ur] \ar[dr] &	&	& 	&	&	&	&	&	&X_1 \ar[ur]	\\ 
&X_2 \ar[dr]	&	&	&	&	&	&	&\tau_\Omega(X_2) \ar[ur] 	&	\\
	&	&X_3	&	& 	&	&	&\tau_\Omega^2(X_3)\ar[ur]	& 	\\
&X_4 \ar[ur] \ar[dr]&	&	&	&	&\tau_\Omega^2(X_4)\ar[ur]	&	&	&	& \\
	&	&X_5	&	&	&\tau_\Omega^3(X_5)\ar[ur]	& 	&	& 	\\
	&\ar[ur]&	&	&\ar[ur]	&	& 	&	& 	\\
\ar@{.}[ur]&	&	&\ar@{.}[ur]	&	&	& 	&	& 	}
\end{xy}$$
This is impossible as it clearly contradicts the results from Corollary \ref{oneinfinitesectionalpathpermodule}. It follows that every full sectional subgraph of $\Gamma$ is of type $A_\infty$ \qed
\bigskip
\\The following example shows that there are left stable components such that for every module $X$ there is an Ext-injective module $I$ such that there is a path from $I$ to $X$. Moreover, it can be seen why the proof for Theorem \ref{subgraphtype}(a) cannot be applied to these components.

\begin{Ex}\label{helical}
Let $A$ be the path algebra of the quiver
$$\begin{xy}
\xymatrix{
  e_1 \ \ \ar@<1ex>[r]^{\alpha} \ar@<-1ex>[r]_{\beta}  & \ \ e_2 \ \  \ar[r]^{\gamma}  & \ \  e_3 \\
}
\end{xy},$$
with the relation $\gamma\beta = 0$. Then the Auslander-Reiten quiver of $A$ contains an infinite left stable component $\Gamma$ such that for all modules $X$ in $\Gamma$ there is a path from $I_3$ to $X$.
\end{Ex}

First of all we write down the Jordan-H\"{o}lder composition series of the projective and injective modules.

 $$\begin{xy}
\xymatrixrowsep{0in}
\xymatrixcolsep{-0.1in}
  \xymatrix{ 
	&	&S_1   	&		&		&				&	&	&					\\	
P_1=\ \ &S_2   	&	& S_2\ \ \ \ \ \ \ \ \ \ P_2= \ \	&S_2	&\ \ \ \ \ \ \ \ \ \ P_3= \ \	&S_3	&	&		\\	
	&S_3    &	&		&	S_3	&				&	&	&
   }
\end{xy}$$
 $$\begin{xy}
\xymatrixrowsep{0in}
\xymatrixcolsep{-0.1in}
  \xymatrix{ 
	&	&		&		&		&	&		&	S_1					\\	
I_1=\ \ &S_1   	& \ \ \ \ \ \ \ \ \ \ I_2= \ \	&S_1	&	&S_1	&\ \ \ \ \ \ \ \ \ \ I_3= \ \	&S_2		\\	
	&	&		&	&	S_2	&		&			&	S_3
   }
\end{xy}$$
We easily verify $\tau(I_3) \cong I_3/S_3$. Then the component of $I_3$ looks the following way:
$$\begin{xy}
\xymatrixrowsep{0.5in}
\xymatrixcolsep{0.5in}

  \xymatrix@!0{ 
& & & & & & & & & & &  \\
& &\cdots& &\tau^3(I_3)\ar@{.>}[ll] \ar[ur] \ar[dr]& &\tau^2(I_3) \ar@{.>}[ll] \ar[dr] \ar[ur]&	&\tau(I_3) \ar@{.>}[ll] \ar[ur]	&	& I_3] \ar[ur] \ar@{.>}[ll]&\\
&\cdots& 	&\tau^2(I_3)\ar@{.>}[ll] \ar[ur] \ar[dr]	& 	&\tau(I_3) \ar@{.>}[ll] \ar[ur]	&	&I_3] \ar[ur] \ar@{.>}[ll]&	&	&		&\\
\cdots&	&\tau(I_3) \ar@{.>}[ll]\ar[ur]	&	&I_3] \ar[ur] \ar@{.>}[ll]&	&	&	&	&	&	&
 }
\end{xy}$$

Clearly, the component must be infinite because all indecomposable projective modules are adjacent in the Auslander-Reiten quiver. Furthermore, we can see that any full sectional subgraph is of type $A_\infty$ in this component. Its shape suggests the following name.

\begin{Def}
\mbox{}
\begin{enumerate}[(a)]
\item We call a connected component $\Gamma$ of the left stable Auslander-Reiten quiver $\Gamma_l$ \textbf{helical} if for every module $X$ there is an Ext-injective module $I$ in $\Gamma$ such that there is a path from $I$ to $X$.
\item Dually, we call a connected component $\Gamma$ of $\Gamma_r$ \textbf{cohelical} if for every module $X$ there is a projective module $P$ in $\Gamma$ and a path from $X$ to $P$.
\end{enumerate}
\end{Def}

Due to the last theorem the following definition makes sense.

\begin{Def}
\mbox{}
\begin{enumerate}[(a)]
\item Let $\Gamma$ be a connected component of $\Gamma_l$. If $\Gamma$ is helical, we define the \textbf{left subgraph type of $\pmb{\Gamma}$} to be $A_\infty$ as any full sectional subgraph is of that type. On the other hand, if $\Gamma$ is not helical, then the \textbf{left subgraph type of $\pmb{\Gamma}$} is defined as the type of a full sectional subgraph $\Sigma$ such that there is no path from an Ext-injective module to any vertex in $\Sigma$. 
\item Dually, let $\Gamma$ denote a connected component of $\Gamma_r$. If $\Gamma$ is cohelical, we define the \textbf{right subgraph type of $\pmb{\Gamma}$} to be $A_\infty$. If $\Gamma$ is not cohelical, then the \textbf{right subgraph type of $\pmb{\Gamma}$} is defined as the type of a full sectional subgraph $\Sigma$ such that there is no path from any vertex in $\Sigma$ to a projective module. 
\item If $\Gamma$ denotes a connected component of the stable Auslander-Reiten quiver $\Gamma_s$, then the \textbf{subgraph type of $\pmb{\Gamma}$} is defined as the type of any full sectional subgraph $\Sigma$ in $\Gamma$.
\end{enumerate}
\end{Def}

It is not hard to see that in a connected component of $\Gamma_s$ subgraph type, left subgraph type and right subgraph type coincide. Let $\Gamma$ denote a connected component of the left stable Auslander-Reiten quiver $\Gamma_l$ of $\Omega$. With some basic matrix calculations we show that $\Gamma$ is finite if its left subgraph type is given by a Dynkin diagram. If $\Gamma$ is finite, it clearly only contains $\tau_\Omega$-periodic orbits. An infinite component must have infinitely many orbits or some non-periodic, infinite, left stable orbits. 
\bigskip
\\We assume $\Gamma$ is infinite, so if $\Gamma$ contains a $\tau_\Omega$-periodic orbit, we know by Theorem \ref{HPRT} that $\Gamma$ is a stable tube. In particular, a full sectional subgraph is of type $A_\infty$ and the left subgraph type is not Dynkin. Secondly, if $\Gamma$ is helical, then every full sectional subgraph is of type $A_\infty$ and the left subgraph type is not Dynkin.

\begin{Def}
\mbox{}
\begin{enumerate}[(a)]
 \item Let $\Gamma$ be a non-helical connected component of $\Gamma_l$ that is not $\tau_\Omega$-periodic and let $\Sigma$ be a full sectional subgraph of $\Gamma$. We say that $\Sigma$ is \textbf{left of all projective and Ext-injective modules} if there is neither a projective module $P$ nor an Ext-injective module $I$ such that there is a path in the whole Auslander-Reiten quiver $\Gamma_\Omega$ from $P$ or $I$ to a module in $\Sigma$ respectively.
 \item Dually, let $\Gamma$ be a non-cohelical connected component of $\Gamma_r$ that is not $\tau_\Omega$-periodic and let $\Sigma$ be a full sectional subgraph of $\Gamma$. We say that $\Sigma$ is \textbf{right of all projective and Ext-injective modules} if there is neither a projective module $P$ nor an Ext-injective module $I$ such that there is a path in the whole Auslander-Reiten quiver $\Gamma_\Omega$ from a module in $\Sigma$ to $P$ or $I$ respectively.
\end{enumerate} 
\end{Def}

Note that a full sectional subgraph left of all projective and Ext-injective modules always exists since $\Gamma$ is not $\tau_\Omega$-periodic but there are only finitely many projective and Ext-injective modules in $\Omega$. Note that if $\Gamma$ is of finite left subgraph type, then it cannot be helical, so the former is a stronger condition on $\Gamma$. Recall that the \textbf{dimension vector} of a module $X$ is the column vector such that its $i$-th entry is given by the number of simple modules $S_i$ in the Jordan-H\"{o}lder decomposition of $X$.

\begin{Def}
\mbox{}
\begin{enumerate}[(a)]
\item Let $\Gamma$ be a connected component of $\Gamma_l$ of finite left subgraph type that is not $\tau_\Omega$-periodic. If $X_1, \ldots, X_n$ denote all modules of a full sectional subgraph left of all projective and Ext-injective modules and $x_1, \ldots, x_n$ their dimension vectors, then $\tau_\Omega$ generates a matrix with its $ij$-th entry being the number of $x_i$ appearing in the dimension vector of $\tau_\Omega(X_j)$ given by $\E_\Omega(X_j)$. We call this matrix a \textbf{translation matrix} of $\Gamma$ and denote it by $M$. 
\item Dually, let $\Gamma$ be a connected component of $\Gamma_r$ of finite right subgraph type that is not $\tau_\Omega$-periodic. If $X_1, \ldots, X_n$ denote all modules of a full sectional subgraph right of all projective and Ext-injective modules and $x_1, \ldots, x_n$ their dimension vectors, then $\tau_\Omega^{-1}$ generates a matrix $M^- = (m_{ij}^-)$ such that $m_{ij}^-$ is the number of $x_i$ appearing in the dimension vector of $\tau_\Omega^{-1}(X_j)$ given by $\E'_\Omega(X_j)$. We call $M^-$ a \textbf{cotranslation matrix} of $\Gamma$.
\end{enumerate}

\end{Def}

\begin{Ex}
Let $A$ be the hereditary path algebra of the quiver 
$$\begin{xy}
  \xymatrix{ 
    &	&e_6 \ar[d] \\
e_1 \ar@<1ex>[dr]\ar@<-1ex>[dr]  &	&e_3 \ar[dr]\ar[d]\ar[dl]    \\
  &	e_2	& e_4 	& e_5.
 	}
\end{xy}
$$ 
We calculate a translation matrix $M = (m_{ij})$ of the preinjective component of $A$.
\end{Ex}

Note that in this algebra we have $I_1 = S_1$ and $I_6 = S_6$. The other injective modules have the following Jordan-H\"{o}lder composition series.
 $$\begin{xy}
\xymatrixrowsep{0in}
\xymatrixcolsep{-0.1in}
  \xymatrix{ 
&	&	&S_6    &				&S_6				&	&S_6	&					&S_6	\\	
I_2=\ \ &S_1^2& &S_3  	&\ \ \ \ \ \ \ \ \ \ I_3= \ \	&S_3	&\ \ \ \ \ \ \ \ \ \ I_4= \ \	&S_3	&\ \ \ \ \ \ \ \ \ \ I_4= \ \		&S_3	\\	
	   & &S_2    	&		&		&				&	&S_4	&					&S_5
   }
\end{xy}$$
For convenience, let us denote $\tau(I_j)$ by $X_j$ for all $1 \leq j \leq 6$. The preinjective component $\Gamma$ looks as follows.
$$\begin{xy}
\xymatrixrowsep{0.5in}
\xymatrixcolsep{0.5in}

  \xymatrix@!0{ 
&\cdots 	& & \tau(X_1) \ar@<1ex>[dr]\ar@<-1ex>[dr] \ar@{.>}[ll]	&	& X_1 \ar@<1ex>[dr]\ar@<-1ex>[dr] \ar@{.>}[ll]&	&I_1] \ar@{.>}[ll]&\\
\cdots&	&\tau(X_2) \ar[dr] \ar@<1ex>[ur]\ar@<-1ex>[ur] \ar@{.>}[ll] &	&X_2 \ar[dr] \ar@<1ex>[ur]\ar@<-1ex>[ur] \ar@{.>}[ll] &	&I_2] \ar@<1ex>[ur]\ar@<-1ex>[ur] \ar[dr] \ar@{.>}[ll]	&	&	\\
&\cdots	& &\tau(X_3) \ar[dr] \ar [ddr] \ar[dddr]\ar[ur] \ar@{.>}[ll]& 	&X_3 \ar[dr] \ar [ddr] \ar[dddr]\ar[ur] \ar@{.>}[ll]&		&I_3] \ar[dddr] \ar@{.>}[ll] &	\\ 
\cdots & &\tau(X_4) \ar[ur] \ar@{.>}[ll]&	&X_4 \ar[ur] \ar@{.>}[ll]&	&I_4] \ar[ur] \ar@{.>}[ll]	\\ 
\cdots&	&\tau(X_5) \ar[uur] \ar@{.>}[ll]	&	&X_5 \ar[uur] \ar@{.>}[ll]	&	&I_5]  \ar[uur] \ar@{.>}[ll] 	&	\\
&	&\cdots	&	& \tau(X_6)  \ar[uuur] \ar@{.>}[ll] 	&	& X_6  \ar[uuur] \ar@{.>}[ll] 	&	& I_6]  \ar@{.>}[ll]	
 }
\end{xy}$$
Clearly, the subquiver consisting of all modules $X_j$ and all arrows between them is a full sectional subgraph of $\Gamma$ that is left of all projective and Ext-injective modules. We denote the dimension vectors of $X_j$ by $x_j$ and the dimension vector of $\tau(X_j)$ by $\tau(x_j)$. It follows from the almost split sequence 
$$\begin{xy}
  \xymatrix{ 
        0 \ar[r]   &	\tau(X_6) \ar[r]		&	X_3  \ar[r]		&	X_6 \ar[r]		&	0	 }
\end{xy}$$
that $\tau(x_6) = x_3 - x_6$. By definition of the translation matrix we get $m_{36} = 1$, $m_{66} = -1$ while $m_{i6} = 0$ for $i = 1,2,4,5$. We similarly obtain $\tau(x_1) = 2x_2 - x_1$ and the entries of the first column of $M$. We conclude from
$$\begin{xy}
  \xymatrix{ 
        0 \ar[r]   &	\tau(X_3) \ar[r]	&	X_2 \oplus X_4 \oplus X_5 \oplus \tau(X_6)  \ar[r]		&	X_3 \ar[r]	&0	 }
\end{xy}$$
that $\tau(x_3) = x_2 + x_4 + x_5 + \tau(x_6) - x_3 =  x_2 + x_4 + x_5 - x_6$ by substituting in $x_3 - x_6$ for $\tau(x_6)$. The other dimension vectors are 
$\tau(x_2) = 2\tau(x_1) + \tau(x_3) - x_2 = 4x_2 - 2x_1 + x_4 + x_5 - x_6$, $\tau(x_4) = \tau(x_3) - x_4 = x_2 + x_5 - x_6$ and $\tau(x_5) = \tau(x_3) - x_5 = x_2 + x_4 - x_6$, which gives rise to 
\begin{center}
\begin{math}
M = $\bordermatrix{ 	&	&	&	&	&	&	\cr
			&  -1	&  -2	&   0	&   0	&   0	&   0	\cr
			&   2	&   4	&   1	&   1	&   1	&   0	\cr
			&   0	&   0	&   0	&   0	&   0	&   1	\cr
			&   0	&   1	&   1	&   0	&   1	&   0	\cr
			&   0	&   1	&   1 	&   1	&   0	&   0	\cr
			&   0	&  -1	&  -1	&  -1	&  -1	&  -1	\cr} $
\end{math}.
\end{center}

\begin{Lemma}
\mbox{}
\begin{enumerate}[(a)]
\item Let $\Gamma$ be a connected component of $\Gamma_l$ that is not $\tau_\Omega$-periodic and is of finite left subgraph type. Let $X_1, \ldots, X_n$ denote the modules of a full sectional subgraph left of all projective and Ext-injective modules and $M = (m_{ij})$ its translation matrix. Then the $ij$-th entry of $M^k$ equals the number of $x_i$ in the dimension vector of $\tau_\Omega^k(X_j)$. 
\item Dually, let $\Gamma$ be a connected component of $\Gamma_r$ that is not $\tau_\Omega$-periodic and is of finite right subgraph type. Let $X_1, \ldots, X_n$ denote the modules of a full sectional subgraph right of all projective and Ext-injective modules and $M^- = (m_{ij}^-)$ its cotranslation matrix. Then the $ij$-th entry of $(M^-)^k$ equals the number of $x_i$ in the dimension vector of $\tau_\Omega^{-k}(X_j)$. 
\end{enumerate}

\end{Lemma}

\proof
\\For convenience, we speak about the number of $x_i$ in $\tau_\Omega(X_j)$ instead of the number of $x_i$ in the dimension vector of $\tau_\Omega(X_j)$. Furthermore, we denote the dimension vector of $\tau_\Omega^l(X_j)$ by $\tau_\Omega^l(x_j)$. We prove the statement by induction on $k$. The statement follows immediately from the definitions for $k = 0$ and $k = 1$.
\bigskip
\\Suppose now the statement is true for $k$. We can write the $ij$-th entry of $M^{k+1}$ as $(y_{i1}, \ldots, y_{in})(m_{1j}, \ldots, m_{nj})^T = \sum_{t=1}^n  y_{it}m_{tj}$, where $(y_{i1}, \ldots, y_{in})$ denotes the $i$-th row vector of $M^k$. By the induction hypothesis $y_{ij}$ is the number of $x_i$ in $\tau_\Omega^k(X_j)$, and by definition $m_{ij}$ is the number of $x_i$ in $\tau_\Omega(X_j)$, which coincides with the number of $\tau_\Omega^k(x_i)$ in $\tau_\Omega^{k+1}(X_j)$ as $\Gamma$ is left stable. It follows that
$$\sum_{t=1}^n  y_{it}m_{tj} = \sum_{t=1}^n (\# \text{ of }x_i \text{ in }\tau_\Omega^k(X_t))(\# \text{ of }\tau_\Omega^k(x_t) \text{ in }\tau_\Omega^{k+1}(X_j)),$$
so the $ij$-entry of $M^{k+1}$ is the number of $x_i$ in $\tau_\Omega^{k+1}(X_j)$.
\qed

\begin{Cor}\label{translationmatrix}
\mbox{}
\begin{enumerate}[(a)]
 \item Let $\Gamma$ be a connected component of $\Gamma_l$ that is not $\tau_\Omega$-periodic and is of finite left subgraph type. Let $X_1, \ldots, X_n$ denote the modules of a full sectional subgraph left of all projective and Ext-injective modules and $M$ its translation matrix. If $v_i$ for $i = 1, \ldots, n$ denotes the row vector with its $t$-th entry $v_{it}$ being the multiplicity of the simple module $S_i$ in the Jordan-H\"{o}lder composition series of $X_t$, then the multiplicity of $S_i$ in $\tau_\Omega^k(X_j)$ is given by $v_i \cdot M^k \cdot e_j$, where $e_j$ is the standard basis column vector where the $j$-th entry is $1$.
 \item Dually, let $\Gamma$ be a connected component of $\Gamma_r$ that is not $\tau_\Omega$-periodic and is of finite right subgraph type. Let $X_1, \ldots, X_n$ denote the modules of a full sectional subgraph right of all projective and Ext-injective modules and $M^-$ its cotranslation matrix. If $v_i$ for $i = 1, \ldots, n$ denotes the row vector with its $t$-th entry $v_{it}$ being the multiplicity of the simple module $S_i$ in the Jordan-H\"{o}lder composition series of $X_t$, then the multiplicity of $S_i$ in $\tau_\Omega^{-k}(X_j)$ is given by $v_i \cdot (M^-)^k \cdot e_j$.
\end{enumerate}

\end{Cor}
\proof
\\For convenience, let us set $(y_{1j}, \ldots, y_{nj})^T = M^k \cdot e_j$, which is the $j$-th column vector of $M^k$. Therefore, $y_{tj}$ is the number of $x_t$ in $\tau_\Omega^k(X_j)$ by the previous lemma. We conclude that
$$\sum_{t=1}^n  v_{it}y_{tj} = \sum_{t=1}^n (\# \text{ of }S_i \text{ in }X_t)(\# \text{ of }x_t \text{ in }\tau_\Omega^k(X_j))$$
and $v_i \cdot M^k \cdot e_j$ equals the number of $S_i$ in $\tau_\Omega^k(X_j)$. \qed
\bigskip
\\Obviously, $M$ and $M^-$ depend on the choice of the full sectional subgraph, but for our purposes it suffices to consider an arbitrary translation matrix, which we choose in a way that the computations are least complicated. If $\Gamma$ is a finite connected component of $\Gamma_l$, then there cannot be a a full sectional subgraph $\Sigma$ such that there is no path in $\Gamma_\Omega$ from a projective module to a module in $\Sigma$ by Theorem \ref{finitecomponent}. So if we try to calculate dimension vectors in finite components via $M$, then there must be an error after finitely many steps. This error can express itself in the following way.

\begin{Lemma}\label{translation}
\mbox{}
\begin{enumerate}[(a)]
\item Let $M$ be a translation matrix of a connected component $\Gamma$ of $\Gamma_l$ that is not $\tau_\Omega$-periodic and of finite left subgraph type. Then there are no standard basis vectors $e_j, e_l$ such that $M^k e_j = - e_l$ for some $k \in \N$ .
\item Dually, let $M^-$ be a cotranslation matrix of a connected component $\Gamma$ of $\Gamma_r$ that is not $\tau_\Omega$-periodic and of finite right subgraph type. Then there are no standard basis vectors $e_j, e_l$ such that $(M^-)^k e_j = - e_l$ for some $k \in \N$.
\end{enumerate}
\end{Lemma}

\proof
\\The modules of a full sectional subgraph left of all projective and Ext-injective modules are again denoted by $X_1, \ldots, X_n$. By Corollary \ref{translationmatrix} the multiplicity of $S_i$ in $\tau_\Omega^k(X_j)$ is given by $v_i \cdot M^k \cdot e_j$, where $v_i$ denotes the row vector with its $t$-th entry being the multiplicity of the simple module $S_i$ in the Jordan-H\"{o}lder composition series of $X_t$. Suppose there are standard basis vectors $e_j, e_l$ and some $k \in \N$ such that $M^k e_j = - e_l$. Since $X_l$ is non-zero, there exists a $v_i$ with its $l$-th entry being greater than $0$. We obtain that the multiplicity of $S_i$ in $\tau_\Omega^k(X_j)$ is
$$v_i \cdot M^k \cdot e_j = v_i \cdot -e_l < 0,$$
which is impossible.

\qed

\begin{Th}\label{Dynkinorbitgraph}
\mbox{}
\begin{enumerate}[(a)]
\item Let $\Gamma$ be a connected component of $\Gamma_l$ such that its left subgraph type is given by a Dynkin diagram. Then $\Gamma$ is finite and, in particular, $\tau_\Omega$-periodic.
\item Dually, if $\Gamma$ is a connected component of $\Gamma_r$ such that its right subgraph type is given by a Dynkin diagram, then $\Gamma$ is finite and $\tau_\Omega$-periodic.
\end{enumerate}

\end{Th}

\proof
\\Firstly, assume $\Gamma$ is of left subgraph type $A_1$, i.e. $\Gamma$ consists of just one $\tau_\Omega$-orbit. Every module in this orbit has a minimal right almost split morphism whose domain consists completely of modules that are not left stable. Hence the only orbit must be $\tau_\Omega$-periodic and $\Gamma$ is finite.
\bigskip
\\Let the left subgraph type of $\Gamma$ now be given by a Dynkin diagram other than $A_1$. Suppose $\Gamma$ is infinite, then it is clearly not $\tau_\Omega$-periodic as $\Gamma$ contains only finitely many $\tau_\Omega$-orbits. So there is a translation matrix $M$ that we can analyze.
\bigskip
\\ 1. $A_n$
\bigskip
\\Let the left subgraph type of $\Gamma$ be given by the Dynkin diagram $A_n$ for $n \geq 2$, i.e a full sectional subgraph left of all projective and Ext-injective modules and its relative Auslander-Reiten translates are given by
$$\begin{xy}
\xymatrixrowsep{0.5in}
\xymatrixcolsep{0.5in}

  \xymatrix@!0{ 
&	& 	&	&	&\tau_\Omega(X_1) \ar[dr] &	& X_1  \ar@{.>}[ll]    \\
&	&	&	&\tau_\Omega(X_2) \ar[ur]&	&X_2 \ar@{.>}[ll] \ar[ur]	&  \\
	&	&	&\ar[ur]	&	&\ar[ur]  \\
	&	&\ar@{.}[ur]	&	&\ar@{.}[ur]	&  \\
	&\tau_\Omega(X_{n-1})\ar[dr]\ar[ur]&	& X_{n-1} \ar@{.>}[ll] \ar[ur]&	&	&  \\
\tau_\Omega(X_n) \ar[ur]&	&X_n \ar@{.>}[ll] \ar[ur]&	&	&	&  }
\end{xy}$$

It is easy to see that the dimension vectors of $\tau_\Omega(X_1), \tau_\Omega(X_2), \ldots ,\tau_\Omega(X_{n-1})$, $\tau_\Omega(X_n)$ are given by $ x_2 - x_1, x_3 - x_1, \ldots , x_n - x_1, - x_1$ respectively.
This yields a translation matrix of $\Gamma$,
\begin{center}
\begin{math}
M_n = $\bordermatrix{ 	&	&	&	&	&	&	  	\cr
			&   -1	&   -1	&\cdots	&   -1	&   -1	&  -1	  	\cr
			&   1	&   	&   	&   	&   	&   	  	\cr
			&   	&   1	&   	&   	&   	&   	  	\cr
			&   	&   	&\ddots	&   	&   	&   	  	\cr
			&   	&   	&   	&   1	&   	&   	 	\cr
			&   	&   	&   	&   	&   1	&   	 	\cr} $
\end{math},
\end{center}
where all blank entries are zero. It is easy to see that $M_n \cdot e_n = -e_1$, which is a contradiction by Lemma \ref{translation}.
\bigskip
\\2. $D_n$
\bigskip
\\If $\Gamma$ is of left subgraph type $D_n$ for $n \geq 4$, a full sectional subgraph left of all projective and Ext-injective modules is obtained by attaching another module $X_n$ to the path of type $A_{n-1}$. We reorder the module names to generate a nice translation matrix.
$$\begin{xy}
\xymatrixrowsep{0.5in}
\xymatrixcolsep{0.5in}

  \xymatrix@!0{ 
	&	&	&	&	&\tau_\Omega(X_{n-1}) \ar[ddr]	&	&X_{n-1} \ar@{.>}[ll]	\\
	&	&	&	&	&\tau_\Omega(X_n) \ar[dr]	&	&X_n\ar@{.>}[ll]		\\
	&	&	&	&\tau_\Omega(X_1) \ar[ur] \ar[uur]&	&X_1 \ar@{.>}[ll] \ar[uur] \ar[ur]&  \\
	&	&	&\ar[ur]	&	&\ar[ur]	&	&  \\
	&	&\ar@{.}[ur]	&	&\ar@{.}[ur]	&	&  \\
	&\tau_\Omega(X_{n-3})\ar[dr]\ar[ur]&		& X_{n-3} \ar@{.>}[ll] \ar[ur]&	  \\
\tau_\Omega(X_{n-2}) \ar[ur]	&	&X_{n-2} \ar@{.>}[ll] \ar[ur]	&	  }
\end{xy}$$

We obtain the multiplicities of $x_i$ in $\tau_\Omega(X_j)$ in the same way as before by considering the almost split sequences ending in $X_j$. They generate translation matrices 

\begin{center}
\begin{math}
M_n = $\bordermatrix{ 	&	&	&	&	&	&	&	&	\cr
			&   1	&   1	&\cdots	&   1	&   1	&   1	&   1	&   1	\cr
			&   1	&   	&	&	&	&	&	&	\cr
			&   	&   1	&   	&   	&   	&   	&   	&   	\cr
			&   	&   	&\ddots	&   	&   	&   	&   	&   	\cr
			&   	&   	&   	&   1	&   	&   	&   	&   	\cr
			&   	&   	&   	&   	&   1	&   	&   	&   	\cr
			&   -1	&   -1	&\cdots	&  -1	&  -1	&  -1	&   -1	&   0	\cr
			&   -1	&   -1	&\cdots	&  -1	&  -1	&  -1	&   0	&  -1	\cr} $
\end{math}.
\end{center}
We observe that $M_n^{n-1} = -\Id$ for n even and
\begin{center}
\begin{math}
M_n^{n-1} = $\bordermatrix{ 	
					&   	&	&   	&   	&   	&   	&   	\cr
					&  -1	&	&	&	&	&	&	\cr
					&	&  -1	&	&	&	&	&	\cr
					&	&	&\ddots	&	&	&	&	\cr
					&	&	&	&  -1	&	&	&	\cr
					&	&	&	&	&  -1	&	&	\cr
					&	&	&	&	&	&   0	&  -1	\cr
					&	&	&	&	&	&  -1	&   0	\cr} $
\end{math}
\end{center}
for n odd. Clearly, in both cases we can find standard basis vectors $e_j$ and $e_l$ such that $M^{n-1} e_j = - e_l$, which again contradicts Lemma \ref{translation}.
\bigskip
\\3. $E_n$
\bigskip
\\Now we consider connected components whose left subgraph types are of Dynkin type $E_n$ for $n = 6, 7, 8$. A full sectional subgraph left of all projective and Ext-injective modules of type $E_n$ is given below.

$$\begin{xy}
\xymatrixrowsep{0.5in}
\xymatrixcolsep{0.5in}

  \xymatrix@!0{ 
  	&	&	&	&	&\tau_\Omega(X_{n-2}) \ar[dddr]	&	&X_{n-2}\ar@{.>}[ll]		&  \\
	&	&	&	&	&	&\tau_\Omega(X_n) \ar[dr]		&	&	X_n \ar@{.>}[ll]	\\
	&	&	&	&	&\tau_\Omega(X_{n-1}) \ar[dr]\ar[ur]	&	&X_{n-1}\ar@{.>}[ll]\ar[ur]	&  \\
	&	&	&	&\tau_\Omega(X_1) \ar[ur]\ar[uuur]&	&X_1\ar@{.>}[ll] \ar[ur]\ar[uuur]&	&  \\
	&	& 	&\ar[ur]	&	&\ar[ur]	&				&  \\
	&	&\ar@{.}[ur] 	&	&\ar@{.}[ur]	&	&	&  \\
	&\tau_\Omega(X_{n-4})\ar[dr]\ar[ur]&		& X_{n-4} \ar@{.>}[ll] \ar[ur]  \\
\tau_\Omega(X_{n-3}) \ar[ur]&		&X_{n-3} \ar@{.>}[ll] \ar[ur]	  }
\end{xy}$$
We again compute the dimension vectors of the relative Auslander-Reiten translates in terms of $x_i$ and obtain the following translation matrices.
\begin{center}
\begin{math}
M_6 = $\bordermatrix{ 	&	&	&	&	&	&	\cr
			&   1	&   1	&   1	&   1	&   1	&   0	\cr
			&   1	&   	&	&	&	&	\cr
			&   	&   1	&   	&   	&   	&   	\cr
			&  -1	&  -1	&  -1	&  -1	&   0	&   0	\cr
			&   0	&   0	&   0 	&   0	&   0	&   1	\cr
			&  -1	&  -1	&  -1	&   0	&  -1	&  -1	\cr} $
\end{math}
\end{center}
\begin{center}
\begin{math}
M_7 = $\bordermatrix{ 	&	&	&	&	&	&	&	\cr
			&   1	&   1	&   1	&   1	&   1	&   1	&   0	\cr
			&   1	&   	&	&	&	&	&	\cr
			&   	&   1	&   	&   	&   	&   	&   	\cr
			&   	&   	&   1	&   	&   	&   	&   	\cr
			&  -1	&  -1	&  -1	&  -1	&  -1	&   0	&   0	\cr
			&   0	&   0	&   0	&   0 	&   0	&   0	&   1	\cr
			&  -1	&  -1	&  -1	&  -1	&   0	&  -1	&  -1	\cr} $
\end{math}
\end{center}
\begin{center}
\begin{math}
M_8 = $\bordermatrix{ 	&	&	&	&	&	&	&	&	\cr
			&   1	&   1	&   1	&   1	&   1	&   1	&   1	&   0	\cr
			&   1	&   	&	&	&	&	&	&	\cr
			&   	&   1	&   	&   	&   	&   	&   	&   	\cr
			&   	&   	&   1	&   	&   	&   	&   	&   	\cr
			&   	&   	&   	&   1	&   	&   	&   	&   	\cr
			&  -1	&  -1	&  -1	&  -1	&  -1	&  -1	&   0	&   0	\cr
			&   0	&   0	&   0	&   0	&   0 	&   0	&   0	&   1	\cr
			&  -1	&  -1	&  -1	&  -1	&  -1	&   0	&  -1	&  -1	\cr} $
\end{math}
\end{center}

Direct calculation shows that $M_8^{15} = -\Id$ and $M_7^9 = - \Id$ . It also turns out that $M_6^6 \cdot e_1 = - e_1$; therefore, we can apply Lemma \ref{translation} to obtain a contradiction in all three cases. Hence $\Gamma$ must be finite and $\tau_\Omega$-periodic if its left subgraph type is given by a Dynkin diagram. \qed

\section{Translation quivers and Riedtmann's structure theorem}\label{treetypesection}

Now we consider the concept of assigning a graph to a stable Auslander-Reiten quiver introduced by Riedtmann, which does not work for non-stable components. The whole section closely follows \cite{R80}. Let $\Gamma$ be a quiver with vertex set $\Gamma_0$ and arrow set $\Gamma_1$ and let $\tau$ be an injective map from a subset of $\Gamma_0$ into $\Gamma_0$. For any vertex $x$ let $x^-$ denote the set of immediate predecessors and $x^+$ the set of immediate successors, i.e.
$$x^- = \{y \in \Gamma_0 | \text{ there is an arrow } y \rightarrow x\},$$
$$x^+ = \{y \in \Gamma_0 | \text{ there is an arrow } x \rightarrow y\}.$$

\begin{Def}
The pair $(\Gamma, \tau)$ is called a \textbf{translation quiver} if
\begin{enumerate}[(a)]
 \item $\Gamma$ neither has loops nor multiple arrows between two vertices.
 \item  Whenever $x \in \Gamma_0$ is such that $\tau(x)$ is defined then $x^- = \tau(x)^+$. 
\end{enumerate}
Moreover, we call $\tau$ the translation of the quiver and say that $\Gamma$ is connected if for all vertices $x$ and $y$ in $\Gamma$ there is a walk between $x$ and $y$.
\end{Def}

\begin{Def}
We say a translation quiver $(\Gamma, \tau)$ is \textbf{locally finite} if $x^+$ and $x^-$ are finite sets for every vertex $x$.
\end{Def}

Many stable Auslander-Reiten quivers of functorially finite resolving subcategories are locally finite translation quivers. However, there are two cases when stable Auslander-Reiten quivers are not translation quivers. The first case are components that contain arrows with non-trivial valuation, i.e. multiple arrows between modules. Secondly, although impossible in $A$-mod, we cannot eliminate the possibility that in a proper functorially finite resolving subcategory of $A$-mod there are irreducible morphisms from a stable module to itself, which would be considered a loop in the translation quiver by Riedtmann's construction.

\begin{Def}
A morphism of quivers $\varphi : \Gamma \to \Gamma'$ assigns to each vertex $x$ of $\Gamma$ a vertex $\varphi(x)$ of $\Gamma'$ and to each arrow $\alpha$ from $x$ to $y$ in $\Gamma$ an arrow $\varphi(\alpha)$ from $\varphi(x)$ to $\varphi(y)$ in $\Gamma'$. We call $\varphi$ a \textbf{morphism of translation quivers} if it commutes with the translation.
\end{Def}

\begin{Def}
A \textbf{directed tree} is a directed graph without loops, multiple arrows or cycles.
\end{Def}

To a directed tree $B$ we associate a translation quiver $\Z B$ as follows. The vertices of $\Z B$ are the pairs $(n,x)$ with $n \in \Z$ and $x$ a vertex of $B$. For each arrow $x \to y$ in $B$ and each $n \in \Z$ we have two arrows $(n,x) \to (n,y)$ and $(n,y) \to (n-1,x)$. The translation is defined via $\tau(n,x) = (n+1,x)$. We consider $B$ to be embedded in $\Z B$ as the subgraph consisting of vertices $(0,x)$ and the arrows connecting them.

\begin{Def}
A group $\Pi$ of automorphisms of a translation quiver $\Gamma$ is said to be \textbf{admissible} if $\{\varphi(n,x)|\varphi \in \Pi\}$ does not contain two vertices with an arrow between them for all vertices $(n,x)$ in $\Gamma$. 
\end{Def}

We then construct the quotient quiver $\Gamma/\Pi$ consisting of vertices $(\Gamma/\Pi)_0 = \Gamma_0/\Pi$ and arrows $(\Gamma/\Pi)_1 = \Gamma_1/\Pi$. Note that $\Gamma/\Pi$ is a again a translation quiver which is stable if $\Gamma$ is stable \cite{R80}.

\begin{Th}\cite[Struktursatz]{R80} \label{Riedtmann}
Let $\Gamma$ be a connected translation quiver, then there is a directed tree $B$ and an admissible group of automorphisms $\Pi \subset Aut (\Z B)$ such that $\Gamma \cong \Z B/\Pi$. The undirected graph $\overline{B}$ associated to $B$ is uniquely determined by $\Gamma$ up to canonical isomorphism, and $\Pi$ is uniquely defined up to conjugation in $\Aut(\Z B)$.
\end{Th}
\proof
\\As it is necessary for further proofs, we give the construction of $B$. We fix a vertex $x$ of $\Gamma$ and consider the set of sectional paths starting in $x$. There is a vertex in $B$ for each sectional path $(x = y_o \to y_1 \to \cdots \to y_n)$. The arrows of $B$ are 
$(x = y_o \to y_1 \to \cdots \to y_n) \longrightarrow (x = y_o \to y_1 \to \cdots \to y_n \to y_{n+1})$. The rest of the proof can be found in \cite{R80}.
\qed

\begin{Def}
Under the assumptions of Theorem \ref{Riedtmann} we call the graph $\overline{B}$ associated to $B$ the \textbf{tree type of $\pmb{\Gamma}$}.
\end{Def}

\begin{Cor}
Let $\Gamma$ be a connected component of the stable Auslander-Reiten quiver $\Gamma_s$ that does not contain loops or multiple arrows, then there is a directed tree $B$ and an admissible group of automorphisms $\Pi \subset Aut (\Z B)$ such that $\Gamma \cong \Z B/\Pi$. The undirected graph $\overline{B}$ associated to $B$ is uniquely determined by $\Gamma$ up to canonical isomorphism, and $\Pi$ is uniquely defined up to conjugation in $\Aut(\Z B)$.
\end{Cor}

\proof
\\Since $\Gamma$ does not contain loops or multiple arrows, it must be a stable translation quiver, so the statement follows immediately from Theorem \ref{Riedtmann}. \qed

\begin{Cor}\label{treetypeandsubgraphtypecoincide}
 Let $\Gamma$ be a connected component of the stable Auslander-Reiten quiver $\Gamma_s$ that does not contain loops or multiple arrows. Then its tree type $\overline{B}$ coincides with its subgraph type $\overline{\Sigma}$.
\end{Cor}
\proof
\\We fix an indecomposable module $X$ and let $\Sigma$ be the subgraph of $\Gamma$ consisting of all sectional paths starting in $X$. Then $\Sigma$ is a full sectional subgraph and the undirected graph $\overline{\Sigma}$ associated to $\Sigma$ coincides with $\overline{B}$ by the construction given in the proof of Theorem \ref{Riedtmann}. \qed
\bigskip
\\We can now proceed preparing the proof of the main result using subadditive functions. Recall the definitions.

\begin{Def}
Let $T$ be a simple-laced graph consisting of edges $T_1$ and vertices $T_0$ such that each vertex has at most finitely many edges attached to it. The \textbf{Cartan matrix} $C(T) = (c_{ij})$ of $T$ is defined as follows. If $i,j \in T_0$, then

$$
c_{ij} =\left\{\begin{array}{cl} 2 \quad&\mbox{\emph{if i = j}}\\
-1 \quad&\mbox{\emph{if there is an edge between i and j}}\\
0 \quad&\mbox{\emph{otherwise.}}\end{array}\right.
$$
\end{Def}

\begin{Def}
A function $n: T_0 \to \Q,\ x \mapsto n_x$ is called \textbf{subadditive} if 
$$\sum_{x \in T_0} c_{xy}n_x \geq 0 \text{ for all }y \in T_0.$$
It is called \textbf{additive} if equality holds for all $y \in T_0$.
\end{Def}

Let us now recall an important result on graphs and subadditive functions.

\begin{Th} \label{additiveDynkin} \cite[Theorem 4.5.8]{B91}
\\Let $T$ be a connected labeled graph and $x \mapsto n_x$ a subadditive function on $T$. If $x \mapsto n_x$ is not additive, then T is a finite Dynkin diagram or $A_\infty$.
\end{Th}

\begin{Th}\label{treetypedynkin}
Suppose $\Gamma$ is a finite, connected component of $\Gamma_s$ without loops or multiple arrows. Then the tree type of $\Gamma$ is a Dynkin diagram and coincides with its left subgraph type.
\end{Th}

\proof
\\Since $\Gamma$ is finite, we know that there are only finitely many orbits in $\Gamma$ and there is an $n \in \N$ such that $\tau_\Omega^n(X) = X$ for all modules $X$ in $\Gamma$. Let $X_1, \ldots, X_n$ denote the modules of an arbitrary orbit, while $Y_{11}, \ldots Y_{1n}, \ldots , Y_{k1}, \ldots Y_{kn}$ name the modules of the adjacent orbits. This means $Y_{1i} \oplus \cdots \oplus Y_{ki}$ is a direct summand of the middle term of $\E_\Omega(X_{i+1})$ and hence
$$l(X_i) + l(X_{i+1}) \geq \sum_{j=1}^k l(Y_{ji}).$$
Clearly, equality holds if and only if all modules in the corresponding almost split sequences are periodic. Adding these equations we obtain 
$$ 2 \sum_{i=1}^n l(X_i) \geq \sum_{j=1}^k \sum_{i=1}^n  l(Y_{ji}).$$
This gives rise to a subadditive function which assigns the sum of the length of $n$ modules $X, \tau_\Omega(X), \ldots, \tau_\Omega^{n-1}(X)$ of an orbit to the corresponding vertex in $\Gamma$'s tree type. $\Gamma$ is finite but does not contain any projective modules, so it cannot be a connected component of $\Gamma_\Omega$ by Theorem \ref{finitecomponent}. Consequently, there is a finite orbit in $\Gamma_\Omega$ adjacent to $\Gamma$ by Lemma \ref{orbit} and thus the aforementioned subadditive function is not additive. By Theorem \ref{additiveDynkin} the tree type of $\Gamma$ must be given by a Dynkin diagram or $A_\infty$. But $A_\infty$ can only be the tree type of an infinite component and hence the tree type of $\Gamma$ is Dynkin.
\bigskip
\\By Lemma \ref{orbit} every Ext-injective module adjacent to any orbit of $\Gamma$ must be contained in a finite orbit. Consequently, $\Gamma$ must also be a connected component of $\Gamma_l$, so left subgraph type, subgraph type and tree type of $\Gamma$ coincide by definition and Corollary \ref{treetypeandsubgraphtypecoincide} respectively. \qed


\section{The main result}

In this section we combine the results on the tree type and the left subgraph type and obtain that $\Omega$ is finite if and only if the left subgraph types of all connected components of the left stable Auslander-Reiten quiver are given by Dynkin diagrams. Firstly, we observe how the left subgraph type behaves if a connected component of $\Gamma_l$ is not a translation quiver. Therefore, we consider the case of multiple arrows between a pair of indecomposable modules in the left stable part of the quiver. It is well known that all algebras of finite representation type do not have multiple arrows in their Auslander-Reiten quivers \cite[VII Proposition 2.2]{ARS95}. However, there are examples of finite Auslander-Reiten quivers of functorially finite resolving subcategories with multiple arrows between modules, where they occur between modules in finite orbits.

\begin{Ex}
Let $A$ be the path algebra of the Kronecker-quiver
$$\begin{xy}
  \xymatrix{ 
    e_2 \ar@<1ex>[r] \ar@<-1ex>[r] &	e_1	    }
\end{xy}.
$$ 
Then there is a functorially finite resolving subcategory of $A$-mod that is finite and contains multiple arrows in its finite component.
\end{Ex}
As $A$ is a hereditary algebra, it has a unique preprojective component, which looks the following way. Note that the projective module $P_i$ is given by $Ae_i$.

$$\begin{xy}
\xymatrixrowsep{0.5in}
\xymatrixcolsep{0.5in}

  \xymatrix@!0{ 
  &	[P_2	 \ar@<1ex>[dr] \ar@<-1ex>[dr]	&	& \tau^{-1}(P_2) \ar@{.>}[ll] \ar@<1ex>[dr] \ar@<-1ex>[dr]	&	& \cdots \ar@{.>}[ll]	\\
[P_1 \ar@<1ex>[ur] \ar@<-1ex>[ur] & 	& \tau^{-1}(P_1) \ar@{.>}[ll] \ar@<1ex>[ur] \ar@<-1ex>[ur]&	& \tau^{-2}(P_1) \ar@{.>}[ll] &	&\cdots \ar@{.>}[ll]  }
\end{xy}$$
 
Since the projective modules are also the standard modules in this order of primitive idempotents, the whole category $\F(\Delta)$ of standard-filtered modules consists only of $P_1$ and $P_2$. Consequently, its Auslander-Reiten quiver is given by 
$$\begin{xy}
\xymatrixrowsep{0.5in}
\xymatrixcolsep{0.5in}

  \xymatrix@!0{ 
  &	[P_2]	 	\\
[P_1] \ar@<1ex>[ur] \ar@<-1ex>[ur]   }
\end{xy}$$
and contains a double arrow, although the subcategory is clearly finite.

\begin{Th}\label{multiplearrows}
\mbox{}
\begin{enumerate}[(a)]
 \item Let $\Gamma$ be a connected component of $\Gamma_l$. If $\Gamma$ contains multiple arrows from a module $X$ to a module $Y$, then $l(\tau_\Omega^n(Y)) \to \infty$ as $n \to \infty$. In particular, $\Gamma$ is infinite.
 \item Dually, if $\Gamma$ is a connected component of $\Gamma_r$ and contains multiple arrows from a module $Y$ to a module $X$, then $l(\tau_\Omega^n(Y)) \to \infty$ as $n \to - \infty$. In particular, $\Gamma$ is infinite.
\end{enumerate}
\end{Th}

\proof
\\Suppose first that $l(\tau_\Omega^n(X)) \leq l(\tau_\Omega^n(Y))$ and $l(\tau_\Omega^{n+1}(Y)) \leq l(\tau_\Omega^n(X))$ for every $n \geq 0$. It follows that there must be some $m \geq 0$ such that all modules $\tau_\Omega^n(X)$ and $\tau_\Omega^k(Y)$ have the same Jordan-H\"{o}lder length for all $n, k \geq m$. Without loss of generality, we can assume that $m = 0$. Since there are multiple arrows from $X$ to $Y$, we have that $\E_\Omega(Y)$ must be of the form
$$\begin{xy}
  \xymatrix{ 
    0 \ar[r]   &   \tau_\Omega(Y) \ar[r] 	&	X \oplus X \ar[r]    &    Y \ar[r]  	&	0	    }
\end{xy}$$
in order to satisfy $l(X \oplus X) = l(Y) + l(\tau_\Omega(Y))$. By Corollary \ref{valuationisinvariantoftranslation} we inductively conclude that for every $n \geq 0$ the almost split sequence $\E_\Omega(\tau_\Omega^n(Y))$ is of the form
$$\begin{xy}
  \xymatrix{ 
0 \ar[r]   &   \tau_\Omega^{n+1}(Y) \ar[r] 	&	\tau_\Omega^n(X) \oplus \tau_\Omega^n(X) \ar[r]    &    \tau_\Omega^n(Y) \ar[r]  	&	0,    }
\end{xy}$$
while $\E_\Omega(\tau_\Omega^n(X))$ looks as follows.
$$\begin{xy}
  \xymatrix{ 
  0 \ar[r]   &   \tau_\Omega^{n+1}(X) \ar[r] 	&	\tau_\Omega^n(Y) \oplus \tau_\Omega^n(Y) \ar[r]    &    \tau_\Omega^n(X) \ar[r]  	&	0}
\end{xy}$$
It follows that the connected component of $X$ and $Y$ in $\Gamma_l$ is of left subgraph type $\widetilde{A}_1$ and has the following shape.
$$\begin{xy}
\xymatrixrowsep{0.5in}
\xymatrixcolsep{0.5in}

  \xymatrix@!0{ 
&  \cdots &	 &\tau_\Omega(X)\ar@{.>}[ll] \ar@<1ex>[dr] \ar@<-1ex>[dr]&	& X \ar@{.>}[ll] \ar@<1ex>[dr] \ar@<-1ex>[dr]	&	& \cdots \ar@{.>}[ll]	\\
\cdots &	 & \tau_\Omega^2(Y)\ar@{.>}[ll] \ar@<1ex>[ur] \ar@<-1ex>[ur] & 	& \tau_\Omega(Y) \ar@{.>}[ll] \ar@<1ex>[ur] \ar@<-1ex>[ur]&	& Y \ar@{.>}[ll] &	&\cdots \ar@{.>}[ll]  }
\end{xy}$$

On the other hand, we know that there must be a non-zero morphism $f \in \rad_\Omega^\infty (P,Y)$ for some projective module $P$. But then $f$ factors over all modules $\tau_\Omega^n(X)$ and $\tau_\Omega^n(Y)$ with $n \geq 0$ and hence must be zero by Lemma \ref{HaradaSai}, which is a contradiction.
Consequently, there is an $n \geq 0$ such that $l(\tau_\Omega^n(X)) > l(\tau_\Omega^n(Y))$ or $l(\tau_\Omega^{n+1}(Y)) > l(\tau_\Omega^n(X))$. Without loss of generality, we can assume that $l(X) > l(Y)$ and there are $k$ arrows from $X$ to $Y$, where $k \geq 2$. There is an almost split sequence
$$\begin{xy}
  \xymatrix{ 
    0 \ar[r]   &   \tau_\Omega (Y) \ar[r] 	&	Z \ar[r]    &    Y \ar[r]  	&	0	    }
\end{xy}.
$$ 
in which $X$ occurs at least twice as a direct summand of $Z$, thus $l(\tau_\Omega (Y)) > l(X)$. Because there are $k$ arrows from $\tau_\Omega (Y)$ to $X$ by Corollary \ref{valuationisinvariantoftranslation}, we can continue inductively and obtain $l(\tau_\Omega^{n+1} (Y)) > l(\tau_\Omega^n(Y))$ for all $n \geq 0$. \qed
\begin{Cor}\label{multiplearrowscontrapositive}
If $\Gamma$ is finite, then its left subgraph type does not contain multiple edges between vertices. 
\end{Cor}

\proof
\\The statement is the contrapositive of Theorem \ref{multiplearrows}. \qed
\bigskip
\\The next step is to consider a left stable component that contains a loop, i.e. a module $X$ such that there is an $\Omega$-irreducible morphism to itself. Such morphisms do exist when the module $X$ is projective and Ext-injective.

\begin{Ex}
Let A be the path algebra
$$\begin{xy}
\xymatrix{
  e_1  \ \ar@(ur,dr)^\beta  \\
}
\end{xy}$$
with the relation $\beta^2 = 0$. Then the subcategory $\Omega = \add(P_1)$ consisting only of copies of the one indecomposable projective module $P_1$ is functorially finite, resolving and there are $\Omega$-irreducible morphisms from $P_1$ to itself.
\end{Ex}

The Auslander-Reiten quiver of $A$-mod is given by 

$$\begin{xy}
\xymatrixrowsep{0.5in}
\xymatrixcolsep{0.5in}

  \xymatrix@!0{ 
& [P_1] \ar[dr] &	&[P_1] \ar[dr] 	&	\\
S_1 \ar[ur]& 	&S_1 \ar[ur] \ar@{.>}[ll]&	&S_1. \ar@{.>}[ll]  }
\end{xy}$$

Therefore, every module in $A$-mod has a right and a left $\Omega$-approximation, so $\Omega$ is functorially finite. Moreover, it is closed under direct summands, extensions, kernels of epimorphisms and contains $_AA$. Hence $\Omega$ is a functorially finite resolving subcategory and, since $P_1$ is the only indecomposable module in $\Omega$, the morphism from $P_1$ to itself factoring over $S_1$ is clearly irreducible in $\Omega$.
\bigskip
\\Although there are to my best knowledge no examples of such morphisms in left stable components, we cannot rule out their existence. However, the existence of such a morphism provides enough information to determine the left subgraph type of the left stable connected component it is contained in.

\begin{Th}\label{irreduciblemorphismtoitself}
Let $X$ be an indecomposable module in $\Omega$ that is not projective or not Ext-injective. If there is an $\Omega$-irreducible morphism $f : X \to X$, then $\tau_\Omega(X) = X$ and the connected components of $\Gamma_l$, $\Gamma_r$ and $\Gamma_s$ containing $X$ coincide and are of subgraph type $A_n$ for some $n <  2^{l(X)}$.
\end{Th}
\proof
\\Let $f: X \to X$ be an irreducible morphism in $\Omega$. Then there exists an $k \in \N$ such that the composition $f^k$ is zero by Lemma \ref{HaradaSai}. Hence if $X$ is not projective then the path obtained by composing $f$ with itself cannot be sectional by Corollary \ref{non-zero}, so we have $\tau_\Omega(X) = X$. On the other hand, if $X$ is projective, we then know by our assumptions that $X$ is not Ext-injective. Therefore, we can apply the same arguments to $\tau_\Omega^{-1}(X)$ and analogously obtain $\tau_\Omega(X) = X$, which contradicts $X$ to be projective. Consequently, $\tau_\Omega(X) = X$ holds under our assumptions and $X$ is neither projective nor Ext-injective, but stable and $\tau_\Omega$-periodic. Moreover, all orbits in the connected components of $\Gamma_l$ and $\Gamma_r$ containing $X$ must in fact be periodic by Lemma \ref{orbit}. We denote the stable connected component of $X$ by $\Gamma$.
\bigskip
\\Assume now the arrow from $X$ to itself in the Auslander-Reiten quiver does not have valuation $(1,1)$. Then $X \oplus X$ is a direct summand of the middle term of $\E_\Omega(X)$, so by Jordan-H\"{o}lder length arguments the arrow's valuation is $(2,2)$ and $X \oplus X$ is in fact the whole middle term of $\E_\Omega(X)$. Consequently, we have a component of the Auslander-Reiten quiver $\Gamma_\Omega$ that contains only $X$, in particular, it does not contain a projective module, which contradicts Theorem \ref{finitecomponent}. Thus the valuation of the arrow from $X$ to itself must be $(1,1)$. 
Let
$$\begin{xy}
  \xymatrix{ 
      X = X_{-1} \ar[r]^f	&	X = X_0 \ar[r]^{f_1} 	&	X_1 \ar[r]^{f_2}     &   \cdots	 \ar[r]	&	X_{m-1} \ar[r]^{f_m}	&	X_m.    }
\end{xy}$$
be a sectional path in the stable component of $X$, then the following holds.
\begin{enumerate}
\item $\tau_\Omega(X_j) = X_j$ for $j = 0, \ldots, m$.
\item $l(X_{j-1}) \geq l(X_{j})$ for $j = 0, \ldots, m$.
\item The almost split sequence $\E_\Omega(X_j)$ contains at most two non-projective indecomposable middle terms for $j = 0, \ldots, m$.
\end{enumerate}
We prove all three statements by induction on $m$. Let $m = 0$, then the first statement holds as we have seen earlier in the proof of this theorem. Moreover, $l(X_{-1}) \geq l(X_0)$ holds as $X_{-1} = X = X_0$. Lastly, we consider the almost split sequence
$$\begin{xy}
  \xymatrix{ 
      0 \ar[r]  	& 	X \ar[r] 	&	X \oplus M \ar[r]     &   X \ar[r]	&	0.	    }
\end{xy}$$
It is easy to see that $l(M) = l(X)$. Suppose now $M$ has two non-projective indecomposable direct summands $M_1$ and $M_2$, then 
%
$$\begin{xy}
\xymatrixrowsep{0.5in}
\xymatrixcolsep{0.5in}

  \xymatrix@!0{ 
 &\cdots&		& \tau_\Omega(M_2) \ar[ddr] \ar@{.>}[ll] &	& M_2 \ar[ddr]  \ar@{.>}[ll]	&	& \cdots \ar@{.>}[ll]	&\\
 &\cdots&		& \tau_\Omega(M_1) \ar[dr] \ar@{.>}[ll] &	& M_1 \ar[dr]  \ar@{.>}[ll]	&	& \cdots \ar@{.>}[ll]	&\\
\cdots&		& X \ar[dr] \ar[ur] \ar[uur] \ar@{.>}[ll]&	& X \ar[dr] \ar[ur] \ar [uur]\ar@{.>}[ll]	&	& X \ar[dr]  \ar@{.>}[ll]	&	&\cdots \ar@{.>}[ll]\\
&\cdots&		& X \ar[dr] \ar[ddr] \ar[ur]\ar@{.>}[ll] &	& X \ar[dr] \ar[ddr] \ar[ur] \ar@{.>}[ll]	&	& X   \ar@{.>}[ll]	&	&\cdots \ar@{.>}[ll] \\
& &\cdots&		& \tau_\Omega(M_1)  \ar[ur] \ar@{.>}[ll]&	& M_1  \ar[ur] \ar@{.>}[ll]	&	& \cdots \ar@{.>}[ll]	&\\
& &\cdots&		& \tau_\Omega(M_2)  \ar[uur] \ar@{.>}[ll]&	& M_2  \ar[uur] \ar@{.>}[ll]	&	& \cdots \ar@{.>}[ll]	&}
\end{xy}$$
is an example of how the Auslander-Reiten quiver could be shaped around $X$. Furthermore, there are two almost split sequences
$$\begin{xy}
  \xymatrix{ 
      0 \ar[r]  	& 	\tau_\Omega(M_i) \ar[r] 	&	X \oplus N_i \ar[r]     &   M_i \ar[r]	&	0.	    }
\end{xy}$$
Since there is an irreducible morphism from $\tau_\Omega(M_i)$ to $X$, $\tau_\Omega(M_i)$ is a direct summand of $M$. Let us assume $\tau_\Omega(M_1)$ is not isomorphic to either $M_1$ or $M_2$. But then $l(M_1) + l(\tau_\Omega(M_1)) < l(X)$ by $\E_\Omega(X)$ and, on the other hand, $l(M_1) + l(\tau_\Omega(M_1)) \geq l(X)$ by $\E_\Omega(M_1)$. Thus $\tau_\Omega(M_1)$ and must be isomorphic to either $M_1$ or $M_2$. We repeat the argument for $M_2$, which leaves us with two cases.
\bigskip
\\The first one is $\tau_\Omega(M_1) = M_1$ and $\tau_\Omega(M_2) = M_2$, which gives us almost split sequences
$$\begin{xy}
  \xymatrix{ 
      0 \ar[r]  	& 	M_i \ar[r] 	&	X \oplus N_i \ar[r]     &   M_i \ar[r]	&	0	    }
\end{xy}$$
and, therefore, $l(X) \leq 2 l(M_i)$ for $i = 1,2$. Putting these inequations into $l(M) = l(X)$ we obtain $2 l(M_1) + 2l(M_2) \leq 2l(X) \leq 2 l(M_1) + 2l(M_2)$ and deduce $2 l(M_i) = l(X)$. Thus $N_i$ must be zero and the whole component of $X$ only consists of $X, M_1$ and $M_2$ but does not contain any projective module, which contradicts Theorem \ref{finitecomponent}. In the second case we have $\tau_\Omega(M_1) = M_2$ and $\tau_\Omega(M_2) = M_1$. We obtain an almost split sequence
$$\begin{xy}
  \xymatrix{ 
      0 \ar[r]  	& 	M_1 \ar[r] 	&	X \oplus N_2 \ar[r]     &   M_2 \ar[r]	&	0	    }
\end{xy}$$
from which we deduce $l(X) = l(M_1) + l(M_2)$ in the same way as in the first case. Hence the component of $X$ is finite again but does not contain projective modules. Due to this contradiction we know that $M$ has at most one indecomposable non-projective direct summand.
\bigskip
\\Suppose now the statements have been proved for $m - 1$. By the induction hypothesis there is an almost split sequence
$$\begin{xy}
  \xymatrix{ 
      0 \ar[r]  	& 	X_{m-1} \ar[r] 	&	X_{m-2} \oplus X_m \oplus P \ar[r]     &   X_{m-1} \ar[r]	&	0	    }
\end{xy}$$
that has two non-projective middle terms $X_{m-2}$ and $X_m$, from which we deduce that $\tau_\Omega(X_m) = X_m$. This proves the first statement. Moreover, the induction hypothesis implies that $l(X_{m-2}) \geq l(X_{m-1})$. Consequently, $l(X_{m-1}) \geq l(P) \oplus l(X_m)$ and, in particular, $l(X_{m-1}) \geq l(X_m)$ hold, which proves the second statement. 
\bigskip
\\Consider now the almost split sequence
$$\begin{xy}
  \xymatrix{ 
      0 \ar[r]  	& 	X_m \ar[r] 	&	X_{m-1} \oplus M \ar[r]     &   X_m \ar[r]	&	0.	    }
\end{xy}$$
Suppose now $M$ has two indecomposable non-projective direct summands $M_1$ and $M_2$, then
$$\begin{xy}
\xymatrixrowsep{0.5in}
\xymatrixcolsep{0.5in}

  \xymatrix@!0{  
\cdots&		& X_{m-1} \ar[dr] \ar@{.>}[ll]&	& X_{m-1} \ar[dr] \ar@{.>}[ll]	&	& X_{m-1} \ar[dr]  \ar@{.>}[ll]	&	&\cdots \ar@{.>}[ll]\\
&\cdots&		& X_m \ar[dr] \ar[ddr] \ar[ur]\ar@{.>}[ll] &	& X_m \ar[dr] \ar[ddr] \ar[ur] \ar@{.>}[ll]	&	& X_m   \ar@{.>}[ll]	&	&\cdots \ar@{.>}[ll] \\
& &\cdots&		& \tau_\Omega(M_1)  \ar[ur] \ar@{.>}[ll]&	& M_1  \ar[ur] \ar@{.>}[ll]	&	& \cdots \ar@{.>}[ll]	&\\
& &\cdots&		& \tau_\Omega(M_2)  \ar[uur] \ar@{.>}[ll]&	& M_2  \ar[uur] \ar@{.>}[ll]	&	& \cdots \ar@{.>}[ll]	&}
\end{xy}$$
is an example of what the Auslander-Reiten quiver could look like around $X_m$. Since $\tau_\Omega(X_m) = X_m$, both $\tau_\Omega(M_1)$ and $\tau_\Omega(M_2)$ must be direct summands of $M$, which is a module whose length is shorter or equal than the length of $X_m$. On the other hand, by the almost split sequences 
$$\begin{xy}
  \xymatrix{ 
      0 \ar[r]  	& 	\tau_\Omega(M_1) \ar[r] 	&	X_m \oplus N_1 \ar[r]     &   M_1 \ar[r]	&	0	    }
\end{xy}$$
and
$$\begin{xy}
  \xymatrix{ 
      0 \ar[r]  	& 	\tau_\Omega(M_2) \ar[r] 	&	X_m \oplus N_2 \ar[r]     &   M_2 \ar[r]	&	0	    }
\end{xy}$$
we know that $l(\tau_\Omega(M_1)) + l(M_1) \geq l(X_m)$ and $l(\tau_\Omega(M_2)) + l(M_2) \geq l(X_m)$. It follows that $2l(M) \geq l(\tau_\Omega(M_1)) + l(M_1) + l(\tau_\Omega(M_2)) + l(M_2) \geq 2l(X_m) \geq 2l(M)$, hence $M \cong M_1 \oplus M_2$ and $l(X_m) = l(M_1) + l(M_2)$. In particular, $\tau_\Omega(M_1) \oplus \tau_\Omega(M_2)$ is isomorphic to $M_1 \oplus M_2$, which leaves us with the same two cases as in the induction start. In the first case there are almost split sequences
$$\begin{xy}
  \xymatrix{ 
      0 \ar[r]  	& 	M_1 \ar[r] 	&	X_m \ar[r]     &   M_1 \ar[r]	&	0	    }
\end{xy}$$
and
$$\begin{xy}
  \xymatrix{ 
      0 \ar[r]  	& 	M_2 \ar[r] 	&	X_m \ar[r]     &   M_2 \ar[r]	&	0	    }
\end{xy}$$
with $2 l(M_1) = 2 l(M_2) = l(X_m)$. In the second case the almost split sequences are given by
$$\begin{xy}
  \xymatrix{ 
      0 \ar[r]  	& 	M_2 \ar[r] 	&	X_m \ar[r]     &   M_1 \ar[r]	&	0	    }
\end{xy}$$
and
$$\begin{xy}
  \xymatrix{ 
      0 \ar[r]  	& 	M_1 \ar[r] 	&	X_m \ar[r]     &   M_2 \ar[r]	&	0.	    }
\end{xy}$$
In both cases the component ends in these orbits, in particular, there are only finitely many modules in the stable component of $X$ and the almost split sequence 
$$\begin{xy}
  \xymatrix{ 
      0 \ar[r]  	& 	X_m \ar[r] 	&	X_{m-1} \oplus M_1 \oplus M_2 \ar[r]     &   X_m \ar[r]	&	0	    }
\end{xy}$$
gives us $l(X_m) = l(X_{m-1})$. Now we show that this cannot hold if there is a projective module in this component. Since $X, X_1, \ldots, X_m, M_1$ and $M_2$ are the only indecomposable non-projective modules in this component, there must be a projective module $P$ such that there is an almost split sequence
$$\begin{xy}
  \xymatrix{ 
      0 \ar[r]  	& 	X_j \ar[r] 	&	X_{j-1} \oplus X_{j+1} \oplus P \ar[r]     &   X_{j-1} \ar[r]	&	0.	    }
\end{xy}$$
Without loss of generality, we can assume that $j$ is maximal with that property. We obtain $l(X_j) > l(X_{j+1})$ by the almost split sequence and the induction hypothesis. As all the almost split sequences $\E_\Omega(X_{j+1}), \ldots, \E_\Omega(X_{m-1})$ have precisely two middle terms given by the original sectional path, it easily follows that $l(X_j) > l(X_{j+1}) > \cdots > l(X_{m-1}) > l(X_m)$, which contradict our previous result. Hence our assumption that $\E_\Omega(X_m)$ has three non-projective middle terms cannot be true, which completes the proof for the third statement. 
\bigskip
\\We deduce that every module in $\Gamma$ has length shorter or equal to $l(X)$ by statements 1 and 2 and is of $\tau_\Omega$-periodicity $1$. In a stable component every full sectional subgraph is of the same type by Lemma \ref{subgraphtype}. Therefore, let $\Sigma$ be an arbitrary full sectional subgraph. By statement 3 we know that for every module $Y$ in $\Sigma$ the almost split sequence 
$$\begin{xy}
  \xymatrix{ 
      0 \ar[r]  	& 	Y \ar[r] 	&	M \ar[r]     &   Y \ar[r]	&	0	    }
\end{xy}$$
contains at most two indecomposable non-projective middle terms, hence every orbit in $\Gamma$ has at most two adjacent orbits in $\Gamma$. It follows that every vertex in the subgraph type of $\Gamma$ has at most two edges connected to it, so the subgraph type of $\Gamma$ is either $A_n$, $A_\infty$ or $A_\infty^\infty$. Suppose it is not $A_n$ for some $n <  2^{l(X)}$, then there is a sectional path 
$$\begin{xy}
  \xymatrix{ 
Y_0 \ar[r]^{g_1}	& Y_1 \ar[r]^{g_2} &    \cdots \ar[r]^{g_{m-1}} 	&	Y_{m-1} \ar[r]^{g_m}     &   Y_m	  }
\end{xy}$$
in $\Gamma$ of length $m \geq 2^{l(X)} - 1$ as each module in $\Sigma$ is stable. Since $l(Y_j) \leq l(X)$ for all $j = 0, \ldots, m$ as a consequence of statements 1 and 2, we conclude that $g_m \cdots g_1 = 0$ by \ref{HaradaSai}. This is a contradiction as, on the other hand, we know $g_m \cdots g_1 \neq 0$ by \ref{non-zero}. Hence the subgraph type of $\Gamma$ is $A_n$ for some $n < 2^{l(X)}$, which completes the proof.
\qed
\bigskip
\\We can now summarize the obtained results to the main theorem.

\begin{Th}\label{maintheorem}
Let $\Omega$ be a functorially finite resolving subcategory. Then the following statements are equivalent.
\begin{enumerate}[(a)]
\item $\Omega$ is finite.
\item The left subgraph type of each connected component of $\Gamma_l$ is given by a Dynkin diagram.
\item The right subgraph type of each connected component of $\Gamma_r$ is given by a Dynkin diagram.
\end{enumerate}
\end{Th}
\proof
\\We only prove the equivalence of the first and second statement since the equivalence of the first and third statement can be shown dually. 
\bigskip
\\Let $\Omega$ be finite and let $\Gamma$ be a connected component of $\Gamma_l$. In particular, $\Gamma$ is stable and does not contain multiple arrows between two fixed indecomposable modules. If there is a module in $\Gamma$ with an irreducible morphism to itself, then the left subgraph type of $\Gamma$ is $A_n$ by Theorem \ref{irreduciblemorphismtoitself}. On the other hand, if none such morphism exists, $\Gamma$ is a translation quiver. By Theorem \ref{treetypedynkin} the tree type of $\Gamma$ is given by a Dynkin diagram and coincides with the left subgraph type.
\bigskip
\\Suppose now that the left subgraph types of all connected components of $\Gamma_l$ are given by Dynkin diagrams and let $\Gamma$ denote such a component. We know that $\Gamma$ is not helical as helical components have left subgraph type $A_\infty$ by definition. Therefore, we can apply Theorem \ref{Dynkinorbitgraph} and obtain that $\Gamma$ is in fact stable and finite. 
\newpage
Suppose now that there are infinitely many stable, finite components $\Gamma_i$. In particular, none of these components contains an Ext-injective module, so every projective module is contained in a finite orbit since the number of non-isomorphic, indecomposable projective modules coincides with the number of non-isomorphic, indecomposable Ext-injective modules. Then, as for every module $X$ there is a projective module $P$ and a surjective morphism from $P$ to $X$, each $\Gamma_i$ must be adjacent to at least one finite orbit in the whole Auslander-Reiten quiver $\Gamma_\Omega$. But there are only finitely many projective modules in $\Gamma_\Omega$ and hence only finitely many finite orbits. Each of these orbits can only be adjacent to finitely many $\Gamma_i$, which is a contradiction. Hence $\Omega$ is finite.\qed
\bigskip
\\In order to show how this theorem generalizes Riedtmann's theorem, we deduce the following.

\begin{Th}\label{Riedtmanngeneralization}
Let $\Omega$ be a functorially finite resolving subcategory such that every projective module in $\Gamma_\Omega$ is contained in a finite orbit and there is no $\Omega$-irreducible morphism $f : X \to X$. Then $\Omega$ is finite if and only if $\Gamma_s$ does not contain multiple arrows and the tree type of every connected component of $\Gamma_s$ is given by a Dynkin diagram. 
\end{Th}
\proof
\\Let $\Gamma$ be a connected component of $\Gamma_s$. If $\Omega$ is finite, then $\Gamma$ does not contain multiple arrows by Lemma \ref{multiplearrows}. Moreover, by our assumptions $\Gamma$ does not contain an arrow $X \to X$, so $\Gamma$ is a stable translation quiver. Then by Theorem \ref{treetypedynkin} the tree type of $\Gamma$ is a Dynkin diagram.
\bigskip
\\On the other hand, suppose that $\Gamma_s$ does not contain multiple arrows and the tree type of each connected component $\Gamma$ of $\Gamma_s$ is a Dynkin diagram. $\Gamma$ is a translation quiver as it neither contains multiple arrows nor an $\Omega$-irreducible morphism $f: X \to X$. By Corollary \ref{treetypeandsubgraphtypecoincide} the subgraph type of $\Gamma$ equals the tree type of $\Gamma$ and, therefore, must be given by a Dynkin diagram. Since the numbers of projective modules and Ext-injective modules in $\Gamma_\Omega$ coincide and every projective module is contained in a finite orbit, every left stable orbit is in fact stable. It follows that $\Gamma_s = \Gamma_l$ and the left subgraph type of each connected component is also given by a Dynkin diagram. Thus $\Omega$ is finite by Theorem \ref{maintheorem}. \qed
\bigskip
\\The following example shows why we need to assume that every projective module is contained in a finite orbit. 
\begin{Ex}
Let $A$ be the path algebra of the quiver
$$\begin{xy}
\xymatrix{
 \ e_1 \ \ar@(ur,dr)^\beta \ar@(ul,dl)_\alpha  \\
}
\end{xy}$$
where all paths of length greater than $1$ are zero. 
\end{Ex}
It is not hard to see that the only indecomposable injective and projective modules have the following Jordan-H\"{o}lder decompositions.
 $$\begin{xy}
\xymatrixrowsep{0in}
\xymatrixcolsep{-0.1in}
  \xymatrix{ 
	&    &S_1&    					&S_1	&	&S_1		\\	
P_1=\ \ &S_1&   	&S_1 \ \ \ \ \ \ \ \ \ \ I_1= \ \	&	&S_1	
   }
\end{xy}$$
Computing the Auslander-Reiten translate of $S_1$ we obtain
 $$\begin{xy}
\xymatrixrowsep{0in}
\xymatrixcolsep{-0.1in}
  \xymatrix{ 
	  &S_1&   	&S_1	&	&S_1		\\	
\tau(S_1)=\ \ &	 &S_1&   	&S_1 	
   }
\end{xy}$$
Therefore, the connected component $\Gamma$ of $\Gamma_A$ containing $I_1, S_1$ and $P_1$ looks as follows.
$$\begin{xy}
\xymatrixrowsep{0.5in}
\xymatrixcolsep{0.5in}

  \xymatrix@!0{ 
\cdots&	& \tau(I_1) \ar@<1ex>[dr] \ar@<-1ex>[dr] \ar@{.>}[ll]&	& I_1] \ar@<1ex>[dr] \ar@<-1ex>[dr]  \ar@{.>}[ll]&	& [P_1 \ar@<1ex>[dr] \ar@<-1ex>[dr]	&	& \tau^{-1}(P_1) \ar@{.>}[ll] &	&\cdots \ar@{.>}[ll]\\
 &\cdots&	&  \tau(S_1) \ar@<1ex>[ur] \ar@<-1ex>[ur] \ar@{.>}[ll] &	& S_1  \ar@<1ex>[ur] \ar@<-1ex>[ur]  \ar@{.>}[ll]&	&\tau^{-1}(S_1) \ar@<1ex>[ur] \ar@<-1ex>[ur]	\ar@{.>}[ll] &	&\cdots \ar@{.>}[ll]}
\end{xy}$$
Note that the orbits of $\Gamma$ are not adjacent to any other orbits since $P_1$ and $I_1$ are the only indecomposable projective and injective modules respectively. Hence the middle term of all almost split sequences in $\Gamma$ consists of two copies of the same module. It follows that the length of modules in the orbit of $I$ is strictly increasing and given by the formula $l(\tau^n(I_1)) = 4(n + 1) - 1$ for all $n \geq 0$. We deduce that $I_1$ and $P_1$ are in different orbits and, in particular, the orbit of $I_1$ is left stable. Consequently, $\Gamma_l$ contains a connected component of the form
$$\begin{xy}
\xymatrixrowsep{0.5in}
\xymatrixcolsep{0.5in}
  \xymatrix@!0{ 
\cdots&	& \tau(I_1) \ar@<1ex>[dr] \ar@<-1ex>[dr] \ar@{.>}[ll]&	& I_1] \ar@<1ex>[dr] \ar@<-1ex>[dr]  \ar@{.>}[ll]\\
 &\cdots&	&  \tau(S_1) \ar@<1ex>[ur] \ar@<-1ex>[ur] \ar@{.>}[ll] &	& S_1  \ar@{.>}[ll]&	&\tau^{-1}(S_1) \ar@{.>}[ll] &	&\cdots \ar@{.>}[ll]}
\end{xy}$$
whose left subgraph type is $\widetilde{A}_1$. On the other hand, the stable part of this connected component consist only of the following orbit.
$$\begin{xy}
\xymatrixrowsep{0.5in}
\xymatrixcolsep{0.5in}
  \xymatrix@!0{ 
\cdots&	&  \tau(S_1)  \ar@{.>}[ll] &	& S_1  \ar@{.>}[ll]&	&\tau^{-1}(S_1) \ar@{.>}[ll] &	&\cdots \ar@{.>}[ll]}
\end{xy}$$
Hence its tree type is given by the Dynkin diagram $A_1$ although the component is infinite.
\\
\\In Riedtmann's original paper \cite{R80} she considered multiple arrows just as single arrows with non-trivial valuation. From that point of view every Auslander-Reiten quiver without loops is a translation quiver. The next example shows why we use a different approach to obtain equivalent statements in Theorem \ref{Riedtmanngeneralization}.
\begin{Ex}
Let $A$ be the path algebra of the quiver
$$\begin{xy}
\xymatrix{
 \ e_1 \ \ar@(ur,dr)^\beta \ar@(ul,dl)_\alpha  \\
}
\end{xy}$$
with relations $\alpha^2 = \beta^2 = \alpha\beta - \beta\alpha = 0$.
\end{Ex}
Note that this algebra is isomorphic to the group algebra of the Klein four-group over a field of characteristic $2$ \cite{B91}.
It is not hard to see that $A$ is selfinjective, i.e. every projective module is also injective, and the only indecomposable projective and injective module has the following Jordan-H\"{o}lder decomposition.

 $$\begin{xy}
\xymatrixrowsep{0in}
\xymatrixcolsep{-0.1in}
  \xymatrix{ 
	&    &S_1&    						\\	
P_1=I_1=\ \ &S_1&   	&S_1	\\
	&	&S_1		
   }
\end{xy}$$
For convenience, let $X$ denote the module $\rad(P_1)$. It follows from \cite[V. Proposition 5.5]{ARS95} that
$$\begin{xy}
\xymatrix{
0 \ar[r]	& X \ar[r]& P_1 \oplus S_1 \oplus S_1 \ar[r]	& \tau^{-1}(X) \ar[r]	&0 
}
\end{xy}$$
is an almost split sequence. Moreover, we compute
 $$\begin{xy}
\xymatrixrowsep{0in}
\xymatrixcolsep{-0.1in}
  \xymatrix{ 
	  &S_1&   	&S_1	&	&S_1		\\	
\tau(S_1)=\ \ &	 &S_1&   	&S_1 	
   }
\end{xy}$$
from which we deduce that the Auslander-Reiten quiver of the connected component containing $P_1$ looks as follows.
$$\begin{xy}
\xymatrixrowsep{0.5in}
\xymatrixcolsep{0.5in}

  \xymatrix@!0{ 
    &	&	&	&	&[P_1] \ar[dr]	&	&	\\
\cdots&	& \tau(X) \ar@<1ex>[dr] \ar@<-1ex>[dr] \ar@{.>}[ll]&	& X \ar[ur] \ar@<1ex>[dr] \ar@<-1ex>[dr]  \ar@{.>}[ll]&	& \tau^{-1}(X) \ar@<1ex>[dr] \ar@<-1ex>[dr] \ar@{.>}[ll]	&	& \tau^{-2}(X) \ar@{.>}[ll] &	&\cdots \ar@{.>}[ll]\\
 &\cdots&	&  \tau(S_1) \ar@<1ex>[ur] \ar@<-1ex>[ur] \ar@{.>}[ll] &	& S_1  \ar@<1ex>[ur] \ar@<-1ex>[ur]  \ar@{.>}[ll]&	&\tau^{-1}(S_1) \ar@<1ex>[ur] \ar@<-1ex>[ur]	\ar@{.>}[ll] &	&\cdots \ar@{.>}[ll]}
\end{xy}$$
Let us denote the subgraph consisting of all stable modules in this component, i.e. every module except for $P_1$, by $\Gamma$. Then $\Gamma$ is a connected component of both $\Gamma_s$ and $\Gamma_l$ and its left subgraph type is $\widetilde{A}_1$, which proves that the component is infinite by Theorem \ref{maintheorem}. If we consider the double arrows as single arrows with valuation, then the tree type of $\Gamma$ is $A_2$. So $\Gamma$ is infinite although its tree type is given by in Dynkin diagram in Riedtmann's point of view. This is why we consider connected components containing multiple arrows separately in Theorem \ref{Riedtmanngeneralization}.

\begin{Cor}\cite[Hauptsatz]{R80}
Let $A$ be an algebra of finite representation type. Then the tree type of each connected component of $\Gamma_s$ is given by a Dynkin diagram.
\end{Cor}
\proof
\\The Corollary is an immediate consequence of Theorem \ref{Riedtmanngeneralization} and Corollary \ref{multiplearrowscontrapositive}. \qed

\newpage

\section{Coray tubes in functorially finite resolving subcategories}

In this section we show that helical components in $\Gamma_\Omega$ are the same as coray tubes. Some of the statements are generalizations of results in \cite{L93} and their proofs closely follow the original proofs in the aforementioned paper.

\begin{Def}
Let $\Gamma$ be a translation quiver. A vertex $x$ of $\Gamma$ is called a \textbf{coray vertex} if there is an infinite sectional path 
$$\begin{xy}
  \xymatrix{ 
\cdots \ar[r]	& x_n \ar[r] &    x_{n-1} \ar[r]	&	\cdots \ar[r] 	&	x_2 \ar[r]     &   x_1 = x	  }
\end{xy}$$
in $\Gamma$ such that for each integer $n > 0$ the path 
$$\begin{xy}
  \xymatrix{ 
x_n \ar[r] &    x_{n-1} \ar[r]	&	\cdots \ar[r] 	&	x_2 \ar[r]     &   x_1 = x	  }
\end{xy}$$
is the only sectional path of length $n$ in $\Gamma$ which ends at $x$. The aforementioned infinite sectional path is called a \textbf{coray} ending in $x$.
\end{Def}

Let $x$ be a coray vertex in a translation quiver $\Gamma$ with a coray as above. For a positive integer $n$ we define the translation quiver $\Gamma[x,n]$ as follows. The vertices of $\Gamma[x,n]$ are those of $\Gamma$ together with pairs $(i,j)$ where $i \geq 1, 1 \leq j \leq n$. The arrows of $\Gamma[x,n]$ are those of $\Gamma$, excluding those ending at $x_i$ other than $x_{i+1} \to x_i$ for $i \geq 1$, together with the following arrows:
\begin{enumerate}
\item $(i + 1, j) \to (i,j)$ for $i \geq 1, 1 \leq j \leq n$,
\item $(i, j+ 1) \to (i + 1, j)$ for $i \geq 1, 1 \leq j \leq n$,
\item $(n + i - 1,1) \to x_i$ for all $i \geq 1$,
\item $y \to (i,n)$ if $y \to x_i$ is an arrow in $\Gamma$ other than $x_{i+1} \to x_i$.
\end{enumerate}

Let $\tau$ be the translation of $\Gamma$. The translation $\tau'$ of $\Gamma[x,n]$ is defined as follows, if $z$ is a vertex of $\Gamma$ and $z \neq x_i$ for all $i \geq 1$ and $\tau(z)$ is defined, then $\tau'(z) = \tau(z)$. Moreover, we set $\tau'(x_i) = (n + i, 1)$ for all $i \geq 1$, $\tau'(i,j) = (i, j + 1)$ for $1 \leq j < n$, $i \geq 1$ and $\tau'(i,n) = \tau (x_i)$ if $\tau(x_i)$ is defined.
\bigskip
\\We inductively  define a translation quiver $\Gamma[x_0,n_0][x_1,n_1] \ldots [x_r,n_r]$, where $n_i$ are positive integers, $x_0$ is a coray vertex of $\Gamma$ and $x_i$ is a coray vertex of $\Gamma[x_0,n_0][x_1,n_1] \ldots [x_{i-1},n_{i-1}]$. Where there is no need to emphasize the coray vertices we say that a translation quiver $\Gamma$ is of the form $\Gamma[n_0,n_1, \ldots ,n_r]$ if $ \Gamma' \cong \Gamma[x_0,n_0][x_1,n_1] \ldots [x_r,n_r]$ for some coray vertices $x_i$. We say $\Gamma'$ is a translation quiver obtained from $\Gamma$ by \textbf{coray insertions}. Dually, one can define the concepts of ray vertices and ray insertions.

\begin{Def}
\mbox{}
\begin{enumerate}[(a)]
\item We call a translation quiver a \textbf{coray tube} if it is obtained from a stable tube by coray insertions.
\item Dually, we call a translation quiver a \textbf{ray tube} if it is obtained from a stable tube by ray insertions.
\end{enumerate}
\end{Def}

\begin{Ex}
Let $A$ be the path algebra of the quiver
$$\begin{xy}
\xymatrix{
  e_1 \ \ \ar@<1ex>[r]^{\alpha} \ar@<-1ex>[r]_{\beta}  & \ \ e_2 \ \  \ar[r]^{\gamma}  & \ \  e_3 \\
}
\end{xy},$$
with the relation $\gamma\beta = 0$. Then $\Gamma_A$ contains a connected component which is a coray tube.
\end{Ex}

We have seen in Example \ref{helical} that the connected component $\Gamma'$ of the injective module $I_3$ looks as follows.
$$\begin{xy}
\xymatrixrowsep{0.5in}
\xymatrixcolsep{0.5in}

  \xymatrix@!0{ 
& & & & & & & & & & &  \\
& &\cdots& &\tau^3(I_3)\ar@{.>}[ll] \ar[ur] \ar[dr]& &\tau^2(I_3) \ar@{.>}[ll] \ar[dr] \ar[ur]&	&\tau(I_3) \ar@{.>}[ll] \ar[ur]	&	& I_3] \ar[ur] \ar@{.>}[ll]&\\
&\cdots& 	&\tau^2(I_3)\ar@{.>}[ll] \ar[ur] \ar[dr]	& 	&\tau(I_3) \ar@{.>}[ll] \ar[ur]	&	&I_3] \ar[ur] \ar@{.>}[ll]&	&	&		&\\
\cdots&	&\tau(I_3) \ar@{.>}[ll]\ar[ur]	&	&I_3] \ar[ur] \ar@{.>}[ll]&	&	&	&	&	&	&
 }
\end{xy}$$

$\Gamma'$ is obtained by coray insertions the following way. Given a stable tube $\Gamma$ of type $\Z A_\infty/\tau$

$$\begin{xy}
\xymatrixrowsep{0.5in}
\xymatrixcolsep{0.5in}

  \xymatrix@!0{ 
&\ar[dr]&	& \ar[dr]	\\
&	&x_4 \ar[dr] \ar[ur] \ar@{.>}[ll]	&	&x_4 \ar[dr] \ar@{.>}[ll]&	& \ar@{.>}[ll] 	\\
&	&	&x_3 \ar[dr] \ar[ur] \ar@{.>}[ll]	& 	&x_3 \ar[dr] \ar@{.>}[ll]&	& \ar@{.>}[ll] 	\\ 
&	&	&	&x_2 \ar[dr] \ar[ur] \ar@{.>}[ll]&	&x_2 \ar[dr] \ar@{.>}[ll]&	& \ar@{.>}[ll] 	\\
& 	& 	&	&	& x_1 \ar[ur] \ar@{.>}[ll]  &	& x_1 \ar@{.>}[ll] 	&	& \ar@{.>}[ll] }
\end{xy}$$
such that $x_1$ is a coray vertex, then $\Gamma[x_1,1]$ 
$$\begin{xy}
\xymatrixrowsep{0.5in}
\xymatrixcolsep{0.5in}

  \xymatrix@!0{ 
& \ar[dr]	&	&\ar[dr]&	& \ar[dr]	\\
&	&x_5 \ar[dr] \ar[ur] \ar@{.>}[ll]	&	&(5,1) \ar[dr] \ar[ur] \ar@{.>}[ll]	&	&x_4 \ar[dr] \ar@{.>}[ll]&	& \ar@{.>}[ll] 	\\
& &	&x_4 \ar[dr] \ar[ur] \ar@{.>}[ll]	&	&(4,1) \ar[dr] \ar[ur] \ar@{.>}[ll]	& 	&x_3 \ar[dr] \ar@{.>}[ll]&	& \ar@{.>}[ll] 	\\ 
& &	&	&x_3 \ar[dr] \ar[ur] \ar@{.>}[ll]	&	&(3,1) \ar[dr] \ar[ur] \ar@{.>}[ll]&	&x_2 \ar[dr] \ar@{.>}[ll]&	& \ar@{.>}[ll] 	\\
& &	& & 	&x_2 \ar[dr] \ar[ur] \ar@{.>}[ll]	&	& (2,1) \ar[ur] \ar[dr] \ar@{.>}[ll]  &	& x_1 \ar@{.>}[ll] 	&	& \ar@{.>}[ll]\\ 
  & &	&	&	&	&x_1 \ar[ur] \ar@{.>}[ll]&	&(1,1)] \ar[ur] \ar@{.>}[ll]&	&  }
\end{xy}$$
is a coray tube which is obviously isomorphic to $\Gamma'$.

\begin{Lemma}\label{shapeofinjectiveend}
Let $\Gamma$ be a connected component of $\Gamma_\Omega$ with trivial valuation.
\begin{enumerate}[(a)]
 \item Assume that $\Gamma$ does not contain projective modules and each module in $\Gamma$ has at most two immediate predecessors. Suppose 
$$\begin{xy}
  \xymatrix{ 
\cdots \ar[r] &		X_{i+1} \ar[r] &		X_i \ar[r]	&	\cdots \ar[r] 	&	X_2 \ar[r]     &   X_1 = X	  }
\end{xy}$$ 
is an infinite sectional path in $\Gamma$ containing an infinite number of arrows with finite global left degree. If $n$ is a positive integer and $\tau_\Omega^{-n}(X_i)$ is defined for some $i \geq 1$, then $\tau_\Omega^{-n}(X_j)$ is also defined for all $j > i$. In particular, if some module $X_i$ is stable, then all modules $X_j$ with $j > i$ are also stable.
 \item Assume that $\Gamma$ does not contain Ext-injective modules and each module in $\Gamma$ has at most two immediate successors. Suppose 
$$\begin{xy}
  \xymatrix{ 
X = X_1 \ar[r] &		X_2 \ar[r]	&	\cdots \ar[r] 	&	X_i \ar[r]     &	X_{i+1} \ar[r] &   \cdots	  }
\end{xy}$$ 
is an infinite sectional path in $\Gamma$ containing an infinite number of arrows with finite global right degree. If $n$ is a positive integer and $\tau_\Omega^n(X_i)$ is defined for some $i \geq 1$, then $\tau_\Omega^n(X_j)$ is also defined for all $j > i$. In particular, if some module $X_i$ is stable, then all modules $X_j$ with $j > i$ are also stable.
\end{enumerate}

\end{Lemma}

\proof
\\Suppose the statement is false, i.e. there are positive integers $i$ and $n$ for which $\tau_\Omega^{-n}(X_i)$ exists, but $\tau_\Omega^{-n}(X_{i+1})$ does not. Without loss of generality, we can assume that $n$ is the least integer for which such an $i$ exists. Hence $\Gamma$ contains a subgraph which looks as follows.
$$\begin{xy}
\xymatrixrowsep{0.5in}
\xymatrixcolsep{0.5in}

  \xymatrix@!0{  
&	&	&	& I] \ar[dr]	&	&  \\
&	&	&	&	& \tau_\Omega^{-n+1}(X_i)\ar[dr]^{f_0}	&			& \tau_\Omega^{-n}(X_i) \ar@{.>}[ll] \\
&	&	&	& \tau_\Omega(Y_0)\ar[ur] \ar[dr]^{f_1}	&	& Y_0 \ar[ur]\ar@{.>}[ll] \\
&	&	& \tau_\Omega(Y_1)\ar[ur]\ar[dr]^{f_2}	& 	& Y_1 \ar[ur]\ar@{.>}[ll] 	\\ 
&	& \tau_\Omega(Y_2)\ar[ur]	&	&  Y_2 \ar[ur]\ar@{.>}[ll]\\
&	\ar[ur]&	& \ar[ur] \\ 
\ar@{.}[ur] &	& \ar@{.}[ur]	 }
\end{xy}$$
We prove by induction that there is a sectional path 
$$\begin{xy}
  \xymatrix{ 
Y_k \ar[r]	&	\cdots \ar[r] 	&	Y_1 \ar[r]     &   Y_0	 \ar[r]	&  \tau_\Omega^{-n}(X_i) }
\end{xy}$$ 
for each length $k \geq 0$ and every irreducible morphism $f_i : \tau_\Omega(Y_{i-1}) \to Y_i$ is a monomorphism. Let 
$$\begin{xy}
  \xymatrix{ 
      0 \ar[r]  	& 	\tau_\Omega^{-n+1}(X_i) \ar[r]^(.6){f_0} 	&	Y_0 \ar[r]     &   \tau_\Omega^{-n}(X_i) \ar[r]	&	0	    }
\end{xy}$$
be the almost split sequence ending in $\tau_\Omega^{-n}(X_i)$, then $f_0$ is a monomorphism. Since $\Gamma$ does not contain projective modules and every module in $\Gamma$ has at most $2$ predecessors, $Y_0$ must be indecomposable, which proves the statement for $k = 0$. Suppose now the statement is true for an arbitrary $k$. Then there is a monomorphism $f_k : \tau_\Omega(Y_{k-1}) \to Y_k$, from which we conclude that the almost split sequence $\E_\Omega(Y_k)$ must have two indecomposable middle terms, $\tau_\Omega(Y_{k-1})$ and $Y_{k+1}$. Then by Lemma \ref{exact} we can construct a short exact sequence
$$\begin{xy}
  \xymatrix{ 
      0 \ar[r]  	& 	\tau_\Omega(Y_k) \ar[r]^{f_{k+1}} 	&	Y_{k+1} \ar[r]     &   \tau_\Omega^{-n}(X_i) \ar[r]	&	0,	    }
\end{xy}$$
which gives us that the irreducible morphism $f_{k+1} : \tau_\Omega(Y_k) \to Y_{k+1}$ is again a monomorphism. As a consequence, there is an infinite sectional path
$$\begin{xy}
  \xymatrix{ 
\cdots \ar[r]	&	Y_k \ar[r]	&	\cdots \ar[r] 	&	Y_1 \ar[r]     &   Y_0 }
\end{xy}$$ 
ending in $Y_0$. However, as we have chosen $n$ to be minimal and $Y_1 \ncong \tau_\Omega^{-n+1}(X_i)$, there is another infinite sectional path
$$\begin{xy}
  \xymatrix{ 
\cdots \ar[r]	&\tau_\Omega^{-n+1}(X_{i+l}) \ar[r]	&\cdots \ar[r] &\tau_\Omega^{-n+1}(X_{i+1}) \ar[r]	&\tau_\Omega^{-n+1}(X_i) \ar[r]     &   Y_0. }
\end{xy}$$ 
Since $\Gamma$ does not contain projective modules, we obtain two infinite sectional paths 
$$\begin{xy}
  \xymatrix{ 
\cdots \ar[r]	&	X_k \ar[r]	&	\cdots \ar[r] 	&	X_{i-1} \ar[r]     &   X_i }
\end{xy}$$ 
and
$$\begin{xy}
  \xymatrix{ 
\cdots \ar[r]	&	\tau_\Omega^n(Y_k) \ar[r]	&	\cdots \ar[r] 	&	\tau_\Omega^n(Y_1) \ar[r]     &   \tau_\Omega^n(Y_0) \ar[r] & X_i }
\end{xy}$$ 
in which by Lemma \ref{twosectionalpathsinfiniteglobaldegree} each arrow has infinite global left degree since $\tau_\Omega^n(Y_0) \ncong X_{i-1}$. It follows that the infinite sectional path
$$\begin{xy}
  \xymatrix{ 
\cdots \ar[r] &		X_i \ar[r]	&	\cdots \ar[r] 	&	X_2 \ar[r]     &   X_1 = X	  }
\end{xy}$$ 
contains only finitely many arrows of finite global left degree, which contradicts our assumptions. \qed

\begin{Th}\label{coray}
Let $\Gamma$ be a connected component of $\Gamma_\Omega$ containing an oriented cycle but no $\tau_\Omega$-periodic module.
\begin{enumerate}[(a)]
 \item If $\Gamma$ does not contain a projective module, then $\Gamma$ is a coray tube.
 \item Dually, if $\Gamma$ does not contain an Ext-injective module, then $\Gamma$ is a ray tube.
\end{enumerate}
\end{Th}

\proof
\\By Theorem \ref{leftstablecomponentstructure} we know that there is an infinite path of type $A_\infty$ in $\Gamma$ which is a full sectional subgraph. Moreover, the theorem says that $\Gamma$ contains only finitely many different orbits and at least one of them is not right stable. 
\bigskip
\\Suppose now that $\Gamma$ contains a stable orbit, then its connected component $\Gamma'$ in the right stable Auslander-Reiten quiver $\Gamma_r$ is of right subgraph type $A_n$ for some $n < s$, where $s$ denotes the number of different $\tau_\Omega$-orbits in $\Gamma$. Then by Theorem \ref{Dynkinorbitgraph} the component $\Gamma'$ is finite and, in particular, $\tau_\Omega$-periodic, which contradicts our assumptions. Hence every orbit contains an Ext-injective module. Let
$$\begin{xy}
  \xymatrix{ 
     \cdots \ar[r] &	\tau_\Omega^{2r}(X_1) \ar[r] &	\tau_\Omega^r(X_s) \ar[r] 	&  \cdots \ar[r] &  \tau_\Omega^r(X_1) \ar[r] & X_s \ar[r] &  \cdots	\ar[r]	&	X_1	    }
\end{xy}$$
be an infinite sectional path in $\Gamma$, which exists by Lemma \ref{infinitesectionalpathinleftstablecomponent}. Then there exist positive integers $n_1, \ldots, n_s$ such that $X_j = \tau_\Omega^{n_j}(I_j)$ for some Ext-injective modules $I_1, \ldots, I_s$. It follows from Lemma \ref{shapeofinjectiveend} that $n_1 \leq n_2 \leq \cdots \leq n_s$. Let $m_1, m_2, \ldots, m_t$ be positive integers such that
$$ n_1 = \cdots = n_{m_1} < n_{m_1+1} = \cdots = n_{m_2} < \cdots < n_{m_{t-1}+1} = \cdots = n_{m_t} = n_s.$$
It then follows that $\Gamma$ as a translation quiver is of the form
$$(\Z A_\infty/\tau^n)[m_1, m_2 - m_1, \ldots , m_t - m_{t-1}],$$
where $\tau$ denotes the translation in $\Z A_\infty$ and $n = r - s$. \qed

\begin{Lemma}\label{coraypath}
Let $n_1, n_2, \ldots , n_t, n_{t+1}$ be positive integers for some $t \geq 0$.
\begin{enumerate}[(a)]
 \item If the coray tube $(\Z A_\infty/\tau^n)[n_1, n_2, \ldots , n_t]$ contains a path from some vertex $x$ to another vertex $y$, then
$(\Z A_\infty/\tau^n)[n_1, n_2, \ldots , n_t, n_{t+1}]$ also contains a path from $x$ to $y$. 
\item Dually, if the ray tube $(\Z A_\infty/\tau^n)[n_1, n_2, \ldots , n_t]$ contains a path from some vertex $x$ to another vertex $y$, then
$(\Z A_\infty/\tau^n)[n_1, n_2, \ldots , n_t, n_{t+1}]$
also contains a path from $x$ to $y$.
 \end{enumerate}
\end{Lemma}
\proof
\\Let $z_1$ be the coray vertex where we insert $n_{t+1}$ and let
$$\begin{xy}
  \xymatrix{ 
     \cdots \ar[r] &	z_{i+1} \ar[r] &	z_i \ar[r] 	&  \cdots \ar[r] &  z_2 \ar[r] & z_1	    }
\end{xy}$$
be the corresponding coray, while 
$$\begin{xy}
  \xymatrix{ 
   \gamma:	& x = x_0 \ar[r] &	x_1 \ar[r] &	\cdots \ar[r] 	&  x_{m-1} \ar[r] &  x_m = y	    }
\end{xy}$$
denotes the path from $x$ to $y$. If $x_{j-1} \to x_j$ is an arrow in $\gamma$ such that either $x_j \neq z_i$ for all $i \geq 1$ or $x_{j-1} \to x_j = z_{i+1} \to z_i$ for some $i \geq 1$, then by definition the same arrows exist in $(\Z A_\infty/\tau^n)[n_1, n_2, \ldots , n_t, n_{t+1}]$. On the other hand, if $x_{j-1} \to x_j$ is an arrow in $\gamma$ such that $x_j = z_i$ but $x_{j-1} \neq z_{i+1}$, then there is path
$$\begin{xy}
  \xymatrix{ 
x_{j-1} \ar[r] &(i,n_{t+1}) \ar[r] &	(i+1 ,n_{t+1} - 1) \ar[r] 	&  \cdots      }
\end{xy}$$
$$\begin{xy}
  \xymatrix{ 
  \cdots \ar[r]  &  (n_{t+1} + i - 2 ,2) \ar[r] & (n_{t+1} + i - 1 ,1) \ar[r] &x_j   }
\end{xy}$$
in $(\Z A_\infty/\tau^n)[n_1, n_2, \ldots , n_t, n_{t+1}]$, which completes the proof.  \qed

\begin{Th}\label{helicalcoray}
Let $\Gamma$ be a connected component of $\Gamma_\Omega$.
\begin{enumerate}[(a)]
\item Suppose $\Gamma$ is left stable. Then $\Gamma$ is helical if and only if it is a coray tube.
\item Dually, if $\Gamma$ is right stable, then $\Gamma$ is cohelical if and only if it is a ray tube.
\end{enumerate}
\end{Th}

\proof
\\Let $\Gamma$ be helical, then by Theorem \ref{helical} it contains an oriented cycle. Moreover, since helical components always contain an Ext-injective module, there cannot be a $\tau_\Omega$-periodic module in $\Gamma$ by Lemma \ref{orbit}. We deduce that $\Gamma$ is a coray tube by Theorem \ref{coray}.
\bigskip
\\Suppose now that $\Gamma$ is a coray tube of the form 
$$(\Z A_\infty/\tau^n)[n_1, n_2, \ldots , n_t].$$
We prove by induction on $t$ that there is a vertex $i$ in $(\Z A_\infty/\tau^n)[n_1, n_2, \ldots , n_t]$ such that $\tau^{-1}(i)$ does not exist and there is a path from $i$ to any vertex $y$ in $(\Z A_\infty/\tau^n)[n_1, n_2, \ldots , n_t]$. If $t = 1$, then we consider a coray tube of the form $(\Z A_\infty/\tau^n)[n_1]$ and let $z_1$ denote the coray vertex of the insertion $n_1$. Since $\Z A_\infty/\tau^n$ is a stable tube, there is a path from $z_1$ to $y$ by Lemma \ref{pathexistence} for every vertex $y$ in $\Z A_\infty/\tau^n$. Additionally, there is also a path from $z_1$ to $y$ in $(\Z A_\infty/\tau^n)[n_1]$ by Lemma \ref{coraypath}. On the other hand, $\tau^{-1}(i)$ does not exist for the vertex $i = (n_1, 1)$ in $(\Z A_\infty/\tau^n)[n_1]$ by construction and there is an arrow $i \to z_1$ in  $(\Z A_\infty/\tau^n)[n_1]$, which gives us a path from $i$ to $y$. 
\bigskip
\\Let $y$ now be a vertex in $(\Z A_\infty/\tau^n)[n_1]$ that has been added in the coray insertion. Hence $y$ is a vertex of the form $(k,l)$ with $k \in \N$ and $1 \leq l \leq n_1$. If 
$$\begin{xy}
  \xymatrix{ 
     \cdots \ar[r] &	z_{i+1} \ar[r] &	z_i \ar[r] 	&  \cdots \ar[r] &  z_2 \ar[r] & z_1    }
\end{xy}$$
is the sectional path for the coray vertex $z_1$, then there is a path
$$\begin{xy}
  \xymatrix{ 
     \tau(z_k) \ar[r] &	(k + 1, n_1) \ar[r] &	(k + 2 , n_1 - 1) \ar[r] 	&  \cdots     }
\end{xy}$$
$$\begin{xy}
  \xymatrix{ 
 \cdots \ar[r] &  (k + n_1 - l , l + 1) \ar[r] & (k + n_1 - l + 1 , l ) \ar[r] & (k + n_1 - l, l ) \ar[r]	& 	\cdots	    }
\end{xy}$$
$$\begin{xy}
  \xymatrix{ 
 \cdots \ar[r] &  (k + 1 , l) \ar[r] & (k , l) = y. }
\end{xy}$$
As we have seen earlier in this proof, there is a path from $i$ to $\tau(z_k)$ that we now compose to a path from $i$ to $y$. This completes the proof for $t = 1$. 
\bigskip
\\Suppose the statement is true for some $t \geq 1$, then for every vertex $y$ in the coray tube $(\Z A_\infty/\tau^n)[n_1, n_2, \ldots , n_t]$ there is path from $i$ to $y$, where $i$ is a vertex such that $\tau^{-1}(i)$ does not exist. Consequently, there exists a path from $i$ to $y$ in $(\Z A_\infty/\tau^n)[n_1, n_2, \ldots , n_t, n_{t+1}]$ by Lemma \ref{coraypath}. If $y = (k,l)$ is a vertex obtained by the coray insertion of $n_{t+1}$ at the coray  
$$\begin{xy}
  \xymatrix{ 
     \cdots \ar[r] &	z_{i+1} \ar[r] &	z_i \ar[r] 	&  \cdots \ar[r] &  z_2 \ar[r] & z_1,	    }
\end{xy}$$
then there is a path
$$\begin{xy}
  \xymatrix{ 
  \gamma: &   \tau(z_k) \ar[r] &	(k+1, n_{t+1}) \ar[r] &	(k + 1 , n_{t+1} - 1) \ar[r] 	&  \cdots     }
\end{xy}$$
$$\begin{xy}
  \xymatrix{ 
 \cdots \ar[r] &  (k + n_{t+1} - l , l - 1) \ar[r] & (k + n_{t+1} - l + 1 , l ) \ar[r] & (k + n_{t+1} - l, l ) \ar[r]	& 	\cdots	    }
\end{xy}$$
$$\begin{xy}
  \xymatrix{ 
 \cdots \ar[r] &  (k + 1 , l) \ar[r] & (k , l) = y. }
\end{xy}$$
By the induction hypothesis we conclude that there exists a path from $i$ to $\tau(z_k)$ in $(\Z A_\infty/\tau^n)[n_1, n_2, \ldots , n_t]$ which by Lemma \ref{coraypath} gives rise to a path from $i$ to $\tau(z_k)$ in $(\Z A_\infty/\tau^n)[n_1, n_2, \ldots , n_t, n_{t+1}]$. Hence we obtain a path from $i$ to $y$ by composing this path with $\gamma$, which proves that the statement is true for all $t \in \N$. Since $\Gamma$ is a coray tube of the form $(\Z A_\infty/\tau^n)[n_1, n_2, \ldots , n_t]$, we deduce that for every module $Y$ there is an Ext-injective module $I$ such that there is a path from $I$ to $Y$, which makes $\Gamma$ a helical component. \qed

\section{Degrees of irreducible morphisms in functorially finite resolving subcategories}

In this section we analyze Auslander-Reiten quivers of functorially finite resolving subcategories using degrees of irreducible morphisms. In particular, we observe under which assumptions we have $l(\tau_\Omega^n(X)) \to \infty$ as $n \to \infty$ for a module $X$ in an infinite connected component. We give a complete description of all infinite connected components such that $l(\tau_\Omega^n(X))$ does not tend to infinity as $n \to \infty$. This section again closely follows \cite{L92}.

\begin{Lemma}\label{degreeequals1}
Let $f : X \to Y$ be an irreducible morphism in a functorially finite resolving subcategory $\Omega$. Then $f$ is a surjective minimal right almost split morphism if and only if $d_\Omega^l(f) = 1$ and all direct summands of $Y$ are not projective.
\end{Lemma}

\proof
\\We first assume that $f$ is a surjective minimal right almost split morphism. Consequently, no direct summand of $Y$ can be projective as $f$ is surjective. Moreover, $Y$ is indecomposable as $f$ is an irreducible minimal right almost split morphism. So we compose $f$ with the corresponding minimal left almost split morphism for $\tau_\Omega(Y)$ and obtain $d_\Omega^l(f) = 1$.
\bigskip
\\On the other hand, let $d_\Omega^l(f) = 1$ and suppose no direct summand of $Y$ is projective. If $Y$ is indecomposable, then $f$ is a surjective minimal right almost split morphism by Corollary \ref{degree1Cor}. In case that $Y$ admits a direct summand $Y_1 \oplus Y_2$ we know that $X$ is indecomposable as $f$ is irreducible. Let $(f_1,f_2)^T : X \to Y_1 \oplus Y_2$ be the corestriction of $f$ to $Y_1 \oplus Y_2$, then $d_\Omega^l((f_1,f_2)^T) = 1$. By Lemma \ref{neighbour} a left neighbour $h$ of $(f_1,f_2)^T$ exists with $d_\Omega^l(h) < 1$, which is impossible. \qed

\begin{Def}
Let $Q$ be a quiver without loops or cycles and let $KQ$ be the hereditary path algebra of $Q$. Suppose $Q$ has $n$ vertices and let $p_1, \ldots ,  p_n$ and $i_1, \ldots, i_n$ denote the dimension vectors of the indecomposable projective and injective $KQ$-modules respectively. Then the \textbf{Coxeter-matrix of $\pmb{Q}$} is the matrix $C$ defined via $C(i_j) = - p_j$.
\end{Def}

Note that $C$ is invertible and that its inverse satisfies $C^{-1}p_j = -i_j$. Since the inverse of the Coxeter-matrix is more convenient for us, we set $C^{-1} = (c_{ij})$. Before we provide a combinatorial description of the entries in $C^{-1}$ in terms of paths in $Q$, we need to establish some more notation. Let $Q$ be a quiver with vertices $e_1, \ldots, e_n$. We say that $e_j$ is a predecessor of $e_i$ if there is a non-trivial path from $e_j$ to $e_i$ in $Q$. If this path consists of only one arrow, then $e_i$ is an immediate predecessor of $e_j$. Successors and immediate successors are defined dually. We denote an arrow from $e_j$ to $e_i$ in $Q$ by $e_j \to e_i$, a path of arbitrary length is denoted by $\gamma: e_j \to e_i$. For example, for a fixed vertex $e_i$ $$\sum_{e_i \to e_j} \sum_{\gamma : e_k \to e_j}$$ is the sum over all paths in $Q$ that end in an immediate successor of $e_i$. Moreover, for any column vector $v$ we denote its $i$-th entry by $(v)_i$. This and other upcoming terminology 
is fixed for the rest of the section.

\begin{Lemma}\label{Coxetercombinatorics}
Let $Q$ be a quiver without loops or cycles. If $Q$ contains vertices $e_1, \ldots, e_n$ and $C^{-1} = (c_{ij})$ denotes the inverse of its Coxeter-matrix, then the entry $c_{ij}$ is given by
$$c_{ij} =  - \# \text{ of paths from $e_i$ to $e_j$} + (\sum_{e_j \to e_k} \# \text{ of paths from $e_i$ to $e_k$})$$
In particular, we have the following.
\begin{enumerate}[(a)]
 \item $c_{ij} = 0$ if there is neither a path from $e_i$ to $e_j$ nor to an immediate predecessor of $e_j$.
 \item $c_{ij} = 0$ if $e_j$ is a non-immediate predecessor of $e_i$.
 \item $c_{ij} = \# \text{ of arrows from $e_j$ to $e_i$}$ if $e_j$ is an immediate predecessor of $e_i$.
  \end{enumerate}
 If $Q$ does not contain multiple arrows between two vertices then we also have the following.
 \begin{enumerate}[(a)]
  \setcounter{enumi}{3}
 \item $c_{ij} = 1$ if $e_j$ is an immediate predecessor of $e_i$.
 \item $c_{ij} = \# \text{ of arrows starting in $e_j$} - 1$ if $e_j = e_i$ or $e_j$ is a successor of $e_i$. 
 \item In particular, if $e_j$ is a successor of $e_i$, we have $c_{ij} = c_{jj}$.
 \end{enumerate}

Dually, if we set $C = (d_{ij})$, then
$$d_{ij} =  - \# \text{ of paths from $e_j$ to $e_i$} + (\sum_{e_k \to e_j} \# \text{ of paths from $e_k$ to $e_i$}).$$

\end{Lemma}

\proof
\\Let $P_1, \ldots, P_n$ and $I_1, \ldots, I_n$ denote the indecomposable projective and injective $KQ$-modules respectively. Alongside with $p_j$ and $i_j$ for dimension vectors of indecomposable projective and injective $KQ$-modules, we denote the dimension vector of the simple module $S_j = P_j/\rad(P_j)$ by $s_j$. Clearly, $c_{ij}$ is given by the $i$-th entry of the column vector of $C^{-1}s_j$. Note that since $KQ$ is an hereditary algebra, $P_j$ is given by all paths that start in $e_j$, while $I_j$ is the module consisting of all paths ending in $e_j$. 
\bigskip
\\Hence the dimension vector of a simple module can be expressed as $s_j = p_j - \rad(p_j) = p_j - \sum_{e_j \to e_k} p_k$. By definition we have $C^{-1}p_k = - i_k$, so the $i$-th entry of this column vector is $-1$ times the number of paths from $e_i$ to $e_k$. Hence 
$$c_{ij} = (C^{-1}s_j)_i = (C^{-1}p_j)_i - \sum_{e_j \to e_k}(C^{-1}p_k)_i = \sum_{e_j \to e_k}(i_k)_i - (i_j)_i$$
$$=  - \# \text{ of paths from $e_i$ to $e_j$} + (\sum_{e_j \to e_k} \# \text{ of paths from $e_i$ to $e_k$}).$$ 
The statements $(a), (b)$ and $(c)$ are special cases of the formula above and follow immediately from the fact that $Q$ contains neither loops nor cycles. Moreover, if $Q$ in addition does not contain multiple arrows, then there is precisely one path from $e_i$ to $e_j$ if and only if $e_i = e_j$ or $e_j$ is a successor of $e_i$. In all other cases there is no path from $e_i$ to $e_j$, from which we deduce statements $(d),(e)$ and $(f)$. \qed
\bigskip
\\Let $\Gamma$ be a connected component of $\Gamma_\Omega$ containing a sectional subgraph $\Sigma$ that consists of modules $M_1, \ldots, M_n$ without multiple arrows between them. Since $\Sigma$ neither contains loops nor cycles, there is a Coxeter-matrix $C = (d_{ij})$ of $\Sigma$ and let $C^{-1} = (c_{ij})$ again denote its inverse. Furthermore, whenever $\Sigma$ occurs in this section, let $m$ denote the column vector $(l(M_1), \ldots , l(M_n))^T$. Then the $i$-th entry of $C^{-1}m$ is given by $(C^{-1}m)_i = \sum_{j = 1}^n c_{ij}l(M_j)$. We split up the sum on the right hand side according to Lemma \ref{Coxetercombinatorics} $(d),(e)$ and obtain
$$(C^{-1}m)_i = \sum_{M_k \to M_i} l(M_k) + \sum_{\gamma: M_i \to M_k}c_{ik} l(M_k)= \sum_{M_k \to M_i} l(M_k) + \sum_{\gamma: M_i \to M_k}c_{kk} l(M_k),$$
where the last equality follows from Lemma \ref{Coxetercombinatorics} $(f)$.
This formula is used frequently in the following proofs.

\begin{Lemma}\label{Coxetersumformula}
Let $M_i$ be a module in $\Sigma$, then we have
\begin{enumerate} [(a)]
 \item $$\sum_{M_i \to M_j}(C^{-1}m)_j =  l(M_i) + \sum_{M_i \to M_j} \sum_{\substack{M_k \to M_j \\   k \neq i}} l(M_k) + \sum_{\gamma: M_i \to M_k} c_{kk}l(M_k)$$
 and
 \item $$\sum_{M_j \to M_i}(Cm)_j =  l(M_i) + \sum_{M_j \to M_i} \sum_{\substack{M_j \to M_k \\   k \neq i}} l(M_k) + \sum_{\gamma: M_k \to M_i} d_{kk}l(M_k).$$
\end{enumerate}

\end{Lemma}
\proof
\\If $M_i$ is a module without successors in $\Sigma$, then the left hand side of the equation is zero, while the right hand side is reduced to $l(M_i) + \sum_{\gamma: M_i \to M_k} c_{kk}l(M_k)$. Since only the trivial path starts in $M_i$, the last sum contains only one summand $c_{ii}l(M_i) = - l(M_i)$ by Lemma \ref{Coxetersumformula} and hence the right hand side is also zero. Let $M_i$ now be any module in $\Sigma$, then we have 
$$\sum_{M_i \to M_j}(C^{-1}m)_j = \sum_{M_i \to M_j}(\sum_{M_k \to M_j} l(M_k) + \sum_{\gamma: M_j \to M_k}c_{kk} l(M_k))$$
$$ = \sum_{M_i \to M_j}\sum_{\substack{M_k \to M_j \\   k \neq i}} l(M_k) + \sum_{M_i \to M_j}(l(M_i) + \sum_{\gamma: M_j \to M_k}c_{kk} l(M_k))$$
$$ = \sum_{M_i \to M_j}\sum_{\substack{M_k \to M_j \\   k \neq i}} l(M_k) + (c_{ii} + 1)l(M_i) + \sum_{M_i \to M_j}\sum_{\gamma: M_j \to M_k}c_{kk} l(M_k))$$
$$= l(M_i) + \sum_{M_i \to M_j} \sum_{\substack{M_k \to M_j \\   k \neq i}} l(M_k) + \sum_{\gamma: M_i \to M_k} c_{kk}l(M_k)$$
\qed

\begin{Cor}\label{Coxeterformula}
For a module $M_i$ in $\Sigma$ we have
\begin{enumerate}[(a)]
 \item $$(C^{-1}m)_i =  - l(M_i) + \sum_{M_j \to M_i} l(M_j) + \sum_{M_i \to M_j}(C^{-1}m)_j - \sum_{M_i \to M_j} \sum_{\substack{M_k \to M_j \\   k \neq i}} l(M_k)$$
 and
 \item $$(Cm)_i =  - l(M_i) + \sum_{M_i \to M_j} l(M_j) + \sum_{M_j \to M_i}(Cm)_j - \sum_{M_j \to M_i} \sum_{\substack{M_j \to M_k \\   k \neq i}} l(M_k).$$
 
\end{enumerate}

\end{Cor}
\proof
\\Again we start by considering the case where $M_i$ does not have successors in $\Sigma$. Then we have
$$ - l(M_i) + \sum_{M_j \to M_i} l(M_j) + \sum_{M_i \to M_j}(C^{-1}m)_j - \sum_{M_i \to M_j} \sum_{\substack{M_k \to M_j \\   k \neq i}} l(M_k)$$
$$ =  - l(M_i) + \sum_{M_j \to M_i} l(M_j) = (C^{-1}m)_i.$$
\\For an arbitrary $M_i$ the statement follows from applying Lemma \ref{Coxetercombinatorics} to the right hand side of the equation, giving 
$$ - l(M_i) + \sum_{M_j \to M_i} l(M_j) + \sum_{M_i \to M_j}(C^{-1}m)_j - \sum_{M_i \to M_j} \sum_{\substack{M_k \to M_j \\   k \neq i}} l(M_k)$$
$$ = \sum_{M_j \to M_i} l(M_j) + \sum_{\gamma: M_i \to M_k} c_{kk}l(M_k) = (C^{-1}m)_i.$$
\qed

\begin{Th}\label{Coxeterinjectivetheorem}
Let $\Gamma$ be a connected component of $\Gamma_\Omega$ containing a sectional subgraph $\Sigma$ that consists of modules $M_1, \ldots, M_n$ without multiple arrows between them. Let $C^{-1} = (c_{ij})$ be the inverse of the Coxeter-matrix $C$ of $\Sigma$. Furthermore, denote the column vectors $(l(M_1), \ldots , l(M_n))^T$, $(l(\tau_\Omega(M_1)), \ldots, l(\tau_\Omega(M_n)))^T$ and $(l(\tau_\Omega^{-1}(M_1)), \ldots, l(\tau_\Omega^{-1}(M_n)))^T$ by $m$, $\tau_\Omega(m)$ and $\tau_\Omega^{-1}(m)$ respectively. 
\begin{enumerate}[(a)]
 \item If none of the $M_i$ is projective, then 
 $$ \tau_\Omega(m) = C^{-1}m + \dim(I),$$
where $I$ is an injective $K\Sigma$-module. 
 \item Dually, if none of the $M_i$ is Ext-injective, then 
 $$ \tau_\Omega^{-1}(m) = Cm + \dim(P),$$
where $P$ is a projective $K\Sigma$-module. 
\end{enumerate}

\end{Th}
\proof
\\We use the same notation as before and prove
$$\tau_\Omega(m)_i = l(\tau_\Omega(M_i)) = (C^{-1}m)_i + \sum_{\gamma : M_i \to M_u} C_u,$$
where $C_u$ are non-negative integers only depending on $u$. Let us define $M_i > M_j$ if there is a path from $M_j$ to $M_i$ and let $M_i$ be a maximal module in this sense. We use induction on this ordering to prove the equation above. Consider the almost split sequence
$$\begin{xy}
  \xymatrix{ 
0 \ar[r] & \tau_\Omega(M_i) \ar[r] &  \bigoplus\limits_{M_j \to M_i}^{\ } M_j \oplus T	\ar[r]	&	M_i \ar[r]	&	0.	    }
\end{xy}$$
Therefore, we have
$$l(M_i) + l(\tau_\Omega(M_i)) = \sum_{M_j \to M_i} l(M_j) + l(T).$$
We reorder the summands in Corollary \ref{Coxeterformula} and obtain
$$l(M_i) = - (C^{-1}m)_i + \sum_{M_j \to M_i}l(M_j)$$
since $M_i$ does not have any successors. Substituting this into the previous equation we obtain
$$- (C^{-1}m)_i + \sum_{M_j \to M_i}l(M_j) + l(\tau_\Omega(M_i)) = \sum_{M_j \to M_i} l(M_j) + l(T).$$
We easily conclude $l(\tau_\Omega(M_i)) = (C^{-1}m)_i + l(T)$ and observe that we want $l(T) = \sum_{\gamma : M_i \to M_u} C_u$. Note that this sum contains only the summand obtained by the trivial path from $M_i$ to $M_i$, hence we set $C_i = l(T)$, which is a non-negative integer. This shows the claim for a maximal $M_i$. Suppose now $M_i$ is a module such that the statement has been proved for all its successors. There is an almost split sequence
$$\begin{xy}
  \xymatrix{ 
  0 \ar[r] & \tau_\Omega(M_i) \ar[r] &  \bigoplus\limits_{M_j \to M_i}^{\ } M_j \bigoplus\limits_{M_i \to M_j}^{\ } \tau_\Omega(M_j) \oplus T_i	\ar[r]	&	M_i \ar[r]	&	0	    }
\end{xy}$$
from which we obtain
$$l(M_i) + l(\tau_\Omega(M_i)) = \sum_{M_j \to M_i} l(M_j) + \sum_{M_i \to M_j} l(\tau_\Omega(M_j)) + l(T_i).$$
We apply the induction hypothesis to all $l(\tau_\Omega(M_j))$ and obtain
$$l(M_i) + l(\tau_\Omega(M_i)) = \sum_{M_j \to M_i} l(M_j) + l(T_i) + \sum_{M_i \to M_j} ((C^{-1}m)_j +  \sum_{\gamma : M_j \to M_u} C_u ). $$
In order to be able to use Corollary \ref{Coxeterformula}, we reorder the summands, add 
$$ \sum_{M_i \to M_j} \sum_{\substack{M_k \to M_j \\   k \neq i}} l(M_k) - \sum_{M_i \to M_j} \sum_{\substack{M_k \to M_j \\   k \neq i}} l(M_k)$$
to the equation and deduce
$$l(\tau_\Omega(M_i)) = -l(M_i) + \sum_{M_j \to M_i} l(M_j) + \sum_{M_i \to M_j} ((C^{-1}m)_j - \sum_{M_i \to M_j}\sum_{\substack{M_k \to M_j \\   k \neq i}} l(M_k)$$
$$ + \sum_{M_i \to M_j}\sum_{\substack{M_k \to M_j \\   k \neq i}} l(M_k) + l(T_i) + \sum_{M_i \to M_j} \sum_{\gamma : M_j \to M_u} C_u $$
$$ = (C^{-1}m)_i + \sum_{M_i \to M_j}\sum_{\substack{M_k \to M_j \\   k \neq i}} l(M_k) + l(T_i) + \sum_{\substack{\gamma : M_i \to M_u \\ i \neq u}} C_u.$$
We set 
$$C_i = \sum_{M_i \to M_j}\sum_{\substack{M_k \to M_j \\   k \neq i}} l(M_k) + l(T_i),$$
which obviously is a non-negative integer, and obtain
$$\tau_\Omega(m)_i = l(\tau_\Omega(M_i)) = (C^{-1}m)_i + \sum_{\gamma : M_i \to M_u} C_u.$$
For convenience, let $v \in \Z^n$ denote the column vector whose $i$-th entry is given by $\sum_{\gamma : M_i \to M_u} C_u$. Consequently, we have
$$\tau_\Omega(m) = C^{-1}m + v,$$
so it remains to show that $v$ is the dimension vector of an injective $K\Sigma$-module. Let $p_{ij}$ denote the number of paths from $M_i$ to $M_j$ in $\Sigma$. Then we have
$$v = \sum_{u=1}^n (p_{1u} C_u, \ldots , p_{nu} C_u)^T = \sum_{u=1}^n C_u (p_{1u}, \ldots , p_{nu})^T, $$
but $(p_{1u}, \ldots , p_{nu})^T$ is precisely the dimension vector of the injective $K\Sigma$-module $I_u$ corresponding to the vertex $M_u$ since $K\Sigma$ is a hereditary algebra. This finishes the proof. \qed
\bigskip
\\For the rest of the section we choose a sectional subgraph $\Sigma$ such that none of its vertices has more than one immediate predecessor, for example by considering a sectional subgraph that consists of all sectional paths starting in a fixed module. As a consequence, the formula 
$$C_i = \sum_{M_i \to M_j}\sum_{\substack{M_k \to M_j \\   k \neq i}} l(M_k) + l(T_i)$$
reduces to $C_i = l(T_i)$, where $T_i$ is again the module that completes the middle term of $\E_\Omega(M_i)$ in the following sense.
$$\begin{xy}
  \xymatrix{ 
  0 \ar[r] & \tau_\Omega(M_i) \ar[r] &  \bigoplus\limits_{M_j \to M_i}^{\ } M_j \bigoplus\limits_{M_i \to M_j}^{\ } \tau_\Omega(M_j) \oplus T_i	\ar[r]	&	M_i \ar[r]	&	0	    }
\end{xy}$$
We deduce that in this case $\tau_\Omega(m) = C^{-1}m$ if and only if $\Sigma$ is a full sectional subgraph or there are some indecomposable projective modules $P_1, \ldots, P_k$ and arrows $M_{i_j} \to P_j$ in $\Gamma$ such that the union of $\Sigma$ with these modules and arrows is a full sectional subgraph.

\begin{Def}
Let $Q$ be a quiver, then a map $\lambda : \Z Q \to \Z$ is called a subadditive function on $\Z Q$ if
$$\lambda (i,k) + \lambda (i, k-1) \geq \sum_{(i,k) \to (j, k')} \lambda (j, k')$$
for all $(i,k) \in \Z Q$. 
\end{Def}

Note that this definition is different from the definition of subadditive functions on simple-laced graphs, which we introduced in Section \ref{treetypesection}. Let $\Gamma$ be a connected component of $\Gamma_l$ containing a sectional subgraph $\Sigma$ consisting of modules $M_1, \ldots, M_n$. We can then define a subadditive function on $\Z\Sigma$ in the following way. First of all we set $\lambda(M_i,k) = l(\tau_\Omega^{k}(M_i))$ for all $M_i$ and $k \in \Z$ such that $\tau_\Omega^{k}(M_i)$ exists. Note that since $\Gamma$ is left stable, $\tau_\Omega^{k}(M_i)$ always exists for $k \geq 0$. If $\tau_\Omega^{k}(M_i)$ does not exist for some $M_i$ and $k < 0$, we set $\lambda(i , k-1) = -\lambda (i,k) + \sum_{(i,k) \to (j, k')} \lambda (j, k')$. Hence we have inductively defined a subadditive function on $\Z \Sigma$.
\bigskip
\\Recall that the \textbf{defect} of a quiver $Q$ is the lowest positive integer $d = d_Q$ such that $C^{-d}x - x$ is in $\rad_q$ for all $x \in \Z^n$, where $q$ is the Euler quadratic form of the underlying graph of $Q$ with undirected edges. The \textbf{defect of a vector} $x$ is the integer $\partial(x)$ such that $C^{-d}x - x = \partial(x) h_q$, where $h_q$ is the generator of $\rad_q$. Note that $\partial : \Z^n \to \Z$ is clearly a morphism of abelian groups. For further detail we refer to \cite[XI.1]{SS07}. Analogously to sectional subgraphs, the type of a quiver $Q$ is the undirected graph associated to $Q$.

\begin{Lemma}\cite[Proposition 4.4]{BB83}\label{subadditivefunctiontendstoinfinity}
Let $Q$ be a quiver of Euclidean type with $n$ vertices and let $\lambda$ be a subadditive function on $\Z Q$. We denote the vector $(\lambda(1,k), \ldots, \lambda(n,k))$ by $\lambda_k$. For convenience, for two vectors $x,y \in \Z^n$ we say $x > y$ if and only if $x_i > y_i$ for all $i = 1, \ldots, n$.
\begin{enumerate}[(a)]
\item If $\lambda_k > 0$ for $0 \leq k < n$ and $\partial(\lambda_0) \geq 0$, then $\lambda_k > 0$ for all $k \geq 0$. Moreover, if $\partial(\lambda_l) > 0$ for some $l \geq 0$, then $\lambda(i,k) \to \infty$ when $k \to \infty$ for all $i = 1, \ldots, n$.
\item Dually, if $\lambda_k > 0$ for $0 \geq k > -n$ and $\partial(\lambda_0) \leq 0$, then $\lambda_k > 0$ for all $k \leq 0$. Moreover, if $\partial(\lambda_l) < 0$ for some $l \leq 0$, then $\lambda(i,k) \to \infty$ when $k \to - \infty$ for all $i = 1, \ldots, n$.
\end{enumerate} 
\end{Lemma}

\begin{Lemma}\cite[XI. Proposition 1.2]{SS07}\label{Coxeterinleftstableksigmacomponent}
Let $M$ be an indecomposable $K\Sigma$-module.
\begin{enumerate}[(a)]
 \item If $M$ is not contained in the preprojective component of $\Gamma_{K\Sigma}$, then we have
$$\dim(\tau^k(M)) = C^{-k}\dim(M)$$
for all $k \geq 0$.
\item  Dually, if $M$ is not contained in the preinjective component of $\Gamma_{K\Sigma}$, then we have
$$\dim(\tau^{-k}(M)) = C^k \dim(M)$$
for all $k \geq 0$. 
\end{enumerate}
\end{Lemma}


\begin{Lemma}\cite[XI. Lemma 1.5]{SS07}\label{Coxeterdelinvariant}
Let $Q$ be a Euclidean quiver with $n$ vertices and let $C^{-1}$ denote the inverse of its Coxeter matrix. Then 
$$\partial(x) = \partial(C^{-1}(x))$$
for all $x \in \Z^n$.
 
\end{Lemma}

Note that the following theorem is a generalization of Theorem \ref{treetypedynkin}.

\begin{Th}\label{Euclideantoinfinity}
Let $\Sigma$ be a sectional subgraph of Euclidean type in a connected component $\Gamma$ of $\Gamma_\Omega$. 
\begin{enumerate}[(a)]
\item If all modules in $\Sigma$ are left stable, then $l(\tau_\Omega^k(X)) \to \infty$ as $k \to \infty$ for each vertex $X$ of $\Sigma$.
\item If all modules in $\Sigma$ are right stable, then $l(\tau_\Omega^k(X)) \to \infty$ as $k \to - \infty$ for each vertex $X$ of $\Sigma$.
\end{enumerate}

\end{Th}

\proof
\\If $\Sigma$ contains multiple arrows, then the statement follows immediately from Theorem \ref{multiplearrows}. Hence we can assume that $\Sigma$ only contains arrows with trivial valuation. Since $\overline{\Sigma}$ is Euclidean, the Auslander-Reiten quiver of $K\Sigma$ has an infinite preinjective component. By Lemma \ref{Coxeterinleftstableksigmacomponent} we have $\dim(\tau^k(I')) = C^{-k}\dim(I')$ for an indecomposable injective $K\Sigma$-module $I'$ and all $k \in \N$, so $C^{-k}\dim(I')$ is always the dimension vector of a $K\Sigma$-module. Therefore, $C^{-d}\dim(I') - \dim(I')$ cannot have only negative entries as all entries in $\dim(I')$ are at most one by choice of $Q$. 
\bigskip
\\If $C^{-d}\dim(I') - \dim(I') = 0$, then by Lemma \ref{Coxeterinleftstableksigmacomponent} we obtain $l(\tau_\Omega^k(I')) \leq p_{I'}$ for all $k \geq 0$ and a constant $p_{I'}$. There are precisely $n$ orbits in the preinjective component of $K\Sigma$ since $n$ denotes the number of vertices in $\Sigma$ and, as we have seen in Lemma \ref{irreducibleboundary}, there is a constant $p$ such that $l(Z') \leq pl(Z)$ for any arrow $Z' \to Z$ in the Auslander-Reiten quiver of $\Omega$. Hence we have $l(X) \leq p^n p_{I'}$ for every module $X$ in the preinjective component of $K\Sigma$. As a result, every morphism from a projective module $P$ to $I'$ must be zero by Lemma \ref{HaradaSai}, which is a contradiction. Thus $\partial(\dim(I'))$ is a positive integer and the defect of any non-zero injective module is also positive.
\bigskip
\\By Theorem \ref{Coxeterinjectivetheorem} and Lemma \ref{Coxeterdelinvariant} we have $\partial(\tau_\Omega(m)) = \partial(C^{-1}m) + \partial(\dim(I)) = \partial(m) + \partial(\dim(I))$ and, consequently, $\partial(\tau_\Omega^{k'}(m)) \geq \partial(\tau_\Omega^k(m))$ if $k' > k$. Suppose that $\partial(\tau_\Omega^k(m)) \leq 0$ for all $k \geq 0$. Hence there must be an integer $r \geq 0$ such that for all $k \geq r$ we have $\partial(\tau_\Omega^k(m)) = \partial(\tau_\Omega^r(m)) \leq 0$. Without loss of generality, we can assume $r = 0$ and that $\Sigma$ is a full sectional subgraph left of all injective and projective modules, which gives us $\tau_\Omega^k(m) = C^{-k}m$ for all $k \geq 0$. 
\bigskip
\\If $\partial(m) < 0$, then there is a $t \in \N$ such that $0 > C^{-t}m = \tau_\Omega^t(m)$, which is impossible as each entry in the latter is given by the length of an $A$-module. Hence $\partial(m) = 0$, then by definition the entries of $C^{-k}m = \tau_\Omega^k(m)$ are bounded for all $k \geq 0$ as $C^{-d}m - m = 0$. We know there is a non-zero morphism $f \in \rad_\Omega^\infty(P, M_1)$ for a projective module $P$. Since $\Sigma$ is a full sectional subgraph left of all injective and projective modules, $f$ must factor over infinitely many modules of the form $\tau_\Omega^k(M_i)$ with $k \geq 0$. So by Lemma \ref{HaradaSai} $f$ must be zero, which is a contradiction. Hence there is a $k \geq 0$ such that $\partial(\tau_\Omega^k(m)) > 0$. Let $\lambda$ be the map that sends a module to its Jordan-H\"{o}lder length. If necessary, we extend $\lambda$ to a subadditive function $\lambda : \Z\Sigma \to \Z$ in the aforementioned way. Then by Lemma \ref{subadditivefunctiontendstoinfinity} $l(\tau_\Omega^k(M_i)) = \lambda(k, M_i) \to \infty$ 
as $k \to \infty$, which 
completes the proof. \qed

\begin{Lemma}\label{nottoinfinity}
Let $\Gamma$ be a connected component of $\Gamma_l$. Assume there is a module $X$ in $\Gamma$ such that $l(\tau_\Omega^n(X))$ does not tend to infinity for $n \to \infty$. Then
\begin{enumerate}[(a)]
\item $l(\tau_\Omega^n(Y))$ does not tend to infinity as $n \to \infty$ for every module $Y$ in $\Gamma$,
\item there is no sectional subgraph of Euclidean type in $\Gamma$ and
\item every arrow $Z \to Y$ in $\Gamma$ has finite global left degree.
\end{enumerate}
Dually, let $\Gamma$ be a connected component of $\Gamma_r$. Assume there is a module $X$ in $\Gamma$ such that $l(\tau_\Omega^n(X))$ does not tend to infinity for $n \to - \infty$. Then
\begin{enumerate}[(a)]
 \setcounter{enumi}{3}
\item $l(\tau_\Omega^n(Y))$ does not tend to infinity as $n \to - \infty$ for every module $Y$ in $\Gamma$,
\item there is no sectional subgraph of Euclidean type in $\Gamma$ and
\item every arrow $Y \to Z$ in $\Gamma$ has finite global right degree.
\end{enumerate}
\end{Lemma}

\proof
\\As we have seen in Lemma \ref{irreducibleboundary}, there is a constant $p$ such that $l(Z') \leq pl(Z)$ for any arrow $Z' \to Z$ in the Auslander-Reiten quiver of $\Omega$. Moreover, by our assumptions there is a constant $p_X$ such that $l(\tau_\Omega^n(X)) \leq p_X$ for infinitely many positive integers $n$. Let $Y$ be a module in $\Gamma$ that is not in the $\tau_\Omega$-orbit of $X$, then we know there is a sectional path of length $k$ from $\tau_\Omega^m(Y)$ to $\tau_\Omega^n(X)$ for some $n, m \in \N$ by Lemma \ref{leftpathexistence} as $\Gamma$ is left stable. If we set $p_Y = p_Xp^k$, we have $l(\tau_\Omega^n(Y)) \leq p_Y$ for infinitely many $n \geq m$, which proves the first statement. From this and Theorem \ref{Euclideantoinfinity} we immediately deduce that $\Gamma$ cannot contain a sectional subgraph of Euclidean type.
\bigskip
\\Let $Y$ and $Z$ be modules in $\Gamma$ such that there is an arrow from $Z$ to $Y$. By our assumptions there exists a constant $p_Y \in \N$ such that there are infinitely many $\tau_\Omega^n(Y)$ such that $l(\tau_\Omega^n(Y)) \leq p_Y$. Let $m \in \N$ be an integer such that $l(\tau_\Omega^n(Y)) \leq p_Y$ for more than $2^{p_Y}$ integers $0 \leq n \leq m$. Consequently, the path 
$$\begin{xy}
  \xymatrix{ 
\tau_\Omega^m(Z) \ar[r]^{f_m} &\tau_\Omega^m(Y) \ar[r]^{g_m}     &   \cdots \ar[r]^{g_2} &\tau_\Omega(Z) \ar[r]^{f_1}	&\tau_\Omega(Y) \ar[r]^{g_1}	&Z \ar[r]^{f_0} &Y  }
\end{xy}$$
is labeled with irreducible morphisms such that $f_0g_1f_1 \cdots g_mf_m = 0$ by Lemma \ref{HaradaSai}. Hence every arrow in $\Gamma$ has finite global left degree. \qed

\newpage

\begin{Cor}\label{Ainftyinftytoinfinity}
\mbox{}
\begin{enumerate}[(a)]
\item Let $\Gamma$ be a connected component of $\Gamma_l$ such that its left subgraph type is $A_\infty^\infty$. Then $l(\tau_\Omega^n(X)) \to \infty$ as $n \to \infty$ for all modules $X$ in $\Gamma$.
\item Dually, if $\Gamma$ is a connected component of $\Gamma_r$ such that its right subgraph type is $A_\infty^\infty$, then $l(\tau_\Omega^n(X)) \to \infty$ as $n \to - \infty$ for all modules $X$ in $\Gamma$.
\end{enumerate}
\end{Cor}

\proof
\\Since $\Gamma$ is left stable and of left subgraph type $A_\infty^\infty$, for every module $X$ in $\Gamma$ there is an $n \in \N$ such that there are two infinite sectional paths 
$$\begin{xy}
  \xymatrix{ 
\cdots \ar[r] &	X_{i+1} \ar[r] &X_i \ar[r] &\cdots \ar[r] &  X_1 \ar[r] & \cdots \ar[r]&  X_0 = X  }
\end{xy}$$
and 
$$\begin{xy}
  \xymatrix{ 
\cdots \ar[r] &	Y_{i+1} \ar[r] &Y_i \ar[r] &\cdots \ar[r] &  Y_1 \ar[r] & \cdots \ar[r]&  Y_0 = X  }
\end{xy}$$
ending in $X$ with $Y_1 \ncong X_1$. Consequently, by Lemma \ref{twosectionalpathsinfiniteglobaldegree} all arrows $X_{i+1} \to X_i$ have infinite global left degree. Hence the Corollary follows immediately from the contrapositive of the previous lemma. \qed

\begin{Lemma}\label{sectionalpathfromorbittoitself}
\mbox{}
\begin{enumerate}[(a)]
 \item Suppose there is a module $X$ in a connected component $\Gamma$ of $\Gamma_l$ such that $l(\tau_\Omega^n(X))$ does not tend to infinity as $n \to \infty$. If there is a sectional path in $\Gamma$ which meets some $\tau_\Omega$-orbit more than once, then $\Gamma$ consists of $\tau_\Omega$-periodic modules.
 \item Dually, let $X$ be a module in a connected component $\Gamma$ of $\Gamma_r$ such that $l(\tau_\Omega^n(X))$ does not tend to infinity as $n \to - \infty$. If there is a sectional path in $\Gamma$ which meets some $\tau_\Omega$-orbit more than once, then $\Gamma$ consists of $\tau_\Omega$-periodic modules.
\end{enumerate}

\end{Lemma}

\proof
\\Let $p$ be the constant from Lemma \ref{irreducibleboundary} such that for all modules $Y$ and $Z$ we have $l(Y) \leq pl(Z)$ if there is an $\Omega$-irreducible morphism from $Y$ to $Z$. By Lemma \ref{nottoinfinity} and our assumptions there is a constant $p_Y$ for every module $Y$ in $\Gamma$ such that $l(\tau_\Omega^n(Y)) \leq p_Y$ for infinitely many positive integers $n$. It follows that for every $r \in \N$ there are infinitely many positive integers $m$ such that $l(\tau_\Omega^{mr}(Y)) \leq p^{2r}p_Y$. Suppose that 
$$\begin{xy}
  \xymatrix{ 
\tau_\Omega^r(Y_0) = Y_n \ar[r] & Y_{n-1} \ar[r] &\cdots \ar[r] &  Y_1 \ar[r]	&Y_0  }
\end{xy}$$
is a sectional path in $\Gamma$ with some $r \geq 0$. Without loss of generality, we may assume that $Y_{n-1}$ and $Y_1$ are not contained in the same $\tau_\Omega$-orbit if $n > 2$. Then there is an infinite path
$$\begin{xy}
  \xymatrix{ 
\cdots \ar[r]  & \tau_\Omega^{2r}(Y_1) \ar[r]  & \tau_\Omega^{2r}(Y_0) \ar[r] & \tau_\Omega^r(Y_{n-1}) \ar[r] & \cdots  }
\end{xy}$$
$$\begin{xy}
  \xymatrix{ 
\cdots \ar[r]& \tau_\Omega^r(Y_1) \ar[r] & \tau_\Omega^r(Y_0) \ar[r] & Y_{n-1} \ar[r] &\cdots \ar[r] &  Y_1 \ar[r]	&Y_0.  }
\end{xy}$$
Since there are infinitely many modules $\tau_\Omega^{mr}(Y)$ on this path with $l(\tau_\Omega^{mr}(Y)) \leq p^{2r}p_Y$, the path cannot be sectional by Corollary \ref{non-zero} and Lemma \ref{HaradaSai}. Hence $Y_{n-1}$ and $Y_1$ belong to the same $\tau_\Omega$-orbit and $n$ is either $1$ or $2$. If $n = 1$, then $Y_0 = \tau_\Omega^{2r-1}(Y_0)$. On the other hand, if $n = 2$, then $r \neq 1$ as the path 
$$\begin{xy}
  \xymatrix{ 
\tau_\Omega^r(Y_0) = Y_2 \ar[r] &  Y_1 \ar[r]	&Y_0  }
\end{xy}$$
is sectional and we obtain $Y_1 = \tau_\Omega^{r-1}(Y_1)$. In both cases $\Gamma$ contains a $\tau_\Omega$-periodic module, and hence $\Gamma$ entirely consists of $\tau_\Omega$-periodic modules by Lemma \ref{orbit} as $\Gamma$ is left stable. For $r < 0$ we prove the statement similarly using the infinite path 
$$\begin{xy}
  \xymatrix{ 
\tau_\Omega^r(Y_0) \ar[r] & Y_{n-1} \ar[r] &\cdots \ar[r] &  Y_1 \ar[r]	&Y_0 \ar[r] & \tau_\Omega^{-r}(Y_{n-1})\ar[r] & \cdots  }
\end{xy}$$
$$\begin{xy}
  \xymatrix{ 
\cdots \ar[r]  & \tau_\Omega^{-r}(Y_1) \ar[r]  & \tau_\Omega^{-r}(Y_0) \ar[r] & \tau_\Omega^{-2r}(Y_{n-1}) \ar[r] & \cdots  }
\end{xy}$$ \qed

\begin{Cor}\label{coraytoinfinity}
\mbox{}
\begin{enumerate}[(a)]
 \item Let $\Gamma$ be a coray tube, then $l(\tau_\Omega^n(X)) \to \infty$ as $n \to \infty$ for all modules $X$ in $\Gamma$.
 \item Let $\Gamma$ be a ray tube, then $l(\tau_\Omega^n(X)) \to \infty$ as $n \to - \infty$ for all modules $X$ in $\Gamma$.
\end{enumerate} 
\end{Cor}

\proof
\\Suppose $l(\tau_\Omega^n(X))$ does not tend to infinity as $n \to \infty$. Since a coray tube is the same as a helical component by Theorem \ref{helicalcoray}, we know there is an oriented cycle from an injective module $I$ to itself in $\Gamma$. Since $\tau_\Omega^n(I) = I$ is only possible for $n = 0$, there is a sectional path from $I$ to $\tau_\Omega^n(I)$ in $\Gamma$ for some $n \geq 0$ by Lemma \ref{leftpathexistence}. Consequently, $\Gamma$ consists of $\tau_\Omega$-periodic modules by the previous lemma, which is a contradiction. \qed

\begin{Lemma}\label{threeleftstablemiddletermsinfinity}
Let $X$ be a non-periodic module in $\Gamma_\Omega$.
\begin{enumerate}[(a)]
 \item If $X$ is left stable and the almost split sequence $\E_\Omega(X)$ admits three left stable middle terms, then $l(\tau_\Omega^n(X)) \to \infty$ as $n \to \infty$.
 \item If $X$ is right stable and the almost split sequence $\E'_\Omega(X)$ admits three right stable middle terms, then $l(\tau_\Omega^n(X)) \to \infty$ as $n \to - \infty$.
\end{enumerate}
\end{Lemma}

\proof
\\Let $\Gamma$ be the connected component of $X$ in $\Gamma_l$, which does not contain $\tau_\Omega$-periodic modules by Lemma \ref{orbit}. Suppose $l(\tau_\Omega^n(X))$ does not tend to infinity as $n \to \infty$, then $\Gamma$ does not contain a sectional subgraph of Euclidean type by Lemma \ref{nottoinfinity}. Moreover, every sectional path in $\Gamma$ meets each $\tau_\Omega$-orbit at most once by Lemma \ref{sectionalpathfromorbittoitself}. Consequently, the left subgraph type of $\Gamma$ must be either $A_\infty$, $A_\infty^\infty$ or $D_\infty$. In all cases there is a infinite sectional path
$$\begin{xy}
  \xymatrix{ 
\cdots \ar[r] & X_{i+1} \ar[r] & X_i \ar[r] & \cdots \ar[r] &  X_1 \ar[r]	& X_0 = \tau_\Omega^r(X)  }
\end{xy}$$
for some $r \geq 0$. Since $\E_\Omega(X)$ has three left stable middle terms, so does $\E_\Omega(\tau_\Omega^r(X))$. Then by Lemma \ref{vmtinfinity} the arrow $\tau_\Omega^r(X_2) \to \tau_\Omega^r(X_1)$ has infinite global left degree, which contradicts Lemma \ref{nottoinfinity}. \qed

\begin{Th}\label{nottoinfinityAinfty}
\mbox{}
\begin{enumerate}[(a)]
 \item Let $\Gamma$ be an infinite connected component of $\Gamma_l$ that contains a module $X$ such that $l(\tau_\Omega^n(X))$ does not tend to infinity as $n \to \infty$. Then the left subgraph type of $\Gamma$ is $A_\infty$.
 \item Dually, if $\Gamma$ is an infinite connected component of $\Gamma_r$ containing a module $X$ such that $l(\tau_\Omega^n(X))$ does not tend to infinity as $n \to - \infty$, then its right subgraph type is $A_\infty$.
\end{enumerate}

\end{Th}

\proof
\\If $\Gamma$ is $\tau_\Omega$-periodic, then it is a stable tube by Theorem \ref{HPRT} and its left subgraph type is $A_\infty$. If $\Gamma$ is not $\tau_\Omega$-periodic, then we know that its left subgraph type is not given by a Dynkin diagram by Theorem \ref{maintheorem}. Moreover, since $l(\tau_\Omega(X))$ does not tend to infinity as $n \to \infty$, $\Gamma$ neither contains a sectional subgraph of Euclidean type nor of type $A_\infty^\infty$ by Lemma \ref{nottoinfinity} and Corollary \ref{Ainftyinftytoinfinity} respectively. The left subgraph type of $\Gamma$ also cannot be $D_\infty$ as $\Gamma$ does not contain a module with three immediate predecessor in $\Gamma$ by Lemma \ref{threeleftstablemiddletermsinfinity}. We deduce that the left subgraph type of $\Gamma$ is $A_\infty$ by process of elimination. \qed
\bigskip
\\Recall that we have already seen in Corollary \ref{coraytoinfinity} that the converse statement does not hold, i.e. there are connected components of $\Gamma_l$ with left subgraph type $A_\infty$ containing a module $X$ such that $l(\tau_\Omega^n(X)) \to \infty$ as $n \to \infty$.

\begin{Cor}\label{finitelymanyorbitsunbounded} 
Let $\Gamma$ be a connected component of $\Gamma_\Omega$. Suppose there are only finitely many $\tau_\Omega$-orbits in $\Gamma$. Then for each positive integer $k$ there are at most finitely many modules of length $k$ in $\Gamma$. 
\end{Cor}

\proof
\\We assume that there is a positive integer $k$ such that there are infinitely many modules of length $k$ in $\Gamma$. Since there are only finitely many $\tau_\Omega$-orbits in $\Gamma$, there exists a module $X$ which is not $\tau_\Omega$-periodic such that $l(\tau_\Omega^n(X)) = k$ for infinitely many integers $n$. Without loss of generality, we assume that there are infinitely many positive integers $n$ with this property. Then by Theorem \ref{nottoinfinityAinfty} we know that $X$ is contained in a connected component of $\Gamma_l$ such that its left subgraph type is $A_\infty$, which is a contradiction. \qed

\begin{Cor}
Suppose there are infinitely many indecomposable modules in $\Omega$ of length $k$ for some $k \in \N$. Then $\Gamma_\Omega$ admits infinitely many $\tau_\Omega$-orbits. 
\end{Cor}

\proof 
\\Suppose the statement is false and $\Gamma_\Omega$ consists of finitely many $\tau_\Omega$-orbits. Then there are clearly only finitely many connected components in $\Gamma_\Omega$, one of which must contain infinitely many indecomposable modules of length $k$. This contradicts Corollary \ref{finitelymanyorbitsunbounded}. \qed
\bigskip
\\We call the following statement the second Brauer-Thrall conjecture for subcategories. If $\Omega$ is infinite, then there are infinitely many positive integers $n_1, n_2, \ldots$ such that for each $i \in \N$ there are infinitely many non-isomorphic indecomposable modules of Jordan-H\"{o}lder length $n_i$ in $\Omega$. As a consequence of the last corollary and Theorem \ref{BrauerThrall1.5}, one could disprove the second Brauer-Thrall conjecture for subcategories by finding an example of an infinite functorially finite resolving subcategory $\Omega$ such that its Auslander-Reiten quiver consists of only finitely many $\tau_\Omega$-orbits. \cite{EMM10} contains an example of an infinite functorially finite resolving subcategory, of which only one connected component with finitely many orbits is known, but it is not clear if this component contains all indecomposable modules in $\Omega$. 

\begin{Th}
Let $\Gamma$ be a connected component of $\Gamma_\Omega$.
\begin{enumerate}[(a)]
 \item If $\Gamma$ is left stable and contains a module $X$ such that $l(\tau_\Omega^n(X))$ does not tend to infinity as $n \to \infty$, then there is a full sectional subgraph
$$\begin{xy}
  \xymatrix{ 
\cdots \ar[r] & X_{i+1} \ar[r] & X_i \ar[r] & \cdots \ar[r] &  X_1 \ar[r]	& X_0}
\end{xy}$$
such that there exists a non-negative integer $n$ with the property that $X_i$ is stable for all $i \geq n$ while $X_j$ is left stable but not right stable for all $0 \leq j < n$.
\item Dually, if $\Gamma$ is right stable and contains a module $X$ such that $l(\tau_\Omega^n(X))$ does not tend to infinity as $n \to - \infty$, then there is a full sectional subgraph
$$\begin{xy}
  \xymatrix{ 
X_0 \ar[r] & X_1 \ar[r] & \cdots \ar[r] &	X_i \ar[r] & X_{i+1} \ar[r] &  \cdots}
\end{xy}$$
such that there exists a non-negative integer $n$ with the property that $X_i$ is stable for all $i \geq n$ while $X_j$ is right stable but not left stable for all $0 \leq j < n$.
\item If $\Gamma$ is stable and contains a module $X$ such that there are infinitely many $\tau_\Omega^n(X)$ with Jordan-H\"{o}lder length at most $k$ for some $k \in \N$, then $\Gamma$ is either a stable tube or isomorphic to $\Z A_\infty$.
\end{enumerate}
\end{Th}

\proof
\\By Theorem \ref{nottoinfinityAinfty} the left subgraph type of $\Gamma$ is $A_\infty$ and $\Gamma$ is either a stable tube or not $\tau_\Omega$-periodic. In the former case the statement is trivial, so we suppose that $\Gamma$ is not $\tau_\Omega$-periodic. Let $\Sigma$ be a full sectional subgraph of $\Gamma$, then $\Sigma$ contains a module $X_0$ that has only one immediate predecessor in $\Gamma$. Since $\Gamma$ is left stable, we obtain a sectional path
$$\begin{xy}
  \xymatrix{ 
\cdots \ar[r] & X_{i+1} \ar[r] & X_i \ar[r] & \cdots \ar[r] &  X_1 \ar[r]	& X_0}
\end{xy}$$
which is a full sectional subgraph by Theorem \ref{subgraphtype}. Suppose now there is a $j \geq 1$ such that $X_j$ is not right stable but $X_{j-1}$ is right stable. Without loss of generality, we can assume that $j$ is minimal with that property. Since $\Gamma$ is a connected component of $\Gamma_\Omega$, we know that
$$\begin{xy}
  \xymatrix{ 
X_{j-1} \ar[r] & X_{j-2} \ar[r] & \cdots \ar[r] &  X_1 \ar[r]	& X_0}
\end{xy}$$
is a full sectional subgraph of a connected component $\Gamma'$ of $\Gamma_r$. Hence the right subgraph type of $\Gamma'$ is $A_j$ and it follows that the component is finite and $\tau_\Omega$-periodic by Theorem \ref{maintheorem}, which is a contradiction. Thus the proof of $(a)$ is completed and $(b)$ can be proved dually. 
\bigskip
\\If $\Gamma$ is stable and contains a module $X$ such that there are infinitely many $\tau_\Omega^n(X)$ with Jordan-H\"{o}lder length at most $k$ for some $k \in \N$, then $l(\tau_\Omega^n(X))$ or $l(\tau_\Omega^{-n}(X))$ does not tend to infinity as $n \to \infty$. So statement $(c)$ follows immediately from $(a)$ or $(b)$.\qed

\end{document}